\documentclass{gtart_h}

\def\ifplaintex{\expandafter\ifx\csname documentclass\endcsname\relax}

\ifplaintex 
\else
\fi

\def\gt{{\mathsurround=0pt\it $\cal G\mskip-2mu$eometry \&\ 
$\cal T\!\!$opology}}        

\def\gtm{{\mathsurround=0pt\it $\cal G\mskip-2mu$eometry \&\ 
$\cal T\!\!$opology $\cal M\mskip-1mu$onographs}}    

\def\gtp{{\mathsurround=0pt\it $\cal G\mskip-2mu$eometry \&\ 
$\cal T\!\!$opology $\cal P\!$ublications}}  

\def\recd{{\small Received:\qua\receiveddate\ifx\reviseddate\relax
\else\qquad Revised:\qua\reviseddate\fi\par}} 

\def\volumenumber#1{\def\thevolumenumber{#1}}
\def\volumeyear#1{\def\thevolumeyear{#1}}
\def\volumename#1{\def\thevolumename{#1}}
\def\papernumber#1{\def\thepapernumber{#1}}
\def\pagenumbers#1#2{\def\startpage{#1}\def\finishpage{#2}}
\def\published#1{\def\publishdate{#1}}
\def\received#1{\def\receiveddate{#1}}
\def\revised#1{\def\reviseddate{#1}}
\def\accepted#1{\def\accepteddate{#1}}

\def\coverauthors#1{\def\thecoverauthors{#1}}
\def\asciiauthors#1{\def\theasciiauthors{#1}}

\def\asciiurl#1{\def\theasciiurl{#1}}
\def\coverauthors#1{\def\thecoverauthors{#1}}
\long\def\asciiabstract#1{\long\def\theasciiabstract{#1}}

\def\shortauthors#1{\def\theshortauthors{#1}}

\let\\\par
\let\thevolumenumber\relax\let\thepapernumber\relax
\let\thevolumeyear\relax\let\startpage\relax
\let\finishpage\relax\let\publishdate\relax\let\receiveddate\relax
\let\reviseddate\relax\let\accepteddate\relax\let\theasciititle\relax
\let\theasciiauthors\relax
\let\theasciiabstract\relax
\let\thecoverauthors\relax
\let\thecoverauthors\relax\let\theerratum\relax\let\theasciiemail\relax
\let\theshortauthors\relax\let\theshorttitle\relax\let\theasciiurl\relax

\def\startpage{1}\def\finishpage{15}\def\thepapernumber{77}

\volumenumber{2}
\volumename{Proceedings of the Kirbyfest}
\volumeyear{1999}

\long\def\maketitlep{   

\count0=\startpage

\gtm\nl        
{\small Volume \thevolumenumber: \thevolumename\nl 
\ifx\theerratum\relax\else Erratum \erratumnumber\nl\fi
Pages \startpage--\finishpage\nl}

\vglue 0.1truein   

{\parskip=0pt\leftskip 0pt plus 1fil\def\\{\par\smallskip}{\ifplaintex\large
\else\Large\fi\bf\thetitle}\par\medskip}   
\vglue 0.05truein 

{\parskip=0pt\leftskip 0pt plus 1fil\def\\{\par}{\sc\theauthors}
\par\medskip}%
 
\vglue 0.03truein 

{\small\leftskip 25pt\rightskip 25pt{\bf Abstract}\stdspace\theabstract

{\bf AMS Classification}\stdspace\theprimaryclass
\ifx\thesecondaryclass\relax\else; \thesecondaryclass\fi\par
{\bf Keywords}\stdspace \thekeywords\par}\vglue 7pt

}   

\font\phead=cmsl9 scaled 950
\font=cmsl9 scaled 1050
\font\pnum=cmbx10 scaled 913
\font\lnum=cmbx10 
\font\pfoot=cmsl9 scaled 950
\font=cmsl9 scaled 1050
\ifplaintex
\headline{\vbox to 0pt{\vskip -4.5mm\line{\small\phead\ifnum
\count0=\startpage ISSN 1464-8997 (on line)
1464-8989 (printed) \hfill {\pnum\folio}\else\ifodd\count0\def\\{ }%
\ifx\theshorttitle\relax\thetitle\else\theshorttitle\fi\hfill{\pnum\folio}
\else\def\\{ and }{\pnum\folio}\hfill\ifx\theshortauthors\relax\theauthors
\else\theshortauthors\fi\fi\fi}\vss}}
\footline{\vbox to 0pt{\vglue 0mm\line{\small\pfoot\ifnum\count0=\startpage
Published \publishdate:\qua\copyright\ \gtp\hfill\else
\gtm, Volume \thevolumenumber\ (\thevolumeyear)\hfill\fi}\vss
}}
\else
\makeatletter
\def\@oddhead{{\small\lhead\ifnum\count0=\startpage ISSN 1464-8997 (on line)
1464-8989 (printed) \hfill {\lnum\number\count0}\else\ifodd\count0
\def\\{ }\ifx\theshorttitle\relax \thetitle \else\theshorttitle\fi\hfill
{\lnum\number\count0}\else\def\\{ and }{\lnum\number\count0}
\hfill\ifx\theshortauthors\relax 
\theauthors\else\theshortauthors\fi\fi\fi}}\def\@evenhead{@oddhead}
\def\@oddfoot{\small\lfoot\ifnum\count0=\startpage Published \publishdate:\qua\copyright\ \gtp\hfill\else
\gtm, Volume \thevolumenumber\ (\thevolumeyear)\hfill\fi}
\def\@evenfoot{@oddfoot}
\makeatother
\fi

\let\maketitlepage\maketitlep
\let\makeshorttitle\maketitlepage
\let\maketitle\maketitlepage

\newwrite\gtoutfile
\long\gdef\makeheadfile{  
{\def\\{, }\def\s{ }
\immediate\openout\gtoutfile head.xxx
\immediate\write\gtoutfile{Proxy-for: \ifx\theasciiauthors\relax
\theauthors\else\theasciiauthors\fi\s<\ifx\theasciiemail\relax\theemail\else\theasciiemail\fi>}
\immediate\write\gtoutfile{\noexpand\\}
\immediate\write\gtoutfile{Authors: \ifx\theasciiauthors\relax
\theauthors\else\theasciiauthors\fi}
{\def\\{ }\immediate\write\gtoutfile{Title: \ifx\theasciititle\relax
\thetitle\else\theasciititle\fi}}
\immediate\write\gtoutfile{Subj-class: GT or SG, GR etc}
\immediate\write\gtoutfile{MSC-class: \theprimaryclass\ifx\thesecondaryclass\relax\else, \thesecondaryclass\fi}
\immediate\write\gtoutfile{Journal-ref: Geom. Topol. Monogr. \thevolumenumber\s
(\thevolumeyear) \startpage-\finishpage}
\immediate\write\gtoutfile{Comments: Published by Geometry and Topology Monographs at}
\immediate\write\gtoutfile{\s\s\s  http://www.maths.warwick.ac.uk/gt/GTMon\thevolumenumber/paper\thepapernumber.abs.html}
\immediate\write\gtoutfile{\noexpand\\}
\immediate\write\gtoutfile{}
\ifx\theasciiabstract\relax
\immediate\write\gtoutfile{\theabstract}\else
\immediate\write\gtoutfile{\theasciiabstract}\fi
\immediate\write\gtoutfile{}
\immediate\write\gtoutfile{\noexpand\\}
\immediate\write\gtoutfile{}
\immediate\closeout\gtoutfile}}  

\def\maketitlepage{\maketitlep\makeheadfile}
\let\makeshorttitle\maketitlepage
\let\maketitle\maketitlepage

\volumenumber{4}
\volumename{Invariants of knots and 3-manifolds (Kyoto 2001)}
\volumeyear{2002}
\papernumber{24}
\pagenumbers{377}{572}
\received{28 December 2001}
\revised{4 December 2002 -- 8 April 2004}
\accepted{8 April 2004}
\published{1 June 2004}

\usepackage{amsmath}
\usepackage{amssymb}
\usepackage{amscd}
\usepackage{pstricks}
\usepackage{graphicx}
\allowdisplaybreaks

\usepackage{makeidx}
\makeindex

\newcommand{\C}{{\Bbb C}}
\newcommand{\cA}{{\cal A}}
\newcommand{\cB}{{\cal B}}
\newcommand{\cC}{{\cal C}}
\newcommand{\cF}{{\cal F}}
\newcommand{\cH}{{\cal H}}

\newcommand{\conn}{_{\mbox{\tiny conn}}}
\newcommand{\DD}{{\Delta\!\Delta}}
\newcommand{\dwn}[1]{\raisebox{-6pt}[0pt][0pt]{#1}}
\newcommand{\e}{{\varepsilon}}
\newcommand{\FF}{{\Bbb F}}

\newcommand{\hLMO}{\hat Z^{\mbox{\scriptsize L{\hspace{-1pt}}M{\hspace{-1pt}}O}}}
\newcommand{\hZ}{{Z}}
\newcommand{\K}{{\Bbb K}}
\newcommand{\LL}{{\Bbb L}}
\newcommand{\LMO}{Z^{\mbox{\scriptsize L{\hspace{-1pt}}M{\hspace{-1pt}}O}}}

\newcommand{\mboxsm}[1]{\mbox{\small #1}}
\newcommand{\mboxss}[1]{\mbox{\scriptsize #1}}
\newcommand{\M}{{\Bbb M}}
\newcommand{\MK}{\M\!\K}
\newcommand{\N}{{\Bbb N}}
\newcommand{\namae}{}
\newcommand{\ora}[1]{\mbox{$\overrightarrow{\mbox{#1}}$}}
\newcommand{\Q}{{\Bbb Q}}
\newcommand{\R}{{\Bbb R}}

\newcommand{\simsubC}[1]{\!\!\underset{{}^{C_{#1}}}{\sim}}
\newcommand{\simsubHL}[1]{\!\!\underset{{}^{H\!L_{#1}}}{\sim}}
\newcommand{\simsubY}[1]{\!\!\underset{{}^{Y_{#1}}}{\sim}}
\newcommand{\subkz}{_{{}_{\mbox{\tiny K{\hspace{-0.7pt}}Z}}}}
\newcommand{\vcA}{\overrightarrow{\cal A}}
\newcommand{\Z}{{\Bbb Z}}
\newcommand{\ZHS}{\mbox{$\Z HS$}}
\newcommand{\ZHSs}{\mbox{$\Z HS$'s}}

\def \ThKeCob {\mbox{\boldmath${\cal C}ob$\unboldmath}}
\def\ThKethrafill{$\mathsurround=0pt \mathord- \mkern-6mu 
\cleaders\hbox{$\mkern-2mu
\mathord- \mkern-2mu$}\hfill \mkern-6mu\mathord\twoheadrightarrow$}
\newcommand {\ThKeonto} [1]{\hbox to #1pt{\ThKethrafill}\,\,\,}
\def \ThKeV {\mbox{\boldmath$\cal V$\unboldmath}}
\def \ThKev {\mbox{\boldmath$\scriptstyle\cal V$\unboldmath}}
\newcommand{\ext}[1] {\mbox{\raisebox{.4ex}{$\bigwedge^{\!#1}$}}\mkern-1mu}
\def \ThKeQc{\mbox{\boldmath${\cal Q}$\unboldmath}}

\newcommand{\pict}[2]{\mbox{$\begin{array}{c}
   \includegraphics[height=#2]{fig-prob/#1.ps}
   \end{array}$}}

\newcommand{\picthw}[3]{\mbox{$\begin{array}{c}
   \includegraphics[height=#2,width=#3]{fig-prob/#1.ps}
   \end{array}$}}

\theoremstyle{plain}
   \newtheorem{thm}{Theorem}[section]

   \newtheorem{conj}[thm]{Conjecture}
   \newtheorem{quest}[thm]{Question}
   \newtheorem{prob}[thm]{Problem}
   \newtheorem{exe}[thm]{Exercise}
\theoremstyle{definition}
   
   \newtheorem*{rem}{Remark}
   \newtheorem*{exm}{Example}
   \newtheorem*{update}{Update}
   
\begin{document}

\makeatletter
\def\@oddhead{{\small\lhead\ifnum\count0=\startpage ISSN 1464-8997 (on line)
1464-8989 (printed) \hfill {\lnum\romannumeral\count0}\else\ifodd\count0
\def\\{ }\ifx\theshorttitle\relax \thetitle \else\theshorttitle\fi\hfill
{\lnum\romannumeral\count0}\else\def\\{ and }{\lnum\romannumeral\count0}
\hfill\ifx\theshortauthors\relax 
\theauthors\else\theshortauthors\fi\fi\fi}}\def\@evenhead{@oddhead}
\def\@oddfoot{\small\lfoot\ifnum\count0=\startpage Published \publishdate:\qua\copyright\ \gtp\hfill\else
\gtm, Volume \thevolumenumber\ (\thevolumeyear)\hfill\fi}
\def\@evenfoot{@oddfoot}
\makeatother

\count0=1

\cl{\Large\bf Problems on invariants~of~knots~and~3-manifolds}
\smallskip\cl{\sc{\rm Edited by} T. Ohtsuki}

\vspace{0.5pc}
\begin{flushright}
\pict{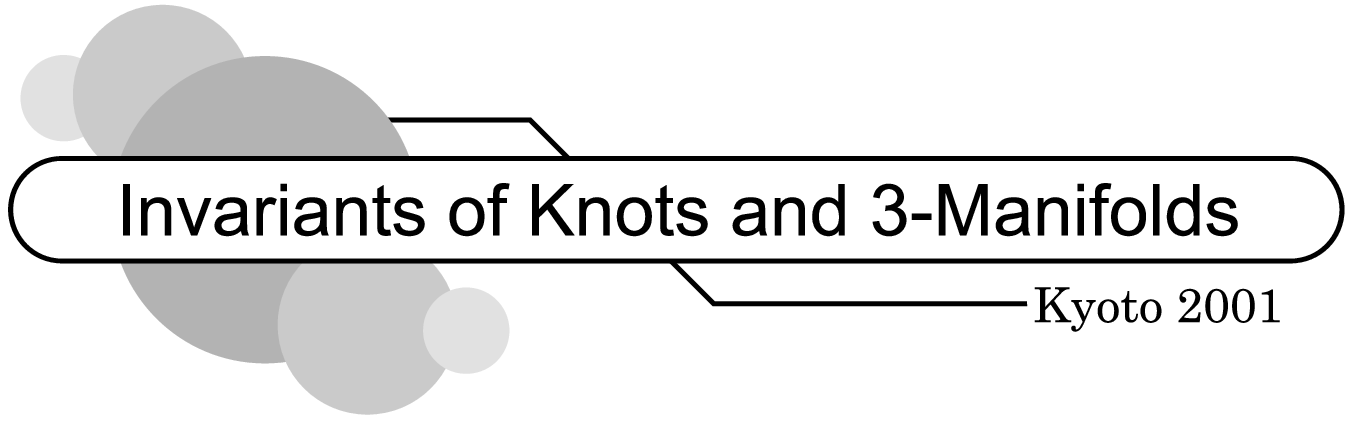}{1.5cm}
\end{flushright}
\vspace{-3pc}

\renewcommand{\thefootnote}{\fnsymbol{footnote}}
\footnotetext[0]{The logo for the workshop and the seminars 
was designed by N. Okuda.}
\renewcommand{\thefootnote}{\arabic{footnote}}

\section*{Preface}
\addcontentsline{toc}{section}{Preface}

{\small
The workshop and seminars on ``Invariants of knots and 3-manifolds''
was held at Research Institute for Mathematical Sciences, Kyoto University
in September 2001.
There were 25 talks in the workshop in September 17--21,
and there were 27 talks in the seminars in the other weeks of September.
Each speaker was requested to give his/her open problems
in a short problem session after his/her talk,
and many interesting open problems were given and discussed
by the speakers and participants in the workshop and the seminars.
Contributors of the open problems were also requested
to give kind expositions of history, background, significance, and/or 
importance of the problems.
This problem list was made by editing these open problems
and such expositions.\footnote{\footnotesize
Open problems on the Rozansky-Witten invariant
were written in a separate manuscript \cite{RoSa02}.
Some fundamental problems are quoted from other problem lists such as
\cite{Jones_prob}, \cite{Kirby}, \cite{Lin_website}, \cite{Morton},
\cite{Morton_Hellas}, \cite[Pages 571--572]{Turaev_book}.}

Since the interaction between geometry and mathematical physics in the 1980s, 
many invariants of knots and 3-manifolds have been discovered and studied.
The discovery and analysis of the enormous number of 
these invariants yielded a new area: 
the study of invariants of knots and 3-manifolds 
(from another viewpoint, the study of the sets of knots and 3-manifolds). 
Recent works have almost completed 
the topological reconstruction of the invariants derived from 
the Chern-Simons field theory,
which was one of main problems of this area.
Further, relations among these invariants have been studied enough well,
and these invariants are now well-organized.
For the future developments of this area,
it might be important to consider 
various streams of new directions;\footnote{\footnotesize
For example, directions related to other areas such as 
hyperbolic geometry via the volume conjecture and 
the theory of operator algebras via invariants arising from 6j-symbols.}
this is a reason why the editor tried to make the problem list expository.
The editor hopes this problem list will clarify 
the present frontier of this area and 
assist readers when considering future directions.

The editor will try to keep up-to-date information on this problem list 
at his web site.\footnote{\footnotesize
\url{http://www.kurims.kyoto-u.ac.jp/~tomotada/proj01/}}
If the reader knows a (partial) solution of any problem in this list,
please let him\footnote{\footnotesize
Email address of the editor is: tomotada@kurims.kyoto-u.ac.jp}
know it.

\vskip 0.5pc

\hspace{24pc}
February, 2003
\newline
\hspace{24pc}
T. Ohtsuki

}

\newpage

{\small\baselineskip 10pt\itemsep 3pt\parskip 4.5pt

\tableofcontents

}

\newpage
\makeatletter
\def\@oddhead{{\small\lhead\ifnum\count0=\startpage ISSN 1464-8997 (on line)
1464-8989 (printed) \hfill {\lnum\number\count0}\else\ifodd\count0
\def\\{ }\ifx\theshorttitle\relax \thetitle \else\theshorttitle\fi\hfill
{\lnum\number\count0}\else\def\\{ and }{\lnum\number\count0}
\hfill\ifx\theshortauthors\relax 
\theauthors\else\theshortauthors\fi\fi\fi}}\def\@evenhead{@oddhead}
\def\@oddfoot{\small\lfoot\ifnum\count0=\startpage Published \publishdate:\qua\copyright\ \gtp\hfill\else
\gtm, Volume \thevolumenumber\ (\thevolumeyear)\hfill\fi}
\def\@evenfoot{@oddfoot}
\makeatother
\count0=377

\title{Problems on invariants of knots and 3-manifolds}
\authors{{\rm Edited by} T. Ohtsuki}      
\shortauthors{T. Ohtsuki (Editor)}      
\asciiauthors{T. Ohtsuki (Editor)}      
\coverauthors{T. Ohtsuki (Editor)}      
\email{tomotada@kurims.kyoto-u.ac.jp}
\urladdr{http://www.kurims.kyoto-u.ac.jp/~tomotada/proj01/}
\asciiurl{http://www.kurims.kyoto-u.ac.jp/ tomotada/proj01/}

\begin{abstract}   
This is a list of open problems on invariants of knots and 3-manifolds
with expositions of their history, background, significance, or importance.
This list was made by editing open problems
given in problem sessions in the workshop and seminars 
on ``Invariants of Knots and 3-Manifolds''
held at Kyoto in 2001.
\end{abstract}
\asciiabstract{%
This is a list of open problems on invariants of knots and 3-manifolds
with expositions of their history, background, significance, or
importance.  This list was made by editing open problems given in
problem sessions in the workshop and seminars on `Invariants of Knots
and 3-Manifolds' held at Kyoto in 2001.}

\primaryclass{20F36, 57M25, 57M27, 57R56}
\secondaryclass{13B25, 17B10, 17B37, 18D10, 20C08, 20G42, 22E99, 41A60, 46L37, 57M05, 57M50, 57N10, 57Q10, 81T18, 81T45}

\maketitle  

\setcounter{section}{-1}

\section{Introduction}
\label{sec.intro}

\renewcommand{\thefootnote}{\fnsymbol{footnote}}
\footnotetext[0]{Chapter \ref{sec.intro} was written by J. Roberts.}
\renewcommand{\thefootnote}{\arabic{footnote}}

The study of quantum invariants of links and three-manifolds has a
strange status within topology. When it was born, with Jones' 1984
discovery of his famous polynomial \cite{Jon85}, it seemed
that the novelty and power of the new invariant would be a wonderful
tool with which to resolve some outstanding questions of
three-dimensional topology. Over the last 16 years, such hopes have
been largely unfulfilled, the only obvious exception being the
solution of the Tait conjectures about alternating knots (see for
example \cite{MeTh93}).

This is a disappointment, and particularly so if one expects the role
of the quantum invariants in mathematics to be the same as that of the
classical invariants of three-dimensional topology. Such a comparison
misses the point that most of the classical invariants were {\em
created} specifically in order to distinguish between things; their
definitions are mainly intrinsic, and it is therefore clear what kind
of topological properties they reflect, and how to attempt to use them
to solve topological problems.

Quantum invariants,\index{quantum invariant} 
on the other hand, should be thought of as having
been {\em discovered}. Their construction is usually indirect (think
of the Jones polynomial,\index{Jones polynomial} 
defined with reference to {\em diagrams} of a
knot) and their existence seems to depend on very special kinds of
algebraic structures (for example, $R$-matrices), whose behaviour is
closely related to three-dimensional combinatorial topology (for
example, Reidemeister moves). Unfortunately such constructions give
little insight into what kind of topological information the
invariants carry, and therefore into what kind of applications they
might have. 

Consequently, most of the development of the subject has taken place
in directions away from classical algebraic and geometric
topology. From the earliest days of the subject, a wealth of
connections to different parts of mathematics has been evident:
originally in links to operator algebras, statistical mechanics, graph
theory and combinatorics, and latterly through physics (quantum field
theory and perturbation theory) and algebra (deformation theory,
quantum group\index{quantum group}
 representation theory).  It is the investigation of
these outward connections which seems to have been most profitable,
for the two main frameworks of the modern theory, that of Topological
Quantum Field Theory and Vassiliev theory (perturbation theory) have
arisen from these.

The TQFT\index{TQFT} 
viewpoint \cite{Atiyah_t} gives a good interpretation of the
cutting and pasting properties of quantum invariants, and viewed as a
kind of ``higher dimensional representation theory'' ties in very well
with algebraic approaches to deformations of representation
categories. It ties in well with geometric quantization theory and
representations of loop groups \cite{Atiyah_g}.
In its physical formulation via 
the Chern-Simons path-integral\index{Chern-Simons!--- path integral}\index{path integral|see {Chern-Simons}}
(see Witten \cite{Witten}),
it even offers a {\em conceptual} explanation of the
invariants' existence and properties, but because this is not
rigorous, it can only be taken as a heuristic guide to the properties
of the invariants and the connections between the various approaches
to them.

The Vassiliev theory\index{Vassiliev invariant} 
(see \cite{BarNatan_V_Top, Kontsevich_F, Thurston_w})
gives geometric
definitions of the invariants in terms of 
integrals over configuration spaces,\index{configuration space!--- integral} 
and also can be viewed as a classification theory, in the
sense that there is a universal invariant, 
the Kontsevich integral\index{Kontsevich invariant}\index{Kontsevich integral|see {Kontsevich invariant}}
 (or
more generally the Le-Murakami-Ohtsuki invariant \cite{LMO}),\index{LMO invariant} 
through which all the other invariants factor. Its drawback is that
the integrals are very hard to work with -- eight years passed between
the definition and calculation \cite{Thurston_w} of the Kontsevich integral
of the {\em unknot}!

These two frameworks have revealed many amazing properties and
algebraic structures of quantum invariants, which show that they are
important and interesting pieces of mathematics in their own right,
whether or not they have applications in three-dimensional
topology. The structures revealed are precisely those which can, and
therefore must, be studied with the aid of three-dimensional pictures
and a topological viewpoint; the whole theory should therefore be
considered as a new kind of algebraic topology specific to three
dimensions.

Perhaps the most important overall goal is simply to {\em really
understand} the topology underlying quantum invariants in three
dimensions: to relate the ``new algebraic topology'' to more classical
notions and obtain good intrinsic topological definitions of the
invariants, with a view to applications in three-dimensional topology
and beyond.

The problem list which follows contains detailed problems in all areas
of the theory, and their division into sections is really only for
convenience, as there are very many interrelationships between
them. Some address unresolved matters or extensions arising from
existing work; some introduce specific new conjectures; some describe
evidence which hints at the existence of new patterns or structures;
some are surveys on major and long-standing questions in the field;
some are purely speculative.

Compiling a problem list is a very good way to stimulate research
inside a subject, but it also provides a great opportunity to ``take
stock'' of the overall state and direction of a subject, and to try to
demonstrate its vitality and worth to those outside the area.  We hope
that this list will do both.

\newpage

\section{Polynomial invariants of knots}

\subsection{The Jones polynomial}

\begin{figure}[ht!]
$$
\begin{array}{ccccc}
\pict{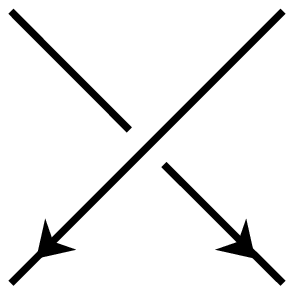}{1.4cm} &\qquad&
\pict{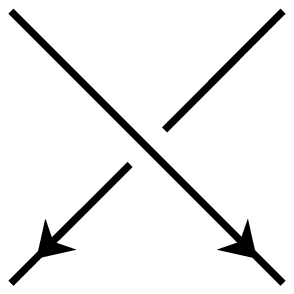}{1.4cm} &\qquad&
\pict{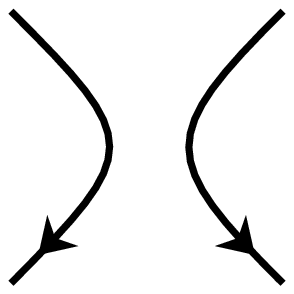}{1.4cm} \\
L_+ && L_- && L_0
\end{array}
$$
\caption{\label{fig.L+L-L0}
Three links $L_+, L_-, L_0$}
\end{figure}

The {\it Kauffman bracket}\index{Kauffman bracket} 
of unoriented link diagrams
is defined by the following recursive relations,
\begin{align*}
& \Big\langle \!\! \pict{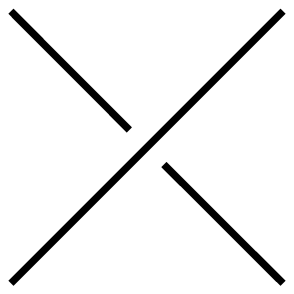}{1.2cm} \!\! \Big\rangle
= A \Big\langle \!\! \pict{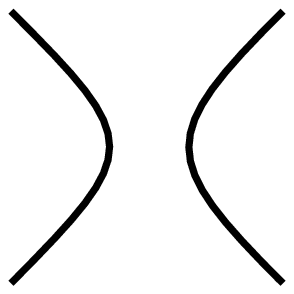}{1.2cm} \!\! \Big\rangle
+ A^{-1} \Big\langle \!\! \pict{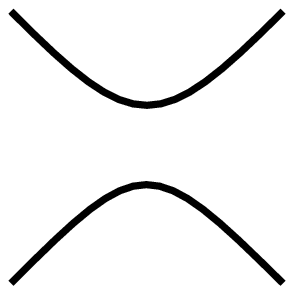}{1.2cm} \!\! \Big\rangle, \\
& \big\langle \pict{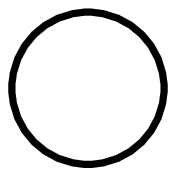}{1cm} D \big\rangle
= (-A^2 - A^{-2}) \langle D \rangle
\qquad \mbox{for any diagram $D$,} \\
& \langle \mbox{the empty diagram $\emptyset$} \rangle = 1,
\end{align*}
where three pictures in the first formula
imply three links diagrams,
which are identical except for a ball,
where they differ as shown in the pictures.
The {\it Jones polynomial}\index{Jones polynomial}
$V_L(t)$
of an oriented link $L$ is defined by
$$
V_L(t) = (-A^2-A^{-2})^{-1} (-A^3)^{-w(D)} \langle D \rangle 
\Big|_{A^2=t^{-1/2}} \in \Z[t^{1/2},t^{-1/2}],
$$
where $D$ is a diagram of $L$,
$w(D)$ is the writhe of $D$,
and $\langle D \rangle$ is the Kauffman bracket of $D$
with its orientation forgotten.
The Jones polynomial is an isotopy invariant of oriented links
uniquely characterized by
\begin{align}
\label{eq.Jones_skein}
& t^{-1} V_{L_+}(t) -t V_{L_-}(t) = (t^{1/2}-t^{-1/2}) V_{L_0}(t), \\
& V_O(t) = 1, \notag
\end{align}
where $O$ denotes the trivial knot, and
$L_+$, $L_-$, and $L_0$ are three oriented links,
which are identical except for a ball,
where they differ as shown in Figure \ref{fig.L+L-L0}.
It is shown, by (\ref{eq.Jones_skein}), that
for any knot $K$, its Jones polynomial $V_K(t)$ belongs to $\Z[t,t^{-1}]$.

\subsubsection{Does the Jones polynomial distinguish the trivial knot?}

\begin{prob}[{\cite[Problem 1]{Jones_prob}}] 
\label{prob.VK=1}
Find a non-trivial knot $K$ with 
$V_K(t)=1$.\index{Jones polynomial!trivial ---}
\end{prob}

\begin{rem}
It is shown by computer experiments that
there are no non-trivial knots with $V_K(t)=1$
up to 17 crossings of their diagrams \cite{DH97},
and up to 18 crossings \cite{Yamada_uJ}.
See \cite{Bigelow_J} (and \cite{Bigelow_website})
for an approach to find such knots
by using representations of braid groups.
\end{rem}

\begin{rem}
Two knots with the same Jones polynomial
can be obtained by mutation.
A {\it mutation}\index{mutation} 
is a relation of two knots,
which are identical except for a ball,
where they differ by $\pi$ rotation of a 2-strand tangle
in one of the following ways (see \cite{APR} for mutations).
$$
\pict{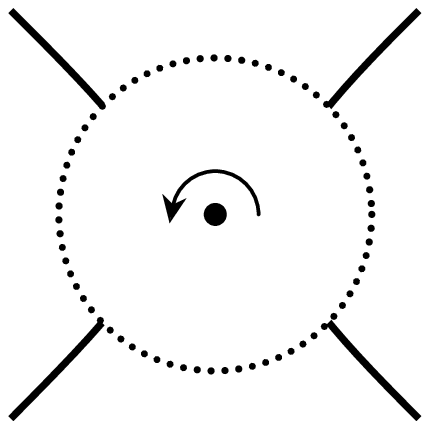}{1.7cm} \qquad
\pict{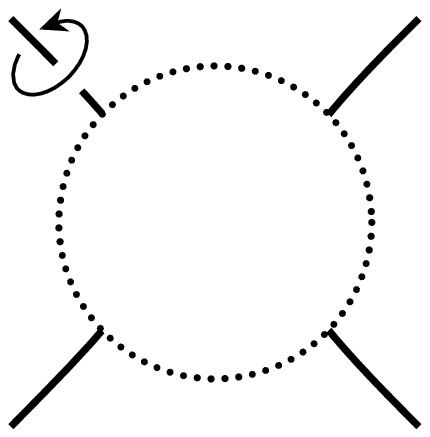}{1.7cm} \qquad
\pict{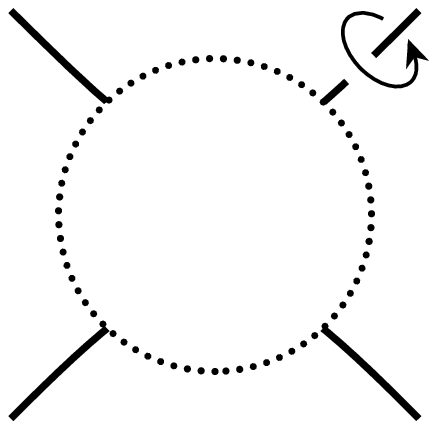}{1.7cm}
$$
For example, the Conway knot and the Kinoshita-Terasaka knot
are related by a mutation.
$$
\pict{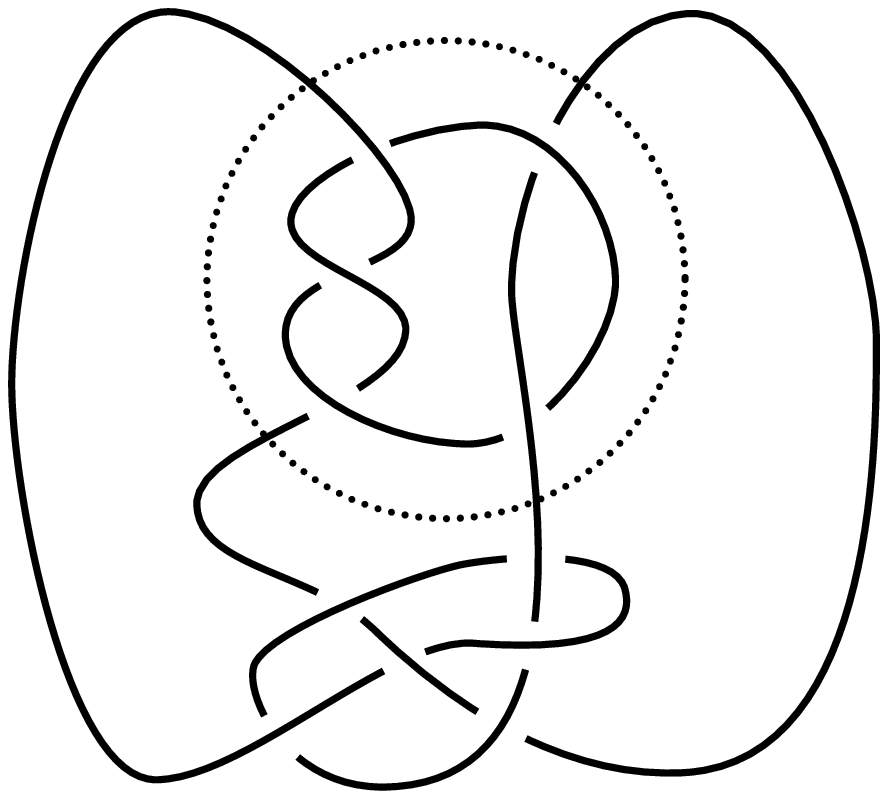}{3cm} \qquad\qquad
\pict{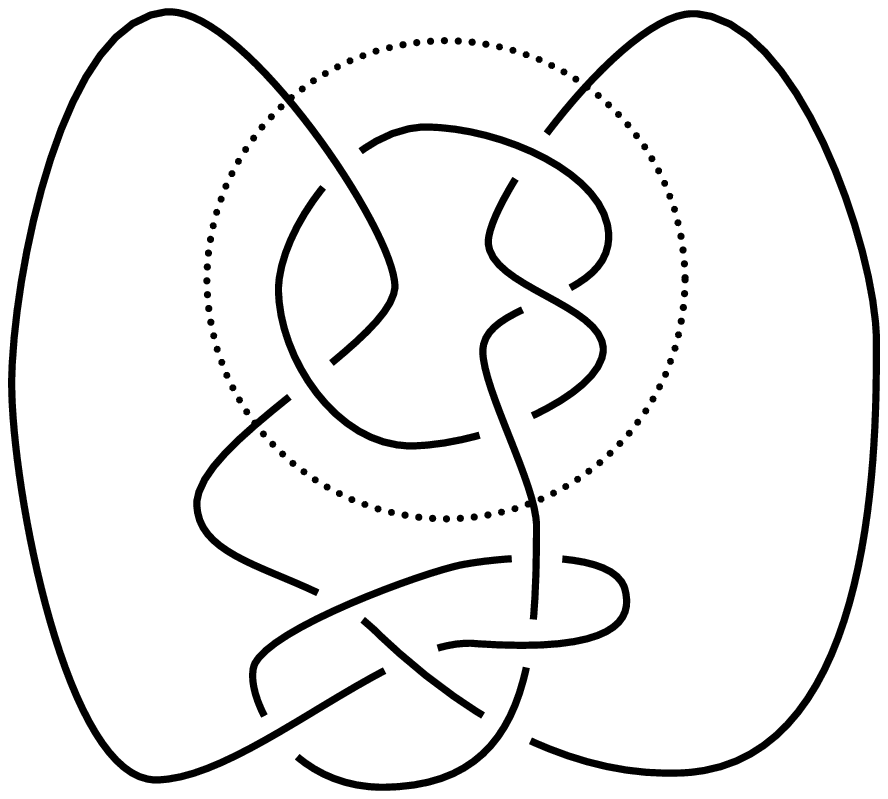}{3cm}
$$
They have the same Jones polynomial,
because their diagrams have the same writhe
and the Kauffman bracket of the tangle shown in the dotted circle
can be presented by
$$
\Big\langle \!\! \pict{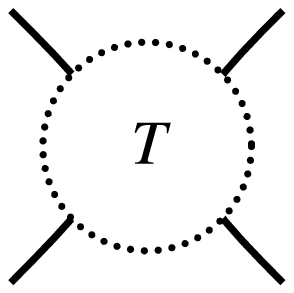}{1.2cm} \!\! \Big\rangle
= x \Big\langle \!\! \pict{c-zero}{1.2cm} \!\! \Big\rangle
+ y \Big\langle \!\! \pict{c-infty}{1.2cm} \!\! \Big\rangle
= \Big\langle \!\! \pict{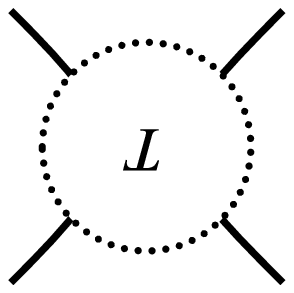}{1.2cm} \!\! \Big\rangle,
$$
with some scalars $x$ and $y$.
\end{rem}

\begin{rem}
The Jones polynomial can be obtained from the Kontsevich invariant\index{Kontsevich invariant!--- and Jones polynomial}
through the weight system $W_{sl_2,V}$
for the vector representation $V$ of $sl_2$
(see, e.g.\ \cite{Ohtsuki_book}).
Problem \ref{prob.VK=1} might be related to the kernel of $W_{sl_2,V}$.
\end{rem}

\begin{rem}
Some links with the Jones polynomial
equal to that of the corresponding trivial links 
are given in \cite{EKT01}.
For example, the Jones polynomial of the following link 
is equal to the Jones polynomial of the trivial 4-component link.
$$
\pict{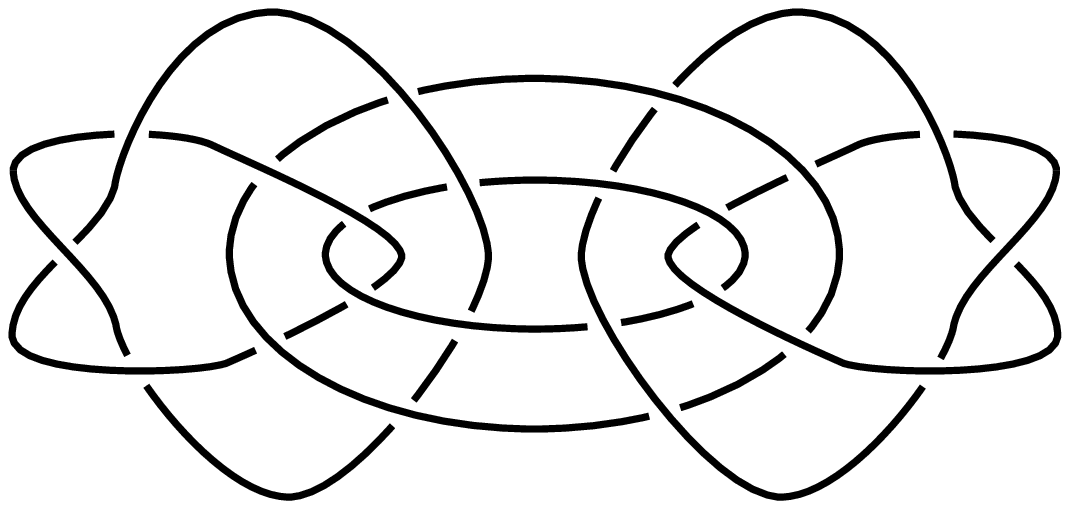}{2.1cm}
$$
\end{rem}

\begin{rem}[\rm (X.-S. Lin \cite{Lin_website})] 
Use Kontsevich integral\index{Kontsevich invariant!--- and Jones polynomial} 
to show the existence of a non-trivial knot with
trivial Alexander-Conway polynomial. This might give us some hints to 
Problem \ref{prob.VK=1}.
\end{rem}

\subsubsection{Characterization and interpretation of the Jones polynomial}

\begin{prob}[{\cite[Problem 2]{Jones_prob}}] 
\label{prob.char_VK}
Characterize those elements of $\Z[t,t^{-1}]$ of 
the form $V_K(t)$.\index{Jones polynomial!image of ---}
\end{prob}

\begin{rem}[\rm \cite{Jones_prob}] 
The corresponding problem for the Alexander polynomial has been solved;
it is known that
a polynomial $f(t) \in \Z[t,t^{-1}]$ is equal to
the Alexander polynomial of some knot $K$
if and only if $f(1)=1$ and $f(t) = f(t^{-1})$.
The formulas $V_K(1)=1$ and
$V_K( \exp \frac{2\pi\sqrt{-1}}{3} ) = 1$
are obtained by the skein relation (\ref{eq.Jones_skein}).
These formulas give weak characterizations of the required elements.
\end{rem}

\begin{rem}[\rm (X.-S. Lin \cite{Lin_website})] 
The Mahler measure (see \cite{Everest} for its exposition) of 
a polynomial $F(x) = a \prod_i (x-\alpha_i) \in \C[x]$
is defined by 
$$
m(F) = \log |a| + \sum_i \log \max \{ 1, |\alpha_i| \}
= \int^1_0 \log | F( e^{2\pi \sqrt{-1} \theta} ) | d \theta.
$$
The Mahler measure can be defined also for a Laurent polynomial similarly. 
Is it true that $m(V_K)>0$ for the Jones polynomial $V_K$ of a knot $K$,
if $K$ is a non-trivial knot?
\end{rem}

\begin{prob}
Find a 3-dimensional topological interpretation 
of the Jones polynomial of links.\index{Jones polynomial!interpretation of ---}
\end{prob}

\begin{rem}
The Alexander polynomial has a topological interpretation
such as the characteristic polynomial
of $H_1 (\widetilde{S^3-K}; \Q)$
of the infinite cyclic cover of the knot complement $S^3-K$,
where $H_1 (\widetilde{S^3-K}; \Q)$ is regarded as
a $\Q[t,t^{-1}]$-module by regarding $t$ as
the action of the deck transformation on $\widetilde{S^3-K}$.
\end{rem}

\begin{rem}
In the viewpoint of mathematical physics,
Witten \cite{Witten} gave a 3-dimensional interpretation
of the Jones polynomial of a link
by a path integral
including a holonomy along the link
in the Chern-Simons field theory.
\end{rem}

\begin{rem}
Certain special values of the Jones polynomial have some interpretations.
The formulas $V_L(1) = (-2)^{\# L -1}$ and 
$V_L( \exp \frac{2\pi\sqrt{-1}}3 ) = 1$
are shown by the skein relation (\ref{eq.Jones_skein}),
where $\# L$ denotes the number of components of $L$.
It is known that
$| V_L(-1) |$ is equal to
the order of $H_1(M_{2,L})$ if its order is finite, and 0 otherwise.
Here, $M_{2,L}$ denotes the double branched cover of $S^3$ branched along $L$.
It is shown, in \cite{Murakami_Arf}, that
$V_L(\sqrt{-1}) = (- \sqrt{2})^{\# L -1} 
(-1)^{{\rm Arf}(L)}$
if $\mbox{Arf}(L)$ exists, and $0$ otherwise.
It is shown, in \cite{LicMil}, that
$V_L( \exp \frac{\sqrt{-1}\pi}3 )
= \pm \sqrt{-1}^{\# L -1} 
\sqrt{-3}^{{\rm dim} H_1(M_{2,L};\Z/3 \Z)}$.
If $\omega$ is equal to a 2nd, 3rd, 4th, 6th root of unity,
the computation of $V_L(\omega)$ can be done in polynomial time
of the number of crossings of diagrams of $L$
by the above interpretation of $V_L(\omega)$.
Otherwise, $V_L(\omega)$ does not have such a topological interpretation,
in the sense that
computing $V_L(\omega)$ of an alternating link $L$
at a given value $\omega$ is $\#{\rm P}$-hard
except for the above mentioned roots of unity
(see \cite{JVW90,Welsh}).
\end{rem}

\begin{prob}[J. Roberts] 
\label{prob.roberts1}
Why is the Jones polynomial\index{Jones polynomial!interpretation of ---} 
a polynomial? 
\end{prob}

\begin{rem}[{\rm (J. Roberts)}] 
A topological invariant of knots should ideally be defined in an
intrinsically $3$-dimensional fashion, so that its invariance under
orientation-preserving diffeomorphisms of $S^3$ is built-in.
Unfortunately, almost all of the known constructions of the Jones
polynomial (via $R$-matrices, skein relations, braid groups or the
Kontsevich integral,\index{Kontsevich invariant!--- and Jones polynomial} 
for example) break the symmetry, requiring the
introduction of an axis (Morsification of the knot) or plane of
projection (diagram of the knot). I believe that the ``perturbative''
construction via configuration space integrals \cite{Thurston}, whose
output is believed to be essentially equivalent to the Kontsevich
integral, is the only known intrinsic construction.

In the definitions with broken symmetry, it is generally easy to see
that the result is an integral Laurent polynomial in $q$ or
$q^{\frac12}$. In the perturbative approach, however, we obtain a
formal power series in $\hbar$, and although we know that it ought to
be the expansion of an integral Laurent polynomial under the
substitution $q=e^\hbar$, it seems hard to prove this directly. A
related observation is that the analogues of the Jones polynomial for
knots in $3$-manifolds other than $S^3$ are {\em not} polynomials, but
merely functions from the roots of unity to algebraic integers. What
is the special property of $S^3$ (or perhaps $\R^3$) which causes this
behaviour, and why does the variable $q$ seem natural only when one
breaks the symmetry?

The typical raison d'etre of a Laurent polynomial is that it is a
character of the circle. (In highbrow terms this is an example of
``categorification'',\index{categorification}
but it is also belongs to a concrete tradition
in combinatorics that to prove that something is a non-negative
integer one should show that it is the dimension of a vector space.)
The idea that the Jones polynomial is related to $K$-theory
\cite{Wil01b} and that it ought to be the $S^1$-equivariant index of
some elliptic operator defined using the special geometry of $\R^3$ or
$S^3$ is something Simon Willerton and I have been pondering for some
time. As for the meaning of $q$, Atiyah suggested the example in
equivariant $K$-theory
\[ K_{SO(3)}(S^2) \cong K_{S^1}(pt) = \Z[q^{\pm 1}],\] in to make the
first identification {\em requires} a choice of axis in $\R^3$. (This
would suggest looking for an $SO(3)$-equivariant $S^2$-family of
operators.)  
\end{rem}


\begin{prob}[J. Roberts]  
\label{prob.roberts2}
Is there a relationship between values of 
Jones polynomials\index{Jones polynomial!special value of ---} 
at roots of unity and branched cyclic coverings of a knot? 
\end{prob}

\begin{prob}[J. Roberts]  
\label{prob.roberts3}
Is there a relationship between the Jones polynomial\index{Jones polynomial!interpretation of ---} 
of a knot and
the counting of points in varieties defined over\index{finite field} 
finite fields? 
\end{prob}

\begin{rem}[{\rm (J. Roberts)}]  
These two problems prolong the ``riff in the key of $q$'': the
amusing fact that traditional, apparently independent uses of that
letter, denoting the number of elements in a finite field, the
deformation parameter $q=e^\hbar$, the variable in the Poincar\'e
series of a space, the variable in the theory of modular forms,
etc. turn out to be related.

The first problem addresses a relationship which holds for the
Alexander polynomial.
For example, the order of the torsion in $H_1$
of the $n$-fold branched cyclic cover equals the product of the values
of the Alexander polynomial at all the $n$th roots of unity. It's hard
not to feel that the variable $q$ has some kind of meaning as a deck
translation, and that the values of the Jones polynomial at roots of
unity should have special meanings. 

The second has its roots in Jones' original formulation of his
polynomial using Hecke algebras. 
The Hecke algebra\index{Hecke algebra} 
$H_n(q)$ is just
the Hall algebra of double cosets of the Borel subgroup inside
$SL(n,\FF_q)$; the famous quadratic relation $\sigma^2 = (q-1) \sigma +
q$ falls naturally out of this. Although the alternative definition of
$H_n(q)$ using generators and relations extends to allow $q$ to be any
complex number (and it is then the roots of unity, at which $H_n(q)$
is not semisimple, which are the obvious special values), it might be
worth considering whether Jones polynomials at prime powers $q=p^s$
have any special properties.

Ideally one could try to find a topological definition of the Jones
polynomial (perhaps only at such values) which involves finite
fields. The coloured Jones polynomials of the unknot are quantum
integers, which count the numbers of points in projective spaces
defined over finite fields; might those for arbitrary knots in $S^3$
count points in other varieties? Instead of counting counting points,
one could consider Poincar\'e polynomials, as the two things are
closely related by the Weil conjectures.

One obvious construction involving finite fields is to count
representations of a fundamental group into a finite group of Lie
type, such as $SL(n, \FF_q)$. Very much in this vein, Jeffrey Sink
\cite{Sin00} associated to a knot a zeta-function formed from the
counts of representations into $SL(2, \FF_{p^s})$, for fixed $p$ and
varying $s$. His hope, motivated by the Weil conjectures, was the idea
that the $SU(2)$ Casson invariant\index{Casson invariant} 
might be related to such counting.
For such an idea to work, it is probably necessary to find some way of
counting representations with signs, or at least to enhance the
counting in some way. Perhaps the kind of twisting used in the
Dijkgraaf-Witten theory \cite{DiWi90} could be used. 
\end{rem}


\begin{prob}[J. Roberts]  
\label{prob.roberts4}
Define the Jones polynomial\index{Jones polynomial!interpretation of ---}
intrinsically using homology of local systems. 
\end{prob}

\begin{rem}[{\rm (J. Roberts)}]  
The Alexander polynomial of a knot can be defined using the
twisted homology of the complement. In the case of the Jones
polynomial, no similar direct construction is known, but the approach
of Bigelow \cite{Big01b} is tantalising. He shows how to construct a
representation of the braid group $B_{2n}$ on the twisted homology of
the configuration space of $n$ points in the $2n$-punctured disc, and
how to use a certain ``matrix element'' of this representation to
obtain the Jones polynomial of a knot presented as a plait. Is there
any way to write the same calculation directly in terms of
configuration spaces of $n$ points in the knot complement, for
example?
\end{rem}


\begin{prob}[J. Roberts]  
\label{prob.roberts5}
Study the relation between the Jones polynomial\index{Jones polynomial!--- and Gromov-Witten theory}\index{Gromov-Witten theory}
and Gromov-Witten theory. 
\end{prob}

\begin{rem}[{\rm (J. Roberts)}]  
The theory of pseudo-holomorphic curves or ``Gromov-Witten
invariants'' has been growing steadily since around 1985, in parallel
with the theory of quantum invariants\index{quantum invariant!--- and Gromov-Witten theory} 
in three dimensional topology.
During that time it has come to absorb large parts of modern geometry
and topology, including symplectic topology, Donaldson/Seiberg-Witten
theory, Floer homology, enumerative algebraic geometry, etc. It is
remarkable that three-dimensional TQFT has remained isolated from it
for so long, but finally there is a connection, as explained in the
paper by Vafa and Gopakumar \cite{GoVa00} (though prefigured by Witten
\cite{Wit95}), and now under investigation by many geometers. The
basic idea is that the HOMFLY polynomial can be reformulated as a
generating function counting pseudo-holomorphic curves in a certain
Calabi-Yau manifold, with boundary condition a Lagrangian submanifold
associated to the knot. (This is the one place where the HOMFLY and
not the Jones polynomial is essential!) The importance of this
connection can hardly be overestimated, as it should allow the
exchange of powerful techniques between the two subjects.
\end{rem}

\begin{figure}[ht!]
\begin{align*}
& \pict{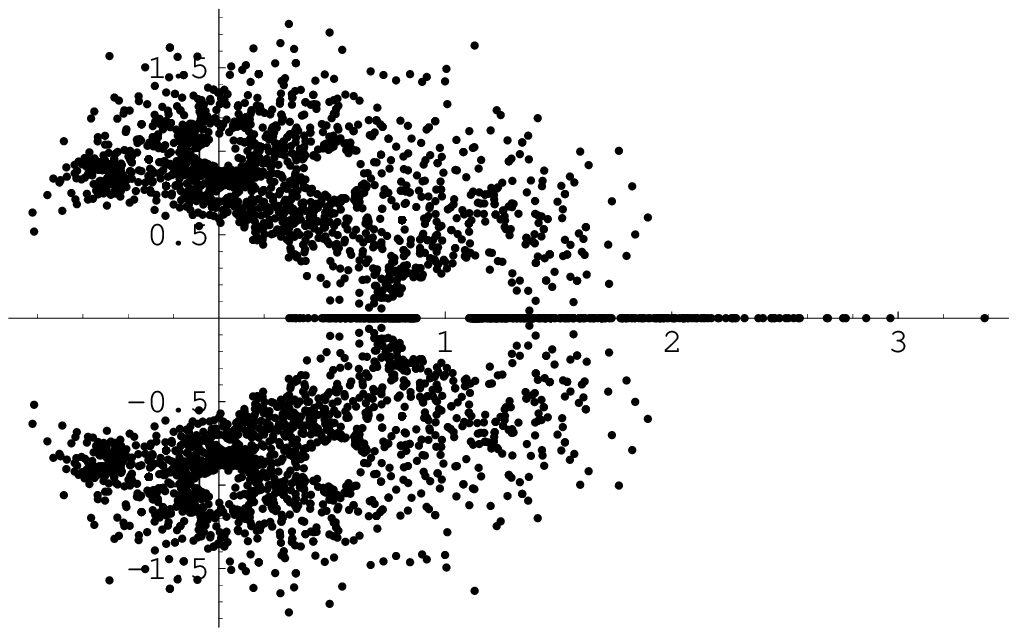}{3cm} \quad
\pict{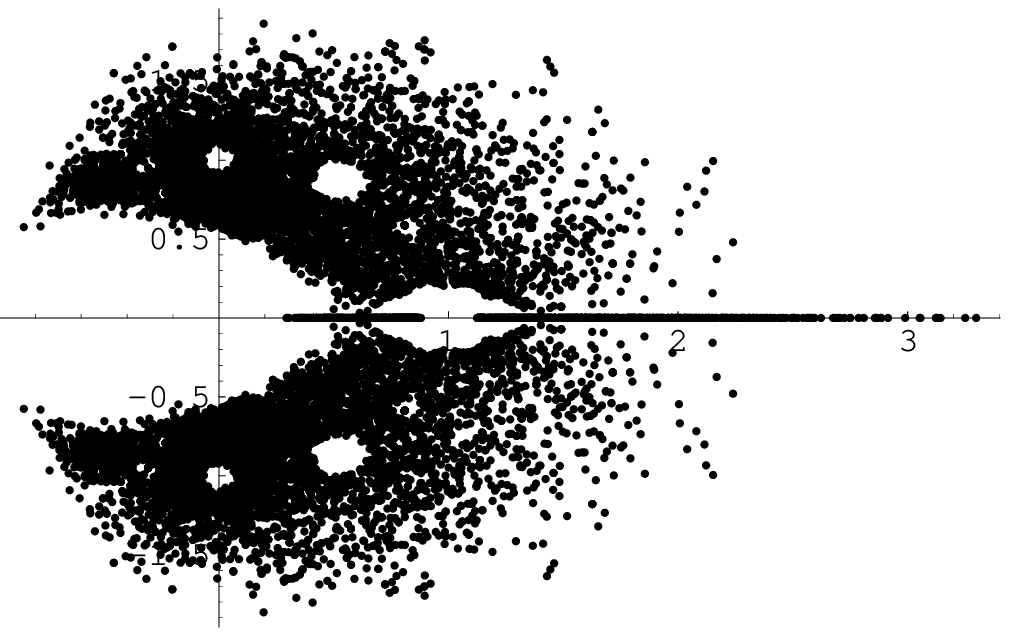}{3cm} \\
& \qquad\qquad \picthw{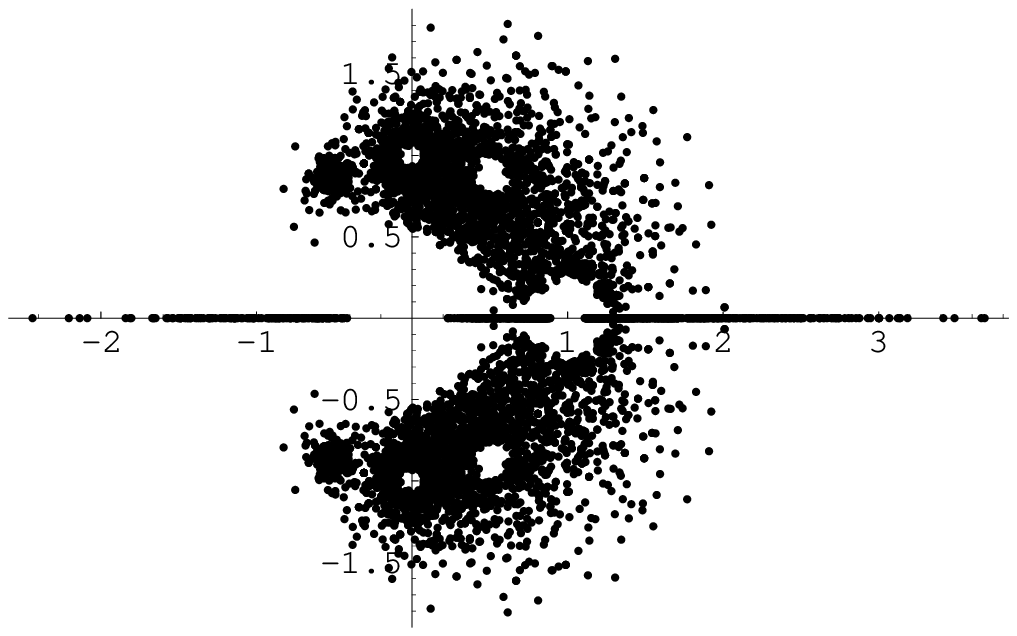}{3cm}{7cm}
\end{align*}
\caption{\label{fig.0V_11_12}
The upper pictures show the distribution of zeros of the Jones polynomial 
for alternating knots of 11 and 12 crossings \cite{Lin_website}.
The lower picture shows the distribution of zeros of the Jones polynomial 
for 12 crossing non-alternating knots \cite{Lin_website}.
See \cite{Lin_website} for further pictures
for alternating knots with 10 and 13 crossings.}
\end{figure}

\begin{figure}[ht!]
\begin{align*}
& \qquad\quad\ \picthw{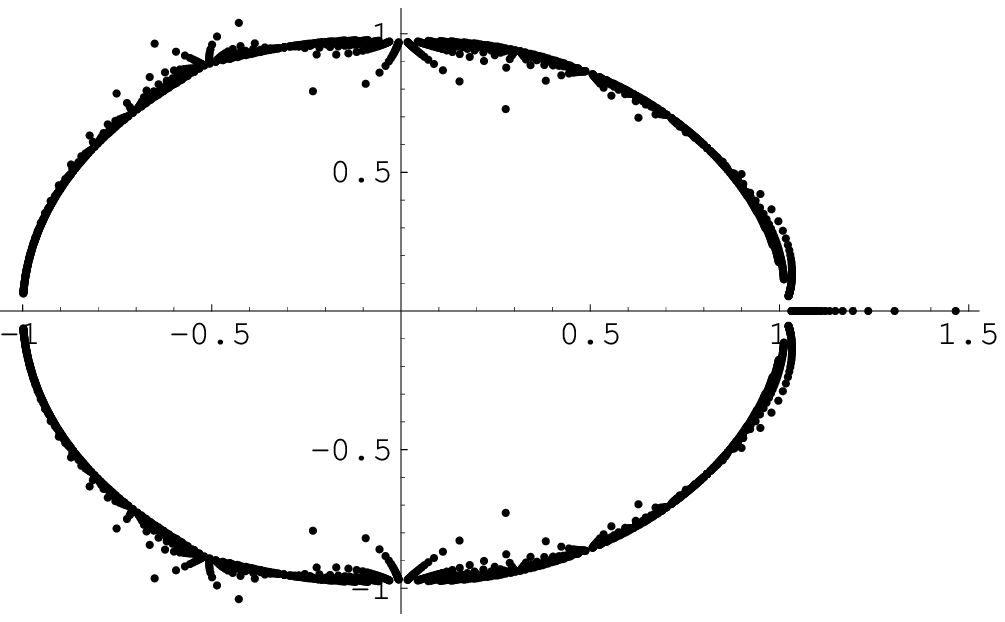}{3.1cm}{3.6cm} \qquad
\picthw{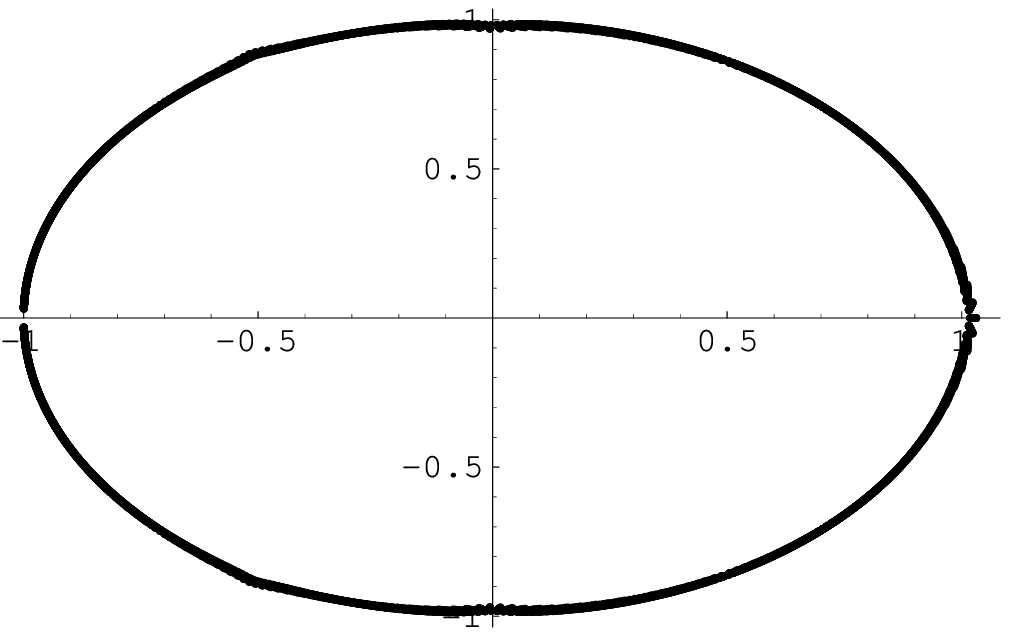}{3cm}{3cm} \\
& \picthw{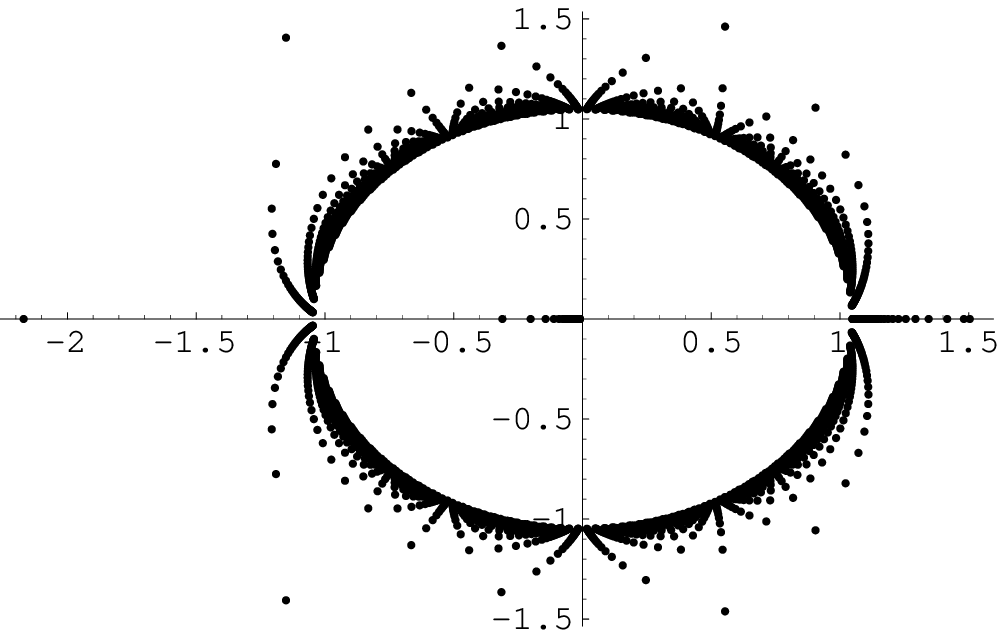}{4cm}{5cm} \qquad
\picthw{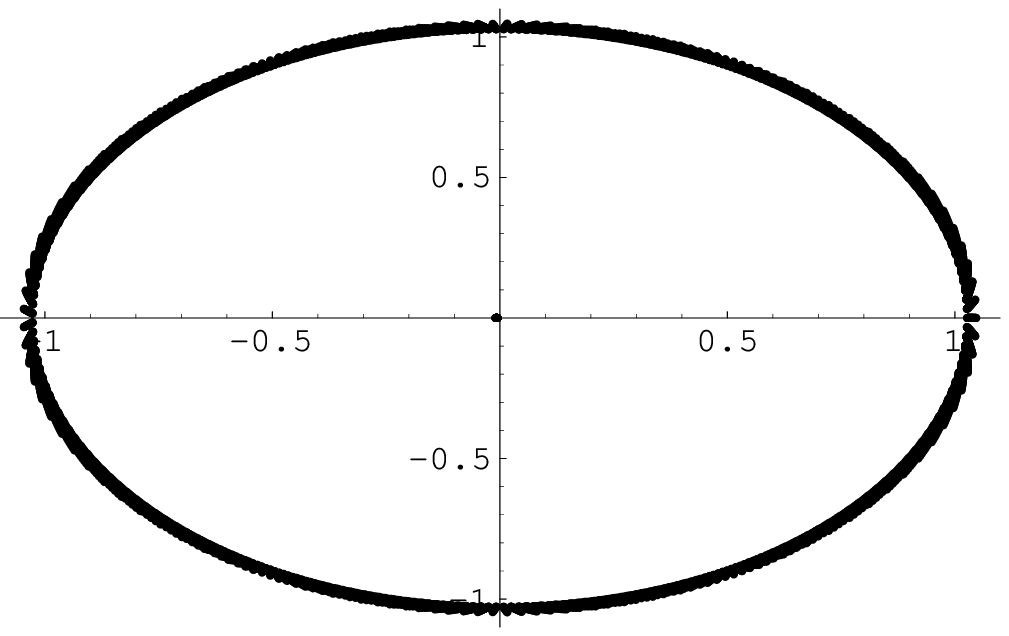}{2.9cm}{2.9cm}
\end{align*}
\caption{\label{fig.0V_tk}
The upper pictures show
the distribution of zeros of the Jones polynomial for 
$n$-twist knots, with $n$ from 1 to 50 and from 51 to 100, respectively
\cite{Lin_website}.
The lower pictures show
the distribution of zeros of the Jones polynomial for 
$(2, 2n-1)$ torus knots, with $n$ from 1 to 50 and 
from 51 to 100, respectively \cite{Lin_website}.
See \cite{Lin_website} for further pictures 
for $(3, 3n+1)$ and $(3,3n+2)$ torus knots.}
\end{figure}

\subsubsection{Numerical experiments}

The following problem might characterize 
the form of the Jones polynomial of knots in some sense.

\begin{prob}[X.-S. Lin] 
\label{prob.desc_0_VK}
Describe the set of zeros of the Jones polynomial\index{Jones polynomial!zeros of ---}
of all (alternating) knots.
\end{prob}

\begin{rem}[{\rm (X.-S. Lin)}] 
The plottings in Figure \ref{fig.0V_11_12} numerically describe
the set of zeros of the Jones polynomial of many knots.
Similar plottings are already published in \cite{WuWang}
for some other infinite families of knots 
for which the Jones polynomial is known explicitly.
See also \cite{ChangShrock} for some other plottings.
\end{rem}

\begin{rem}[{\rm (X.-S. Lin)}] 
It would be a basic problem 
to look into the zero distribution of the family of polynomials 
with bounded degree such that coefficients are all integers 
and coefficients sum up to 1, 
and compare it with the zero distribution of the Jones polynomial 
on the collection of (alternating) knots with bounded crossing number. 
The paper \cite{OdlyzkoPoonen} discusses
the zero distribution of the family of polynomials 
with 0,1 coefficients and bounded degree.
It is particularly interesting to compare the plotting shown in this paper
with the plottings 
in Figures \ref{fig.0V_11_12} and \ref{fig.0V_tk}
for the zeros of the Jones polynomials.
\end{rem}

\begin{figure}[ht!]
$$
\rotatebox{-90}{\pict{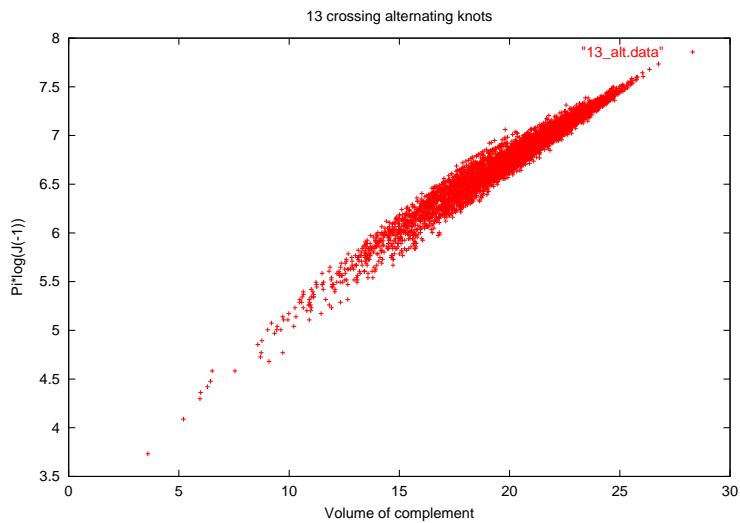}{10cm}}
$$
\caption{\label{fig.vol_Jones1}
The distribution of pairs of 
the hyperbolic volume of knot complements and $\pi \log V_K(-1)$
for alternating knots with 13 crossings \cite{Dunfield_website}.}
\end{figure}

\begin{figure}
\begin{align*}
& \rotatebox{-90}{\pict{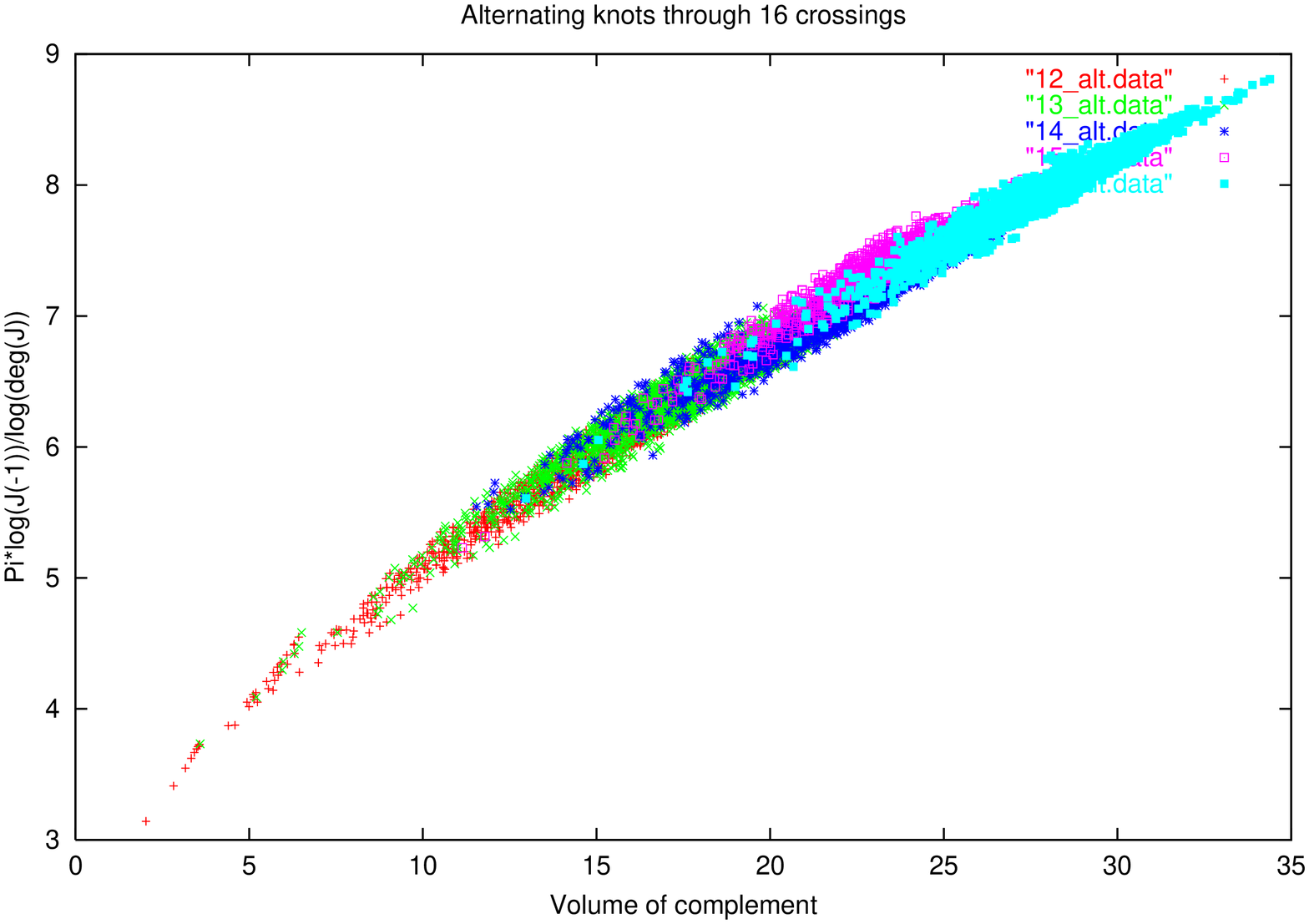}{10cm}} \\
& \rotatebox{-90}{\pict{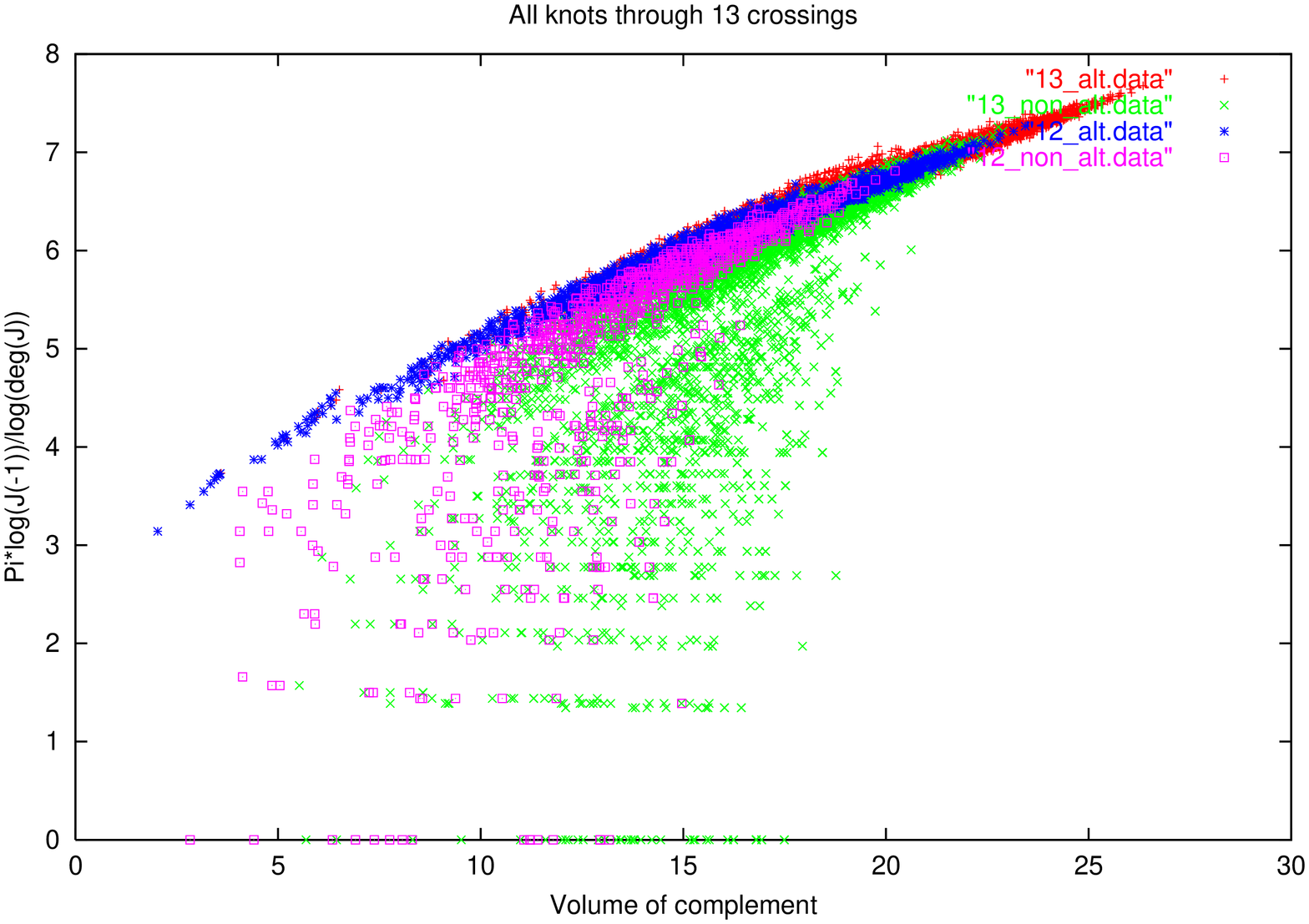}{10cm}}
\end{align*}
\caption{\label{fig.vol_Jones2}
The distributions of pairs of the hyperbolic volume of knot complements 
and $\pi \log V_K(-1) / \log \mbox{\rm deg} V_K(t)$.
The upper picture is for all alternating knots with 12 and 13 crossings 
and samples of alternating knots with 14, 15, and 16 crossings, 
and the lower picture is for all knots with 13 or fewer crossings
\cite{Dunfield_website}.}
\end{figure}

\begin{prob}[N. Dunfield] 
\label{prob.vol_VK-1}
Find the relationship between 
the hyperbolic volume\index{hyperbolic volume!--- and $\log V_K(-1)$} 
of knot complements
and $\log V_K(-1)$
(resp. $\log V_K(-1)/\log \mbox{\rm deg} V_K(t)$).
\end{prob}

\begin{rem}[{\rm (N. Dunfield \cite{Dunfield_website})}] 
    $V_K(-1)$ is just $\Delta_K(-1)$, which is the order of the torsion
  in the homology of the double cover of $S^3$ branched over $K$.
  $\log V_K(-1)$ is one of the first terms of the volume conjecture
  (Conjecture \ref{conj.vol_conj}).  Figure \ref{fig.vol_Jones1}
  suggests that for alternating knots with a fixed number of
  crossings, $\log V_K(-1)$ is almost a linear function of the volume.
 
  Figure \ref{fig.vol_Jones2} suggests that there should be an
  inequality
  $$
  \frac{\log V_K(-1)}{\log \mbox{\rm deg} V_K(t)} < a \cdot \mbox{\rm vol}(S^3-K) + b
  $$
  for some constants $a$ and $b$.  For 2-bridge knots, Agol's work
  on the volumes of 2-bridge knots \cite{Agol_volume} can be used to
  prove such an inequality with $a = b = 2/v_3$ (here, $v_3$ is the
  volume of a regular ideal tetrahedron).
\end{rem}

\subsubsection{Categorification of the Jones polynomial}

Khovanov \cite{Khovanov_J,Khovanov_t}
defined certain homology groups of a knot 
whose Euler characteristic is equal to the Jones polynomial,
which is called the {\it categorification}\index{categorification} 
of the Jones polynomial.
See also \cite{BarNatan_Kh} for an exposition of it.

\begin{prob}
Understand Khovanov's categorification\index{categorification} 
of the Jones polynomial.\index{Jones polynomial!categorification of ---}
\end{prob}

\begin{prob}
Categorify other knot polynomials.
\end{prob}

\begin{rem}[\rm (M. Hutchings)] 
There does exist a categorification of the Alexander polynomial, or more
precisely of $\Delta_K(t)/(1-t)^2$,\index{Alexander polynomial!categorification of ---}
where $\Delta_K(t)$ denotes the (symmetrized) Alexander polynomial 
of the knot $K$.  
It is a kind of Seiberg-Witten Floer homology of
the three-manifold obtained by zero surgery on $K$.  
One can regard it as $\Z \times \Z/2\Z$ graded, 
although in fact the column whose Euler characteristic gives
the coefficient of $t^k$ is relatively $\Z/2k\Z$ graded.
\end{rem}

\subsection{The HOMFLY, Q, and Kauffman polynomials}

The {\it skein polynomial}\index{skein polynomial} 
(or the 
{\it HOMFLY polynomial})\index{HOMFLY polynomial|see {skein polynomial}}
$P_L(l,m) \in \Z[l^{\pm1},m^{\pm1}]$
of an oriented link $L$
is uniquely characterized by
\begin{align*}
& l^{-1} P_{L_+}(l,m) - l P_{L_-}(l,m) = m P_{L_0}(l,m), \\
& P_O(l,m) = 1, \notag
\end{align*}
where $O$ denotes the trivial knot, and
$L_+$, $L_-$, and $L_0$ are three oriented links,
which are identical except for a ball,
where they differ as shown in Figure \ref{fig.L+L-L0}.
For a knot $K$, $P_K(l,m) \in \Z[l^{\pm2},m]$.
The {\it Kauffman polynomial}\index{Kauffman polynomial}
$F_L(a,z) \in \Z[a^{\pm1},z^{\pm1}]$
of an oriented link $L$
is defined by
$F_L(a,z) = a^{-w(D)} [D]$
for an unoriented diagram $D$ presenting $L$ (forgetting its orientation),
where $[D]$ is uniquely characterized by
\begin{align*}
& \left[ \pict{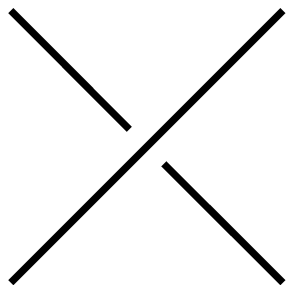}{1cm} \right] 
+ \left[ \pict{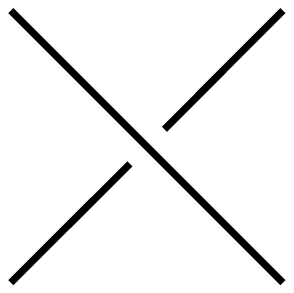}{1cm} \right]
= z \left( \left[ \pict{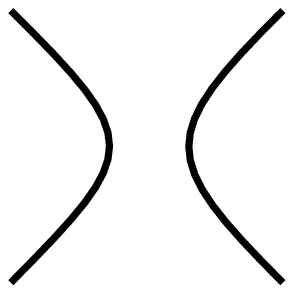}{1cm} \right]
+ \left[ \pict{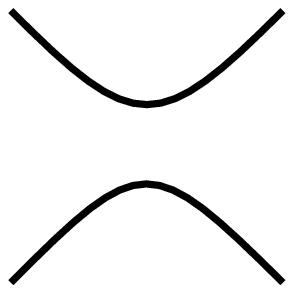}{1cm} \right] \right) \\*
& \left[ \pict{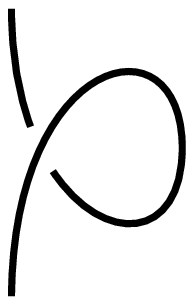}{1cm} \right] 
= a \left[ \ \pict{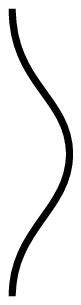}{1cm} \ \right], \\*
& [O] = 1.
\end{align*}
For a knot $K$, $F_K(a,z) \in \Z[a^{\pm1},z]$.
The {\it Q polynomial}\index{Q polynomial} 
$Q_L(x) \in \Z[x^{\pm1}]$ of an unoriented link $L$
is uniquely characterized by
\begin{align*}
& Q \Big( \pict{bigelow/bmw1}{1cm} \Big)
+ Q \Big( \pict{bigelow/bmw2}{1cm} \Big)
= x \left( Q \Big( \pict{bigelow/bmw3}{1cm} \Big)
+ Q \Big( \pict{bigelow/bmw4}{1cm} \Big) \right) \\*
& Q(O) = 1.
\end{align*}
It is known that
\begin{align*}
V_L(t) &= P_L(t,t^{1/2}-t^{-1/2}) = F_L(-t^{-3/4},t^{1/4}+t^{-1/4}), \\*
\Delta_L(t) &= P_L(1,t^{1/2}-t^{-1/2}), \\*
Q_L(z) &= F_L(1,z), 
\end{align*}
where $\Delta_L(t)$ denotes the Alexander polynomial of $L$.\index{Alexander polynomial}
The variable $m$ of $P_L(l,m)$ is called the {\it Alexander variable}.
See, e.g.\ \cite{Kawauchi_book,Lickorish_book},
for details of this paragraph.

Let the {\it span} of a polynomial denote 
the maximal degree minus the minimal degree of the polynomial.

\begin{prob}[A. Stoimenow] 
\label{prob.VK_span}
Does the Jones polynomial\index{Jones polynomial!--- of given span}
 $V$ admit only finitely many values
of given span? 
What about the Q polynomial\index{Q polynomial!--- of given span} 
or the skein,\index{skein polynomial!--- of given span}  
Kauffman polynomials\index{Kauffman polynomial!--- of given span}   
(when fixing the span in both variables)?
\end{prob}

\begin{rem}[\rm (A. Stoimenow)] 
It is true for the skein polynomial when bounding the canonical
genus (for which the Alexander degree of the skein polynomial 
is a lower bound by Morton), in particular it is true for the skein
polynomial of homogeneous links \cite{Cromwell_h}.
It is true for the Jones, Q and 
Kauffman F
polynomial of alternating links (for F more generally for adequate
links). One cannot bound the number of different links, at least
for the skein and Jones polynomial, 
because Kanenobu \cite{Kanenobu_s} gave
infinitely many knots with the same skein polynomial.
\end{rem}

\begin{prob}[A. Stoimenow] 
\label{prob.QKunit}
Why are the unit norm complex numbers $\alpha$ for which the value
$Q_K(\alpha)$\index{Q polynomial!values of ---}
has maximal norm statistically concentrated around
$e^{11\pi \sqrt{-1}/25}$\,?
\end{prob}

\begin{rem}[\rm (A. Stoimenow)] 
The maximal point of
$| Q_K(e^{2 \pi \sqrt{-1} t}) |$ for $t \in [0,1)$
is statistically concentrated around $t = 11/50$.
This was revealed by an experiment in an attempt to estimate the
asymptotical growth of the coefficients of the $Q$ polynomial.
There seems no difference in the behaviour of alternating and non-%
alternating knots.
\end{rem}

\begin{prob}[M. Kidwell, A. Stoimenow] 
\label{prob.PWK=FK}
Let K be a non-trivial knot,
and let $W_K$ be a Whitehead double of K.
Is then\index{skein polynomial!--- of Whitehead double}\index{Kauffman polynomial!--- and Whitehead double}  
$$
\deg_m P_{W_K}(l,m) = 2\deg_z F_K(a,z) +2 \, ?
$$
\end{prob}

\begin{rem}[\rm (A. Stoimenow)] 
It is true for $K$ up to 11 crossings. 
$\deg_m P_{W_K}(l,m)$ is independent on the twist of $W_K$
if it is $>2$ by a simple skein argument.
\end{rem}

\begin{update}
Gruber \cite{Gruber} showed that,
if $K$ is a prime alternating knot and 
$W_K$ is its untwisted Whitehead double, then  
$\deg_m P_{W_K}(l,m) \le 2\deg_z F_K(a,z) +2$.
\end{update}

\begin{prob}[E. Ferrand, A. Stoimenow] 
\label{prob.sign_PL_FL}
Is for any alternating link $L$,\index{signature!--- of link}\index{skein polynomial!--- and signature}\index{Kauffman polynomial!--- and signature}  
\[
\sigma(L)\ge \mbox{min\,deg}_l \big(P_L(l,m) \big)
\ge \mbox{min\,deg}_a \big( F_L(a^{-1},z) \big)\,?
\]
\end{prob}
 
\begin{rem}[\rm (A. Stoimenow)] 
The second inequality is conjectured by Ferrand \cite{sto.Fer} 
(see also comment on Problem \ref{stoi.prm}), and related to estimates
of the Bennequin numbers of Legendrian knots. As for the first
inequality, by Cromwell \cite{Cromwell_h}
we have $\mbox{min\,deg}_l \big( P_L(l,m) \big)\le 1-\chi(L)$ and
classically $\sigma(L) \le 1-\chi(L)$. 
\end{rem}

\begin{prob}[A. Stoimenow] 
\label{prob.Conway_cr}
If $\nabla_{k}$ is the coefficient of $z^k$ in 
the Conway polynomial\index{Conway polynomial!--- and crossing number}
and $c(L)$ is the crossing number\index{crossing number} 
of a link $L$, is then
\[
\big|\nabla_{k}(L)\big|\le\frac{c(L)^k}{2^k\,k!}\,?
\]
\end{prob}

\begin{rem}[\rm (A. Stoimenow)] 
The inequality is non-trivial only for $L$ of $k+1,k-1,\dots$
components. It is also trivial for $k=0$, easy for $k=1$ ($\nabla_1$
is just the linking number of 2 component links) and proved
by Polyak-Viro \cite{PoVi_C} 
for {\em knots} and $k=2$. There are constants $C_k$ with
\[
\big|\nabla_{k}(L)\big|\le C_k\,{c(L)^k}\,,
\]
following from the proof 
(due to \cite{BarNatan_pp,Stanford_cV} for knots, 
due to \cite{Stoimenow_LWconj} for links)
of the Lin-Wang conjecture \cite{LiWa96} for links, but
determining $C_k$ from the proof is difficult. Can the inequality
be proved by Kontsevich-Drinfel'd, say at least for knots, using
the description of the weight systems of $\nabla$ of
Bar-Natan and Garoufalidis \cite{BG_MMR}? 
More specifically, one can ask
whether the $(2,n)$-torus links (with parallel orientation)
attain the maximal values of $\nabla_k$. One can also ask about
the shape of $C_k$ for other families of Vassiliev invariants,
like $\frac{d^k}{dt^k} V_L(t) \big|_{t=1}$.
\end{rem}

\begin{prob}[A. Stoimenow] 
\label{stoi.prm}
Does $\mbox{min\,deg}_a \big( F_L(a^{-1},z) \big)\le 1-\chi(L)$
hold
for any link $L$?\index{Kauffman polynomial!--- and genus}
If $u(K)$ is the unknotting number of a knot $K$,
does\break $\mbox{min\,deg}_a \big( F_K(a^{-1},z) \big)\le 2u(K)$
hold
for any knot $K$?\index{Kauffman polynomial!--- and unknotting number} 
\end{prob}

\begin{rem}[\rm (A. Stoimenow)] 
For the common lower bound of $2u$ and $1-\chi$ for knots, 
$2g_s$, there is a 15 crossing knot $K$ 
with $2g_s(K) < \mbox{min\,deg}_a \big( F_K(a^{-1},z) \big)$.
Morton \cite{Morton} conjectured long ago that 
$1-\chi(L)\ge \mbox{min\,deg}_l \big( P_L(l,m) \big)$.
There are recent counterexamples, but only of 19 to 21 crossings.
Ferrand \cite{sto.Fer} observed that very often 
$\mbox{min\,deg}_l \big( P_K(l,m) \big) \ge 
\mbox{min\,deg}_a \big( F_K(a^{-1},z) \big)$
(he conjectures it in particular always to hold for alternating $K$), 
so replacing\break `$\mbox{min\,deg}_a \big( F(a^{-1},z) \big)$' 
for `$\mbox{min\,deg}_l \big( P_K(l,m) \big)$'
enhances the difficulty of Morton's problem 
(the counterexamples are no longer such).
\end{rem}

\subsection{The volume conjecture}

In \cite{bene.K2} R. Kashaev defined a series of invariants
$\langle L \rangle_N \in \C$
of a link $L$ for $N=2,3,\cdots$ 
by using the quantum dilogarithm.\index{dilogarithm!quantum ---}
In \cite{Kashaev_hv} he observed, by formal calculations, that
$$
2 \pi \cdot \lim_{N \to \infty}
\frac{\log \langle L \rangle_N}{N} = \mbox{\rm vol} (S^3-L)
$$
for $L = K_{4_1}, K_{5_2}, K_{6_1}$,
where $\mbox{\rm vol}(S^3 - L)$ denotes the hyperbolic volume of $S^3 -L$.
Further, he conjectured that
this formula holds for any hyperbolic link $L$.
In 1999, H. Murakami and J. Murakami \cite{vol.MM} proved that
$\langle L \rangle_N = J_N(L)$ for any link $L$,
where $J_N(L)$ denotes the $N$-colored\index{Jones polynomial!colored ---}
Jones polynomial\footnote{\footnotesize 
This is the invariant obtained as the quantum invariant
 of links
associated with the $N$-dimensional irreducible representation
of the quantum group $U_q(sl_2)$.}
of $L$ evaluated at $e^{2\pi\sqrt{-1}/N}$.

\begin{conj}[The volume conjecture, \cite{Kashaev_hv,vol.MM}] 
\label{conj.vol_conj}
For any knot $K$,\index{volume conjecture!--- for knots}
\begin{equation}
\label{eq.vol_conk}
2 \pi \cdot \lim_{N \to \infty} \frac{\log | J_N(K) |}{N} =
v_3  || S^3 - K ||, 
\end{equation}
where $|| \cdot ||$ denotes the 
simplicial volume\index{simplicial volume!--- and colored Jones polynomial}  
and $v_3$ denotes the 
hyperbolic volume\index{hyperbolic volume!--- and colored Jones polynomial}  
of the regular ideal tetrahedron.
\end{conj}

\begin{rem}
For a hyperbolic knot $K$, (\ref{eq.vol_conk}) implies that
$$
2 \pi \cdot \lim_{N \to \infty} \frac{\log | J_N(K) |}{N} 
= \mbox{\rm vol}(S^3 - K).
$$
\end{rem}

\begin{rem}[{\rm \cite{vol.MM}}] 
Both sides of (\ref{eq.vol_conk}) behave well
under the connected sum and the mutation of knots.
Namely, 
\begin{align*}
& || S^3 - (K_1 \# K_2) || = || S^3 - K_1 || + || S^3 - K_2 ||, \\*
& J_N(K_1 \# K_2) = J_N(K_1) J_N(K_2),
\end{align*}
and $J_N(K)$ and $||S^3-K||$ do not change under a mutation
of $K$.
For details see \cite{vol.MM} and references therein.
\end{rem}

\begin{rem}
The statement of the volume conjecture for a link $L$
should probably be the same statement as (\ref{eq.vol_conk})
replacing $K$ with $L$.
It is necessary to assume that $L$ is not a split link,
since $J_N(L) = 0$ for a split link $L$
(then, the left hand side of (\ref{eq.vol_conk}) does not make sense).
\end{rem}

\begin{exm}
It is shown \cite{KaTi00} that
for a torus link $L$
$$
\lim_{N \to \infty} \frac{\log \langle L \rangle_N}{N} = 0,
$$
which implies that (\ref{eq.vol_conk}) is true for torus links.
\end{exm}

\begin{rem}
Conjecture \ref{conj.vol_conj} has been proved
for the figure eight knot $K_{4_1}$
(see \cite{Hitoshi_aY} for an exposition).
However,
we do not have a rigorous proof
of this conjecture
for other hyperbolic knots so far.
We explain its difficulty below,
after a review of a proof for $K_{4_1}$.
\end{rem}

We sketch a proof of 
Conjecture \ref{conj.vol_conj}
for the figure eight knot $K_{4_1}$;
for a detailed proof see \cite{Hitoshi_aY}.
It is known that
\begin{equation}
\label{eq.JNK41}
J_N(K_{4_1}) = \sum_{n=0}^{N-1} (q)_n (q^{-1})_n,
\end{equation}
where we put $q = e^{2\pi\sqrt{-1}/N}$ and
$$
(q)_n = (1-q)(1-q^2)\cdots(1-q^n), \qquad (q)_0 = 1.
$$
As $N$ tends to infinity fixing $n/N$ in finite,
the asymptotic behaviour of the absolute value of $(q)_n$ 
is described by
\begin{align*}
\log |(q)_n| 
& =  \sum_{k=1}^n \log \Big( 2 \sin \frac{\pi k}{N} \Big)
= \frac{N}{\pi} \int^{n\pi/N}_0 \log ( 2 \sin t) d t + O(\log N) \\
& = - \frac{N}{2\pi} 
\mbox{Im} \big( \mbox{Li}_2(e^{2\pi n\sqrt{-1}/N})\big) + O(\log N),
\end{align*}
where $\mbox{Li}_2$ denotes the 
{\it dilogarithm function}\index{dilogarithm!--- function}
defined on $\C - \{ x \in \R \ | \ x > 1 \}$ by
$$
\mbox{Li}_2(z) = \sum_{n=1}^\infty \frac{z^n}{n^2} 
= - \int_0^z \frac{\log(1-s)}{s} d s.
$$
Noting that each summand of (\ref{eq.JNK41}) is real-valued,
we have that
$$
J_N(K_{4_1}) = \sum_{0 \le n < N} \exp \Big( \frac{N}{2\pi} \mbox{Im} \big( 
\mbox{Li}_2(e^{-2\pi n\sqrt{-1}/N})-\mbox{Li}_2(e^{2\pi n\sqrt{-1}/N}) \big)
+ O(\log N) \Big).
$$
The asymptotic behaviour of this sum
can be described by
the maximal point $z_0$ of
$\mbox{Im} \big( \mbox{Li}_2 (1/z) - \mbox{Li}_2(z) \big)$
on the unit circle $\big\{ z \in \C \ \big| \ |z| = 1 \big\}$.
In fact this $z_0$ is a critical point of
$\mbox{Li}_2 (1/z) - \mbox{Li}_2(z)$ in $\C$,
and hence
$\mbox{Im} \big( \mbox{Li}_2 (1/z_0) - \mbox{Li}_2(z_0) \big)$ 
gives the hyperbolic volume of $S^3 - K_{4_1}$.
Therefore, the conjecture holds in this case.

Next, we sketch a formal argument 
toward Conjecture \ref{conj.vol_conj}
for the knot $K_{5_2}$.
Following \cite{Kashaev_hv}, we have that
$$
J_N(K_{5_2})
= \sum_{0 \le m \le n < N} \frac{(q)_n^2}{(q)_m^\ast} q^{-m(n+1)},
$$
where the asterisk implies the complex conjugate.
By applying the formal approximation\footnote{\footnotesize 
It might be difficult to justify this approximation in a usual sense,
since the argument of $(q)_n$,
given by $(q)_n = |(q)_n| \cdot q^{-n(n+1)/2} (-\sqrt{-1})^n$,
changes discretely and quickly near the limit.}
\begin{align}
\label{eq.qn_app}
(q)_n & \underset{?}{\sim} \exp \Big( \frac{N}{2\pi\sqrt{-1}} \big(
\mbox{Li}_2(1) - \mbox{Li}_2(e^{2\pi n\sqrt{-1}/N}) \big) \Big), \\*
(q)_n^\ast & \underset{?}{\sim} \exp \Big( \frac{N}{2\pi\sqrt{-1}} \big(
\mbox{Li}_2(e^{-2\pi n\sqrt{-1}/N}) - \mbox{Li}_2(1) \big) \Big), \notag
\end{align}
we have that
\begin{align*}
J_N(K_{5_2})
\underset{?}{\sim} \sum_{0 \le m \le n < N}
& \exp \Big( \frac{N}{2\pi\sqrt{-1}} \big( \frac{\pi^2}{2} 
-2 \mbox{Li}_2(e^{2\pi n\sqrt{-1}/N}) \\*
& \qquad\qquad
- \mbox{Li}_2(e^{-2\pi m\sqrt{-1}/N})
+ \frac{2\pi n}{N} \frac{2\pi m}{N} \big) \Big).
\end{align*}
Further, by formally replacing\footnote{\footnotesize 
It might be seriously difficult to justify this replacement,
since there is a large parameter $N$ in the power of the summand,
which exponentially contributes the summand.}
the sum with an integral
putting $t = n/N$ and $s = m/N$, we have that
\begin{align}
J_N(K_{5_2}) &\underset{??}{\sim}
N^2 \int_{0 \le s \le t \le 1}
 \exp \frac{N}{2\pi\sqrt{-1}} \big( \frac{\pi^2}{2}
- 2 \mbox{Li}_2 (e^{2\pi\sqrt{-1}t}) \notag \\*
&\qquad\qquad\qquad\qquad\qquad - \mbox{Li}_2(e^{-2\pi\sqrt{-1}s}) 
+ 2 \pi t \cdot 2 \pi s \big) d s d t \label{eq.qn_app2} \\*
&= - \frac{N^2}{4 \pi^2}
\int \exp \frac{N}{2\pi\sqrt{-1}} \big( \frac{\pi^2}{2}
- 2 \mbox{Li}_2 (z) - \mbox{Li}_2( \frac{1}{w} ) - \log z \log w \big) 
\frac{d w}{w} \frac{d z}{z}, \notag
\end{align}
where the second integral is over the domain
$\big\{ (z,w) \in \C^2 \ \big| \ |z|=|w|=1, \ 0 \le \mbox{arg}(w) 
\le \mbox{arg}(z) \le 2\pi \big\}$,
and the equality is obtained by putting
$z=e^{2\pi\sqrt{-1}t}$ and $w=e^{2\pi\sqrt{-1}s}$.
By applying the saddle point method\footnote{\footnotesize 
The saddle point method in multi-variables is not established yet.
This might be a technical difficulty.}
the asymptotic behaviour might be described by a critical value of 
\begin{equation}
\label{eq.Li2-52}
\frac{\pi^2}{2}
- 2 \mbox{Li}_2 (z) - \mbox{Li}_2( \frac{1}{w} ) - \log z \log w.
\end{equation}
Since a critical value of this function gives a hyperbolic volume
of $S^3 - K_{5_2}$,
this formal argument suggests 
Conjecture \ref{conj.vol_conj} for $K_{5_2}$.

It was shown by Yokota \cite{vol.Y},
following ideas due to Kashaev \cite{bene.K2} and 
Thurston \cite{Thurston_grenoble}, 
that the hyperbolic volume of the complement of any hyperbolic knot $K$
is given by a critical value of such a function as (\ref{eq.Li2-52}),
which is obtained from a similar computation of $J_N(K)$ as above.

\begin{prob}
Justify the above arguments rigorously.
\end{prob}

\begin{rem}
The asymptotic behaviour of $J_N(K)$
might be described by using
quantum invariants\index{quantum invariant!--- of 3-manifold!asymptotic behaviour of ---} 
 of $S^3-K$.
We have some ways to compute 
the asymptotic behaviour of such a quantum invariant, say,
when $K$ is a fibered knot
(in this case, $S^3-K$ is homeomorphic to 
a mapping torus of a homeomorphism of a punctured surface),
and when we choose a simplicial decomposition of
(a closure of) $S^3-K$.
For details, see remarks of Conjecture \ref{conj.vol_olim}. 
\end{rem}

The following conjecture is a complexification of 
the volume conjecture (Conjecture \ref{conj.vol_conj}).

\begin{conj}
\label{conj.vol_CS}
{\rm (H. Murakami, J. Murakami, M. Okamoto, T. Takata, 
\newline Y.~Yokota \cite{MMOTY})}
\ 
For a hyperbolic link $L$,\index{volume conjecture!--- for links} 
$$
2 \pi \sqrt{-1} \cdot \lim_{N \to \infty} \frac{\log J_N(L)}{N}
= {\rm CS}(S^3-L) + \sqrt{-1} \mbox{\rm vol}(S^3-L)
$$
for an appropriate choice of a branch of the logarithm,
where CS and vol denote
the Chern-Simons invariant\index{Chern-Simons!--- invariant!--- and colored Jones polynomial}
and the\index{hyperbolic volume!--- and colored Jones polynomial}  
hyperbolic volume respectively.
Moreover,
\begin{equation}
\label{eq.limJN1Jn}
\lim_{N \to \infty}
\frac{J_{N+1}(L)}{J_N(L)}
= \exp \Big( \frac{1}{2\pi\sqrt{-1}} 
\big({\rm CS}(S^3-L) + \sqrt{-1} \mbox{\rm vol}(S^3-L) \big) \Big).
\end{equation}
\end{conj}

\begin{rem}
It is shown \cite{MMOTY}, by formal calculations
(such as (\ref{eq.qn_app}) and (\ref{eq.qn_app2})), that
Conjecture \ref{conj.vol_CS} is ``true'' for 
$K_{5_2}, K_{6_1}, K_{6_3}, K_{7_2}, K_{8_9}$
and the Whitehead link.
\end{rem}

\begin{rem}
The statement for non-hyperbolic links
should probably be the same statement,
replacing $\mbox{\rm vol}(S^3-L)$ with $v_3 || S^3-L ||$.
Note that, if $L$ is not hyperbolic,
then it is also a problem (see Problem \ref{prob.topCS}) 
to find an appropriate definition of ${\rm CS}(S^3-L)$,
which might be given by (\ref{eq.limJN1Jn}).
It is necessary to assume that $L$ is not a split link,
since $J_N(L) = 0$ for a split link $L$.
\end{rem}

\begin{rem}[{\rm (H. Murakami)}] 
Zagier \cite{Zagier_V} gave a conjectural presentation
of the asymptotic behaviour of the following sum,
$$
J_N(K_{3_1}) = \sum_{k=0}^{N-1} (q)_k 
\underset{N \to \infty}{\sim}
\exp \Big( - \frac{\pi\sqrt{-1}}{12} ( N-3 + \frac{1}{N}) \Big) N^{3/2}
+ \sum_{k \ge 0} \frac{b_k}{k!} \big( - \frac{2\pi\sqrt{-1}}{N} \big)^k
$$
for some series $b_k$.
This suggests that
$\lim \frac{\log J_N(K_{3_1})}{N}$ should be $-\pi\sqrt{-1}/12$.
\end{rem}

\begin{prob}[H. Murakami] 
\label{prob.CS_torus_knot}
For a torus knot $K$, 
calculate ${\rm CS}(S^3-K)$ (giving an appropriate definition of it)
and calculate $\lim \frac{\log J_N(K)}{N}$
(fixing an appropriate choice of a branch of the logarithm).
\end{prob}

\newpage

\section{Finite type invariants of knots}
\label{sec.Vinv}

Let $R$ be a commutative ring with $1$ such as $\Z$ or $\Q$.
We denote by $\K$ the set of isotopy classes of oriented knots.
A {\it singular knot} is an immersion of $S^1$ into $S^3$
whose singularities are transversal double points.
We regard singular knots as in $R \K$ 
by removing each singularity linearly by
$$
\pict{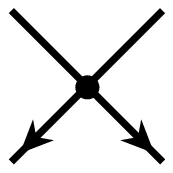}{1cm}
= \pict{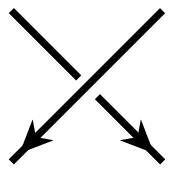}{1cm}
- \pict{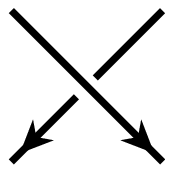}{1cm}.
$$
Let $\cF_d(R\K)$ denote the submodule of $R \K$ spanned by 
singular knots with $d$ double points,
regarding them as in $R \K$.
Then, we have a descending series of submodules,
$$
R \K = \cF_0(R \K)\ \supset\ \cF_1(R \K)\ \supset\ 
\cF_2(R \K)\ \supset\ \cdots.
$$
An $R$-homomorphism $v: R \K \to R$
(or, a homomorphism $\Z \K \to A$ for an abelian group $A$)
is called a {\it Vassiliev invariant}\index{Vassiliev invariant} 
(or a {\it finite type invariant})\index{finite type invariant!--- of knots|see {Vassiliev invariant}}
of degree $d$
if $v |_{\cF_{d+1}(R\K)} = 0$.
See \cite{BarNatan_website}
for many references of Vassiliev invariants.

\begin{figure}[ht!]
$$
\begin{array}{rcl}
\mboxsm{The AS relation}&:&
\pict{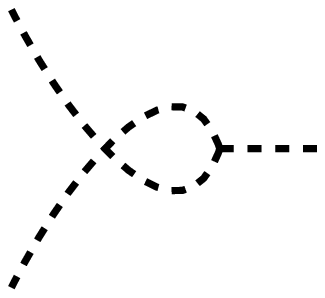}{50pt} = - \pict{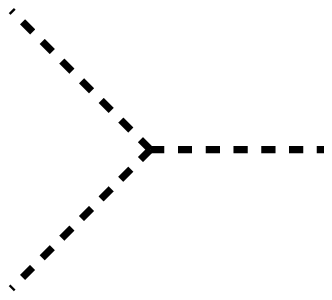}{50pt} \\
\mboxsm{The IHX relation}&:&
\pict{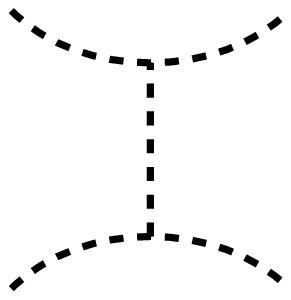}{50pt} = \pict{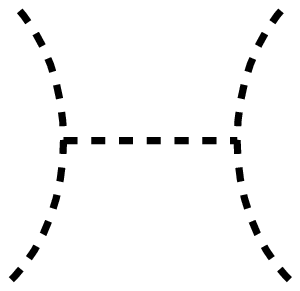}{50pt} - \pict{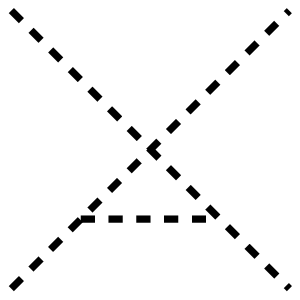}{50pt} \\
\mboxsm{The STU relation}&:&
\pict{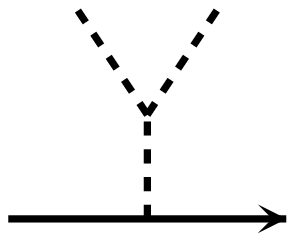}{50pt} = \pict{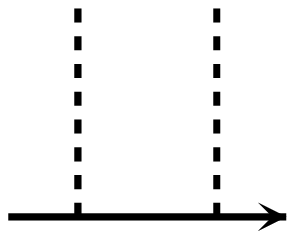}{50pt} - \pict{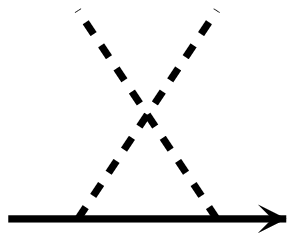}{50pt} \\
\mboxsm{The FI relation}&:&
\quad \pict{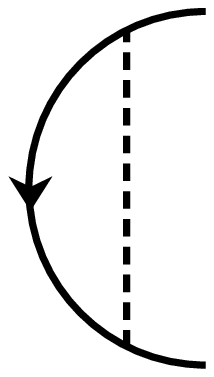}{1.3cm} = 0
\end{array}
$$
\caption{\label{fig.AS-IHX-STU} The AS, IHX, STU, and FI relations}
\end{figure}

A trivalent vertex of a graph is called {\it vertex-oriented} 
if a cyclic order of the three edges
around the trivalent vertex is fixed.
A\index{Jacobi diagram} 
{\it Jacobi diagram\/}\footnote{\footnotesize 
A Jacobi diagram is also called a {\it web diagram}
or a {\it trivalent diagram} in some literatures.
In physics this is often called a {\it Feynman diagram}.}
on an oriented 1-manifold $X$ is
the manifold $X$ together with
a uni-trivalent graph such that
univalent vertices of the graph are distinct points on $X$ and 
each trivalent vertex is vertex-oriented.
The {\it degree} of a Jacobi diagram
is half the number of univalent and trivalent vertices
of the uni-trivalent graph of the Jacobi diagram.
We denote by $\cA(X;R)$ the module over $R$ 
spanned by Jacobi diagrams on $X$
subject to the AS, IHX, and STU relations
shown in Figure \ref{fig.AS-IHX-STU},
and denote by $\cA(X;R)^{(d)}$ the submodule of $\cA(X;R)$ 
spanned by Jacobi diagrams of degree $d$.
There is a canonical surjective homomorphism
\begin{equation}
\label{eq.AS1toFdFd1}
\cA(S^1;R)^{(d)}/\mbox{FI} \to \cF_d(R\K)/\cF_{d+1}(R\K),
\end{equation}
where FI is the relation shown in Figure \ref{fig.AS-IHX-STU}.
This map is known to be an isomorphism when $R=\Q$
(due to Kontsevich).
For a Vassiliev invariant $v : R \K \to R$ of degree $d$,
its weight system\index{weight system} 
$\cA(S^1;R)^{(d)}/\mbox{FI} \to R$ is defined
by the map (\ref{eq.AS1toFdFd1}).

\newpage

\subsection{Torsion and Vassiliev invariants}

Let $R$ be a commutative ring with $1$, say $\Z/n\Z$.
Then, 
$\Q$-, $\Z$-, $R$-valued Vassiliev invariants and their weight systems
and the Kontsevich invariant
form the following commutative diagram.
\begin{center}
\begin{picture}(380,160)
\put(242,150){\vector(-1,0){160}}
\put(120,155){\footnotesize Kontsevich invariant}
\put(0,80){
$
\begin{CD}
\cA(S^1;\Q)/\mbox{FI} @. @. \K \\
@VV{\mbox{\footnotesize proj}}V @. @VVV \\
\cA(S^1;\Q)^{(d)}/\mbox{FI} @>{\mbox{\footnotesize isom}}>>
\cF_d(\Q\K)/\cF_{d+1}(\Q\K) @>{\subset}>>
\Q\K/\cF_{d+1}(\Q\K) @>>> \Q \\
@AA{\cdot \otimes \Q}A @AA{\cdot \otimes \Q}A @A{\cdot \otimes \Q}AA @A{\cdot \otimes \Q}AA \\
\cA(S^1;\Z)^{(d)}/\mbox{FI} @>{\mbox{\footnotesize surj}}>>
\cF_d(\Z\K)/\cF_{d+1}(\Z\K) @>{\subset}>>
\Z\K/\cF_{d+1}(\Z\K) @>>> \Z \\
@VV{\mbox{\footnotesize proj}}V @VV{\mbox{\footnotesize proj}}V @V{\mbox{\footnotesize proj}}VV @V{\mbox{\footnotesize proj}}VV \\
\cA(S^1;R)^{(d)}/\mbox{FI} @>{\mbox{\footnotesize surj}}>>
\cF_d(R\K)/\cF_{d+1}(R\K) @>{\subset}>>
R\K/\cF_{d+1}(R\K) @>>> R
\end{CD}
$}
\end{picture}
\end{center}
Here, the right horizontal maps are derived from Vassiliev invariants,
and the compositions of horizontal maps are their weight systems.

\begin{conj}[{\cite[Problem 1.92 (N)]{Kirby}}] 
\label{conj.KK_tf}
$\cF_d(\Z\K) / \cF_{d+1}(\Z\K)$ is 
torsion free\index{torsion!--- free} 
for each $d$.
\end{conj}

\begin{rem}[\rm (see {\cite[Remark on Problem 1.92 (N)]{Kirby}})] 
Goussarov has checked the conjecture for $d \le 6$.
It has been checked that
$\cF_d(\Z\K) / \cF_{d+1}(\Z\K)$ 
has no 2-torsion for $d \le 9$ by Bar-Natan,
and for $d \le 12$ in \cite{Kneissler}.
\end{rem}

\begin{rem}
If this conjecture was true, 
then $\Z$-valued and $\Q$-valued Vassiliev invariants 
carry exactly the same information about knots.
Moreover, any $(\Z/ n\Z)$-valued Vassiliev invariants
would be derived from $\Z$-valued Vassiliev invariants.
\end{rem}

\begin{conj}
\label{conj.AS1_tf}
$\cA(S^1;\Z)$ is\index{Jacobi diagram!space of ---!torsion of ---}
torsion free.\index{torsion!--- free} 
\end{conj}

\begin{rem}[{\rm (T. Stanford)}] 
This conjecture would imply Conjecture \ref{conj.KK_tf} 
because of the Kontsevich integral.
However, it is possible that 
there is torsion in $\cA(S^1,\Z)^{(d)}$ which is in the kernel
of the map (\ref{eq.AS1toFdFd1}).
\end{rem}

\begin{conj}[{X.-S. Lin \cite{Lin_website}}] 
\label{conj.torsion_ws}
Let $R$ be a commutative ring with 1, say $\Z/2\Z$.
Every weight system\index{weight system!realization by Vassiliev invariant}
$\cA(S^1;R)^{(d)}/\mbox{FI} \to R$
is induced by some Vassiliev invariant 
$R \K \to R$.\index{Vassiliev invariant!weight system of ---}
\end{conj}

\begin{rem}
If the map (\ref{eq.AS1toFdFd1}) 
is an isomorphism and
$\cF_d(R\K)/\cF_{d+1}(R\K)$
is a direct summand of $R\K/\cF_{d+1}(R\K)$,
then this conjecture is true 
(see the diagram at the beginning of this section).
\end{rem}

\begin{rem}
When $R=\Q$, this conjecture is true,
since the composition of the Kontsevich invariant
 and a weight system
gives a Vassiliev invariant, which induces the weight system.
If the Kontsevich invariant with coefficients in $R$ would be constructed
(see Problem \ref{prob.Kinv_ff}),
this conjecture would be true.
\end{rem}

\begin{figure}[ht!]
$$
\pict{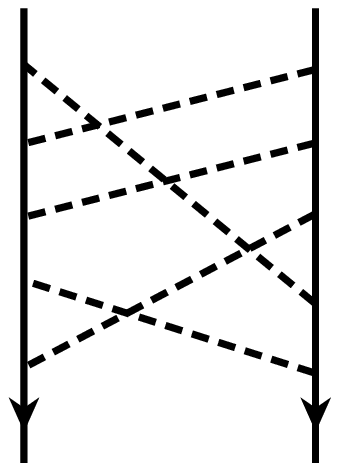}{2.2cm}
\ - \ \pict{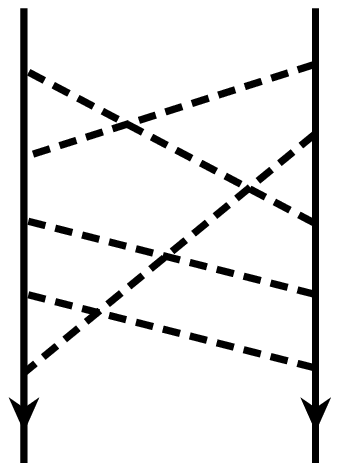}{2.2cm}
$$
\caption{\label{fig.DK2torsion}
A 2-torsion element in $\cA(\downarrow\downarrow;\Z)$
due to Dogolazky--Kneissler}
\end{figure}

\begin{rem}[\rm (T. Stanford)] 
The chord diagram module $\cA(\downarrow\downarrow,\Z)$
corresponds to finite-type invariants\index{finite type invariant!--- of string links} 
of two-strand string links.  
Jan Kneissler and Ilya Dogolazky (see \cite{Dogolazky}) showed that  
there is a 2-torsion element in $\cA(\downarrow\downarrow,\Z)^{(5)}/\mbox{FI}$
(see Figure \ref{fig.DK2torsion}).
I have done recent calculations (to be written up soon) which show that 
there is no $\Z/2\Z$-valued invariant of string links 
corresponding to this torsion element.  
Thus there is a $\Z/2\Z$ weight system $\cA(\downarrow\downarrow,\Z/2\Z)/\mbox{FI} \to \Z/2\Z$ 
which is not induced by a $\Z/2\Z$-valued finite-type invariant.  
So for string links, Conjecture 2.1 is false.
\end{rem}

\noindent
{{\namae}(T. Stanford)}\quad
Note that the Kontsevich integral
 works (for rational invariants)
for string links just as well as for knots.  
Since this calculation shows that 
there is no $\Z/2\Z$ Kontsevich integral for string links,
it suggests that there is no $\Z/2\Z$ Kontsevich integral for knots.


\begin{quest}[T. Stanford] 
\label{quest.2tor_cA}
The Dogolazky-Kneissler 2-torsion\index{torsion}  
element in $\cA(\downarrow\downarrow,\Z)$ 
(see Figure \ref{fig.DK2torsion})
can be embedded into a chord diagram in $\cA(S^1,\Z)$ in many ways.  
Such an embedding will always produce an element $x \in \cA(S^1,\Z)$ 
with $2x = 0$.  
Is it possible to produce such an $x$ which is nontrivial?  
If so, this would give a counterexample to Conjecture \ref{conj.torsion_ws}. 
\end{quest}

\subsection{Do Vassiliev invariants distinguish knots?}

\begin{conj}
\label{conj.V_dist_knots}
Vassiliev invariants distinguish oriented 
knots.\index{Vassiliev invariant!strength of ---}
(See Conjecture \ref{conj.K_dist_knots} for an equivalent statement
of this conjecture.)
\end{conj}

\begin{rem}
Two knots with the same Vassiliev invariant
up to an arbitrarily given degree
can be obtained; see \cite{Ohyama_sim}
and Goussarov-Habiro theory \cite{Gou_kg,Gou_fn,Habiro_GT}.
Hence, we need infinitely many Vassiliev invariants
to show this conjecture.
\end{rem}

\begin{prob}
\label{prob.VdistKfromO}
Does there exists a non-trivial oriented knot 
which can not be distinguished from the trivial knot\index{Vassiliev invariant!strength of ---}
by Vassiliev invariants?
(See Problem \ref{prob.KdistKfromO} for an equivalent problem.)
\end{prob}

\begin{rem}
The volume conjecture (Conjecture \ref{conj.vol_conj})
suggests that the answer is no; see \cite{vol.MM}.
\end{rem}

\begin{conj}[{see \cite[Problem 1.89 (B)]{Kirby}}]
\label{conj.Vdist_rev}
For any oriented knot $K$,
no Vassiliev invariants distinguish $K$ from $-K$.\index{Vassiliev invariant!strength of ---}
(See Conjecture \ref{conj.Kdist_rev} for an equivalent statement
of this conjecture.)
\end{conj}

\begin{rem}[\rm ({\cite[Remark on Problem 1.89]{Kirby}})] 
The first oriented knot which is different from its 
reverse\index{reversed knot} 
is $8_{17}$.
It is known that no Vassiliev invariants of degree $\le 9$
can distinguish a knot from its reverse.
\end{rem}

\begin{rem}
This conjecture is reduced to the problem
to find $D \in \cA(S^1)$ with $D \ne -D$,
where $-D$ is $D$ with the opposite orientation of $S^1$.
If such a $D$ existed,
the conjecture would fail.
Such a $D$ has not been known so far.
\end{rem}

\begin{rem}
Kuperberg \cite{Kuperberg_d} showed that 
all Vassiliev invariants either distinguish all oriented knots, 
or there exist prime, unoriented knots which they do not distinguish.
\end{rem}

\subsection{Can Vassiliev invariants detect other invariants?}

\noindent 
{{\namae}(T. Stanford)}\quad 
Let $h_G(K)$ be the number of homomorphisms
from the fundamental group of the complement of a knot $K$
to a finite group\index{finite group} 
$G$.
This is not a Vassiliev invariant \cite{Eisermann_nV}.
$h_{{\frak S}_3}(K)$ of the 3rd symmetric group ${\frak S}_3$
is presented by the number of 3-colorings of $K$, 
and $h_{D_5} (K)$ of the dihedral group $D_5$ of order 10
is presented by the number of 5-colorings of $K$.  
These are determined by the Jones and Kauffman polynomials, respectively 
(see the remark of Problem \ref{prob.skeinmod_7col}),  
and therefore are determined by invariants of finite type.  
In fact, by the usual power-series
expansions of the Jones and Kauffman polynomials,
we see that 
$h_{{\frak S}_3}$
and $h_{D_5}$ are the (pointwise) limits of respective
sequences of finite-type invariants.

\begin{quest}[T. Stanford] 
\label{quest.hG_Vinv}
Can we approximate $h_G$ by 
Vassiliev invariants\index{Vassiliev invariant!detectability by ---} 
for other $G$
than dihedral groups?
\end{quest}

\begin{rem}[{\rm (T. Stanford)}] 
It is known (due to W. Thurston) that
knot groups are residually finite.
So if $h_G$ can be approximated by finite-type
invariants for all finite groups $G$, then
Vassiliev invariants would distinguish the unknot.
\end{rem}

\begin{rem}[{\rm (T. Stanford)}] 
If $p$ is a prime, then there exists a nontrivial $p$-coloring
of a knot $K$, and hence a nontrivial representation of the
fundamental group of $K$ 
into the dihedral group $D_p$ of order $2p$, if and only if
$\Delta_K(-1)$ is divisible by $p$.
Thus the Alexander polynomial contains information about $h_{D_p}$, though
it may not determine $h_{D_p}$ completely.  
Suppose that $G$ is a finite, non-abelian group,
not isomorphic to $D_p$.
Even if we cannot approximate $h_G$ by finite-type
invariants, it would at least be interesting to know
whether finite-type invariants provide any information
at all about $h_G$.
\end{rem}

\begin{rem}
Let $h_X(K)$ denote
the number of homomorphisms
from the knot quandle of a knot $K$
to a finite quandle $X$.
The number $h_G(K)$ can be presented by the sum of $h_X(K)$
for subquandles $X$ of the conjugation quandle of $G$.
In this sense, it is a refinement of Question \ref{quest.hG_Vinv}
to approximate $h_X$ of finite quandles $X$
by Vassiliev invariants.
It is known \cite{Inoue} that
$h_X(K)$ for certain Alexander quandles $X$
can be presented by the $i$th Alexander polynomial of $K$.
\end{rem}

\begin{prob}[{X.-S. Lin \cite{Lin_website}}] 
\label{prob.sign_Vinv}
Is the knot signature\index{signature!--- of knot} 
the limit of a sequence of\index{Vassiliev invariant!detectability by ---}
Vassiliev invariants?
\end{prob}

\begin{rem}
It is known \cite{Dean} that
the signature of knots is not a Vassiliev invariant.
\end{rem}

\subsection{Vassiliev invariants and crossing numbers}

Let $v_2$ and $v_3$ be
$\R$-valued Vassiliev invariants of degree 2 and 3 respectively
normalized by the conditions that
$v_2(\overline{K})=v_2(K)$
and $v_3(\overline{K})=-v_3(K)$
for any knot $K$ and its mirror image $\overline{K}$
and that $v_2(K_{\overline{3}_1})=v_3(K_{\overline{3}_1})=1$
for the right trefoil knot $K_{\overline{3}_1}$.
They are primitive Vassiliev invariants,
and the image of $v_2 \times v_3$ is equal to 
$\Z \times \Z \subset \R \times \R$.

\begin{prob}[N. Okuda \cite{Okuda}] 
\label{prob.Okuda_fish}
Describe the set\index{Vassiliev invariant!--- and crossing number}\index{crossing number} 
\begin{equation}
\label{eq.fish_set}
\Big\{ \big( \frac{v_2(K)}{n^2}, \frac{v_3(K)}{n^3} \big) 
\in \R \times \R \ \Big| \ 
\mbox{$K$ has a knot diagram with $n$ crossings}
\Big\}.
\end{equation}
\end{prob}

\begin{figure}[ht!]
\begin{center}
\begin{picture}(250,200)
\put(0,95){\pict{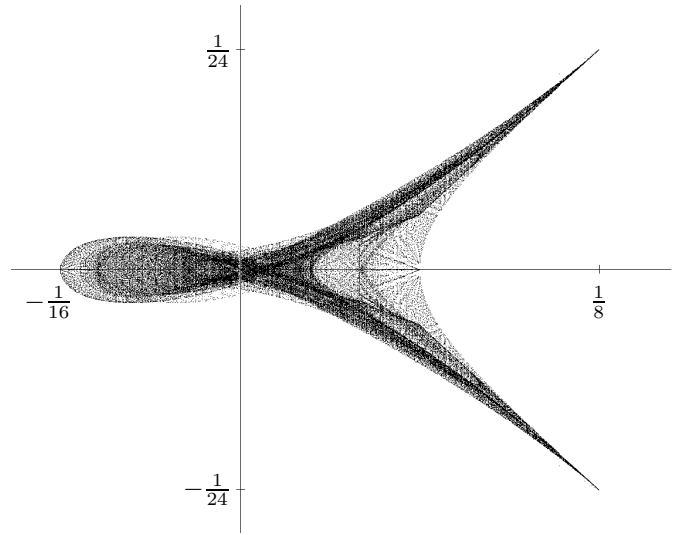}{7cm}}
\put(10,85){\small $-\frac{1}{16}$}
\put(224,85){\small $\frac{1}{8}$}
\put(78,180){\small $\frac{1}{24}$}
\put(70,14){\small $-\frac{1}{24}$}
\end{picture}
\end{center}
\caption{\label{fig.fish}
The plottings of the set (\ref{eq.fish_set}) 
for some family of knots \cite{Okuda}}
\end{figure}

\begin{rem}
Willerton \cite{Wil_fish} observed that
the set of $(v_2(K),v_3(K))$ for knots $K$
with a (certain) fixed crossing number gives a fish-like graph.
This fish-like graph is discussed in \cite{DaLeLi01}
from the point of view of the Jones polynomial.
\end{rem}

\begin{rem}[{\rm (N. Okuda)}] 
It is shown by Okuda \cite{Okuda} 
(the right inequality of (\ref{eq.v2_n}) is due to \cite{PoVi_C}) 
that, 
if a knot $K$ has a diagram with $n$ crossings, then
\begin{align}
\label{eq.v2_n}
& - \Big\lfloor \frac{n^2}{16} \Big\rfloor \le v_2(K) \le
\Big\lfloor \frac{n^2}{8} \Big\rfloor, \\*
\label{eq.v3_n}
& |v_3(K)| \le \Big\lfloor \frac{n(n-1)(n-2)}{15} \Big\rfloor,
\end{align}
where $\lfloor x \rfloor$ denotes
the greatest integer less than or equal to $x$.
It follows that the set (\ref{eq.fish_set}) is included in the rectangle
$[ -1/16, 1/8 ] \times [ -1/15, 1/15 ]$.
It is a problem to describe the smallest domain including this set.
The plottings in Figure \ref{fig.fish}
numerically describe the set (\ref{eq.fish_set})
for a large finite subset of a certain infinite family of knots.
Okuda \cite{Okuda} identified the boundary
of the domain including this set for this infinite family of knots.
This boundary is given by curves
parameterized by some polynomials of degree 2 (for the $v_2$-coordinate)
and of degree 3 (for the $v_3$-coordinate).
Further, the points $(1/8,\pm1/24)$ are the limits of the points 
given by the $(2,n)$ torus knot and its mirror image.
The point $(-1/16,0)$ is the limit of the points given by the knots
$$
\begin{picture}(100,70)
\put(0,27){\pict{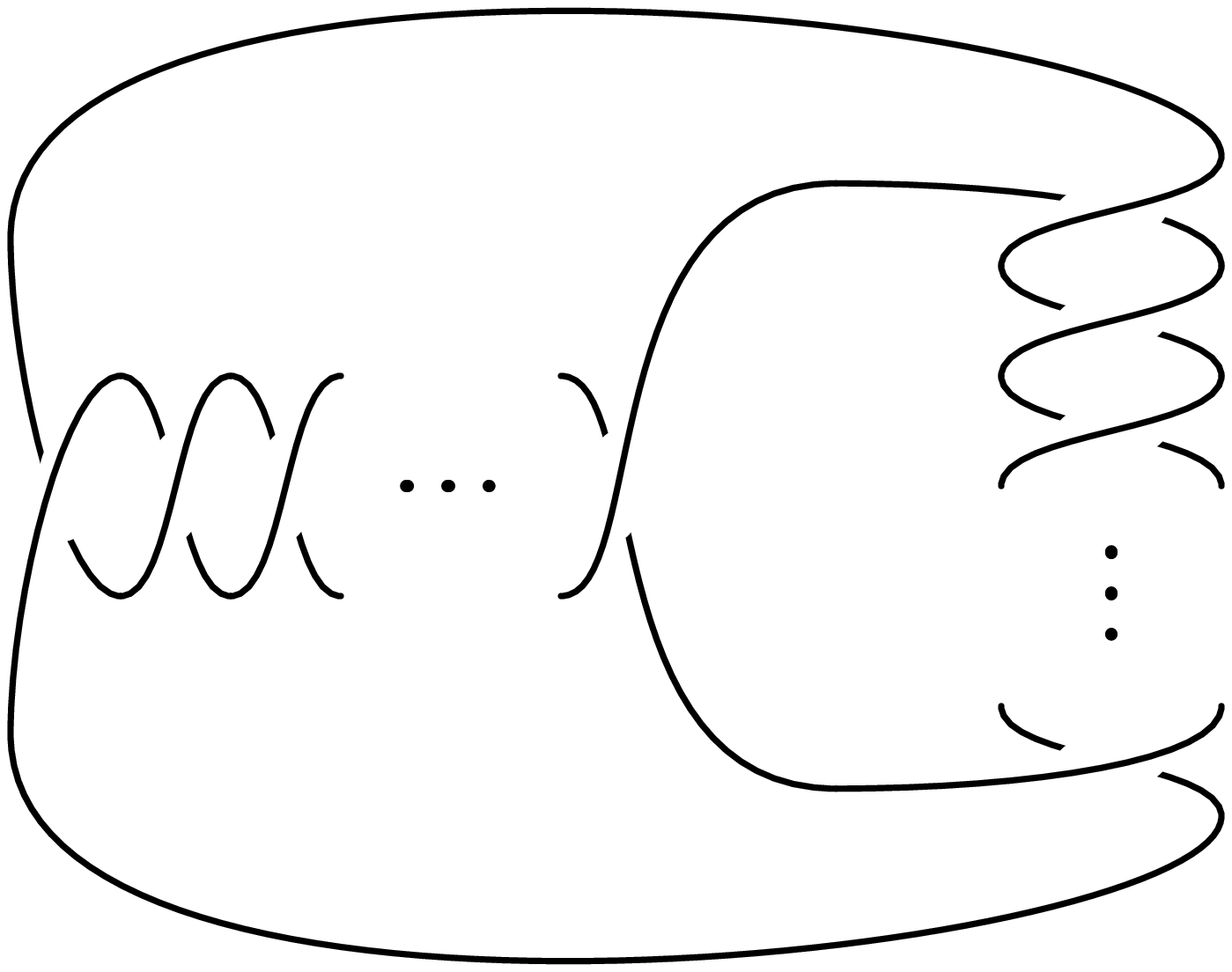}{2.5cm}}
\put(22,47){\small $n/2$}
\put(100,30){\small $n/2$}
\end{picture}
$$
for $n$ divisible by 4, where each twisting part has $n/2$ crossings.
These knots gives the bounds of (\ref{eq.v2_n}),
while the inequality in (\ref{eq.v3_n}) might not be best possible
(see Conjecture \ref{conj.bound_v3} below).
\end{rem}

\begin{rem}[{\rm (O. Viro)}] 
The experimental data (in Figure \ref{fig.fish})
suggest that there might exists an additional
invariant(s) which together with $v_2$, $v_3$, and $n$ satisfy an
algebraic equation(s) such that the set (\ref{eq.fish_set}) is the projection
of the algebraic set defined by the equation(s).
\end{rem}

\begin{conj}[S. Willerton \cite{Wil_fish}] 
\label{conj.bound_v3}
Let $v_3$ be as above.
If a knot $K$ has a diagram with $n$ crossings,\index{crossing number}  
then
$$
|v_3(K)| \le \Big\lfloor \frac{n(n^2-1)}{24} \Big\rfloor.
$$
\end{conj}

\begin{rem}
It is shown in \cite{Wil_fish} that
the $(2,n)$ torus knot gives the equality of this formula.
\end{rem}


\subsection{Dimensions of spaces of Vassiliev invariants}

We denote by $\cA(S^1;R)\conn$ the submodule of $\cA(S^1;R)$
spanned by Jacobi diagrams with connected uni-trivalent graphs.
As a graded vector space
$\cA(S^1;\Q)$ is isomorphic to
the symmetric tensor algebra of $\cA(S^1;\Q)\conn$.
A Vassiliev invariant $v$ is called {\it primitive}
if $v( K_1 \# K_2 ) = v(K_1) + v(K_2)$
for any oriented knots $K_1$ and $K_2$.
The degree $d$ subspace of $\cA(S^1;\Q)\conn$ 
is dual to the $d$th graded vector space 
for $\Q$-valued primitive Vassiliev invariants.

\begin{prob}
\label{prob.dimV}
Determine the dimension\index{Vassiliev invariant!dimension of the space of ---}
of the space of primitive Vassiliev invariants of each degree $d$.
Equivalently,
determine the dimension of the 
space\index{Jacobi diagram!space of ---!dimension of ---}  
$\cA(S^1;\Q)\conn^{(d)}$ for each $d$.
\end{prob}

\begin{table}[ht!]
\begin{center}
\begin{tabular}{|c||c|c|c|c|c|c|c|c|c|c|c|} \hline
\small $d$ &\small 0 &\small 1 &\small 2 &\small 3 &\small 4 &\small 5 &\small 6 &\small 7 &\small 8 &\small 9 &\small 10
\\ \hline\hline
dim $\cA(S^1)\conn^{(d)}$ &\small 
0 &\small 1 &\small 1 &\small 1 &\small 2 &\small 3 &\small 5 &\small 8 &\small 12 &\small 18 &\small 27 
\\ \hline
\small dim $\cA(S^1)^{(d)}$ &\small 
1 &\small 1 &\small 2 &\small 3 &\small 6 &\small 10 &\small 19 &\small 33 &\small 60 &\small 104 &\small 184 
\\ \hline
\small dim $\cA(S^1)^{(d)}/\mbox{FI}$ &\small 
1 &\small 0 &\small 1 &\small 1 &\small 3 &\small 4 &\small 9 &\small 14 &\small 27 &\small 44 &\small 80 
\\ \hline
\end{tabular}
\end{center}
\begin{center}
\begin{tabular}{|c||c|c|c|c|} \hline
\small $d$ &\small 11 &\small 12 &\small 13 &\small 14\\ \hline\hline
\small dim $\cA(S^1)\conn^{(d)}$ 
&\small 39 &\small 55 &\small $\ge78$ &\small $\ge108$ \\ \hline
\small dim $\cA(S^1)^{(d)}$ 
&\small 316 &\small 548 &\small $\ge932$ &\small $\ge1591$\\ \hline
\small dim $\cA(S^1)^{(d)}/\mbox{FI}$ 
&\small 132 &\small 232 &\small $\ge384$ &\small $\ge 659$\\ \hline
\end{tabular}
\end{center}
\caption{Some dimensions given in \cite{Broadhurst,Kneissler}} 
\end{table}

\begin{rem}
The dimension of $\cA(S^1;\Q)\conn^{(d)}$ can partially be computed as follows.

Let $\cal B$ denote the vector space over $\Q$
spanned by vertex-oriented uni-trivalent graphs
subject to the AS and IHX relations,
and let $\cB\conn^{(d)}$ denote the subspace of $\cB$
spanned by connected uni-trivalent graphs with $2d$ vertices.

It is known that $\cA(S^1;\Q)\conn^{(d)}$ is isomorphic to $\cB\conn^{(d)}$
by (\ref{eq.A=B=B}).
Let $\cB\conn^{(d,u)}$ be the subspace of $\cB\conn^{(d)}$
spanned by uni-trivalent graphs
with $u$ univalent vertices
(hence, with $2d-u$ trivalent vertices),
and $\beta_{d,u}$ its dimension.
Then, the dimension of $\cA(S^1;\Q)\conn^{(d)}$
is presented by  $\sum_{u \ge 2} \beta_{d,u}$.

Bar-Natan \cite{BarNatan_scV} 
gave a table of $\beta_{d,u}$ for $d \le 9$
and for some other $(d,u)$ by computer.

The series of $\beta_{k,k}$ is given as follows.
The direct sum $\oplus_k \cB\conn^{(k,k)}$ is isomorphic to 
the polynomial ring $\Q[x^2]$
as a graded vector space by (\ref{eq.loop0});
in other words, it is spanned by ``wheels''.
Hence, the series of $\beta_{k,k}$ is presented by 
the following generating function,
$$
\sum_{k \ge 0} \beta_{k,k} t^k = \frac1{1-t^2}.
$$
That is, 
$\beta_{k,k} = 1$ if $k$ is even, and $0$ otherwise.

The series of $\beta_{k+1,k}$ is given as follows.
The direct sum $\oplus_k \cB\conn^{(k+1,k)}$ is isomorphic
to $\Q[\sigma_2, \sigma_3^2]$
as a graded vector space by (\ref{eq.loop1}),
where $\sigma_i$ denotes
the $i$-th elementary symmetric polynomial in some variables.
Hence, the series of $\beta_{k+1,k}$ is presented by 
the following generating function,
$$
\sum_{k \ge 0} \beta_{k+1,k} t^k = \frac1{(1-t^2)(1-t^6)}.
$$

The series of $\beta_{k+2,k}$ is presented by
$$
\sum_{k \ge 0} \beta_{k+2,k} t^k = \frac1{(1-t^2)(1-t^4)(1-t^6)},
$$
since $\oplus_k \cB\conn^{(k+2,k)}$
is isomorphic, as a graded vector space,
to $\Q[\sigma_2, \sigma_3^2, \sigma_4]$
with elementary symmetric polynomials in some variables by (\ref{eq.loop2}).

It is conjectured \cite{Dasbach_III} that
the series of $\beta_{k+3,k}$ would be presented by
$$
\sum_{k \ge 0} \beta_{k+3,k} t^k 
\stackrel{?}{=}
\frac{1+t^2+t^8-t^{10}}{(1-t^2)(1-t^4)(1-t^6)(1-t^{10})}.
$$

It has been shown that
$\beta_{d,u} = 0$ for $d \le 9$ and for $d \le u+2$.
However, it is conjectured yet for other $(d,u)$.

A conjecture of
a two-variable generating function for the series of $\beta_{d,u}$
with two parameters $d$ and $u$
is given in \cite{Broadhurst}.
\end{rem}

\begin{table}[ht!]
\begin{center}
\begin{tabular}{|c||c|c|c|c|c|c|c||c|} \hline
\small $\beta_{d,u}$ &\small $u=2$ &\small $u=4$ &\small $u=6$ &\small $u=8$ &\small $u=10$ &\small $u=12$ &\small $u=14$ &\small total \\ 
\hline\hline
\small $d=1$ &\small  1 &&&&&&&\small  1 \\ \hline
\small $d=2$ &\small  1 &&&&&&&\small  1 \\ \hline
\small $d=3$ &\small  1 &&&&&&&\small  1 \\ \hline
\small $d=4$ &\small  1 &\small  1 &&&&&&\small  2 \\ \hline
\small $d=5$ &\small  2 &\small  1 &&&&&&\small  3 \\ \hline
\small $d=6$ &\small  2 &\small  2 &\small  1 &&&&&\small  5 \\ \hline
\small $d=7$ &\small  3 &\small  3 &\small  2 &&&&&\small  8 \\ \hline
\small $d=8$ &\small  4 &\small  4 &\small  3 &\small  1 &&&&\small  12 \\ \hline
\small $d=9$ &\small  5 &\small  6 &\small  5 &\small  2 &&&&\small  18 \\ \hline
\small $d=10$ &\small  6 &\small  8 &\small  8 &\small  4 &\small  1 &&&\small  27 \\ \hline
\small $d=11$ &\small  8 &\small  10 &\small  11 &\small  8 &\small  2 &&&\small  39 \\ \hline
\small $d=12$ &\small  9 &\small  13 &\small  15 &\small  12 &\small  5 &\small  1 &&\small  55 \\ \hline
\small $d=13$ &\small $\ge11$ &\small $\ge16$ &\small $\ge20$ &\small $\ge18$ &\small $\ge10$ &\small 3 &&\small  $\ge78$ \\ \hline
\small $d=14$ &\small $\ge13$ &\small $\ge19$ &\small $\ge25$ &\small $\ge26$ &\small $\ge17$ &\small 7 &\small 1 &\small  $\ge108$ \\ \hline
\end{tabular}
\end{center}
\caption{\label{tbl.bdu}
   A table of $\beta_{d,u}$ \cite{Broadhurst,Kneissler}}
\end{table}

\begin{rem}
An asymptotic evaluation of a lower bound of 
$\mbox{dim} \cA(S^1)\conn^{(d)}$ was given in \cite{ChDu99};
$\mbox{dim} \cA(S^1)\conn^{(d)}$ grows at least as $d^{\log d}$
when $d \to \infty$.
Further, it was improved in \cite{Dasbach_III};
$\mbox{dim} \cA(S^1)\conn^{(d)}$
grows at least as
$e^{c \sqrt{d}}$ for any $c < \pi \sqrt{2/3}$ when $d \to \infty$.
\end{rem}

\begin{rem}
Upper bounds of $\mbox{dim} \cA(S^1)\conn^{(d)}$
were given
$\mbox{dim} \cA(S^1)\conn^{(d)} \le (d-1)!$ in \cite{ChDu94}
and 
$\mbox{dim} \cA(S^1)\conn^{(d)} \le (d-2)!/2$ (for $d > 5$)
in \cite{NgSt99}.
Stoimenow \cite{Stoimenow98}
introduced the number $\xi_d$
of ``regular linearized chord diagrams'',
and showed that
$\mbox{dim} \cA(S^1)^{(d)}/\mbox{FI} \le \xi_d$.
Further, he showed that
$\xi_d$ is asymptotically at most $d! / 1.1^d$,
which was improved by
$d! / ( 2 \ln 2 + o(1) )$ in \cite{BoRi00}.
Furthermore, Zagier \cite{Zagier_V} showed that
\begin{equation}
\label{eq.Zag}
\sum_{n=0}^\infty (1-q)(1-q^2) \cdots (1-q^n)
= \sum_{d=0}^\infty \xi_d (1-q)^d
\in \Z[[1-q]],
\end{equation}
and that
$$
\xi_d \sim \frac{d! \sqrt{d}}{(\pi^2/6)^d} 
\big( C_0 + \frac{C_1}{d} + \frac{C_2}{d^2} + \cdots \big)
$$
with $C_0 = 12 \sqrt{3} \pi^{-5/2} e^{\pi^2/12} \approx 2.704$,
$C_1 \approx -1.527$, $C_2 \approx -0.269$.
It follows that
the asymptotic growth of
$\mbox{dim}\cA(S^1)^{(d)}/\mbox{FI}$
is at most $O( d! \sqrt{d} (\pi^2/6)^{-d} )$.
\end{rem}

\begin{table}[ht!]
\begin{center}
\begin{tabular}{|c||c|c|c|c|c|c|c|c|c|c|c|} \hline
\small $d$ &\small 0 &\small 1 &\small 2 &\small 3 &\small 4 &\small 5 
&\small 6 &\small 7 &\small 8 &\small 9 &\small 10 
\\ \hline\hline
\small dim $\cA(S^1)^{(d)}/\mbox{FI}$ &\small 
1 &\small 0 &\small 1 &\small 1 &\small 3 &\small 4 &\small 9 &\small 14 
&\small 27 &\small 44 &\small 80 
\\ \hline
\small $\xi_d$ &\small  1 &\small  1 &\small  2 &\small  5 &\small  15 
&\small  53 &\small  217 &\small  1014 &\small  5335 &\small  31240 
&\small  201608 
\\ \hline
\end{tabular}
\end{center}
\caption{Upper bounds $\xi_d$ of $\mbox{dim}\cA(S^1)^{(d)}/\mbox{FI}$
(see \cite{Stoimenow98})
}
\end{table}

\subsection{Milnor invariants}

\noindent
{{\namae}(T. Stanford)}\quad
Fix $k$, and consider $k$-strand string links.  Let $V_n$ be
the subspace of rational-valued finite-type invariants of
order $\le n$ (of $k$-strand string links).  Let $M_n
\subset V_n$ be the subspace of Milnor invariants and
products of Milnor invariants.\index{Milnor invariant}
It is known that in general
$M_n$ is a proper subspace of $V_n$.  

\begin{quest}[T. Stanford] 
\label{quest.Milnorinv_Vinv}
Does $M_n$ have an
interesting complementary space in $V_n$?  Consider, for
example, the space $N_n \subset V_n$ of invariants $v$ with
the property that $v(L) = 0$ for any string link $L$ such
that $\pi_1 (B^3 - L)$ is free.  Is $N_n$ nontrivial?  Do
$N_n$ and $M_n$ together span $V_n$?
\end{quest}

Here is some background and motivation.

When considering finite-type invariants of string links, the
first ones that come to mind are the Milnor invariants.
These were defined by Milnor \cite{Milnor_lg}
in 1954 as numbers associated
to links.  They are not quite invariants of links, in the
usual sense, because of some indeterminacy.  They are,
however, well-defined as invariants of string links, and
this point of view was taken by Habegger and Lin \cite{HaLi_lh}.  
After Vassiliev's work appeared, Bar-Natan \cite{BarNatan_hsl}
and Lin \cite{Lin_ps} showed
(independently) that the Milnor invariants are finite-type
invariants.  Habegger and Masbaum \cite{HaMa_KM}
showed that on the chord
diagram level, the Milnor invariants (including products of
Milnor invariants) are exactly the ones 
that vanish on Jacobi diagrams that contain internal loops, and also that the
Milnor invariants are the only rational-valued finite-type
invariants of string links which are also concordance
invariants.

String links may have local knots in the strands, and such
knots are not detected by Milnor invariants.  If a string
link $L$ has local knots, then $\pi_1 (B^3 - L)$ is not
free.  Hence the question as to whether finite-type invariants
can show that the complement of a string link is not free.

\medskip

\noindent
{{\namae}(M. Polyak)}\quad
Let us review the constructions of Milnor $\overline\mu$-invariant
in \cite{Cochran}.
For a $n$-component link $L = L_1 \cup \cdots \cup L_n$,
regard the homotopy class of $L_n$ as in 
$\pi_1\big( S^3 - (L_1 \cup \cdots \cup L_{n-1}) \big)$,
and write it in terms of meridians $m_1, \cdots, m_{n-1}$
of $L_1, \cdots, L_{n-1}$.
Consider its Magnus expansion putting $m_i = 1 + X_i$
for non-commutative variables $X_i$.
Then, Milnor's $\overline\mu$-invariant
$\overline{\mu}_{i_1 \cdots i_k, n}(L)$
is defined to be the coefficient of $X^{i_1} \cdots X^{i_k}$
in the expansion,
which is an invariant under the assumption that
the lower $\overline\mu$-invariants vanish.
For example, $\mu_{1,2}$ is equal to 
the linking number $\mbox{lk}(L_1,L_2)$ of $L_1$ and $L_2$.
Further, 
if $\mu_{i,j}(L)=0$ for any $i,j$, 
then $\mu_{12,3}(L) = \mbox{lk}(L_{12},L_3)$,
where $L_{12}$ denotes the link
which is the intersection of Seifert surfaces of $L_1$ and $L_2$.
In general,
under the vanishing assumption of the lower $\overline\mu$-invariants,
$\overline{\mu}_{12\cdots n-1, n}(L) = \mbox{lk}(L_{12\cdots n-1}, L_n)$
where $L_{12\cdots k}$ (for $k=2,3,\cdots, n-1$)
denotes the link which is the intersection of Seifert surfaces
of $L_{12\cdots k-1}$ and $L_k$.

\begin{prob}[M. Polyak] 
\label{prob.top_pres_muinv}
Milnor's $\overline\mu$-invariants of string links\index{Milnor invariant}
can be defined similarly as above (see \cite{Polyak_sG}).
Find a topological presentation of a $\overline\mu$-invariant
of string links
(not assuming the vanishing of the lower $\overline\mu$-invariants).
\begin{itemize}
\item[\rm(1)]
Show that $\mbox{lk}(L_{12\cdots n-1},L_n)$
is well-defined in an appropriate sense.
\item[\rm(2)]
Identify it with $\overline{\mu}_{12\cdots n-1, n}(L)$.
\end{itemize}
\end{prob}

\subsection{Finite type invariants of virtual knots}

\begin{figure}[htp]
\begin{align*}
\mboxsm{The 6T relation}: & 
\pict{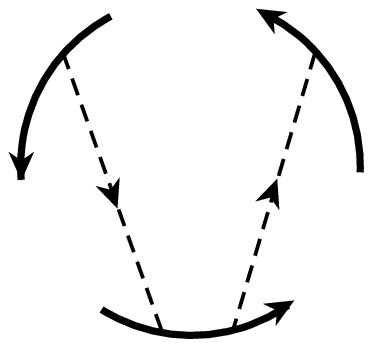}{1.5cm} + \pict{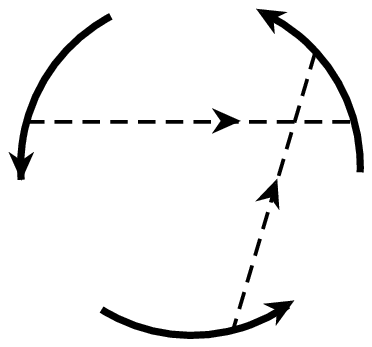}{1.5cm} + \pict{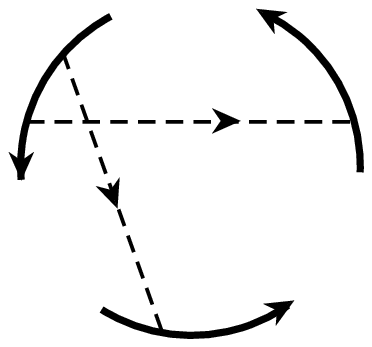}{1.5cm} \\
& \qquad\qquad
= \pict{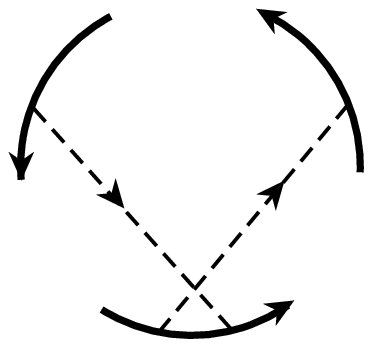}{1.5cm} + \pict{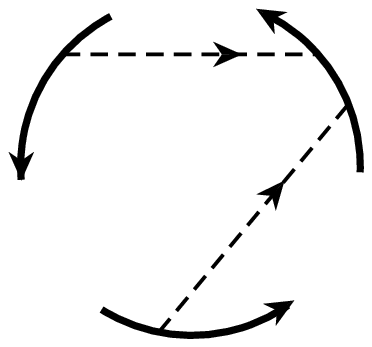}{1.5cm} + \pict{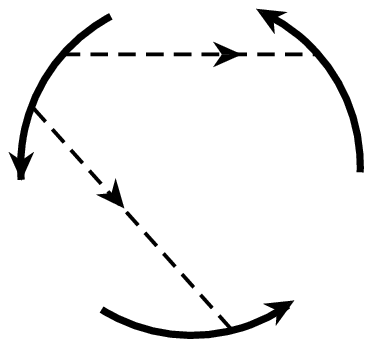}{1.5cm} \\
\mboxsm{The \ora{FI} relation}: &
\pict{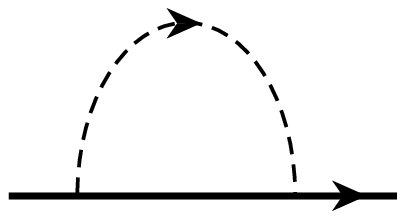}{1.05cm} = 0 = 
\pict{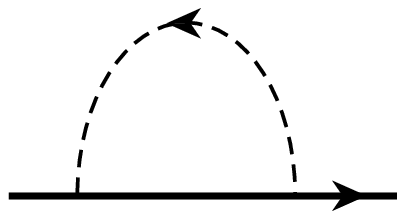}{1.05cm} \\
\mboxsm{The weak \ora{FI} relation}: &
\pict{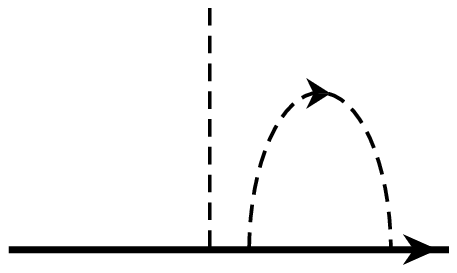}{1.4cm} \!\!=\!\! \pict{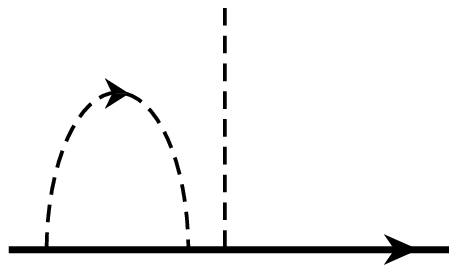}{1.4cm}, \\
& \pict{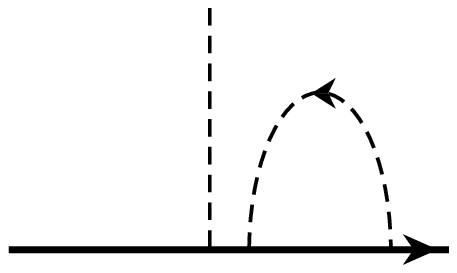}{1.4cm} \!\!=\!\! \pict{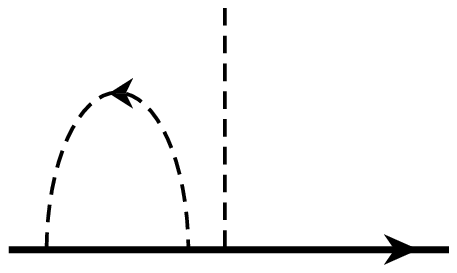}{1.4cm} \\
\mboxsm{The \ora{AS} relation}: & 
\pict{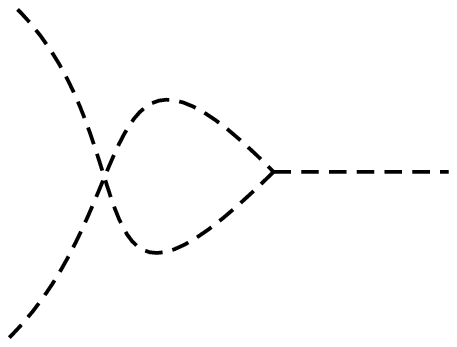}{1.7cm} = - \pict{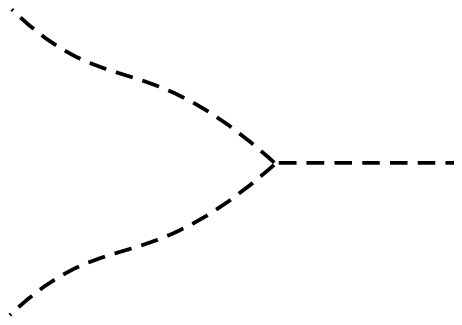}{1.4cm} \\
\mboxsm{The \ora{IHX} relation}: & 
\pict{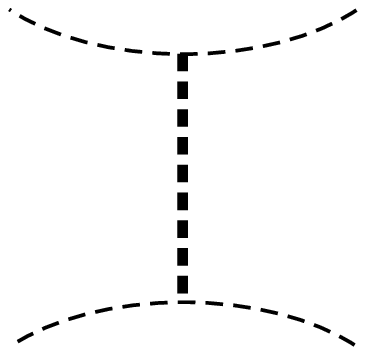}{1.5cm} = \pict{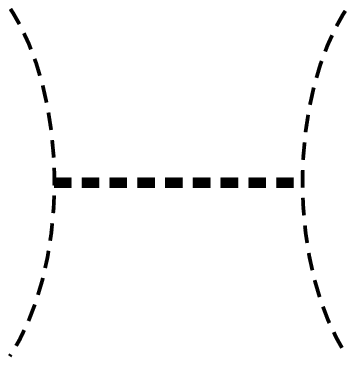}{1.5cm} - \pict{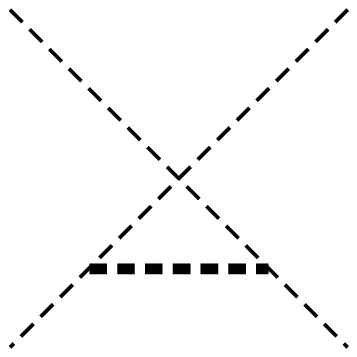}{1.5cm} \\
\mboxsm{The \ora{STU} relation}: & 
\pict{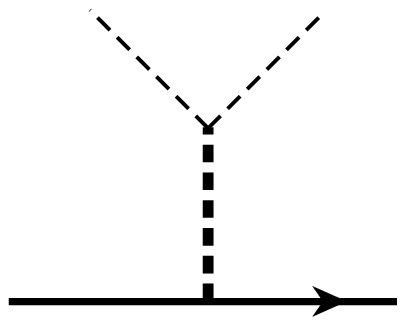}{1.5cm} = \pict{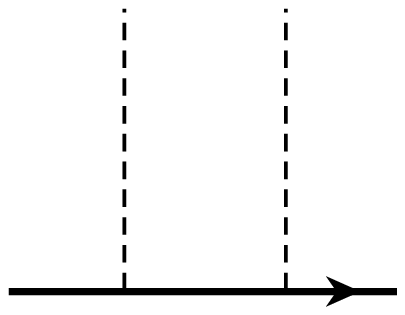}{1.5cm} - \pict{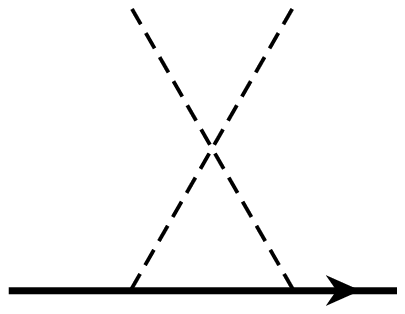}{1.5cm} 
\end{align*}
\caption{\label{fig.oASIHXSTU}
The 6T and the oriented FI, AS, IHX, and STU relations.
Here, a thick dashed line implies the sum of the two orientations,
and corresponding thin dashed lines of pictures in the same formula
have the same (arbitrarily given) orientation.}
\end{figure}

A {\it virtual knot}\index{virtual knot} 
\cite{Kauffman_vk}
is defined by a knot diagram with virtual crossings modulo Reidemeister moves.
Finite type invariants\index{finite type invariant!--- of virtual knots} 
of virtual knots were studied in \cite{GPV},
where their weight systems\index{weight system}  
are defined on 
the space $\vcA(X;R)/\ora{FI}$ of arrow diagrams.
Here an {\it arrow diagram}\index{arrow diagram}
\cite{Polyak_ad}
is a chord diagram with oriented chords,
and $\vcA(X;R)$ denotes the module over a commutative ring $R$
spanned by arrow diagrams on $X$
subject to the 6T relation,
and $\ora{FI}$ denotes the oriented FI relation
(see Figure \ref{fig.oASIHXSTU} for these relations).
It is known \cite{Polyak_ad} that $\vcA(X;R)$ is isomorphic to the module
spanned by acyclic oriented Jacobi diagrams on $X$
subject to the relations
$$
\pict{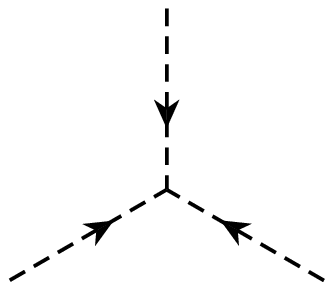}{1.3cm} = 0 = \pict{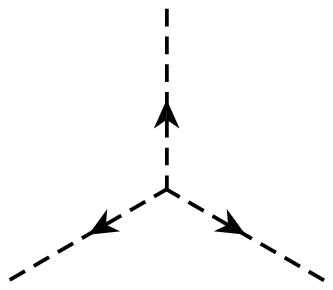}{1.3cm}
$$
and the \ora{AS}, \ora{IHX}, and \ora{STU} relations
(see Figure \ref{fig.oASIHXSTU}).

%

\begin{prob}
Let $I$ denote an oriented interval.
\begin{itemize}
\item[\rm(1)]
Determine the dimensions of\index{arrow diagram!space of ---!dimension of ---}
$\vcA(S^1;\Q)^{(d)}$ and $\vcA(I;\Q)^{(d)}$ for each $d$.
\item[\rm(2)]
Determine the dimensions of 
$\vcA(S^1;\Q)^{(d)}/\ora{FI}$ and $\vcA(I;\Q)^{(d)}/\ora{FI}$ for each $d$.
\item[\rm(3)]
Determine the dimensions of\index{arrow diagram!space of ---!dimension of ---} 
$\vcA(S^1;\Q)^{(d)}/\mbox{(weak \ora{FI})}$ and \newline
$\vcA(I;\Q)^{(d)}/\mbox{(weak \ora{FI})}$ for each $d$.
\end{itemize}
\end{prob}

\begin{rem}
It is shown by elementary computation 
that $\vcA(S^1,\Q)^{(2)}/\mbox{FI} = 0$
and that
$\vcA(I,\Q)^{(2)}/\mbox{FI}$ is a 2-dimensional vector space spanned by
$$
\pict{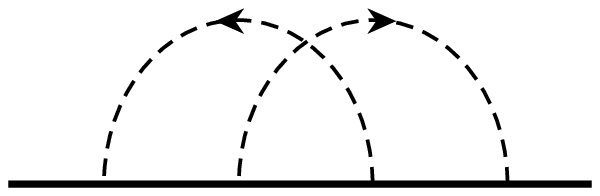}{1.1cm}, \qquad\qquad
\pict{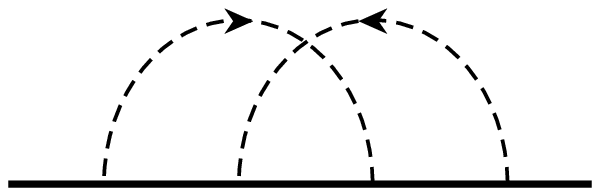}{1.1cm}.
$$
Note that the dimensions of
$\vcA(S^1;\Q)^{(d)}$ and $\vcA(I;\Q)^{(d)}$
differ unlike the unoriented case.
\end{rem}

\begin{rem}
Constructive weight systems\index{weight system!--- from Lie bialgebra}
on $\vcA(X;R)$ can be defined
by using Lie bialgebras
(see, e.g.\ \cite{Drinfeld_ICM86,EiKa96}, for Lie bialgebras),
where the weight systems of the following diagrams
$$
\pict{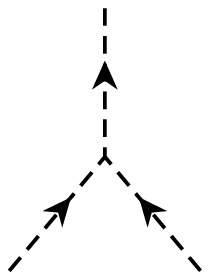}{1.8cm}\ 
\begin{CD} {\frak g} \\ @AAA \\ {\frak g} \otimes {\frak g} \end{CD}
\qquad\qquad
\pict{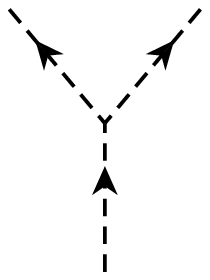}{1.8cm}\ 
\begin{CD} {\frak g} \otimes {\frak g} \\ @AAA \\ {\frak g} \end{CD}
$$
are defined to be the bracket and the co-bracket of 
a Lie bialgebra $\frak g$.
Such weight systems are helpful
when we estimate lower bounds of the dimensions of the spaces $\cA(X;R)$.
\end{rem}

\begin{conj}[M. Polyak] 
\label{conj.inj_A2vA}
The following two maps are injective,\index{arrow diagram!space of ---}
\begin{align*}
& \cA(I)^{(d)} \longrightarrow \vcA(I)^{(d)} \\
& \cA(I)^{(d)}/\mbox{FI} \longrightarrow \vcA(I)^{(d)} / \ora{FI},
\end{align*}
where they are defined by
$$
\pict{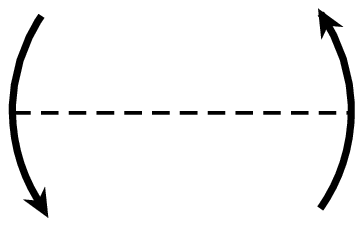}{0.9cm} \longmapsto
\pict{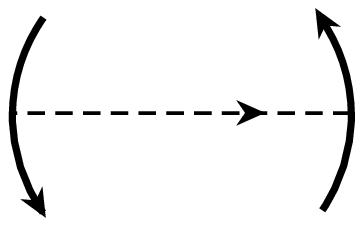}{0.9cm} + \pict{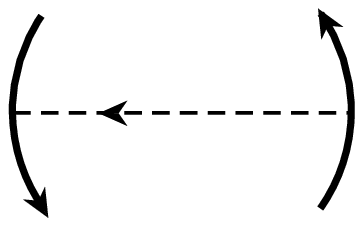}{0.9cm}.
$$
\end{conj}

\begin{rem}
If these maps are injective,
then weight systems on $\cA(I)^{(d)}$ and $\cA(I)^{(d)}/\mbox{FI}$
would be detected by weight systems on 
$\vcA(I)^{(d)}$ and $\vcA(I)^{(d)}/\ora{FI}$;
in other words, the upper rightward map in the following diagram
would be surjective.
$$
\begin{CD}
\left\{ \begin{array}{l} \mbox{degree $d$ weight systems} \\
   \mbox{for long virtual knots} \end{array} \right\}
& @>>> & 
\left\{ \begin{array}{l} \mbox{degree $d$ weight systems} \\
   \mbox{for classical knots} \end{array} \right\} \\
@AAA && @AAA \\
\left\{ \begin{array}{l} \mbox{degree $d$ finite type invariants} \\
   \mbox{for long virtual knots} \end{array} \right\}
& @>>> & 
\left\{ \begin{array}{l} \mbox{degree $d$ Vassiliev invariants} \\
   \mbox{for classical knots} \end{array} \right\}
\end{CD}
$$
Hence, this conjecture follows from Conjecture \ref{conj.Vinv_ckvk} below,
which implies that the lower rightward map in the above diagram is surjective.
\end{rem}

\begin{conj}{\rm\cite{GPV}}\qua 
\label{conj.Vinv_ckvk}
Every Vassiliev invariant of classical knots\index{Vassiliev invariant!extension for virtual knots}
can be extended to 
a finite type invariant\index{finite type invariant!--- of virtual knots}  
of long virtual knots.\index{virtual knot!finite type invariant of ---} 
(See also Problem \ref{prob.Kinv_vk}.)
\end{conj}

\subsection{Finite type invariants derived from local moves}

One aspect of the study of knot invariants
is the study of the set of knots.
A local move and finite type invariants derived from it
might give an approach of this study.

\begin{figure}[ht!]
$$
\begin{array}{ll}
\mboxsm{A crossing change ``$\times$''}: & 
\pict{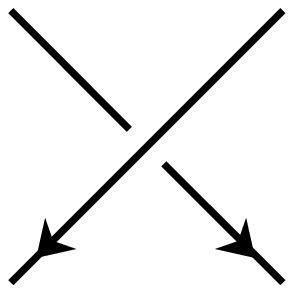}{1.2cm} \longleftrightarrow \pict{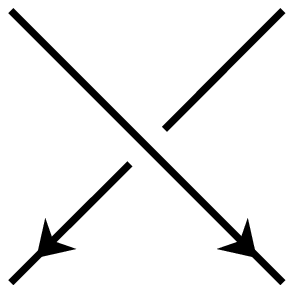}{1.2cm} \\
\mboxsm{A double crossing change ``$\times\!\times$''}: & 
\pict{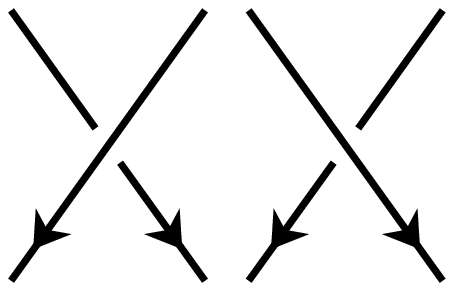}{1.2cm} \longleftrightarrow \pict{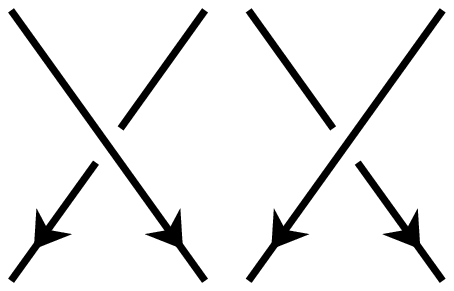}{1.2cm} \\
\mboxsm{A $\#$ move}: & 
\pict{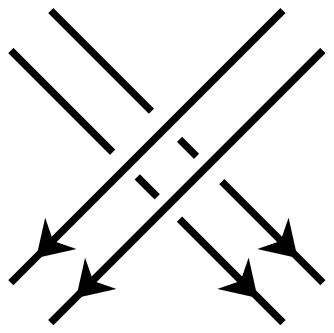}{1.45cm} \longleftrightarrow \pict{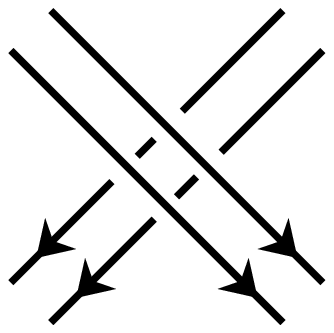}{1.45cm} \\
\mboxsm{A pass move}: & 
\pict{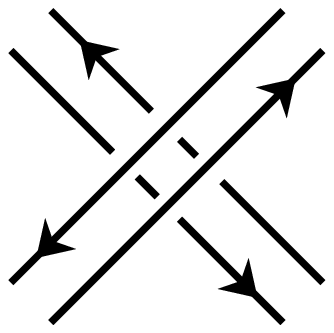}{1.45cm} \longleftrightarrow \pict{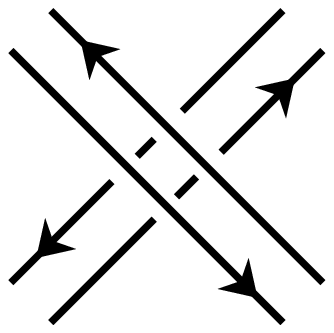}{1.45cm} \\
\mboxsm{A $\Delta$ move}: & 
\pict{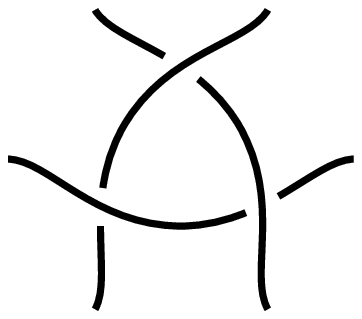}{1.5cm} \longleftrightarrow \pict{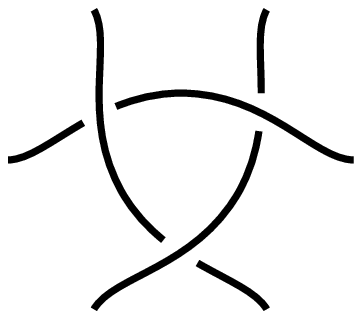}{1.5cm} \\
\mboxsm{A doubled delta move $\DD$}: &
\pict{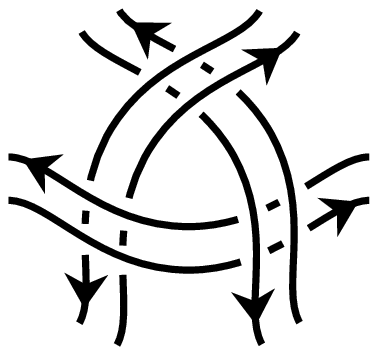}{1.55cm} \longleftrightarrow \pict{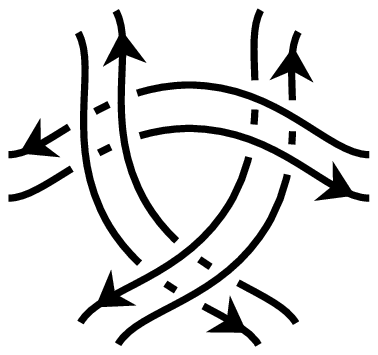}{1.55cm} \\
\mboxsm{An $n$-gon move}: & 
\pict{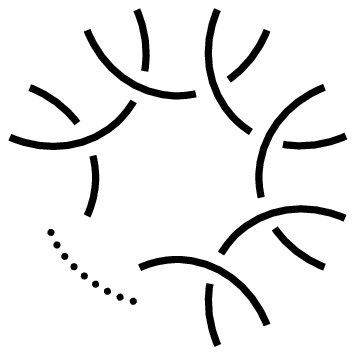}{1.55cm} \longleftrightarrow \pict{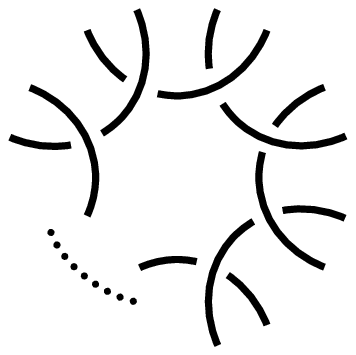}{1.55cm} 
\end{array}
$$
\caption{\label{fig.local_move} Some local moves among oriented knots.
The strands of both sides of a $\Delta$ move and an $n$-gon move
have any orientations
such that corresponding strands from opposite sides of the moves
are oriented in the same way.
Each side of an $n$-gon move has $n$ strands.}\index{local move}\index{local move!$\Delta$ move}\index{local move!doubled delta move} 
\end{figure}

A {\it local move}\index{local move}
 is a move between two knots,
which are identical except for a ball,
where they differ as shown in both sides of a move
in Figure \ref{fig.local_move}.
Let $R$ be a commutative ring with $1$,
and $\K$ the set of isotopy classes of oriented knots, as before.
For a local move $\frak m$,
we define $\cF_d(R\K,{\frak m})$ as follows.
Let $K$ be an oriented knot
with $d$ disjoint balls $B_1,B_2,\cdots,B_d$
such that $K$ is as shown in one side of $\frak m$ in each $B_i$.
For any subset $S \subset \{ 1,2,\cdots,d \}$,
we denote by $K_S$
the knot obtained from $K$ by applying $\frak m$
in each $B_i$ for $i \in S$.
We define $\cF_d(R\K,{\frak m})$ 
to be the submodule of $R \K$ spanned by
\begin{equation}
\label{eq.sumKS}
\sum_S (-1)^{\# S} K_S
\end{equation}
for any $K$ with $d$ balls,
where $\# S$ denotes the number of elements of $S$, and
the sum runs over all subsets $S$ of $\{ 1,2,\cdots,d \}$.
Then, we have a descending series of submodules,
$$
R \K = \cF_0(R \K,{\frak m}) \ \supset\ \cF_1(R \K,{\frak m}) 
\ \supset\ \cF_2(R \K,{\frak m}) \ \supset\ \cdots.
$$
Note that $\cF_d(R\K) = \cF_d(R\K,\times)$
for a crossing change ``$\times$''.
An $R$-homomor\-ph\-ism $v: R \K \to R$
is called 
a {\it finite type invariant of 
$\frak m$-degree $d$},\index{finite type invariant!--- by local move}\index{local move!finite type invariant by ---}
or an $\frak m$ {\it finite type invariant of degree $d$},
if $v |_{\cF_{d+1}(R\K,{\frak m})} = 0$.

It is a fundamental problem of finite type invariants
to calculate the corresponding graded spaces,
which would enable us to identify finite type invariants in some sense.

\begin{prob}
Calculate\index{finite type invariant!--- by local move}\index{local move!finite type invariant by ---}
$\cF_d(\Z\K,{\frak m})/\cF_{d+1}(\Z\K,{\frak m})$,
letting $\frak m$ be a local move such as
\begin{itemize}
\item[\rm(1)]
a $\#$ move,
\item[\rm(2)]
a pass move,
\item[\rm(3)]
a $\Delta$ move,\index{local move!$\Delta$ move}
\item[\rm(4)]
an $n$-gon move.
\end{itemize}
\end{prob}

\begin{rem}
It is known that
crossing change, double crossing change,
$\#$ move (see \cite{Murakami_sharp}),
$\Delta$ move (see \cite{MuNa_Delta_move}),
$n$-gon move (see \cite{Aida_ngon})
are unknotting operations,
i.e.\
any oriented knot can be related to the trivial knot
by a sequence of isotopies and each of these moves.
Hence,
$\cF_0(\Z\K,{\frak m})/\cF_1(\Z\K,{\frak m}) \cong \Z$
for these moves $\frak m$.

It is known \cite{Kauffman_formal} that
Arf invariant gives the bijection
$$
\{ \mbox{knots} \}/\mbox{(pass move)} \longrightarrow \Z/2\Z.
$$
Hence,
$\cF_0(\Z\K,\mbox{pass move})/\cF_1(\Z\K,\mbox{pass move}) \cong \Z \oplus \Z$.
\end{rem}

\begin{rem}
$\Delta$ finite type invariants were introduced in \cite{Mellor_ft};
see also \cite{Stanford_delta}.
\end{rem}

\begin{rem}[\rm (K. Habiro)] 
The following relations hold,
\begin{align*}
\cF_{2d}(\Z\K,\times) \ & \supset \ 
\cF_{d}(\Z\K,\Delta) \ \supset \ 
\cF_{3d}(\Z\K,\times), \\
\cF_{d}(\Z\K,\times) \ & \supset \ 
\cF_{d}(\Z\K,\#) \ \supset \ 
\cF_{d}(\Z\K,\Delta).
\end{align*}
These relations imply that
$\frak m$ finite type invariants are Vassiliev invariants, and
Vassiliev invariants are $\frak m$ finite type invariants,
for $\frak m = \#, \Delta$.
Further, the rank of 
$\cF_d(\Z\K,{\frak m})/\cF_{d+1}(\Z\K,{\frak m})$ is finite 
for these $\frak m$.
\end{rem}

\begin{rem}
For the Kontsevich invariant
$\hZ$ (introduced in Chapter \ref{sec.Kinv}),
we have that
$$
\hZ \Big( \!\!\! \pict{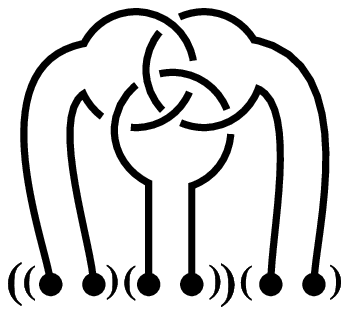}{1.8cm} \!\!\! \Big) - 
\hZ \Big( \pict{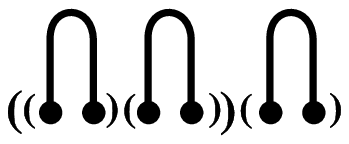}{0.9cm} \Big) 
= \pict{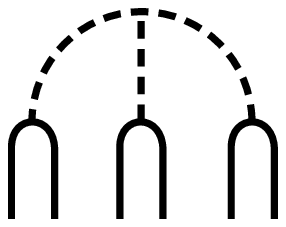}{1.5cm}
+ \left(
\begin{array}{l} \mboxsm{terms of} \\ \mboxsm{higher degrees} \end{array}
\right),
$$
where two tangles in the left hand side are related by a $\Delta$ move.
Hence, the image of 
$$
\cF_d(\Q\K,\Delta) \longrightarrow
\cF_{2d}(\Q\K) \longrightarrow
\cF_{2d}(\Q\K)/\cF_{2d+1}(\Q\K) \cong \cA(S^1;\Q)^{(2d)}
$$
is equal to the subspace of $\cA(S^1;\Q)^{(2d)}$
spanned by Jacobi diagrams on $S^1$
whose uni-trivalent graphs are
disjoint unions of $d$ dashed Y graphs.
\end{rem}

\begin{rem}
Finite type invariants derived from a double crossing change
were introduced in \cite{Appleboim_f},
to study finite type invariants of links with a fixed linking matrix.
For knots,
they are equal to Vassiliev invariants,
that is,
$\cF_d(\Z\K;\times\!\times) = \cF_d(\Z\K,\times)$.
\end{rem}

\begin{figure}[ht!]
$$
\mboxsm{The surgery on} \!\!\! \!\!\pict{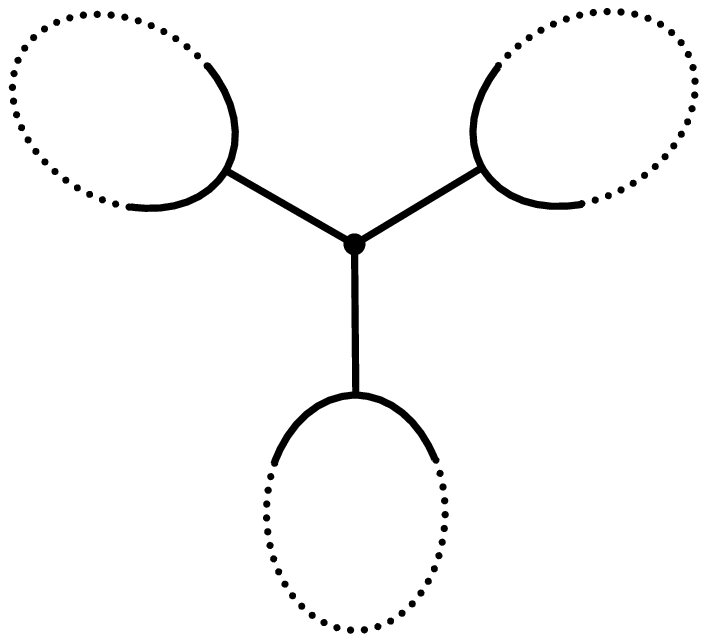}{2.7cm} \!\!\! \!\!
\mboxsm{ is defined to be }
\mboxsm{the surgery along} \!\!\! \!\!\pict{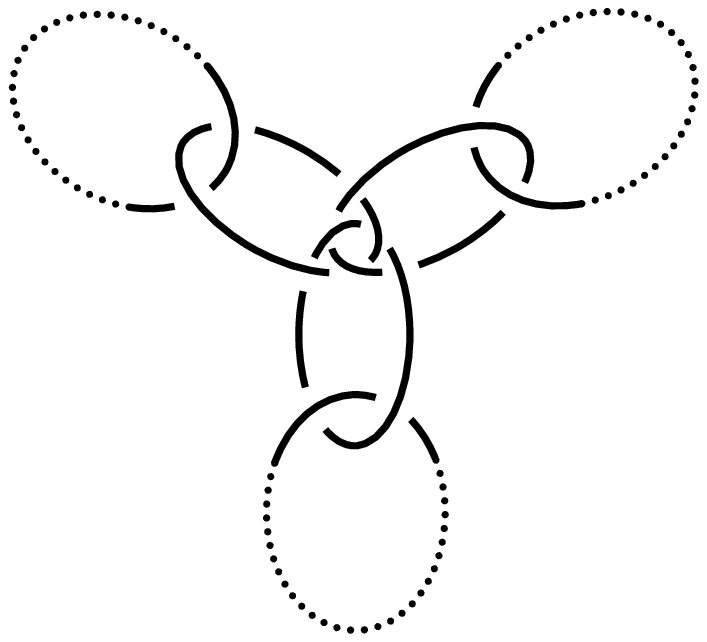}{2.7cm}\!\!\! \!\!.
$$
\caption{\label{fig.surgeryY}
Definition of the surgery on a Y graph.
Dotted lines imply strands possibly knotting and linking.
Three circles (partially dotted) in the left picture
are called {\it leaves}.}
\end{figure}

\noindent
{{\namae}(Y. Ohyama)}\quad
In the case all arcs in a $\Delta$ move are contained
in the same component, it is called a {\it self $\Delta$ move}.
If two links can be transformed into each other by a
finite sequence of self $\Delta$ moves,\index{local move!$\Delta$ move} 
they are said to be
$\Delta$ link homotopic.

\begin{prob}[Y. Ohyama] 
\label{prob.D_link_htpc}
Find necessary and sufficient conditions for two
$\mu$-component links ($\mu>2$) to be $\Delta$ link homotopic.
\end{prob}

\begin{rem}[\rm (Y. Ohyama)] 
For a $\mu$-component link $K=K_{1}\cup K_{2}\cup\ldots\cup K_{\mu}$,
let $\delta_{1}=a_{\mu -1}(K)$ and
$\delta_{2}=a_{\mu +1}(K)-a_{\mu -1}(K)\times (\sum_{i=1}^{\mu}
a_{2}(K_{i})$
for the coefficient $a_{i}(K)$ of the term $z^{i}$ in the Conway
polynomial of $K$.

It is known \cite{Matveev87,MuNa_Delta_move} that
two knots (or links) can be transformed into each other by a finite
sequence of $\Delta$ moves if and only if they have the same number
of components, and, for properly chosen orders and orientations,
they have the same linking numbers between the corresponding
components.
In particular, if two links are $\Delta$ link homotopic,
then their $\delta_1$ coincide.
Further, it is known \cite{ohyama.Nakanishi} that
if two $\mu$ component links are $\Delta$ link homotopic, then
their $\delta_{2}$ coincide.
These are necessary conditions of this problem.

Moreover, for 2-component links, a pair of $\delta_{1}$ and $\delta_{2}$ 
is a faithful invariant of $\Delta$ link homotopy. 
Namely, for two 2-component links, they are $\Delta$ link homotopic
if and only if their $\delta_{1}$ and $\delta_{2}$ coincide
\cite{ohyama.NaOh}.
This gives a required condition of this problem for 2-components links.
\end{rem}

\subsection{Loop finite type invariants}

The {\it loop-degree} of a Jacobi diagram on $S^1$
is defined to be half of the number given by
the number of trivalent vertices minus the number of univalent vertices
of the uni-trivalent graph of the Jacobi diagram.
The filtration of $\cA(S^1)$ given by loop-degrees
is related to a filtration of $\Q\K$
through the Kontsevich invariant.
The theory of the corresponding filtration in $\Z(\MK)$ (given below)
is developed in \cite{GR00}
(noting that this definition also appears
in the September 1999 version of \cite{Kricker_b}).

We denote by $\MK$
the set of pairs $(M,K)$ such that
$M$ is an integral homology 3-sphere and 
$K$ is an oriented knot in $M$.
Consider a move between two pairs $(M,K)$ and $(M',K')$ in $\MK$
such that $(M',K')$ is obtained from $(M,K)$
by surgery on a Y graph (see Figure \ref{fig.surgeryY})
embedded in $M-K$ whose leaves have linking number zero with $K$.
We call this move a {\it loop move}.\index{local move!loop move}
Finite type invariants of degree $d$ derived from a loop move
by (\ref{eq.sumKS})
are called {\it loop finite type invariants of degree $d$},
or {\it finite type invariants of 
loop-degree $d$}.\index{finite type invariant!loop ---}\index{local move!finite type invariant by ---!loop finite type invariant}
We denote the corresponding submodule of $R(\MK)$ by
$\cF_l( \Z (\MK); \mbox{loop})$.

A doubled delta move\index{local move!doubled delta move} 
$\DD$ (see Figure \ref{fig.local_move})
was introduced by Naik-Stanford \cite{NaSt99}
as a move characterizing S-equivalence classes;
two knots are {\it S-equivalent} 
if they are indistinguishable by Seifert matrices.
A doubled delta move $\DD$
can be presented by a surgery on such a Y graph as above.
Thus, 
we have the map
$\cF_l( R \K; \DD) \to \cF_l( R (\MK); \mbox{loop} )$,
taking a knot $K$ to $(S^3,K) \in \MK$.
Hence, a loop finite type invariant\index{finite type invariant!loop ---} 
gives a $\DD$ finite type invariant.

\begin{prob} 
\label{prob.grad_loop_fil}
Let $R$ be a commutative ring with $1$, say, $\Z$ or $\Q$.
\begin{itemize}
\item[\rm(1)]
Describe the spaces
$\cF_l(R (\MK); \mbox{loop}) / \cF_{l+1}(R (\MK); \mbox{loop})$.
\item[\rm(2)]
Describe the spaces
$\cF_l(R \K; \DD) / \cF_{l+1}(R \K; \DD)$.
\item[\rm(3)]
Describe the image of the above map 
$\cF_l( R \K; \DD) \to \cF_l( R (\MK); \mbox{loop} )$.
\end{itemize}
\end{prob}

\begin{rem}[\rm (A. Kricker)] 
It follows by a short argument from \cite{NaSt99} and \cite{Matveev87}
that the following map taking a pair $(M,K)$ to a Seifert matrix of $K$ in $M$
is bijective,\index{local move!loop move}
\begin{equation}
\label{eq.MK/loop}
\MK / \mbox{(loop move)} 
\stackrel{=}{\longrightarrow}
\{ \mbox{S-equivalence classes of Seifert matrices} \}.
\end{equation}
(This implies that
$K$ and $K'$ are related by a sequence of doubled delta moves
if and only if $(S^3,K)$ and $(S^3,K')$ are related by
a sequence of loop moves.)
Hence, 
$\cF_0(\Z(\MK);\mbox{loop}) / \cF_1(\Z(\MK);\mbox{loop})$
is isomorphic to the module over $\Z$
freely spanned by S-equivalence classes.
Moreover, by (\ref{eq.MK/loop}), we have that $\Z(\MK) = \oplus_s \Z(\MK_s)$, 
where the sum runs over all S-equivalence classes $s$.
Further,
{\small$$
\cF_l(\Z(\MK);\mbox{loop}) / \cF_{l+1}(\Z(\MK);\mbox{loop}) 
=\bigoplus_s \cF_l(\Z(\MK_s);\mbox{loop}) / \cF_{l+1}(\Z(\MK_s);\mbox{loop}).
$$}%
Hence, the problem (1) splits into problems
of describing the direct summands on the right hand sides:
describe the spaces
$$\cF_l(\Z(\MK_s);\mbox{loop}) / \cF_{l+1}(\Z(\MK_s);\mbox{loop})$$
for each S-equivalence class $s$.
For the S-equivalence class $u$ including the unknot,
$\cF_l(\Q(\MK_u);\mbox{loop}) / \cF_{l+1}(\Q(\MK_u);\mbox{loop})$
is isomorphic to
$\cA^{\Z[t^{\pm1}]}(\emptyset;\! \Q)^{({\rm loop}\,l)}$
by the map (\ref{eq.lAtoB2}) of the loop expansion of the Kontsevich invariant
(see also \cite{GR00});
for the definition of the space 
$\cA^{\Z[t^{\pm1}]}(\emptyset; \Q)^{({\rm loop}\ l)}$
see Section \ref{sec.loopKinv}.
\end{rem}

\begin{rem}
A surgery on a Y graph in the definition of loop finite type invariants
lifts to a surgery of the infinite cyclic cover of the knot complement,
which does not change its homology.
Hence, it is shown, topologically, that
all coefficients of the Alexander polynomial are 
finite type invariants of loop-degree $0$.

It follows that
all coefficients of the Alexander polynomial are 
finite type invariants of $\DD$-degree $0$.
It can also be shown from the fact that
the Alexander polynomial can be defined by the Seifert matrix of a knot,
which is unchanged by finite type invariants of $\DD$-degree 0
as shown in \cite{NaSt99}.

The Alexander polynomial is universal
among Vassiliev invariants
which are of finite type of $\DD$-degree $0$;
more precisely,
$\log \Delta_K(e^\hbar)$
as a power series of $\hbar$
is universal
among $\Q$-valued primitive Vassiliev invariants
which are of finite type of $\DD$-degree $0$.
An equivalent statement has been shown in \cite{MO_ftS},
using Vassiliev invariants of S-equivalence classes of Seifert matrices.
\end{rem}

\begin{rem}
As shown in \cite{NaSt99} we have a bijection,
$$
\{ \mbox{knots} \}/ \DD 
\stackrel{=}{\longrightarrow}
\{ \mbox{S-equivalence classes} \},
$$
by taking a knot to its S-equivalence class.
Hence,
$\cF_0(\Z\K;\DD) / \cF_1(\Z\K;\DD)$
is isomorphic to the module over $\Z$
freely spanned by S-equivalence classes.
\end{rem}

\begin{rem}[\rm (A. Kricker)] 
The dual space of
$$
\frac{\cF_l(\K\otimes\Q;\DD) \cap \cF_d(\K\otimes\Q;\times)}
{\big((\cF_{l+1}(\K\otimes\Q;\DD) \cap \cF_d(\K\otimes\Q;\times)\big) 
+ \big(\cF_l(\K\otimes\Q;\DD) \cap \cF_{d+1}(\K\otimes\Q;\times)\big)}
$$
is isomorphic to 
the subspace of $\cB$ spanned by connected uni-trivalent graphs
of degree $d$ and of loop-degree $l$,
i.e.\
the space $\cB\conn^{(d,d-l)}$
in the notation given in a remark in Problem \ref{prob.dimV}.
\end{rem}

(A. Kricker)\quad
Let $\MK$ denote
the set of pairs $(M,K)$ such that
$M$ is an integral homology 3-sphere and 
$K$ is an oriented knot in $M$, as before.
A {\it mod $p$ loop move}\index{local move!mod $p$ loop move} 
in $\MK$ is defined to be 
a surgery on a Y graph (see Figure \ref{fig.surgeryY})
such that each leaf has linking number 0 modulo $p$ with the knot.
We consider the question:
what are the mod $p$ loop move equivalence classes of knots?

To state the conjecture below,
we give some notation.
Consider a pair $(M,K)$ 
of an integral homology 3-sphere $M$ and a knot $K$ in $M$.
Let $\Sigma^p_{(M,K)}$ be the $p$-fold branched cyclic cover
of $(M,K)$, and assume that $\Sigma^p_{(M,K)}$ is 
a rational homology 3-sphere.
Observe that there is an action of $\Z/p\Z$ on the
homology group $H_1(\Sigma^p_{(M,K)};\Z)$ 
(induced from the covering transformations). 
Observe also that the linking pairing on the torsion of 
$H_1(\Sigma^p_{(M,K)};\Z)$ (which is the whole group)
is invariant under the action of $\Z/p\Z$.
Here, the linking pairing on the torsion of $H_1(N;\Z)$
of a 3-manifold $N$ is the map
$\mbox{Tor}\big( H_1(N;\Z) \big) \otimes \mbox{Tor}\big( H_1(N;\Z) \big) 
\to \Q/\Z$
taking $\alpha \otimes \beta$ to
$1/n$ times the algebraic intersection of $F$ and $\beta$,
where $F$ is a compact surface bounding $n \alpha$ 
for some non-zero integer $n$.

\begin{conj}[A. Kricker] 
\label{conj.p_loop_move}
Take $(M_1,K_1)$ and $(M_2,K_2)$ of the above sort. 
Then, there exists a $(\Z/p\Z)$-equivariant isomorphism
$\phi: H_1(\Sigma^p_{(M_1,K_1)};\Z)$ $\to H_1(\Sigma^p_{(M_2,K_2)};\Z)$
preserving the linking pairing\index{linking pairing}
if and only if
$(M_1,K_1)$ is equivalent to $(M_2,K_2)$
by a finite sequence of mod $p$ loop moves.
\end{conj}

\begin{rem}[\rm (A. Kricker)] 
The case of $p=1$ would recover Matveev's theorem \cite{Matveev87}:
two closed 3-manifolds $M$ and $N$ are equivalent 
by a finite sequence of surgeries on Y graphs
if and only if there is an isomorphism
$H_1(M;\Z) \to H_1(N;\Z)$
preserving the linking pairing on the torsion.

Also, the limit as $p \to \infty$
should recover a theorem due to Naik-Stanford \cite{NaSt99}:
two knots are equivalent by a finite sequence of loop moves
if and only if they have isometric Blanchfield pairings. 
(Recall that the Blanchfield pairing is the equivariant linking pairing
on the universal cyclic cover.)
\end{rem}

\subsection{Goussarov-Habiro theory for knots}

Related to Vassiliev invariants of knots,
equivalence relations among knots have been studied by 
Goussarov \cite{Gou_kg,Gou_fn} and Habiro \cite{Habiro_GT},
which is called the Goussarov-Habiro theory for knots.\index{Goussarov-Habiro theory!--- for knots}
These equivalence relations are helpful for us
to study structures of the set of knots.

The\index{equivalence relation!$C_d$-equivalence}\index{equivalence relation!$d$-equivalence}  
{\it $C_d$-equivalence\/}\footnote{\footnotesize 
The $C_d$-equivalence is also called
the {\it $(d-1)$-equivalence} (due to Goussarov) in some literatures.}
($d=1,2,3,\cdots$)
among oriented knots is the equivalence relation
generated by either of the following relations,
\begin{itemize}
\item[\rm(1)] 
$C_{d}$-move,\index{local move!$C_d$-move} 
i.e.\
surgery along a tree clasper with $d$ trivalent vertices
whose leaves are disc-leaves \cite{Habiro_GT},
\item[\rm(2)]
relation on a certain collection of $d$ crossing changes
(Goussarov's $(d-1)$-equivalence) \cite{Goussarov_CJ,Goussarov_e,Gou_fn},
\item[\rm(3)]
surgery by an element in the $d$th group in the lower central series 
of pure braid group \cite{Stanford_Vpb},
\item[\rm(4)]
capped grope cobordism of class $d$ \cite{CoTe00}.
\end{itemize}
It is known that
these relations generate the same equivalence relation among knots.
The $C_d$-equivalence is defined among links, string links, $\cdots$,
in the same way.

It is known \cite{Habiro_GT} that
there exists a natural surjective homomorphism
\begin{equation}
\label{eq.AtoKsim}
\cA(S^1;\Z)\conn^{(d)} \longrightarrow 
\{ K \ \simsubC{d} O \} / \simsubC{d+1}
\end{equation}
such that the tensor product of this map and $\Q$ is an isomorphism,
where $O$ denotes the trivial knot.
In particular, 
$\{ K \ \simsubC{d} O \} / \simsubC{d+1}$ forms
an abelian group
with respect to the connected sum of knots,
and hence, so does $\{ \mbox{knots} \} / \simsubC{d+1}$.

\begin{conj}
The map (\ref{eq.AtoKsim}) is an isomorphism.
\end{conj}

This conjecture might be reduced to Conjecture \ref{conj.AS1_tf} and
the following conjecture.

\begin{conj}
$\{ K \ \simsubC{d} O \} / \simsubC{d+1}$ is 
torsion free\index{torsion!--- free}  
for each $d$.
\end{conj}

\begin{rem}
Conjecture \ref{conj.AS1_tf} implies this conjecture,
since the surjective homomorphism (\ref{eq.AtoKsim}) 
gives a $\Q$-isomorphism.
\end{rem}

It is known \cite{Gou_kg,Stanford_Vpb,Habiro_GT} that
two knots $K$ and $K'$ are $C_d$-equivalent if and only if
$v(K) = v(K')$ for any $A$-valued Vassiliev invariant
$v$ of degree $< d$
for any abelian group $A$.
In fact, a natural quotient map
$\{ \mbox{knots} \} \to \{ \mbox{knots} \} / \simsubC{d}$
is a Vassiliev invariant of degree $< d$, 
which classifies $C_d$-equivalence classes of knots.

\begin{conj}[K. Habiro \cite{Habiro_GT}, {\rm see also} {\cite[``Theorem 5'']{Gou_fn}}] 
\label{conj.Cd_string_link}
Two $m$-strand string links $L$ and $L'$ are 
$C_d$-equivalent\index{equivalence relation!$C_d$-equivalence}  
if and only if $v(L) = v(L')$
for any $A$-valued 
finite type invariant\index{finite type invariant!--- of string links} 
$v$ of degree $< d$
for any abelian group $A$.
\end{conj}

\begin{rem}[\rm (M. Polyak)]
The corresponding assertion for links does not hold;
note that $\{ \mbox{links} \} / \simsubC{d}$ 
does not (naturally) form a group.
Recall that $\{ \mbox{knots} \} / \simsubC{d}$ forms an abelian group,
which guarantees the corresponding assertion for knots,
as mentioned above.
The set $\{ \mbox{$m$-strand string links} \} / \simsubC{d}$ 
forms a group with respect to the composition of string links,
though it is not abelian.
\end{rem}

\begin{prob}[M. Polyak] 
\label{prob.GH_vknot}
Establish the Goussarov-Habiro theory\index{Goussarov-Habiro theory!--- for virtual knots} 
for virtual knots.\index{virtual knot!Goussarov-Habiro theory for ---} 
\end{prob}

\begin{rem}
Polyak suggested that the following moves,
$$
\pict{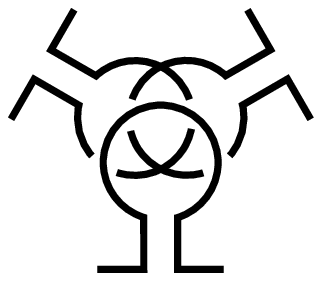}{1.6cm} \longleftrightarrow
\pict{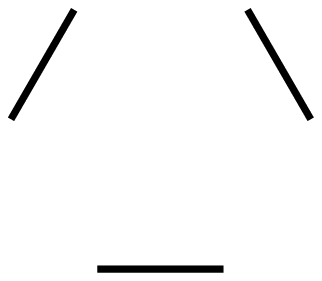}{1.6cm} \longleftrightarrow
\pict{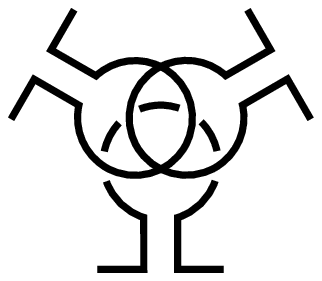}{1.6cm},
$$
(which appear in \cite{GPV})
might play a similar role as the $C_2$-move plays among knots.
They are related to the following diagrams respectively,
$$
\pict{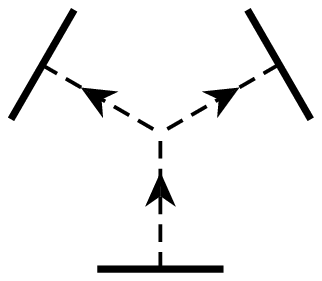}{1.6cm}, \qquad
\pict{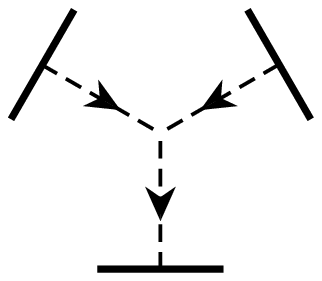}{1.6cm}.
$$
Further, Habiro suggested that the move,
$$
\pict{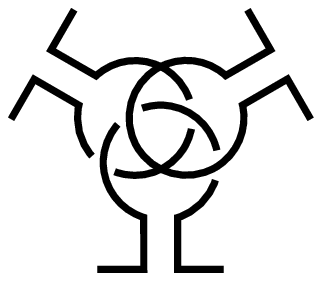}{1.6cm} \longleftrightarrow
\pict{w/vkmove2}{1.6cm},
$$
should be added to the above moves.
It is a problem to define a sequence of equivalence relations
among virtual knots (an extension of the $C_d$-move)
which induces finite type invariants of virtual knots.
Are there surjective homomorphisms from certain modules
of arrow graphs 
(oriented Jacobi diagrams)\index{Jacobi diagram!oriented ---|see {arrow diagram}}
to the graded sets derived from such equivalence relations?
\end{rem}

\noindent
{{\namae}(K. Habiro)}\quad
We denote by $\MK$
the set of pairs $(M,K)$ such that
$M$ is an integral homology 3-sphere and 
$K$ is an oriented knot in $M$.
The {\it $H\!L_d$-equivalence}\index{equivalence relation!homology $d$-loop equivalence}   
({\it homology $d$-loop equivalence})
in $\MK$ is the equivalence relation
generated by either of the following relations,
\begin{itemize}
\item[\rm(1)]
surgery on a tree clasper with $d$ trivalent vertices 
with null-homologous leaves, 
\item[\rm(2)]
surgery on a graph clasper with $d$ trivalent vertices 
with null-homologous leaves,
\item[\rm(3)]
surgery by an element of the $d$th lower central series subgroup of
the Torelli group of compact connected surfaces 
embedded in a null-homologous way.
\end{itemize}
Here, ``null-homologous'' means null-homologous in knots complements.
These relations generate the same equivalence relation in $\MK$.

\begin{prob}[K. Habiro] 
\label{prob.MK_HL}
Describe the abelian group \newline
$\{ (M,K) \simsubHL{d} (S^3,\mbox{unknot}) \} / \! \simsubHL{d+1}$
for each $d$.
\end{prob}

\begin{rem}[\rm (K. Habiro)] 
Two pairs $(M,K)$ and $(M',K')$ in $\MK$
are $H\!L_d$-equivalent
if and only if 
$v(M,K) = v(M',K')$
for any $A$-valued 
loop finite type invariant\index{finite type invariant!loop ---} 
$v$ of loop degree $<d$
for any abelian group $A$.
Thus, the $H\!L$-equivalence gives 
the Goussarov-Habiro theory\index{Goussarov-Habiro theory!--- for loop finite type invariant} 
for loop finite type invariants.

The {\it homotopy $d$-loop equivalence}\index{equivalence relation!homotopy $d$-loop equivalence}  
is defined by
using ``null-homotopic leaves''
instead of ``null-homologous leaves''
in the definition of the $H\!L_d$-equivalence.
These equivalences might be related to
the rational $Z$ invariant $Z^{rat}$.
The homotopy loop equivalence relates
(\ZHS, boundary link) to (\ZHS, boundary link).
A high loop-degree part of $Z^{rat}$ might
be invariant under the homotopy loop equivalence.

The quotient set $\MK/\simsubHL{1}$ can be identified with
the commutative monoid of S-equivalence classes of Seifert matrices.
(See a remark of Problem \ref{prob.grad_loop_fil}.)

Define the equivalence relation $H\!L_d'$ among knots in $S^3$
to be the equivalence relation
generated by
surgery on a tree clasper with $d$ trivalent vertices
with null-homologous leaves in the complement of a knot
such that at least one leaf bounds a disc
with zero intersection number with the knot.
Then, there exists a split exact sequence,
$$
\{ \mbox{knots in $S^3$} \} / \!\!\underset{{}^{H\!L'_{d}}}{\sim} 
\longrightarrow
\MK / \simsubHL{d} \longrightarrow \{ \ZHSs \} / \simsubY{d},
$$
where the first map takes a knot $K$ to $(S^3,K)$
and the second map is the map forgetting knots.

A refinement of Problem \ref{prob.MK_HL} is to consider
the graded sets of the double sequence given by 
the $C_d$-equivalence and the $H\!L_n$-equivalence.
\end{rem}

\subsection{Other problems}

\noindent
{{\namae}(D. Bar-Natan)}\footnote{\footnotesize This part is a quotation from 
\newline\url{http://www.ma.huji.ac.il/~drorbn/Misc/Nullstellensatz/}
}\quad
Is there a Hilbert's Nullstellensatz for finite type invariants of links? 

Let $k$ be an algebraically closed field 
and let $I$ be an ideal in the polynomial ring $k[x_1, \cdots,x_n]$. 
The Hilbert Nullstellensatz\index{Nullstellensatz}
(see e.g.\ \cite{Eisenbud}) says
that the ideal of polynomials in $k[x_1,\cdots,x_n]$ 
that vanish on the variety 
defined by the common zeros of all polynomials in $I$ is the radical of $I$. 

\begin{prob}[D. Bar-Natan] 
\label{prob.HN_Vinv}
Is there a similar statement for 
finite type invariants\index{finite type invariant!--- of links} 
of links? 
Let $I$ be an ideal in the algebra $V$ of finite type invariants of links.
Let $Z$ be the set of links that are annihilated by all members of $I$, 
and let $J$ be the ideal in $V$ of all invariants that vanish on $Z$. 
Clearly, $J$ always contains the radical of $I$. 
Are they always equal? 
\end{prob}

\begin{exm}[{\rm (D. Bar-Natan)}] 
Let $I$ be the ideal generated by linking numbers. 
In this case, $Z$ is the set of algebraically split links. 
Is it true that every finite type invariant 
that vanishes on algebraically split links 
is a sum of multiples of linking numbers? 
I believe it is true, and I believe it follows
from the results of Appleboim \cite{Appleboim_f}, 
but I'm afraid Appleboim's paper is incomplete and 
while I believe it I cannot vouch for its validity. 
\end{exm}

\begin{rem}[{\rm (D. Bar-Natan)}] 
One may also ask, ``what is the Zariski closure of a given set of links?''. 
I believe that in the light of the paragraphs above the
meaning of this question should be clear. 
I know of at least one interesting example: 
In \cite{Ng} Ng shows that the Zariski closure of the set of
ribbon knots is the set of knots whose Arf invariant vanishes. 
\end{rem}

\noindent
{\bf Is the similarity index of two different knots finite?} \newline
{{\namae}(M.-J. Jeong, C.-Y. Park)}

\noindent
K. Habiro and T. Stanford independently showed that  for each positive
integer $n$, two knots $K$ and $L$ have the same values for any
Vassiliev invariants of type $<n$ if and only if they are $LCS_n$-equivalent. 
Y. Ohyama introduced triviality index of knots and
K. Taniyama extended this to the similarity index of links;
see \cite{Ohyama_sim}.
Ohyama showed that if two knots are $n$-similar then they have the same value
for any Vassiliev invariants of type $<n$. It is not difficult to see
that two knots are $n$-similar if they are $LCS_n$-equivalent.
D. Bar-Natan gave a problem whether Vassiliev invariants can distinguish
all of knots or not. This problem is equivalent to the problem, whether
the similarity index of any two different knots have finite similarity
index. We will give a new criterion to calculate the similarity index of
knots and, based on this, raise problems to calculate similarity index.
For example, for two given knots, which knot invariants will give the
best upper bound to calculate the similarity index of knots, along our
above new result? As a partial problem, can we show that  the triviality
index of a non-trivial knot is finite by using our results?

\medskip
\noindent
{\bf Polynomial invariants and Vassiliev invariants} \newline
{{\namae}(M.-J. Jeong, C.-Y. Park)}\index{Vassiliev invariant!--- and polynomial invariants}

\noindent
In 1993, J.S. Birman and X.-S. Lin \cite{BiLi_V} showed that, after a
suitable change of variables, each coefficient of the Jones, HOMFLY and
Kauffman polynomial is a Vassiliev invariant. So we can obtain various
Vassiliev invariants from the derivatives of knot polynomials.

In 2001, by using some specific kinds of tangles, we gave two operations
$~\bar{}~$ and $^*$ operations to get new polynomial invariants from a
given Vassiliev invariant. These new polynomial invariants are also
Vassiliev invariants. So we can obtain various Vassiliev invariants from
the coefficients of these polynomial invariants.

Let $V_n$ be the space of Vassiliev invariants of degrees $\leq n$. For
$A_n \subset V_n$, let $(A_n)$ be the set of Vassiliev invariants
obtained from $A_n$ by using finite numbers of $~\bar{}~$ and $^*$
operations repeatedly.

\begin{prob}[M.-J. Jeong, C.-Y. Park] 
\label{quest.An_Vn}
Find a minimal finite subset $A_n$ of $V_n$ such that span$(A_n)=V_n$.
\end{prob}

\newpage

\section{The Kontsevich invariant}
\label{sec.Kinv}

The framed Kontsevich invariant\index{Kontsevich invariant}
$\hZ(L) \in \cA(\sqcup^l S^1;\Q)$
of an oriented framed link $L$ with $l$ components is defined
by using monodromy along solutions of the formal version of the KZ equation.
Forgetting its framing,
the Kontsevich invariant $\hZ(L)$ of an oriented link $L$
is defined in $\cA(\sqcup^l S^1;\Q)/\mbox{FI}$.
The Kontsevich invariant is universal among quantum invariants
in the sense that
the quantum $({\frak g},R)$ invariant recovers from the Kontsevich invariant
through the weight system substituting 
a Lie algebra $\frak g$ and its representation $R$ into Jacobi diagrams.
Moreover, the Kontsevich invariant is universal among Vassiliev invariants
in the sense that\index{Vassiliev invariant!universal ---}
each coefficient of the Kontsevich invariant is a Vassiliev invariant
and any Vassiliev invariant can be presented by
a linear sum of coefficients of the Kontsevich invariant.

\subsection{Calculation of the Kontsevich invariant}

\begin{prob}
\label{prob.cal_ZK}
For each oriented knot $K$,\index{Kontsevich invariant!calculation of ---}
calculate the Kontsevich invariant $\hZ(K)$ for all degrees.
\end{prob}

\begin{rem}
For each $d$
the degree $d$ part of $\hZ(K)$ is a Vassiliev invariant.
Hence, it is algorithmically possible to calculate it in a finite procedure.
It is a problem to calculate $\hZ(K)$ for all degrees.
\end{rem}

\begin{rem}
D. Bar-Natan, T. Le, and D. Thurston \cite{BaLeTh02}
gave the following presentation of 
the Kontsevich invariant\index{Kontsevich invariant!--- of the trivial knot}
 of the trivial knot $O$,
\begin{equation}
\label{eq.Kinv_O}
\log_\sqcup \hZ (O) 
= \frac12 \log \frac{ \sinh (x/2) }{ x/2 },
\end{equation}
where $x$ is an element in $\cB$ (see (\ref{eq.QxB0})),
and $\cB$ is a space isomorphic to $\cA(S^1)$ (see (\ref{eq.A=B=B})).
The Kontsevich invariant of a cable knot of a knot $K$ can be calculated 
by applying a cabling formula \cite{BaLeTh02}
to the Kontsevich invariant of $K$.
The Kontsevich invariant of the connected sum of knots
is given by the connected sum of the Kontsevich invariant of the knots.
Hence, we can calculate the Kontsevich invariant of knots
obtained from the trivial knot
by finite sequences of cabling and connected sum.
To calculate the Kontsevich invariant of other knots
in a combinatorial way, 
we probably need an associator,
whose combinatorial direct presentation for all degrees
is not known yet (see Problem \ref{prob.present_ass}).
\end{rem}

\subsection{Does the Kontsevich invariant distinguish knots?}

\begin{conj}
\label{conj.K_dist_knots}
The Kontsevich invariant distinguishes oriented knots.\index{Kontsevich invariant!strength of ---}
(See Conjecture \ref{conj.V_dist_knots} for an equivalent statement
of this conjecture.)
\end{conj}

\begin{rem}
Kuperberg \cite{Kuperberg_d} showed that 
all finite type invariants either distinguish all oriented knots, 
or there exist prime, unoriented knots which they do not distinguish.
\end{rem}

\begin{prob}
\label{prob.KdistKfromO}
Does there exists a non-trivial oriented knot $K$ 
such that $\hZ(K)$ $= \hZ(O)$
for the trivial knot $O$?
(See Problem \ref{prob.VdistKfromO} for an equivalent problem.)
\end{prob}

\begin{conj}
\label{conj.Kdist_rev}
$\hZ(K) = \hZ(-K)$ for any oriented knot $K$,
where $-K$ denotes $K$ with the opposite orientation.
(See Conjecture \ref{conj.Vdist_rev} for an equivalent statement
of this conjecture.)
\end{conj}

\subsection{Characterization and interpretation of the Kontsevich invariant}

The space $\cA(S^1)$ is an algebra
with the product given by connected sum of Jacobi diagrams on $S^1$.
Since the Kontsevich invariant $\hZ(K)$ of a knot $K$
is group-like in $\cA(S^1)$,
its logarithm $\log \hZ (K)$ belongs to $\cA(S^1)\conn$,
where $\cA(S^1)\conn$ denotes the vector subspace of $\cA(S^1)$
spanned by Jacobi diagrams on $S^1$ with connected uni-trivalent graphs.

\begin{prob}
Characterize those elements of $\hat\cA(S^1)\conn$ of the form $\log \hZ (K)$, 
or those elements of $\cB\conn$ of the form $\log_\sqcup \hZ(K)$.\index{Kontsevich invariant!image of ---}
\end{prob}

\begin{rem}
If the Kontsevich invariant was injective,
this problem would be a step of the classification problem of knots.
It is known (see, for example, \cite{Ohtsuki_book}) that
those elements of $\cA(S^1)\conn^{(\le d)}$ of 
the form of the degree $\le d$ part of $\log\hZ(K)$
forms a lattice, which is isomorphic to the lattice in $\cA(S^1)\conn$
spanned by Jacobi diagrams over $\Z$,
and that the coefficients of $\log \hZ(K)$ are invariants
which are independent to each other.
Hence, it would be meaningful to characterize the form of infinite sums
of coefficients of $\log \hZ(K)$, resp. $\log_\sqcup \hZ(K)$.

$W_{{\frak g},R}\big( \hZ(K) \big)$ is a polynomial in $q^{\pm 1/2N}$
for any simple Lie algebra $\frak g$ and its representation $R$,
where $N$ is the determinant of the Cartan matrix of $\frak g$
(see \cite{Le_isq}),
since it is equal to the quantum $({\frak g},R)$ invariant of $K$.
This somehow characterizes the form of $\hZ(K)$.

The loop expansion characterizes the infinite sum of
subsequences of $\log_\sqcup \hZ(K)$ in each loop-degrees;
see (\ref{eq.RsRC0}), (\ref{eq.RsRC1}), and (\ref{eq.RsRC2})
in the cases of low loop-degrees.
Since the image of the Kontsevich invariant is a countable set,
there should be more restrictive properties.
\end{rem}

\begin{prob}[J. Roberts]  
\label{prob.roberts_tpK}
Give a good topological construction of the Kontsevich integral.\index{Kontsevich invariant!interpretation of ---}
\end{prob}

\begin{rem}[{\rm (J. Roberts)}]  
The Kontsevich integral is, in my opinion, the deepest part of
the existing theory of quantum invariants,
and it has two
(conjecturally) equivalent formulations, each with its mysteries.

(a). In Kontsevich's original formulation of his integral, the part
relating to {\em braids} is reasonably well-understood: it can be
described using configuration spaces of points in the plane, the
Knizhnik-Zamolodchikov equation, $1$-minimal models in rational
homotopy theory, Chen's iterated integrals and Magnus expansions. The
fact that this actually extends to a {\em knot} invariant does not
seem to appear naturally in these pictures, however. Passing from
braids to (Morsified) knots suggests thinking about configuration
spaces of varying numbers of points in the plane, and allowing some
kind of annihilation and creation of pairs. Is there some way to
utilise such spaces? (A related question is Problem \ref{drinfeldprob}.)

(b). In the perturbative integral formulation, the diagrammatic power
series is introduced as a formal device for keeping track of which
linear combinations of the individual (non-invariant) coefficient
integrals give give knot invariants. It isn't really clear from this
point of view why this series should turn out to have good properties
such as multiplicativity, Kricker/Rozansky rationality, etc. Is there
an ``all-in-one'' definition?
\end{rem}

\subsection{The Kontsevich invariant in a finite field}

\begin{prob}
\label{prob.Kinv_ff}
Construct the Kontsevich invariant\index{Kontsevich invariant!--- in a finite field}\index{finite field}
(i.e.\ a universal Vassiliev invariant)\index{Vassiliev invariant!universal ---!--- in a finite field}
with coefficients in a finite field.
\end{prob}

\begin{rem}
If we could find a solution $(R,\Phi)$ 
of the pentagon and hexagon relations with coefficients in a finite field,
such a solution would give
a combinatorial construction of the Kontsevich invariant
with coefficients in that field.
In this case
we can not put 
$R = \exp \Big( 
\begin{picture}(25,10)\put(-3,-3){\pict{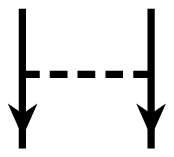}{0.6cm}}\end{picture} 
/ 2 \Big)$
unlike the case of $\Q$ coefficients,
because $p^{-1}$ of the order $p$ of the field
appears in the expansion of the exponential.
\end{rem}

\subsection{The Kontsevich invariant in arrow diagrams}

\begin{conj}[D. Bar-Natan, A. Haviv] 
\label{conj.arrow_Ktriv}
$$
\iota \big( \hZ(O) \big) = \mbox{\rm closure} \left( \exp \Big( \frac12 \big( 
\pict{cd/FIo1}{1.05cm} \!\!-\!\! \pict{cd/FIo2}{1.05cm} \big) \Big) \right),
$$
where $\hZ(O)$ denotes the Kontsevich invariant of the trivial knot
(see \cite{BGRT_w})\index{Kontsevich invariant!--- in arrow diagrams}\index{arrow diagram!Kontsevich invariant in ---}\index{Kontsevich invariant!--- of the trivial knot}
and $\iota$ is the map of Conjecture \ref{conj.inj_A2vA}.
\end{conj}

\begin{rem}[\rm (D. Bar-Natan, A. Haviv)] 
This conjecture is true in any semi-simple Lie algebra.
\end{rem}

\begin{prob}[M. Polyak] 
\label{prob.Kinv_vk}
Construct the ``Kontsevich invariant''\index{Kontsevich invariant!--- of virtual knots}
(i.e.\ a universal 
finite type invariant)\index{finite type invariant!--- of virtual knots!universal ---}\index{virtual knot!finite type invariant of ---!universal ---} 
of virtual knots in $\vcA(I)$.
(See also Conjecture \ref{conj.Vinv_ckvk}.)
\end{prob}

\begin{rem}[\rm (M. Polyak)] 
It is shown by Goussarov (see \cite{GPV}) that
there exists a Gauss diagram formula
for any Vassiliev invariant of classical knots.
His proof is an algorithmical proof,
assuming the existence of such a Vassiliev invariant,
and does {\it not} give a new proof of Kontsevich theorem
``any weight system can be integrated to an invariant of knots''.
It would be nice
to have a new direct combinatorial proof,
which would imply Kontsevich theorem.
Then, it would work for virtual knots.
\end{rem}

\begin{rem}[\rm (M. Polyak)] 
It is known (see, for example, \cite{Ohtsuki_book}) that
quantum invariants
of knots can be defined
by using quasi-triangular quasi-Hopf algebras
with associators $\Phi$.
When $\Phi = 1$,
such definition can naturally extend for virtual knots.
However, when $\Phi \ne 1$
(as in the combinatorial definition of the Kontsevich invariant
of classical knots),
this extension does not work.
\end{rem}

\begin{prob}[D. Thurston] 
\label{prob.csi_vknot}
Construct a series of 
configuration space integrals\index{configuration space!--- integral!--- for virtual knots}  
whose value is in $\vcA(I)$
so that it gives all 
finite type invariants\index{finite type invariant!--- of virtual knots}\index{virtual knot!finite type invariant of ---}
of virtual knots.
\end{prob}

\begin{rem}[\rm (D. Thurston)] 
A technical difficulty is to kill the hidden strata
of the configuration spaces
(see also Problem \ref{prob.khs}).
A way to kill a hidden strata is to use an involution on the strata,
but, in this case, such an involution takes
the following left diagram to the right diagram,
$$
\pict{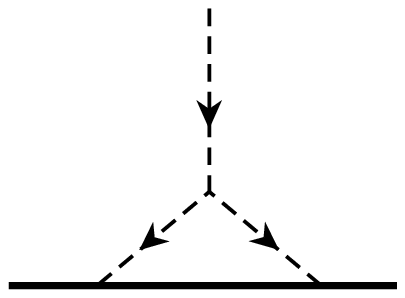}{1.7cm}, \qquad
\pict{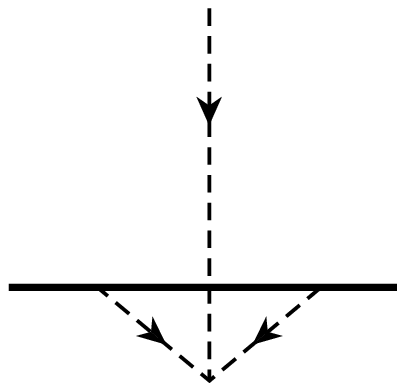}{2.3cm},
$$
where the right diagram is equal to 0 by definition,
while the left one is not necessarily equal to 0.
\end{rem}

\noindent
{{\namae}(M. Polyak)}\quad
Each of the following three approaches gives all 
Vassiliev invariants.\index{Vassiliev invariant!--- by configuration space}
\begin{itemize}
\item
Construction of the Kontsevich invariant
using monodromy along solutions of the KZ equation.
\item
Configuration space integrals
motivated by perturbative Chern-Simons theory.
\item
Gauss diagram formulas,
which count configurations of crossings of knot diagrams.
\end{itemize}
The invariants derived from these three approaches
are expected to be naturally equivalent 
in the following sense.\footnote{\footnotesize 
S. Poirier \cite{Poirier_csi} showed the equivalence between 
the invariants derived from the first and second approaches,
under the assumption of the vanishing of anomaly,
by comparing these invariants for quasi-tangles
(see Question \ref{quest.Kinv=CS}).}{\,}\footnote{\footnotesize 
D. Thurston suggests that
Etingof--Kazhdan R matrices \cite{EiKa96} might be helpful
to relate the invariants derived from the first and third approaches.}
The integral of the second approach gives
an integral presentation of the mapping degree
of a certain map on a configuration space,
and it is shown in the degree 2 case \cite{PoVi_C} that
the invariants of the first and third approaches
can be obtained
by localizing the integral presentation 
with respect to appropriate volume forms on the target space.
A technical difficulty to show this
in a general degree
is to compute the localization on the ``hidden strata'';
it is a part of the boundary of a configuration space,
whose contribution to the derivative of the integral
is killed by an involution on the strata.

\begin{prob}[M. Polyak] 
\label{prob.khs}
Find another way to kill the 
hidden strata,\index{configuration space!--- integral!hidden strata of ---} 
so that the above three approaches
can naturally present the mapping degree of the same map.
\end{prob}

\subsection{The Chern-Simons series of configuration space integrals}
\label{sec.CSseries}

\renewcommand{\thefootnote}{\fnsymbol{footnote}}
\footnotetext[0]{Section \ref{sec.CSseries} was written by C. Lescop.}
\renewcommand{\thefootnote}{\arabic{footnote}}

\begin{quest}[C. Lescop] 
\label{quest.Kinv=CS}
Is the Kontsevich integral\index{Kontsevich invariant!--- by configuration space integral}
of a (zero-framed) knot
equal to the Chern-Simons series of 
configuration space integrals\index{configuration space!--- integral!Kontsevich invariant by ---} 
of the same knot (with Gauss integral 0)?
\end{quest}

The (normalized) Chern-Simons series of configuration space integrals
is a universal Vassiliev knot invariant 
that admits a natural and beautiful symmetric definition that will be given
below before describing the present situation
of this question that was first raised by Kontsevich in \cite{Kontsevich_V}.

\medskip

In 1833, Carl Friedrich Gauss defined the first example of a 
{\em configuration space integral\/}\index{configuration space!--- integral}  
for an oriented two-component link. 
Let us formulate his definition in a modern 
language. Consider an embedding 
$$L: S^1_1 \sqcup S^1_2 \hookrightarrow {\R}^3$$
of the disjoint union of two circles $S^1=\{z \in {\C} \;\;\mbox{s.t.}\;\; |z|=1 \}$ 
into ${\R}^3$. With an element $(z_1,z_2)$ of $S^1_1 \times S^1_2$ that will be called a {\em configuration\/},\index{configuration space}  
we may associate the oriented direction
$\Psi((z_1,z_2))$ of the vector $\overrightarrow{L(z_1)L(z_2)}$. $\Psi((z_1,z_2)) \in S^2$.
Thus, we have associated a map $$\Psi: S^1_1 \times S^1_2 \longrightarrow S^2$$ from a compact oriented 2-manifold to another one
with our embedding.
This map has an integral degree $\mbox{\rm deg}(\Psi)$ that can be defined in several equivalent ways.
For example, it is the number of preimages of a regular value of $\Psi$ counted with signs that can easily be computed from a regular diagram of our two-component link as
$$\mbox{\rm deg}(\Psi)= \sharp \begin{pspicture}[.2](0,0)(.6,.45)
\psline{->}(0.45,0)(.05,.4)
\rput(.25,0){\tiny 1}
\rput(.55,0){\tiny 2}
\psline[border=1pt]{->}(.05,0)(0.45,.4)
\end{pspicture} - \sharp \begin{pspicture}[.2](0,0)(.6,.45)
\rput(.25,0){\tiny 2}
\rput(.55,0){\tiny 1}
\psline{->}(.05,0)(0.45,.4)
\psline[border=1pt]{->}(0.45,0)(.05,.4)
\end{pspicture}=  \sharp \begin{pspicture}[.2](0,0)(.6,.45)
\psline{->}(0.45,0)(.05,.4)
\rput(.25,0){\tiny 2}
\rput(.55,0){\tiny 1}
\psline[border=1pt]{->}(.05,0)(0.45,.4)
\end{pspicture} - \sharp \begin{pspicture}[.2](0,0)(.6,.45)
\rput(.25,0){\tiny 1}
\rput(.55,0){\tiny 2}
\psline{->}(.05,0)(0.45,.4)
\psline[border=1pt]{->}(0.45,0)(.05,.4)
\end{pspicture}\ . $$
It can also be defined as the following 
{\em configuration space integral\/}\index{configuration space!--- integral} 
$$\mbox{\rm deg}(\Psi)=\int_{S^1 \times S^1}\Psi^{\ast}(\omega)$$
where $\omega$ is the homogeneous volume form on $S^2$ such that $\int_{S^2}{\omega}=1$.
It is obvious that this integral degree, that depends continuously on our embedding, is an isotopy invariant; and the reader has recognized that $\mbox{\rm deg}(\Psi)$
is nothing but the linking number of the two components of $L$.

We can again follow Gauss 
and associate the following similar {\em Gauss integral\/}
$I(K)$
to a $C^{\infty}$ embedding $K: S^1 \hookrightarrow {\R}^3$.
Here, we consider the configuration space $C =S^1 \times ]0,2\pi[$, 
and the map $$\Psi: C \longrightarrow S^2$$ 
that maps $(z_1,\eta)$ to the oriented direction of 
$\overrightarrow{K(z_1)K(z_1e^{i\eta})}$, and we set 
$$
I(K)=\int_{C}\Psi^{\ast}(\omega).
$$
This Gauss integral is NOT an isotopy invariant, 
and it can be seen as an exercise that 
it takes any real value on any given isotopy class of knots.

However, we can follow Guadagnini, Martellini and Mintchev 
and associate configuration space integrals to our embedding $K$ 
and to any Jacobi diagram on the circle
$\Gamma$ without small loop like $\begin{pspicture}[.2](0,0)(.6,.4)
\psline[linestyle=dashed,dash=2pt 1pt]{-*}(0.05,.2)(.25,.2)
\psline[linestyle=dashed,dash=2pt 1pt]{-}(.25,.2)(.4,.05)(.55,.2)(.4,.35)(.25,.2)
\end{pspicture}$. 
A {\em configuration\/}\index{configuration space}  
of such a diagram is an embedding $c$ 
of the set $U \cup T$ of its vertices into $\R^3$ 
whose restriction to the set $U$ of univalent vertices factors 
through the knot embedding $K$ so that the factorization 
induces the cyclic order of $U$. Denote the set of these configurations by
$C(K;\Gamma)$. $C(K;\Gamma)$ is an open submanifold of 
$(S^1)^U \times (\R^3)^T$.
Denote the set of dashed edges of $\Gamma$ by $E$, 
and fix an orientation for these edges. 
Then we can define the map
$\Psi:C(K;\Gamma) {\longrightarrow} \left(S\sp{2}\right)\sp{E}$
whose projection to the $S^2$ factor indexed by an edge from a vertex $v_1$ 
to a vertex $v_2$ is the direction of $\overrightarrow{c(v_1)c(v_2)}$.
This map $\Psi$ is again a map between two orientable manifolds 
that have the same dimension, 
namely the number of dashed half-edges of $\Gamma$, 
and we can write the 
{\em configuration space integral:}\index{configuration space!--- integral} 
$$
I(K;\Gamma)=\int_{C(K;\Gamma)}\Psi\sp{\ast}(\Lambda\sp{E}\omega).
$$
For example, if $\theta$ denotes the Jacobi diagram 
$\begin{pspicture}[.2](-.05,0)(.35,.4)
\psline[linestyle=dashed,dash=2pt 1pt]{*-*}(0,.2)(.3,.2)
\pscircle(.15,.2){.15}
\end{pspicture}$, 
then $I(K;\theta)=I(K)$.
Bott and Taubes have proved that this integral is convergent 
\cite{Bott-Taubes}. 
Thus, this integral is well-defined up to sign. 
In fact, an orientation of the trivalent vertices of $\Gamma$ provides 
$I(K;\Gamma)$ with a well-defined sign\footnote{\footnotesize 
Since $S^2$ is equipped with its standard orientation, 
it is enough to orient $C(K;\Gamma) \subset (S^1)^U \times (\R^3)^T$ 
in order to define this sign. 
This will be done by providing the set of the natural coordinates of 
$(S^1)^U \times (\R^3)^T$ with some order up to an even permutation. 
This set is in one-to-one correspondence with 
the set of dashed half-edges of $\Gamma$, 
and the vertex-orientation of the trivalent vertices provides 
a natural preferred such one-to-one correspondence 
up to some (even!) cyclic permutations of three half-edges 
meeting at a trivalent vertex.
Fix an order on $E$, then the set of half-edges becomes ordered by
(origin of the first edge, endpoint of the first edge, 
origin of the second edge,
\dots, endpoint of the last edge), and this order orients $C(L;\Gamma)$.
As an exercise, check that the sign of $I(K;\Gamma)[\Gamma]$ 
does depend neither on our choices
nor on the vertex orientation of $\Gamma$.}
such that the product $I(K;\Gamma)[\Gamma] \in {\cal A}(S^1;\R)$ 
does not depend on the vertex orientation of $\Gamma$.

Now, the {\em perturbative expansion of the Chern-Simons theory
for knots in $\R^3$} is the following sum 
running over all the Jacobi diagrams without small loops 
and without vertex orientation:
$$
Z_{\rm CS}(K)=
\sum{ \frac{I(K;\Gamma)}{\sharp \mbox{Aut}\Gamma}[\Gamma]} \;\; 
\in \cA(S^1;\R)
$$
where $\sharp \mbox{Aut}\Gamma$ is the number of automorphisms of $\Gamma$ as 
a uni-trivalent graph whose univalent vertices are cyclically ordered, but
without vertex-orientation for the trivalent vertices.
The degree one part of $Z_{\rm CS}$ is $\frac{I(K;\theta)}{2}$ 
and therefore $Z_{\rm CS}$ is not invariant under isotopy. 
However, the evaluation\footnote{\footnotesize 
Actually, this evaluation is equal to 
$Z_{\rm CS}(K)\exp(-\frac{I(K;\theta)}{2}\alpha)$ for any representative $K$,
 where $\alpha \in {\cal A}([0,1];\R)$ is the Bott and Taubes anomaly.} 
of $Z_{\rm CS}$ at representatives of knots with null Gauss integral 
is an isotopy invariant that is a universal Vassiliev invariant of knots 
\cite{Bott-Taubes,Altschuler-Freidel,Thurston,Poirier_csi}.
Now, the still open question raised by Kontsevich in \cite{Kontsevich_V} is:
{\em Is the Kontsevich integral\index{Kontsevich invariant!--- by configuration space integral}
 of a zero framed representative of a knot $K$
equal to the above series of 
configuration space integrals\index{configuration space!--- integral!Kontsevich invariant by ---}  
of a representative of $K$ with Gauss integral 0?}

This question has been reduced by Sylvain Poirier \cite{Poirier_csi} 
to the computation of the following constant in 
${\cal A}(S^1; \R)={\cal A}([0,1];\R)$ 
that is called the Bott and Taubes {\em anomaly}.\index{anomaly} 
In order to define the anomaly, 
replace the above knot $K$ by a straight line $D$, 
and consider a Jacobi diagram $\Gamma$ on the oriented line. 
Define $C(D;\Gamma)$ and $\Psi$ as before.
Let $\hat{C}(D;\Gamma)$ be the quotient of $C(D;\Gamma)$ 
by the translations parallel to $D$ and by the positive homotheties, 
then $\Psi$ factors through $\hat{C}(D;\Gamma)$ that has two dimensions less.
Now, allow $D$ to run among all the oriented lines 
through the origin of $\R^3$ and define $\hat{C}(\Gamma)$ as the total space
of the fibration over $S^2$ 
where the fiber over the direction of $D$ is $\hat{C}(D;\Gamma)$. 
$\Psi$ becomes a map between two smooth oriented\footnote{\footnotesize 
$\hat{C}(\Gamma)$ carries a natural smooth structure 
and can be oriented as follows: 
orient $C(D;\Gamma)$ as before, orient $\hat{C}(D;\Gamma)$ so that 
$C(D;\Gamma)$ is locally homeomorphic to the oriented product
(translation vector of the oriented line, ratio of homothety) 
$\times \hat{C}(D;\Gamma)$ and orient $\hat{C}(\Gamma)$ 
as the local product $\mbox{base} \times \mbox{fiber}$.}
manifolds of the same dimension. Then we can again define 
$$
I(\Gamma)=\int_{\hat{C}(\Gamma)}\Psi\sp{\ast}(\Lambda\sp{E}\omega).
$$
Now, the anomaly\index{anomaly}  
is the following sum running over all Jacobi diagrams on the oriented lines (again without vertex-orientation and without small loop):
$$
\alpha=\sum{ \frac{I(\Gamma)}{\sharp \mbox{Aut}\Gamma}[\Gamma]} \;\; 
\in {\cal A}([0,1];\R).
$$
Its degree one part is 
$$
\alpha_1= \begin{pspicture}[0.2](0,0)(.4,.8)
\psline{->}(0.05,0.05)(0.05,0.75)
\pscurve[linestyle=dashed,dash=2pt 1pt]{*-*}(0.05,.15)(.3,.35)(.05,.5)
\end{pspicture}. 
$$
It is not hard to see that for any integer $n$, $\alpha_{2n}=0$.
In \cite{Poirier_csi}, Sylvain Poirier proved that if all $\alpha_i$ vanish
for $i \geq 2$, then the answer to the above Kontsevich question is YES, and he
computed $\alpha_3=0$. He also computed $\alpha_5=0$ with the help of Maple.
In \cite{Lescop_dK}, it is proved that $\alpha$ is a combination of diagrams 
with two univalent vertices. Poirier also gave an equivalent definition of the
anomaly that allows one to see that, for any $i>1$,  $\alpha_i$ is a combination of diagrams 
with at least 6 univalent vertices.

As a corollary, all coefficients of the HOMFLY polynomial properly normalized that are Vassiliev invariants of degree less than seven can be explicitly written as 
combinations of the above configuration space integrals. A positive answer to the Kontsevich question would allow one to express any canonical Vassiliev invariant as an explicit combination of the above configuration space integrals.

G. Kuperberg and D. Thurston have constructed a universal finite type
invariant for homology spheres as a series of configuration space integrals
similar to the above Chern-Simons series in \cite{Kuperberg-Thurston}. 
Their construction yields two natural questions
that are stated in Question \ref{quest.KT_LMO}.

\subsection{Associators}

An {\it associator}\index{associator} 
$\Phi$ is defined to be 
an invertible group-like element 
in $\cA ( \downarrow\downarrow\downarrow ; \C)$
satisfying that $\e_2 \Phi = 1 \in \cA ( \downarrow\downarrow ;\C)$
and the following relations,
\begin{align*}
& \pict{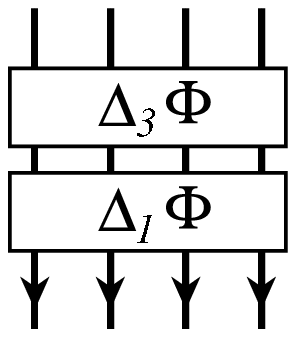}{1.8cm} = \pict{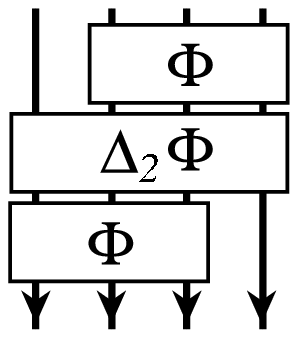}{1.8cm}, \\
& \pict{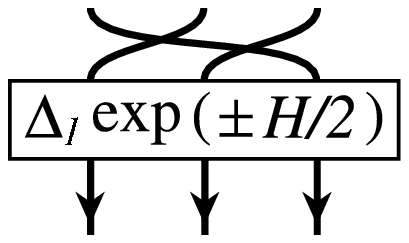}{1.2cm} = \pict{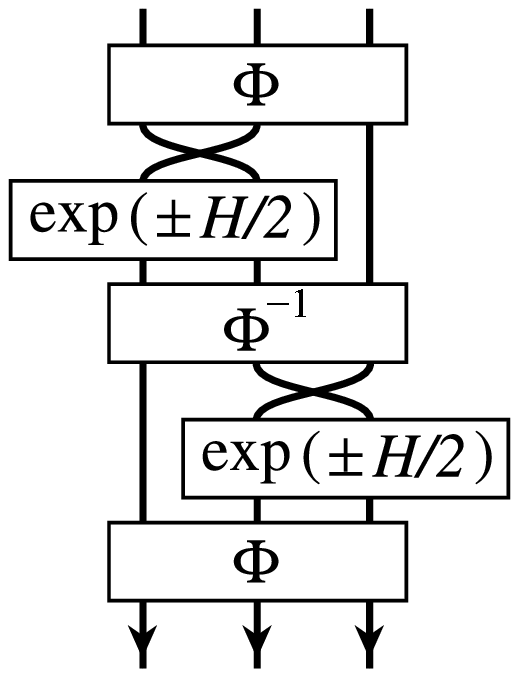}{3.5cm} 
\mbox{ where we put } H = \pict{w/w56gon5}{0.6cm}.
\end{align*}
Here,
$\Delta_i$ and $\e_i$ are
the comultiplication and the counit
acting on the $i$-th solid line;
see \cite{BarNatan_n} for these notations.
An associator is derived from a Drinfel'd series $\varphi(A,B)$ by
\begin{equation}
\label{eq.Phi=phi}
\Phi = \varphi \big(
\pict{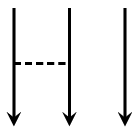}{1cm}, 
\pict{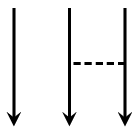}{1cm} \big)
\in \cA ( \downarrow\downarrow\downarrow ;\C),
\end{equation}
where a {\it Drinfel'd series}
is an invertible group-like power series $\varphi(A,B)$
of non-commutative indeterminates $A$ and $B$
satisfying certain relations.

The Drinfel'd associator is given as follows.
We consider the differential equation 
\begin{equation}
\label{eq.Gz}
G'(z) = \frac{1}{2 \pi \sqrt{-1}}
\Big( \frac{A}{z} + \frac{B}{z-1} \Big) G(z),
\end{equation}
for an analytic function $G$ of the variable $z$,
where $G(z)$ belongs to the formal power series ring
$\C \langle\!\langle A,B \rangle\!\rangle$
of non-commutative indeterminates $A$ and $B$.
There exists unique solutions 
$G_{\scriptscriptstyle (\bullet\bullet)\bullet}$ and 
$G_{\scriptscriptstyle \bullet(\bullet\bullet)}$ 
of the above differential equation of the forms
\begin{align*}
G_{\scriptscriptstyle (\bullet\bullet)\bullet}(z) 
&= f(z) z^{A/2\pi\sqrt{-1}} \\
G_{\scriptscriptstyle \bullet(\bullet\bullet)}(z) 
&= g(1-z) (1-z)^{B/2\pi\sqrt{-1}} 
\end{align*}
where $f(z)$ and $g(z)$ are analytic functions
with $f(0) = g(0) = 1 \in \C \langle\!\langle A,B \rangle\!\rangle$
defined in a neighborhood of $0 \in \C$.
The power series 
$\varphi\subkz(A,B) \in \C \langle\!\langle A,B \rangle\!\rangle$
is defined by
$G_{\scriptscriptstyle (\bullet\bullet)\bullet} 
= G_{\scriptscriptstyle \bullet(\bullet\bullet)} \varphi\subkz(A,B)$.
The associator derived from $\varphi\subkz(A,B)$ by (\ref{eq.Phi=phi})
is called the {\it Drinfel'd associator}.

\begin{prob}
\label{prob.present_ass}
Find a combinatorial direct presentation of 
an associator\index{associator!presentation of ---}  
for all degrees,
in particular, an associator with rational coefficients.
\end{prob}

\begin{rem}
We still do not have
a combinatorial direct presentation of any associator for all degrees.
This implies that we still do not know a combinatorial direct presentation
of the Kontsevich invariant\index{Kontsevich invariant!calculation of ---}
 of each knot for all degrees
(except for the trivial knot);
see Problem \ref{prob.cal_ZK} and its remarks.
Bar-Natan \cite{BarNatan_n} showed a combinatorial degree-by-degree proof
of the existence of solutions of the defining relations of 
a pair $(R,\Phi)$.
Our definition of $\Phi$ follows from the defining relations
when $R$ is given by $\exp\big(\frac12 \pict{w/w56gon5}{0.6cm} \big)$.
\end{rem}

\begin{rem}
The only associator 
whose coefficients can be directly presented for all degrees so far
is the Drinfel'd associator.
We can present all degrees of the Drinfel'd associator
by a limit of iterated integrals (see (\ref{eq.phiKZ_int}))
of by multiple zeta functions (see (\ref{eq.phiKZ_mz})).
It is known \cite{LM_H} that
all associators are related to each other by ``twists'',
which are some actions of symmetric elements in 
$\cA(\downarrow\downarrow; \C)$ on associators.
\end{rem}

\begin{rem}
$\varphi\subkz (A,B)$ is presented by the following limit ,
\begin{equation}
\label{eq.phiKZ_int}
\varphi\subkz (A,B) 
= \lim_{\e \to 0} \e^{- B/2\pi\sqrt{-1}} G_\e ( 1-\e) \e^{A/2\pi\sqrt{-1}},
\end{equation}
where we regard $\e^x$ as 
$$
\e^{x} = \exp( x \log \e )
= 1 + x \log \e + x^2 \frac{ (\log\e)^2 }{2} + \cdots.
$$
Further,
$G_\e$ is a solution of (\ref{eq.Gz}) given by
$$
G_\e ( 1-\e)
= 1 +
\sum_{m=1}^\infty 
\int_{\e \le t_1 \le \cdots \le t_m \le 1-\e}
w(t_m) \cdots w(t_1) dt_1 \cdots dt_m,
$$
putting
$$
w(t) = \frac{1}{2\pi\sqrt{-1}}
\Big( \frac{A}{t} + \frac{B}{t-1} \Big).
$$
\end{rem}

\begin{rem}
In \cite{LM_H}, 
$\varphi\subkz(A,B)$ is presented by
\begin{align}
\varphi\subkz(A,B) = 1 + 
\sum_{l=1}^{\infty}
\sum_{\scriptsize \mbox{\bf a}, \mbox{\bf b}, \mbox{\bf p}, \mbox{\bf q}}
(-1)^{\scriptsize | \mbox{\bf b}| + | \mbox{\bf p} |}
& \eta(\mbox{\bf a} + \mbox{\bf p}, \mbox{\bf b} + \mbox{\bf q})
\binom{\mbox{\bf a} + \mbox{\bf p}}{\mbox{\bf p}}
\binom{\mbox{\bf b} + \mbox{\bf q}}{\mbox{\bf q}} \notag \\*
& \times B^{\scriptsize |\mbox{\bf q}|} 
(A,B)^{\scriptsize (\mbox{\bf a}, \mbox{\bf b})}
A^{\scriptsize | \mbox{\bf p} |},
\label{eq.phiKZ_mz}
\end{align}
where the second sum runs over
$\mbox{\bf a}, \mbox{\bf b}, \mbox{\bf p}, \mbox{\bf q}$
such that
the sum of their length is equal to $l$
and entries of them are non-negative integers.
Here, 
the notations are given by
\begin{align*}
& \eta(\mbox{\bf a}, \mbox{\bf b})
= \zeta(\underbrace{1,1,\cdots,1}_{a_1-1}, b_1+1, 
\underbrace{1,1,\cdots,1}_{a_1-1}, b_2+1,\cdots, 
\underbrace{1,1, \cdots,1}_{a_l-1},b_l+1), \\
& | \mbox{\bf a} | = a_1 + a_2 + \cdots + a_l, \\
& \binom{\mbox{\bf a}}{\mbox{\bf b}} = 
\binom{a_1}{b_1} \binom{a_2}{b_2} \cdots \binom{a_l}{b_l}, \\
& (A,B)^{\scriptsize (\mbox{\bf a}, \mbox{\bf b})} 
= A^{a_1} B^{b_1} \cdots A^{a_l} B^{b_l}.
\end{align*}
for $\mbox{\bf a}= (a_1, \cdots, a_l)$ and $\mbox{\bf b}= (b_1, \cdots, b_l)$,
where the multiple zeta function is defined by
$$
\zeta(a_1,a_2,\cdots,a_k) 
= \sum_{n_1 < n_2 < \cdots < n_k \in \N} 
n_1^{-a_1} n_2^{-a_2} \cdots n_k^{-a_k}.
$$
In particular,
$$
\varphi\subkz(A,B) 
= 1 + \frac{1}{24} [A,B]
- \frac{\zeta(3)}{(2 \pi \sqrt{-1})^3} ([A,[A,B]] + [B,[A,B]])
+ \left( \begin{array}{l} \mboxsm{terms of} \\ \mboxsm{degree $\ge 4$}
   \end{array} \right).
$$
\end{rem}

\begin{rem}
In \cite{BarNatan_n},
an associator with rational coefficients is given in low degrees by
\begin{align*}
& \log \varphi (A,B) = 
\frac{[A,B]}{48} - \frac{ 8 [A, [A, [A, B]]] + [A, [B, [A, B]]]}{11520} \\*
& + \frac{ [A,[A,[A,[A,[A,B]]]]] }{60480}
+ \frac{ [A,[A,[A,[B,[A,B]]]]] }{1451520} 
+ \frac{ 13[A,[A,[B,[B,[A,B]]]]] }{1161216} \\*
& + \frac{ 17[A,[B,[A,[A,[A,B]]]]] }{1451520} 
+ \frac{ [A,[B,[A,[B,[A,B]]]]]}{1451520} \\*
& - (\mbox{interchange of $A$ and $B$}) \\*
& + (\mbox{terms of degree $\ge 8$}).
\end{align*}
\end{rem}

\begin{prob}[J. Roberts]  
\label{drinfeldprob} 
Construct a rational Drinfel'd associator in
the context of rational homotopy theory.\index{associator!--- in rational homotopy theory}  
\end{prob}

\begin{rem}[{\rm (J. Roberts)}]  
The theory of $1$-minimal models provides a representation of the
pure braid group, which is the fundamental group of the configuration
space of distinct ordered points in $\C$, the ``pure braid space'' for
short. This is the representation coming from the Kontsevich integral.
A better way to describe it is as a representation of the
fundamental {\em groupoid} of the pure braid space, using ``basepoints
at infinity'' described by associations (bracketings) of the
points. In this picture, the Drinfel'd associator is the image of a
certain path which changes the basepoint. Is there a theory of
$1$-minimal models for fundamental groupoids which gives a
straightforward construction of a (rational-valued) associator, as an
alternative to the tricky iterative procedures of \cite{BarNatan_n}?
\end{rem}

\subsection{Graph cohomology}
\label{sec.graph_cohom}

\renewcommand{\thefootnote}{\fnsymbol{footnote}}
\footnotetext[0]{Section \ref{sec.graph_cohom} was written by J. Roberts.}
\renewcommand{\thefootnote}{\arabic{footnote}}

\begin{prob}[J. Roberts]  
\label{prob.roberts6}
What is graph cohomology\index{graph cohomology} 
the cohomology of? 
\end{prob}

\begin{rem}[{\rm (J. Roberts)}]  
In the theory of quantum knot invariants
 such as the Jones
polynomial, the topology and algebra (in this case, the group $SU(2)$)
are entangled somewhat confusingly. Passing to the theory of finite
type invariants, they become separated: there is a purely topological
part (the Kontsevich integral of a knot) and a purely algebraic part
(the weight system associated to $SU(2)$) whose intermediary is the
space of Jacobi diagrams.

Viewing this space as (part of) Kontsevich's {\em graph (co)homology}
\cite{Kontsevich_F}, we see that quantum invariants arise from a pairing
between elements of graph cohomology and homology. But what actually
is this cohomology? A good geometric interpretation of it might lead
to better understanding of the topological and algebraic constructions
involving it, and their composite. 

Most of the intuition about graph cohomology has been built up from
the algebraic side: it has been portrayed primarily as a kind of
universal invariant theory for Lie algebras. Vogel has pursued this
idea the furthest, but he also showed \cite{Vogel_alg} that not all weight
systems come from classical Lie algebras. In fact, the work of
Rozansky and Witten \cite{RoWi97} and Kapranov \cite{Kap99}
demonstrates that compact holomorphic symplectic manifolds can be used
instead of Lie algebras to define Vassiliev weight systems, and this
gives quite a different perspective on graph cohomology, which Simon
Willerton and I have been studying \cite{Rob01}.

In a similar vein, Bar-Natan, Le and Thurston \cite{Thurston_w} have proved
the so-called ``wheeling conjectures'', diagrammatic generalisations
of the Duflo isomorphism of Lie theory. Their theorem is far too
striking for a purely combinatorial interpretation to be
satisfactory. Does it have a geometric interpretation?

Kontsevich \cite{Kontsevich_F} has given three topological interpretations of
graph cohomology. The first is that it is the twisted cohomology of
``outer space'', the classifying space of the group of outer
automorphisms of a free group. This is analogous to the fact that a
certain complex of {\em fatgraphs} gives the cohomology of the moduli
space of Riemann surfaces. The answer is unsatisfying because the
natural geometric model for the classifying space is, unlike the
Riemann moduli space, not a smooth orbifold, and if we are seeking
geometric constructions underlying the various kinds of diagrammatic
operations we encounter, smoothness would seem to be an essential
property. Is there is a better model?

A second approach comes from configuration spaces of points in
$\R^3$. The complex of graphs (with distinguished legs) maps to the de
Rham complex of configuration spaces, and gives a model for its
cohomology. This kind of viewpoint was exploited by Kontsevich (and
Taubes, and Axelrod and Singer) in defining the perturbative
invariants of $3$-manifolds, and by Bott and Taubes \cite{Bott-Taubes} for
knots.
\end{rem}

In this context, Lie algebra 
weight systems\index{weight system!--- and configuration space}
are functionals on the
cohomology of the configuration spaces,\index{configuration space!homology of ---}  
and might be thought of as
homology classes, or even cycles. Hence the following problem, posed
by Raoul Bott:

\begin{prob}[R. Bott]  
\label{prob.bott}
Give a geometric construction of these homology classes 
coming from Lie algebras. 
\end{prob}

The third and currently best interpretation of graph cohomology is
that it is the cohomology of an infinite-dimensional Lie algebra of
formal Hamiltonian vector fields. Kontsevich uses this to explain (and
vastly generalise) Rozansky-Witten weight systems in terms of
Gelfand-Fuchs cohomology. Can this interpretation be employed on the
topological rather than algebraic side? In other words, is there a
construction involving knots and algebras of formal vector fields
which yields\index{Kontsevich invariant!interpretation of ---} 
the Kontsevich integral?

\subsection{The loop expansion of the Kontsevich invariant}
\label{sec.loopKinv}

The loop expansion is the series of the rational presentations
of the Kontsevich invariant in loop-degrees.\index{Kontsevich invariant!loop expansion of ---}\index{loop expansion|see {Kontsevich invariant}}
It was conjectured by Rozansky \cite{Rozansky_rc}.
The existence of such rational presentations has been proved by 
Kricker \cite{Kricker_lK}
(though such a rational presentation itself is not necessarily
a knot invariant in a general loop degree).
Further, Garoufalidis and Kricker \cite{GK_rat}
defined a knot invariant in any loop degree, 
from which such a rational presentation can be deduced.

We have three isomorphic algebras
\begin{equation}
\label{eq.A=B=B}
\cA(S^1) \cong \cB \cong \cB_\sqcup,
\end{equation}
where the first isomorphism
is the formal Poincare-Birkhoff-Witt isomorphism,
and $\cB$ has the product structure
related, by the isomorphism, to the product structure of $\cA(S^1)$
given by connected sum.
Further, the second isomorphism is the wheeling isomorphism \cite{BGRT_w}
between $\cB$ and $\cB_\sqcup$,
where $\cB_\sqcup$ is $\cB$ as a space
and has the product
given by the disjoint union of uni-trivalent graphs.

We denote by $\cB\conn$ the vector subspace of $\cB_\sqcup$
spanned by connected uni-trivalent graphs,
and denote by $\cB\conn^{({\rm loop}\ l)}$
the vector subspace of $\cB\conn$ spanned by
connected uni-trivalent graphs of loop-degree $l$,
where the {\it loop-degree} of a uni-trivalent graph
is defined to be half of the number given by
the number of trivalent vertices minus the number of univalent vertices.
Then, 
$$
\cB\conn = \bigoplus_{l=0}^\infty \cB\conn^{({\rm loop}\ l)}.
$$
Each $\cB\conn^{({\rm loop}\ l)}$
can be presented by using the polynomial rings in $H^1(G)$
for trivalent graphs $G$ of loop-degree $l$
subject to $\mbox{Aut}(G)$ and the AS and IHX relations.
We will present $\cB\conn^{({\rm loop}\ l)}$
for $l=0,1,2$ in this way,
to state the loop expansion in these loop-degrees.

When $l=0$, we have the map
\vspace{-1pc}
\begin{equation}
\label{eq.QxB0}
\Q[x] \longrightarrow
\cB\conn^{({\rm loop}\ 0)}, \qquad
x^n \longmapsto \pict{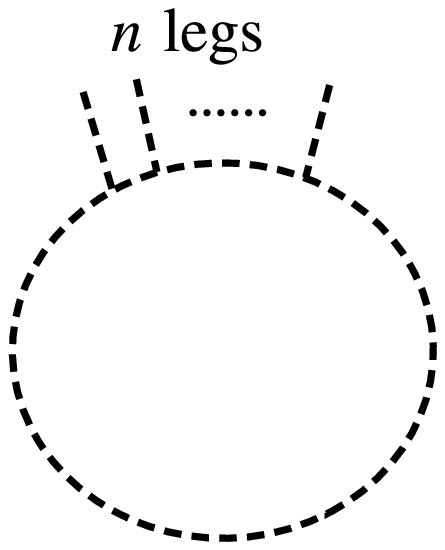}{2.8cm},
\end{equation}
\vspace{-1.2pc}

\noindent
regarding $x$ as a basis of $H^1(\mbox{circle})$.
Since the orientation-reversing automorphism of $S^1$
takes $x^n$ to $-x^n$ by the AS relation,
the above map deduces the following isomorphism,
\begin{equation}
\label{eq.loop0}
\cB\conn^{({\rm loop}\ 0)} \cong \Q[x^2].
\end{equation}
For a knot $K$,
\begin{equation}
\label{eq.RsRC0}
\big( \log_\sqcup \hZ (K) \big)^{({\rm loop}\ 0)} 
= \frac12 \log \frac{ \sinh (x/2) }{ x/2 }
- \frac12 \log \Delta_K(e^x),
\end{equation}
where $\log_\sqcup$ is the logarithm in $\cB_\sqcup$
regarding $\hZ(K)$ as in $\cB_\sqcup$,
and the left hand side is the summand of 
$\log_\sqcup \hZ(K) \in \cB\conn$ in $\cB\conn^{({\rm loop}\ 0)}$.
This development follows from the theory of \cite{BG_MMR}.
See also \cite{Kricker_lK,GK_rat} (and references therein) 
for a recent direct calculation.

When $l=1$, we have the map
\vspace{-1pc}
$$
\Q[x_1,x_2,x_3]
\longrightarrow
\cB\conn^{({\rm loop}\ 1)},
\qquad
x_1^{n_1} x_2^{n_2} x_3^{n_3} \longmapsto \pict{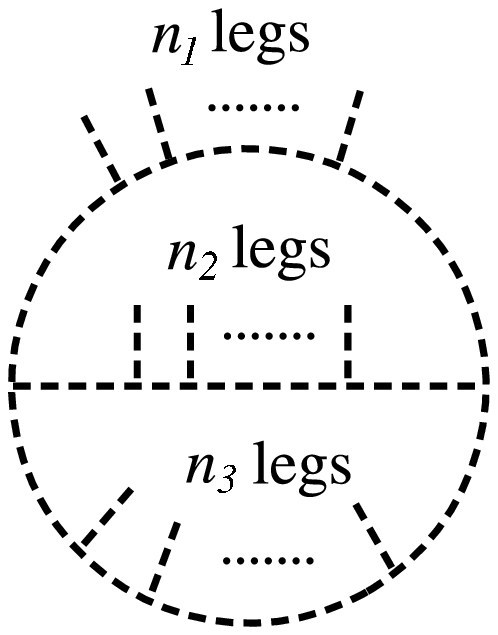}{3.1cm},
$$
\vspace{-1.2pc}

\noindent
regarding $H^1(\mbox{$\theta$-graph})$
as the vector space spanned by $x_1$, $x_2$, $x_3$
subject to the relation $x_1 + x_2 + x_3 = 0$.
Since $\mbox{Aut}(\mbox{$\theta$-graph}) \cong {\frak S}_2 \times {\frak S}_3$,
the above map deduces
\begin{align}
\cB\conn^{({\rm loop}\ 1)}
&\cong\Q[x_1,x_2,x_3] \mbox{\Large $/$} 
   ({\frak S}_2 \times {\frak S}_3,\ x_1+x_2+x_3=0) \notag \\*
&\cong\Big( \Q[x_1,x_2,x_3] \mbox{\Large $/$} (x_1+x_2+x_3=0) \Big)
^{{\frak S}_2 \times {\frak S}_3} \notag \\*
&\cong\Big( \Q[\sigma_1,\sigma_2,\sigma_3] \mbox{\Large $/$} (\sigma_1=0) \Big)
^{({\rm even})}
\cong \Q[\sigma_2, \sigma_3^2], \label{eq.loop1}
\end{align}
where $\sigma_i$ denotes
the $i$-th elementary symmetric polynomial in $x_1$, $x_2$, and $x_3$.
(To compute $\cB\conn^{({\rm loop}\ 1)}$ in a precise argument,
we must also consider the space of ``dumbbell diagram'' with legs.
Since this space is injectively mapped to 
the right hand side of the above formula,
we omit its computation here.)
For a knot $K$ there exists a polynomial $P^\theta_K(t_1,t_2,t_3)$, 
called {\it the 2-loop polynomial},\index{2-loop polynomial}
satisfying that
\begin{equation}
\label{eq.RsRC1}
\big( \log_\sqcup \hZ (K) \big)^{({\rm loop}\ 1)} 
= \frac{ P^\theta_K( e^{x_1}, e^{x_2}, e^{x_3} )}
       {\Delta_K(e^{x_1}) \Delta_K(e^{x_2}) \Delta_K(e^{x_3})}.
\end{equation}
The 2-loop polynomial $P^\theta_K(t_1,t_2,t_3)$ in $t_1, t_2, t_3$
satisfying $t_1 t_2 t_3 =1$
is uniquely determined by each knot $K$.
It is an invariant of $K$
satisfying that \newline
$P^\theta_K(t_i^{\pm1}, t_j^{\pm1}, t_k^{\pm1}) = P^\theta_K(t_1,t_2,t_3)$
for any signs and any $\{ i,j,k \} = \{ 1,2,3 \}$.

\begin{prob}
\label{prob.P_K}
Find a topological construction of 
the 2-loop polynomial\index{2-loop polynomial!topological construction of ---} 
$P^\theta_K$.
\end{prob}

\begin{rem}
As in (\ref{eq.RsRC0})
the loop-degree 0 part of the Kontsevich invariant
is presented by the Alexander polynomial,
which can be constructed from 
the homology of the infinite cyclic cover of the knot complement.
It is shown, in \cite{GR00}, that
the ``first derivative'' of the 2-loop polynomial
is given in terms of linking functions
associated to the infinite cyclic cover of the knot complement.
It is expected \cite{GR00} that
the 2-loop polynomial would be described
in terms of invariants
of the infinite cyclic cover of the knot complement.
\end{rem}

\begin{rem}
A table of the 2-loop polynomial for knots with up to 7 crossings
is given by Rozansky \cite{Rozansky_rc2}.
See also a computer program \cite{Rozansky_website},
which calculates the 2-loop polynomial of each knot.
For example,
\begin{align*}
12 P^\theta_{3_1}(t_1,t_2,\frac1{t_1t_2})
&= -t_1^2 t_2 + t_1^2, \\
12 P^\theta_{4_1}(t_1,t_2,\frac1{t_1t_2})
&= 0, \\
12 P^\theta_{5_1}(t_1,t_2,\frac1{t_1t_2})
&= 2 t_1^4 t_2^2 - 2 t_1^4 t_2 + 2 t_1^4 - t_1^2 t_2 + t_1^2.
\end{align*}
\end{rem}
The 2-loop polynomial for the torus knots is calculated
independently by March\'e \cite{Marche_K} and Ohtsuki \cite{Ohtsuki_c2lp}.

The following problem is a step to Problem \ref{prob.P_K}.

\begin{prob}[A. Kricker] 
\label{prob.pres_Theta}
Let $K_T$ be the knot obtained from a tangle $T$
as shown in Figure \ref{fig.T,K_T}.
Find a presentation of 
the 2-loop polynomial\index{2-loop polynomial!--- of knots of genus 1} 
$P^\theta_{K_T}$
of $K_T$ by using the Kontsevich invariant $\hZ(T)$ of $T$.\index{Kontsevich invariant!--- and 2-loop polynomial}
\end{prob}

\begin{rem}[{\rm (A. Kricker)}] 
$P^\theta_{K_T}$ might be 
presented by the degree $\le 3$ part of $\hZ(T)$.
Generalize the presentation
$\Delta_K(t) = \mbox{det} ( t^{1/2} S - t^{-1/2} S^T )$
of the Alexander polynomial $\Delta_K(t)$
by a Seifert matrix $S$ of $K$.
\end{rem}

\begin{figure}[ht!]
$$
\begin{array}{ccc}
\pict{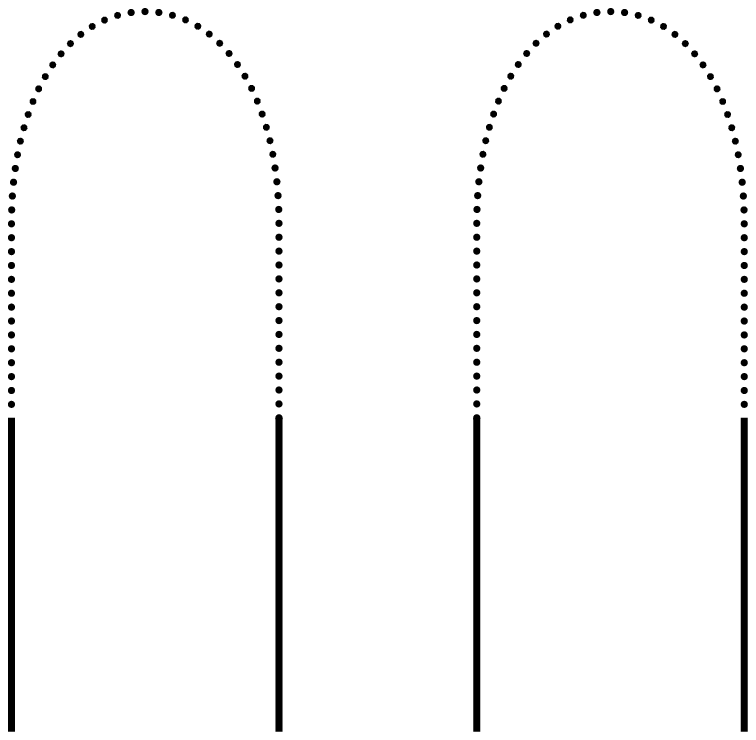}{3.3cm} & \quad & \pict{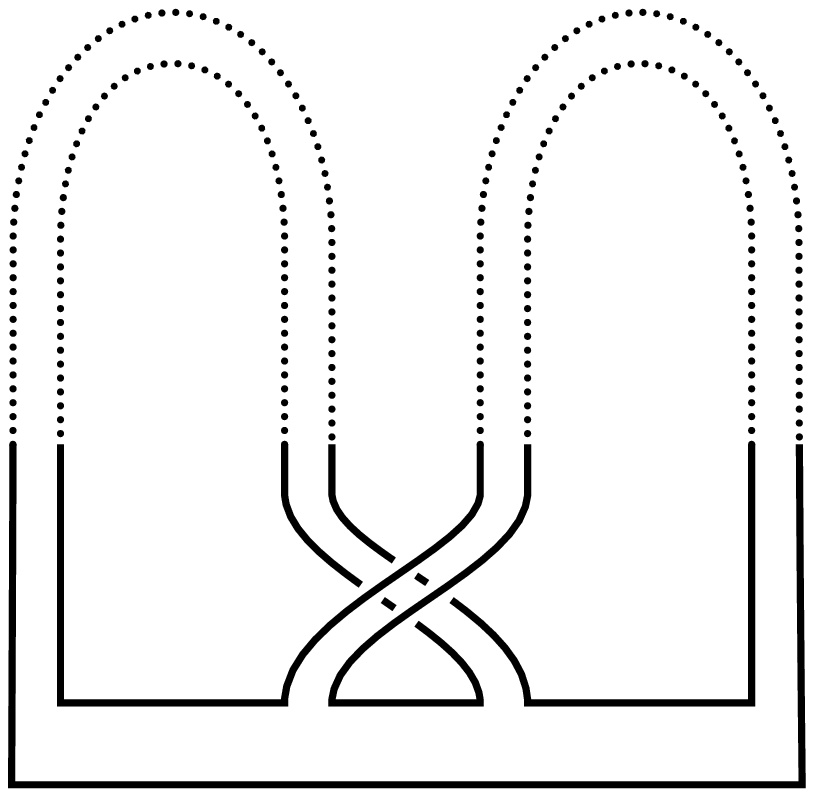}{3.5cm} \\
T && K_T
\end{array}
$$
\caption{\label{fig.T,K_T}
The knot $K_T$ is obtained from the 2-parallel of a 2-strand tangle $T$
by adding the tangle depicted in solid lines in the right picture.
The dotted lines imply strands possibly knotted and linked in some fashion.}
\end{figure}

When $l=2$, we have the map
\vspace{-1pc}
$$
x_1^{n_1} x_2^{n_2} \cdots x_6^{n_6}
\longmapsto \pict{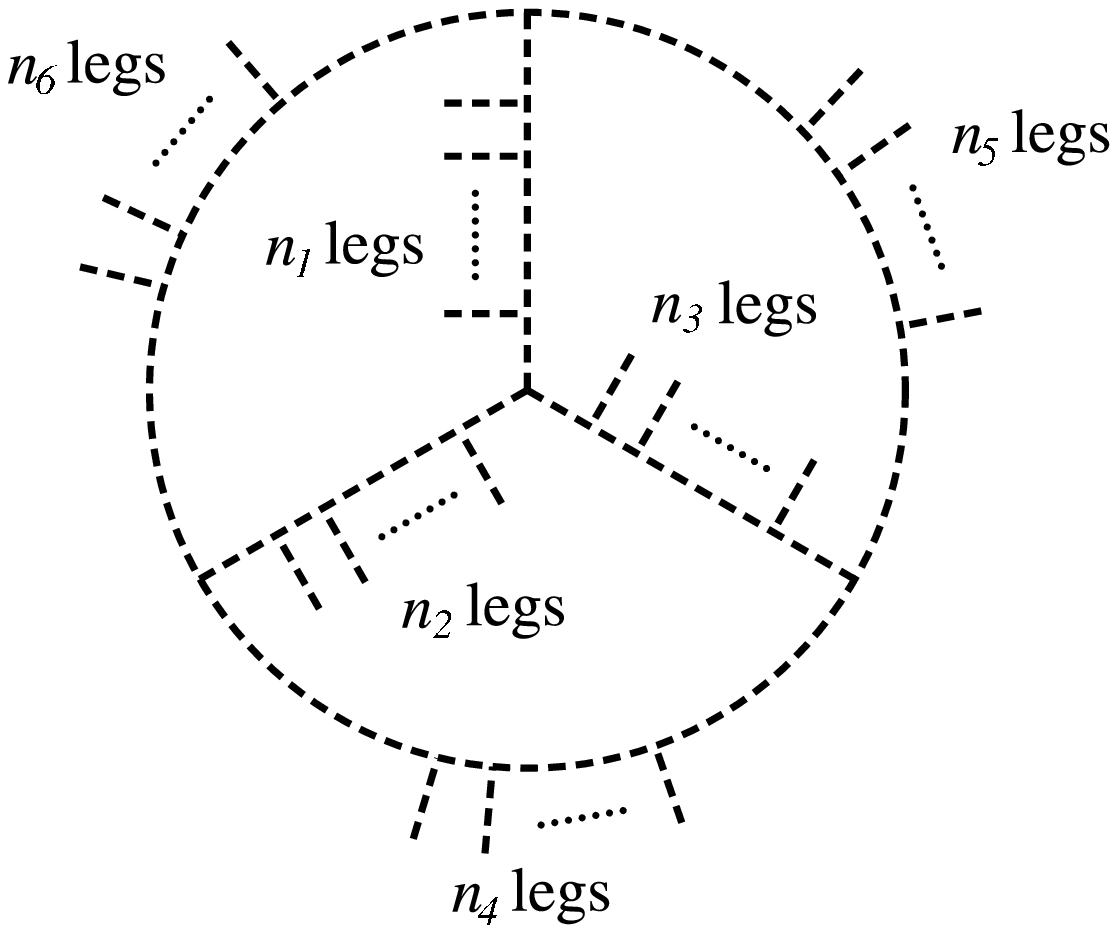}{4.5cm},
$$
\vspace{-1.2pc}

\noindent
which deduces the following isomorphism,
$$
\begin{array}{r}
\cB\conn^{({\rm loop}\ 2)} \cong
\Q[x_1,x_2,\cdots,x_6] \mbox{\Large $/$}
({\frak S}_4, 
x_1+x_2+x_3 = x_1+x_6-x_5 = 0 \ \  \\
x_2+x_4-x_6 = x_3+x_5-x_4 = 0). \end{array}
$$
Corresponding to faces of a tetrahedra, we put
$y_1 =x_1-x_2-x_6$, \newline
$y_2 =x_2-x_3-x_4$,
$y_3 =x_3-x_1-x_4$,
and $y_4 =x_4+x_5+x_6$.
Then,
\begin{align*}
\cB\conn^{({\rm loop}\ 2)} 
& \cong
\Q[y_1,y_2,y_3,y_4]
\mbox{\Large $/$}
({\frak S}_4, y_1+y_2+y_3+y_4 = 0) \\*
& \cong
\Big( \Q[y_1,y_2,y_3,y_4]
\mbox{\Large $/$}
(y_1+y_2+y_3+y_4 = 0) \Big)^{{\frak S}_4},
\end{align*}
where the action of $\tau \in {\frak S}_4$ takes
a polynomial $p(y_1,y_2,y_3,y_4)$ to \newline
$(\mbox{sgn}\tau) p(y_{\tau(1)}, y_{\tau(2)}, y_{\tau(3)}, y_{\tau(4)})$.
Hence,
\begin{equation}
\label{eq.loop2}
\cB\conn^{({\rm loop}\ 2)} \cong
\Big( \Q[\sigma_2,\sigma_3,\sigma_4] \Big)^{({\rm even})}
\cong \Q[\sigma_2, \sigma_3^2, \sigma_4],
\end{equation}
where $\sigma_i$ is
the $i$-th elementary symmetric polynomial in $y_1$, $y_2$, $y_3$, and $y_4$.
(To compute $\cB\conn^{({\rm loop}\ 2)}$ in a precise argument,
we need some more computations,
which are omitted here.)
For a knot $K$ there exists a polynomial $P'_K(t_1,t_2,\cdots, t_6)$ 
satisfying that
\begin{equation}
\label{eq.RsRC2}
\big( \log_\sqcup \hZ (K) \big)^{({\rm loop}\ 2)} 
= \frac{ P'_K( e^{x_1}, e^{x_2}, \cdots, e^{x_6} )}
       {\Delta_K(e^{x_1}) \Delta_K(e^{x_2}) \cdots \Delta_K(e^{x_6})}.
\end{equation}
$P'_K( e^{x_1}, e^{x_2}, \cdots, e^{x_6} )$
is uniquely determined by a knot $K$
(hence, is an invariant of $K$)
in the completion of $\Q[\sigma_2, \sigma_3^2, \sigma_4]$.

\begin{prob}
Find a topological construction of the polynomial $P'_K$ given above.
\end{prob}

\begin{figure}[ht!]
\begin{align*}
& \pict{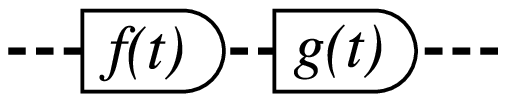}{0.6cm} = \pict{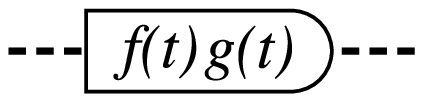}{0.6cm} \\
& \pict{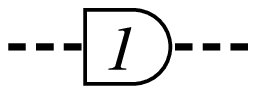}{0.55cm} = \pict{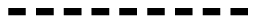}{0.08cm} \\
& \pict{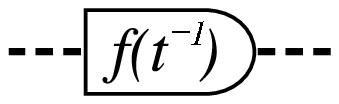}{0.6cm} = \pict{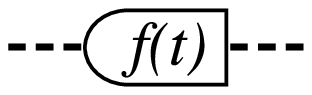}{0.57cm} \\
& \pict{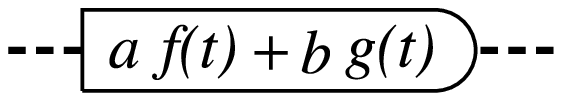}{0.6cm} = a \pict{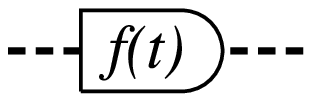}{0.6cm} + b \pict{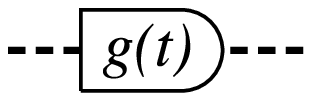}{0.6cm}
\end{align*}
\caption{\label{fig.ml_rel} The multi-linear relations.
Here, $f(t), g(t) \in S$, and $a$, $b$ are scalars.}
\end{figure}

\begin{figure}[ht!]
$$
\pict{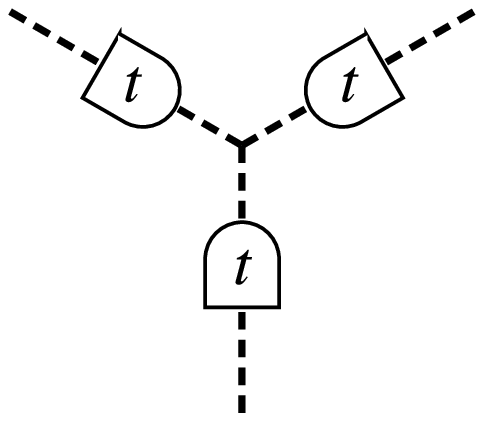}{2.5cm}
= \pict{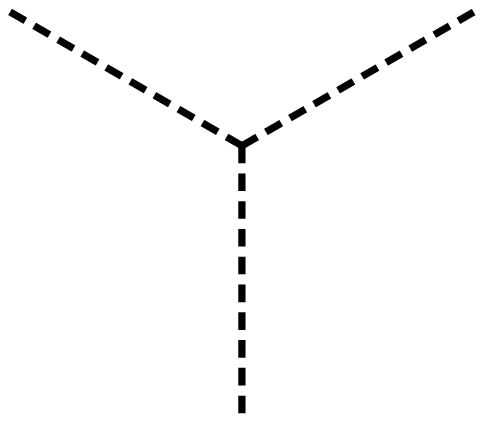}{2.5cm}
$$
\caption{\label{fig.push_rel} The push relation}
\end{figure}

The loop expansion in a general loop-degree is described as follows.
Let $R$ be a field, say $\Q$, and let $S$ be a subring of $R(t)$
which is invariant under the involution $t \mapsto t^{-1}$, 
where $t$ is an indeterminate.
A {\it labeled Jacobi diagram on $\emptyset$} is 
a vertex-oriented trivalent graph,
whose edges are labeled by pairs of local orientations and elements of $S$.
We define $\cA^S(\emptyset;R)$ to be the vector space over $R$
spanned by labeled Jacobi diagrams on $\emptyset$
subject to the AS, IHX, multilinear, and push relations
(see Figures \ref{fig.ml_rel} and \ref{fig.push_rel}).
The {\it loop-degree} of a labeled Jacobi diagram
is half the number of trivalent vertices of the Jacobi diagram.
For a polynomial $A(t)$ with $A(1)=1$ and $A(t)=A(t^{-1})$,
we have a map
\begin{equation}
\label{eq.lAtoB}
\cA^{\Q[t^{\pm1},1/A(t)]}(\emptyset;\Q) \longrightarrow \cB,
\end{equation}
defined by
$$
\pict{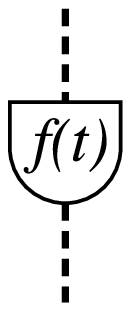}{2cm} \quad \longmapsto \quad c_0 \pict{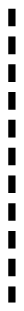}{2cm}
+ c_1 \pict{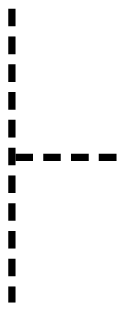}{2cm} + c_2 \pict{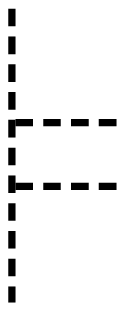}{2cm} + \cdots
+ c_n 
\begin{picture}(35,20)\put(0,5){\pict{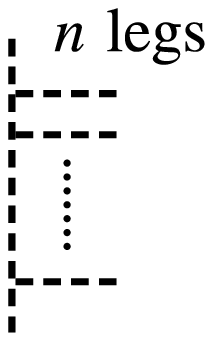}{2.15cm}}\end{picture} + \cdots,
$$
where
$f(t) \in \Q[t^{\pm1},1/A(t)]$ is written
$f(e^h) = \sum_{k=0}^\infty c_k h^k$.
In particular, the map
\begin{equation}
\label{eq.lAtoB2}
\cA^{\Q[t^{\pm1}]}(\emptyset;\Q) \longrightarrow \cB
\end{equation}
is defined by
$$
\pict{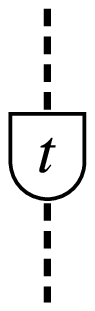}{2cm} \longmapsto \pict{ld/t2}{2cm}
+ \pict{ld/t3}{2cm} + \frac12 \pict{ld/t4}{2cm} + \cdots
+ \frac{1}{n!} 
\begin{picture}(35,20)\put(0,5){\pict{ld/t5}{2.15cm}}\end{picture} + \cdots.
$$
The loop expansion of the Kontsevich invariant
is described by the rational Z invariant
$Z^{\rm rat}(K) \in \cA^{\Q[t^{\pm1},1/\Delta_K(t)]}(\emptyset;\Q)$
which is taken to $\log_\sqcup \hZ(K)$
by the map (\ref{eq.lAtoB}).
In particular, when $\Delta_K(t)=1$, 
$Z^{\rm rat}(K) \in \cA^{\Q[t^{\pm1}]}(\emptyset;\Q)$.
(The existence of $Z^{\rm rat}(K)$ has been shown in \cite{Kricker_lK},
and the canonicality of $Z^{\rm rat}(K)$ has been shown in \cite{GK_rat}.)

\begin{prob}
Find a topological construction of 
the loop-degree $l$ part of 
the rational Z invariant\index{rational Z invariant}
$Z^{\rm rat}(K) \in \cA^{\Q[t^{\pm1},1/\Delta_K(t)]}(\emptyset;\Q)$
of a knot $K$, 
for each $l$.
\end{prob}

\begin{prob}
Find a basis of the space 
$\cA^{\Q[t^{\pm1},1/A(t)]}(\emptyset;\Q)^{({\rm loop}\ l)}$, 
for each $l$,
where $A(t)$ is a polynomial with $A(1)=1$ and $A(t)=A(t^{-1})$.
In particular, find a basis of the space
$\cA^{\Q[t^{\pm1}]}(\emptyset;\Q)^{({\rm loop}\ l)}$.
\end{prob}


\begin{conj}{\rm\cite{Rozansky_rc,GK_rat}}\qua
\label{conj.hair_map_inj}
The map (\ref{eq.lAtoB}) is injective.
In particular, the map (\ref{eq.lAtoB2}) is injective.
\end{conj}

\begin{rem}
If this conjecture is true,
$Z^{\rm rat}(K)$ is determined by the Kontsevich invariant.
\end{rem}

\subsection{The Kontsevich invariant of links in $\Sigma \times [0,1]$}
\label{sec.KinFxI}

\renewcommand{\thefootnote}{\fnsymbol{footnote}}
\footnotetext[0]{Section \ref{sec.KinFxI} was written by T. Kohno.}
\renewcommand{\thefootnote}{\arabic{footnote}}

Let $\Sigma$ be a closed oriented surface.
We denote by ${\mathcal A}_{\Sigma}$ the algebra of chord diagrams on $\Sigma$.
It is defined to be the vector space over $\C$ spanned by the homotopy
classes of continuous maps from chord diagrams to $\Sigma$ modulo 4T relations.

\begin{prob}[T. Kohno] 
\label{kohno_prob.1}
Construct explicitly a universal invariant of 
finite type\index{finite type invariant!--- of links in $\Sigma\times I$}
for links in $\Sigma \times [0, 1]$ with values in ${\mathcal A}_{\Sigma}$.
\end{prob}

\noindent
In the case of genus 0 the above problem is solved by 
Kontsevich integral.
In higher genus case a suggestion for a construction of a universal invariant
was given by Deligne at Oberwolfach meeting 1995.
In the case of a punctured surface the problem
was solved by Andersen, Mattes and Reshetikhin.

Let $G$ be a simple Lie group and
${\mathcal M}^G(\Sigma)$ the moduli space of $G$ flat connections on $\Sigma$.
The space of smooth functions on ${\mathcal M}^G(\Sigma)$ denoted by
$C({\mathcal M}^G(\Sigma))$ has a structure of a Poisson algebra 
coming from a symplectic structure on ${\mathcal M}^G(\Sigma)$.
The algebra ${\mathcal A}_{\Sigma}$ has also
a Poisson algebra structure (see \cite{kohno_AMR96}).
If each component of ${\mathcal A}_{\Sigma}$ is colored
by a representation of $G$, then
there is a natural Poisson algebra homomorphism
$$
\tau : {\mathcal A}_{\Sigma}
\rightarrow
C({\mathcal M}^G(\Sigma)).
$$
Problem \ref{kohno_prob.1} is related to the following problem.

\begin{prob}[T. Kohno] 
\label{prob.qPalg}
Give a deformation quantization\index{deformation quantization} 
of the Poisson algebra ${\mathcal A}_{\Sigma}$
which descends to a deformation quantization of $C({\mathcal M}^G(\Sigma))$.
\end{prob}

\noindent
The above problem will give a new insight on quantization of 
${\mathcal M}^G(\Sigma)$.
It would also be  interesting to investigate
a relation to the geometric quantization of ${\mathcal M}^G(\Sigma)$.

\begin{prob}[T. Kohno] 
\label{prob.qCMG}
Clarify the relation between 
a deformation quantization\index{deformation quantization}  
of $C({\mathcal M}^G(\Sigma))$ 
at a special parameter and 
the space of conformal blocks\index{conformal block} 
in WZW models.
\end{prob}

\begin{prob}[T. Kohno] 
\label{prob.ker_tau}
Determine the image and the kernel of the above map $\tau$.
\end{prob}

The space of conformal blocks in WZW model is defined
as the space of coinvariant tensors in the following way.
Let $p_1, \cdots, p_n$ be marked points on $\Sigma$ and $H_1, \cdots, H_n$ be
representations of the affine Lie algebra $\widehat{\frak g}$.
The space of conformal blocks is defined to be the set of linear forms
$$
\phi : H_1 \otimes \cdots \otimes H_n \longrightarrow \C
$$
invariant under the action of meromorphic functions 
with values in ${\frak g}$ with poles at most at $p_1, \cdots, p_n$, 
where the action is defined by the Laurent expansion at these points.
There is a twisted version of the above
construction, where the above meromorphic functions are
replaced by meromorphic sections of a ${\frak g}$ local system.

\begin{prob}[T. Kohno] 
\label{prob.hol_cb}
Compute the holonomy of 
the space of conformal blocks\index{conformal block} 
of the twisted WZW model.
In particular, determine the action of 
the braid group of $\Sigma$\index{braid group!--- of $\Sigma$}
on the space of conformal blocks for each $G$ flat connection on $\Sigma$.
\end{prob}

There is also a notion of the algebra of
chord diagrams on $n$ strings with
horizontal chord
on $\Sigma$, which we shall denote by
${\mathcal A}_n(\Sigma)$.

\begin{prob}[T. Kohno] 
\label{prob.PnS2AnS}
Let $P_n(\Sigma)$\index{braid group!--- of $\Sigma$}
denote the pure braid group of $\Sigma$ with $n$ strings. 
Does there exist an injective multiplicative homomorphism
$$
\theta : P_n(\Sigma) \rightarrow {\mathcal A}_n(\Sigma)
$$
defined over $\Q$?
\end{prob}

\newpage

\section{Skein modules}
\label{sec.skein_module}

\renewcommand{\thefootnote}{\fnsymbol{footnote}}
\footnotetext[0]{
The original version of Chapter \ref{sec.skein_module} 
was written by J. H. Przytycki.
It was revised by T. Ohtsuki following suggestions given by the referee.
Based on it, Przytycki wrote this chapter.}
\renewcommand{\thefootnote}{\arabic{footnote}}

Skein module is an algebraic object associated to a manifold, 
usually constructed as
a formal linear combination of embedded (or immersed) submanifolds,
modulo locally defined relations.
In a more restricted setting\index{skein module}
 a {\it skein module\/}\footnote{\footnotesize
Alexander first wrote down the skein relation for his polynomial.
Conway rediscovered the relation and placed in the abstract setting
of ''linear skein". He predicted the corresponding skein module for
a tangle.
General skein modules of 3-manifolds were first considered 
in 1987 by Przytycki and Turaev independently 
\cite{prz.P-4}, \cite{prz.Tu-1}.}
is a module associated to a 3-dimensional manifold, by considering
linear combinations of links in the manifold, modulo properly chosen
(skein) relations. It is a main object of the {\it algebraic topology
based on knots}. In the choice of relations one takes
into account several factors:
\begin{enumerate}
\item [(i)]  Is the module we obtain accessible (computable)?
\item [(ii)]
How precise are our modules in distinguishing 3-manifolds and links
in them?
\item [(iii)] Does the module reflect topology/geometry of a 3-manifold
(e.g.\ surfaces in a manifold, geometric decomposition of a manifold)?
\item [(iv)]
Does the module admit some additional structure (e.g.\ filtration,
gradation, multiplication, Hopf algebra structure)? Is it leading
to a Topological Quantum Field Theory (TQFT)\index{TQFT!--- and skein module}
by taking a finite dimensional quotient?
\end{enumerate}

One of the simplest skein modules is a $q$-deformation of the first
homology group of an oriented 3-manifold $M$, denoted by ${\cal S}_2(M;q)$.
It is based on the skein relation (between oriented framed links in $M$):
\vspace{-0.2pc}
$\pict{skein+}{0.7cm} = q \pict{skein0}{0.7cm}$; 
it also satisfies the framing relation 
\vspace{0.3pc}
$\begin{picture}(20,10)
\put(-5,8){\rotatebox{180}{\pict{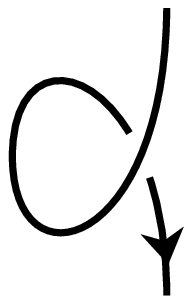}{0.9cm}}}\end{picture}
 = q 
\begin{picture}(17,10)
\put(0,8){\rotatebox{180}{\pict{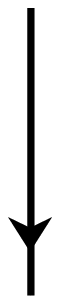}{0.9cm}}}\end{picture}$, 
where the diagrams in each formula imply framed links,
which are identical except in a ball,
where they differ as shown in the diagrams.
Already this simply defined skein module ``sees"
non-separating surfaces in $M$. These surfaces are responsible for
torsion part of the skein module \cite{prz.P-10}.

There is more general pattern: most of analyzed skein modules reflect
various surfaces in a manifold.

The best studied skein modules use skein relations
which worked successfully in the classical
knot theory (when defining polynomial invariants of links in $\R^3$).

\subsection{The Kauffman bracket skein module}

Let $M$ be an oriented 3-manifold,
and put $R=\Z[A^{\pm1}]$.
The {\it Kauffman bracket skein module}\index{skein module!Kauffman bracket ---}\index{Kauffman bracket!--- skein module|see {skein module}}  
$S_{2,\infty}(M)$ of $M$
is defined to be
the $R$ module spanned by
unoriented framed links in $M$ (including the empty link)
subject to the relations
\begin{align*}
& \!\! \pict{pos-cr}{1.2cm} \!\! 
= A  \!\! \pict{c-zero}{1.2cm} \!\! 
+ A^{-1} \!\! \pict{c-infty}{1.2cm} \!\! , \\
& \pict{O}{1cm} = -A^2 - A^{-2}, 
\end{align*}
where three diagrams in the first formula
imply three framed links,
which are identical except in a ball,
where they differ as shown in the diagrams.
The Kauffman bracket gives an isomorphism between 
$S_{2,\infty}(S^3)$ and $R$.
Thus, $S_{2,\infty}(M)$ is a generalization of the Jones polynomial
(in its Kauffman bracket interpretation).
The Kauffman bracket skein module
is best understood among the Jones type skein modules.
It can be interpreted as a quantization of the co-ordinate ring
of the character variety of $SL(2,\C)$ representations of the
fundamental group of the manifold $M$, 
\cite{prz.Bu-2,prz.P-S-0,prz.B-F-K,prz.P-S}.

\begin{prob}
Calculate\index{skein module!Kauffman bracket ---}  
$S_{2,\infty}(M)$ for each oriented 3-manifold $M$.
Find a convenient methodology to calculate it.
\end{prob}

\begin{rem}
It is known that
$S_{2,\infty}(L(p,q))$ of the lens space $L(p,q)$\index{lens space}
is a free $R$ module with $[p/2]+1$ generators \cite{HoPr_lens},
and that
$S_{2,\infty}(S^1 \times S^2) \cong 
R \oplus \bigoplus_{i=1}^\infty R/(1-A^{2i+4})$
\cite{HoPr_S1xS2}.
The Kauffman bracket skein modules are also calculated for
$I$-bundles over surfaces \cite{HoPr_dl,prz.P-4},
the exteriors of $(2,n)$ torus knots \cite{Bullock_tk},
and Whitehead manifolds \cite{HoPr_Wh}.
A connected sum formula is given in \cite{Prz_conn_sum}.
Skein modules at the 4th roots of unity are calculated in \cite{Sikora_4r}.
It is shown in \cite{Lofaro_MV} that
$S_{2,\infty}(M_1 \cup_F M_2$)
for orientable 3-manifolds $M_1$ and $M_2$ with a common boundary $F$
is expressed as a quotient module of a direct sum of tensor products of 
relative skein modules of $M_1$ and $M_2$.
\end{rem}


\begin{prob}[J. Przytycki] 
\label{prob.prz.it_tS2}
Incompressible tori and 2-spheres in $M$ yield 
torsion\index{torsion}   
in $S_{2,\infty}(M)$\index{skein module!Kauffman bracket ---} 
\cite{Prz_fund}.
It is a question of fundamental importance whether other surfaces
can yield torsion as well.
\end{prob}

\begin{conj} 
If every closed incompressible surface in $M$ is parallel to $\partial M$, 
then $S_{2,\infty}(M)$\index{skein module!Kauffman bracket ---}  
is torsion free.\index{torsion!--- free}  
\end{conj}

\begin{rem}
The Kauffman bracket skein module of
the 3-manifold obtained by an integral surgery along the trefoil knot
is finitely generated if and only if
the 3-manifold contains no essential surface \cite{Bullock_tre_sur}.
\end{rem}

The test case for the conjecture is the manifold 
$M=F_{0,3} \times S^1$, where $F_{0,3}$ is a 2-sphere with 3 holes,
because it contains immersed $\pi_1$-injective torus.

\begin{prob}[J. Przytycki]
\label{prob.prz.cS2}
Compute\index{skein module!Kauffman bracket ---}  
$S_{2,\infty}(F_{0,3} \times S^1)$.
\end{prob}

\begin{prob}
Let $F$ be a surface and $I$ an interval.
Describe the algebra 
$S_{2,\infty}(F \times I)$.\index{skein module!Kauffman bracket ---} 
\end{prob}

\begin{rem}
$S_{2,\infty}(F \times I)$ is an algebra (usually noncommutative). 
It is finitely generated
algebra for a compact $F$ \cite{prz.Bu-1}, 
and has no zero divisors \cite{prz.P-S}.
The center of the algebra is generated by boundary 
components of $F$ \cite{prz.B-P,prz.P-S}.
\end{rem}

\begin{prob}
Calculate the skein homology\index{skein homology}
based on the Kauffman bracket skein relation.\index{Kauffman bracket} 
\end{prob}

\begin{rem}
The skein homology were introduced in \cite{BuFrKa_sh}
(see also \cite{KaPrSi_sh}).
\end{rem}


\begin{prob}
We define the {\it $sl_3$ skein module}\index{skein module!$sl_3$ ---}  
$S^{sl_3}(M)$
of an oriented 3-manifold $M$
by the defining relations of the $sl_3$ linear skein 
\cite{Kuperberg_G2,OhYa_SU3}.
Calculate $S^{sl_3}(M)$ of each 3-manifold $M$.
\end{prob}

\begin{rem}
The quantum $sl_3$ invariant
 of links gives an isomorphism
between $S^{sl_3}(S^3)$ and the coefficient ring;
see, e.g.\ \cite{Ohtsuki_book}.
Thus, $S^{sl_3}(M)$ gives a generalization of 
the quantum $sl_3$ invariant of links.
\end{rem}

\subsection{The Homflypt skein module}

Let $M$ be an oriented 3-manifold,
and put $R=\Z[v^{\pm1},z^{\pm1}]$.
The {\it Homflypt skein module}\index{skein module!Homflypt ---} 
$S_3(M)$ of $M$
is defined to be the $R$ module spanned by oriented links in $M$ 
subject to the relation
$$
v^{-1} \pict{skein+}{1.2cm} - v \pict{skein-}{1.2cm} = z \pict{skein0}{1.2cm},
$$
where three diagrams in the formula imply three oriented links,
which are identical except in a ball,
where they differ as shown in the diagrams.
The Homflypt polynomial gives an isomorphism between $S_3(S^3)$ and $R$.
The Homflypt skein modules generalize skein modules
based on Conway relation which were hinted by Conway. 
$S_3(M)$ is related to the algebraic set of
$SL(n,\C)$ representations of the
fundamental group of the manifold $M$ \cite{prz.Si}.

\begin{prob}
Calculate\index{skein module!Homflypt ---}  
$S_3(M)$ for each oriented 3-manifold $M$.
Find a convenient methodology to calculate it.
\end{prob}

\begin{rem}
It is known that
$S_3(F \times I)$ is an infinitely generated free module \cite{prz.P-6},
and that
$S_3(S^1 \times S^2)$ is isomorphic to the direct sum of $R$
and an $R$-torsion module \cite{GiZh_Hs_S1xS2}.
A connected sum formula is given in \cite{GiZh_conn_sum}.
\end{rem}

\begin{prob}
Let $F$ be a surface and $I$ an interval.
Describe the algebra $S_3(F \times I)$.\index{skein module!Homflypt ---} 
\end{prob}

\begin{rem}
$S_3(F \times I)$ is a Hopf algebra 
(usually neither commutative nor co-commutative)
\cite{prz.Tu-2,prz.P-6}. 
$S_3(F\times I)$ is a free module (as mentioned above)
and can be interpreted as a quantization 
\cite{prz.Tu-1,prz.H-K,prz.Tu-2,prz.P-5}.
\end{rem}



\subsection{The Kauffman skein module}

Let $M$ be an oriented 3-manifold,
and put $R=\Z[a^{\pm1},x^{\pm1}]$.
The {\it Kauffman skein module}\index{skein module!Kauffman ---}  
$S_{3,\infty}(M)$ of $M$
is defined to be
the $R$ module spanned by
unoriented framed links in $M$ subject to the relations
\begin{align}
\label{eq.Ksm_rel1}
& \pict{bigelow/bmw1}{1cm} 
+ \pict{bigelow/bmw2}{1cm} 
= x \left( \pict{bigelow/bmw3}{1cm} 
+ \pict{bigelow/bmw4}{1cm} \right), \\*
\label{eq.Ksm_rel2}
& \pict{bigelow/RI1}{1cm} 
= a \ \pict{bigelow/RI2a}{1cm} \ , 
\end{align}
where the diagrams in each formula
imply framed links,
which are identical except in a ball,
where they differ as shown in the diagrams.

\begin{prob}
Calculate\index{skein module!Kauffman ---}   
$S_{3,\infty}(M)$ for each oriented 3-manifold $M$.
Find a convenient methodology to calculate it.
\end{prob}

\begin{rem}
$S_{3,\infty}(F \times I)$ is known to be a free module.
The case of $F$ being a torus was solved by Hoste, Kidwell and Turaev.
It is calculated in \cite{Lieberum_sc} for a surface $F$ with boundary.
$S_{3,\infty}(S^1 \times S^2)$ is calculated in \cite{ZhLu_S1xS2}.
A connected sum formula is given in \cite{Zhong_conn_sum}.
\end{rem}

\begin{prob}
Calculate the higher skein modules based on the Kauffman skein relation
$W_i^{3,\infty}(M)$ and $\hat W^{3,\infty}(M)$
(see below for their definitions).
\end{prob}

\begin{rem}
The higher skein modules were introduced in \cite{prz.P-7}.
They are discussed (in the case of the Conway skein triple) in 
\cite{Rong,Rong-Lick} and 
\cite{AnTu_hsm,AnTu_hsm2}.
In the case of the Kauffman
skein relation, definitions are as follows:  
Let $R\LL$ denote the free $R$ module spanned by
the ambient isotopy classes of
unoriented framed links in an oriented 3-manifold $M$ 
modulo the framing relation (\ref{eq.Ksm_rel2}),
where $R=Z[a^{\pm 1},x^{\pm 1}]$.
We regard singular links with a finite number of double points
as elements in $R\LL$
by replacing a double point with the difference of
the two sides of (\ref{eq.Ksm_rel1}).
We introduce a (singular links) filtration 
$R\LL = C_0 \supset C_1 \supset C_2 \supset C_3 \supset \cdots$,
where the module $C_i$ is generated by singular links with $i$ double points.
We define the $i$th higher Kauffman skein module
as: $W_i^{3,\infty}(M) = R\LL/C_{i+1}$ and the completed
higher Kauffman skein module, $\hat W^{3,\infty}(M)$, 
as the completion of $R\LL$ with respect to the filtration $\{C_i\}$.
\end{rem}



\begin{prob}
Construct invariants of 3-manifolds
via a linear skein theory based on the 
Kauffman skein module.\index{skein module!Kauffman ---}  
\end{prob}

\begin{rem}
It is known that
quantum invariants
 of 3-manifolds can be constructed via linear skein theories
based on the Kauffman bracket skein modules (see \cite{Lickorish_book}) 
and the Homflypt skein modules \cite{Yokota_SUNskein}.
\end{rem}

\begin{update}
Beliakova and Blanchet have done this \cite{BeBl01}.
\end{update}

\subsection{The $q$-homotopy skein module}

Let $M$ be an oriented 3-manifold,
and put $R=\Z[q^{\pm1},z]$.
The {\it $q$-homotopy skein module}\index{skein module!homotopy ---} 
$HS^q(M)$ of $M$
is defined to be the $R$ module spanned by oriented links in $M$ 
subject to the link homotopy relation
\vspace{-0.4pc}
$\pict{skein+}{0.8cm} \!=\! \pict{skein-}{0.8cm}$
for self-crossings
and the skein relation
\vspace{-0.4pc}
$q^{-1} \pict{skein+}{0.8cm} \!- q \pict{skein-}{0.8cm} 
\!= z \pict{skein0}{0.8cm}$
for ``mixed crossings'',
i.e.\
we assume that the two strings of
\vspace{-0.2pc}
$\pict{skein+}{0.8cm}$
(or $\pict{skein-}{0.8cm}$)
of the skein relation belong to different components of the link.

We have an isomorphism between $HS^q(S^3)$
and $\Z[q^{\pm1},t,z]$,
regarding $t^k$ as the trivial link with $k$ components,
and this isomorphism is given by the linking numbers \cite{prz.P-11}.

\begin{prob}
Calculate\index{skein module!homotopy ---}  
$HS^q(M)$ for each 3-manifold $M$.
\end{prob}

\begin{rem}
$HS^q(F \times I)$ is a quantization \cite{prz.H-P-1,prz.Tu-2,prz.P-11}, 
and as noted by Kaiser it can be almost completely understood using
singular tori technique of X.-S. Lin.
$HS^q(M)$ is free if and only if
$\pi_1(M)$ is abelian and $2 b_1(M) = b_1(\partial M)$ \cite{Kaiser_p}.
\end{rem}



\subsection{The $(4,\infty)$ skein module}

\begin{figure}[ht!]
$$ 
\pict{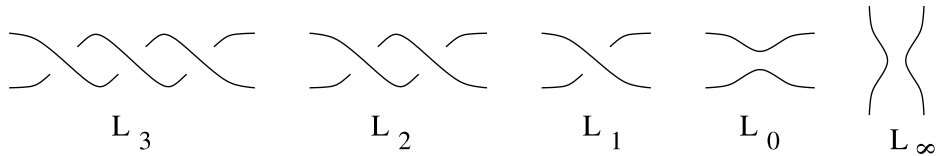}{2cm} 
$$
\caption{\label{p_fig.1}
$L_3, \cdots, L_0, L_\infty$ are framed links
which are identical except in a ball, 
where they differ as shown in the pictures.
Links $L_k$ for $k=4,5,\cdots$ are similarly defined.}
\end{figure}

We generalize the Kauffman bracket and Kauffman skein modules
by considering the general, unoriented skein relation
$b_0L_0 + b_1L_1 + \cdots + b_{n-1}L_{n-1} + b_{\infty}L_{\infty}$
(see Figure \ref{p_fig.1}).
The first new case to analyze, $n=4$, is described in this section.
We call it the $(4,\infty)$ skein module and denote by 
${\cal S}_{4,\infty}(M;R)$.
This problem is very interesting even for $M=S^3$.

The definitions are as follows.
Let $M$ be an oriented 3-manifold, 
${\cal L}_{fr}$ the set of unoriented framed links in $M$ 
(including the empty knot, $\emptyset$) and
$R$ any commutative ring with unity.  
We fix $a,b_0,b_3$ to be invertible elements in $R$ and
fix $b_1,b_2,b_{\infty}$ to be elements of $R$.
Then we define the 
{\it $(4,\infty)$ skein module}\index{skein module!$(4,\infty)$ ---}  
as:
${\cal S}_{4,\infty}(M;R) = R{\cal L}_{fr}/I_{(4,\infty)}$,
where $I_{(4,\infty)}$ is the submodule of $R{\cal L}_{fr}$
generated by the following two relations:
\begin{align*}
\mbox{the $(4,\infty)$ skein relation: } \quad
& b_0L_0  + b_1L_1 + b_2L_2 + b_3L_3 + b_{\infty}L_{\infty} = 0, \\
\mbox{the framing relation: } \quad\quad\ \ \, 
& L^{(1)} = a L,
\end{align*}
where $L_0, \cdots, L_\infty$ are framed links
which are identical except in a ball, 
where they differ as shown in Figure \ref{p_fig.1},
and $L^{(1)}$ denotes a link 
obtained from $L$ by adding $+1$ framing to some component of $L$.

\begin{prob}[J. Przytycki] 
\label{p_4.1}
\mbox{}
\begin{enumerate}
\item[\rm(i)]
Find generators of ${\cal S}_{4,\infty}(S^3,R)$.\index{skein module!$(4,\infty)$ ---}  
\item[\rm(ii)] 
For which parameters of the $(4,\infty)$ skein and framing relations, 
trivial links are linearly independent in ${\cal S}_{4,\infty}(S^3;R)$?
\item[\rm(iii)] 
For which parameters of the $(4,\infty)$ skein and framing relations, 
the trivial knot is not representing
a torsion\index{torsion}   
element of ${\cal S}_{4,\infty}(S^3,R)$?
\end{enumerate}
\end{prob}

A generalization of the Montesinos-Nakanishi conjecture \cite{prz.P-Ts}
said that\break ${\cal S}_{4,\infty}(S^3,R)$ is generated by
trivial links and that 
the $(4,\infty)$ skein module (suitably defined) for $n$-tangles 
is generated by $\prod_{i=1}^{n-1} (3^i +1)$ certain basic $n$-tangles.
This would give a generating set for the $(4,\infty)$ skein module 
of $S^3$ or $D^3$ with $2n$ boundary points (for $n$-tangles). However,
the Montesinos-Nakanishi 3-move conjecture has been disproved by
M.Dabkowski and J.H.Przytycki in February 2002 
\cite{DaPr_3m} and \cite{Prz_smd}.
Therefore $\prod_{i=1}^{n-1} (3^i +1)$ is only the lower bound
for the number of generators.

In \cite{prz.P-Ts} we extensively analyze the possibilities
that trivial links are linearly independent;
if $b_{\infty} = 0$, then this may happen only if $b_0b_1=b_2b_3$.
These leads to the following conjecture (cases (1)--(2)):

\begin{conj}[J. Przytycki{\rm , see \cite{Morton_Hellas}}] 
\label{3.2}
\mbox{}
\begin{enumerate}
\item [\rm(1)] There is a polynomial invariant of unoriented links,
$P_1(L) \in \Z[x,t]$ which satisfies:
\begin{enumerate}
\item [\rm(i)] Initial conditions: $P_1(T_n) = t^n$, where $T_n$ is a trivial
link of $n$ components.
\item [\rm(ii)] Skein relation $P_1(L_0) + x P_1(L_1) - x P_1(L_2) - P_1(L_3)=0$
where $L_0,L_1,L_2,L_3$ is a standard, unoriented skein quadruple
($L_{i+1}$ is obtained from $L_{i}$ by a right-handed half twist on
two arcs involved in $L_{i}$; compare Figure \ref{p_fig.1}.)
\end{enumerate}
\item [\rm(2)] There is a polynomial invariant of unoriented framed links,
$P_2(L) \in \Z[A^{\pm 1},t]$ which satisfies:
\begin{enumerate}
\item [\rm(i)] Initial conditions: $P_2(T_n) = t^n$,
\item [\rm(ii)] Framing relation: $P_2(L^{(1)})=-A^3P_2(L)$ where $L^{(1)}$ is
obtained from a framed link $L$ by a positive half twist on its framing.
\item [\rm(iii)] Skein relation: $P_2(L_0) + A(A^2 + A^{-2})P_2(L_1) +
(A^2 + A^{-2})P_2(L_2) + A P_2(L_3)=0$.
\end{enumerate}
\item [\rm(3)] There is a rational function invariant of unoriented framed links,
$P_3(L) \in \Z[a^{\pm 1},x,y, (x+y+xy+y^2)^{-1}]$ which satisfies:
\begin{enumerate}
\item [\rm(i)] Initial conditions: $P_3(T_n) =
(\frac{ -a^{3}(x+y+x y +x^2) + a^7(x+y+1)^2 -a^{-1} }
{x + y + x y+y^2})^{n-1}$,
\item [\rm(ii)] Framing relation: $P_3(L^{(1)})=a P_3(L)$,
\item [\rm(iii)] Skein relation: $P_3(L_0) + a x P_3(L_1) +
a^2y P_3(L_2) - a^3(x+y+1)P_3(L_3)=0$.
\end{enumerate}
\item [\rm(4)] 
The invariant predicted in (1) (respectively (2) and (3))
is not uniquely defined (if it exists).
\end{enumerate}
\end{conj}

Note that a solution to (3) becomes a solution to (1) 
under the substitution $a=1$, $x=-y$ 
and that a solution to (3) becomes a solution to (2) under the substitution 
$a=-A^3$, $x=-1-A^{-4}$, $y=A^{-4}+ A^{-8}$.
As for the uniqueness of (4),
note that all such invariants agree 
on trivial links and therefore they agree on the space spanned by trivial 
links in the related cubic skein module.

\begin{figure}[ht!]
$$
\pict{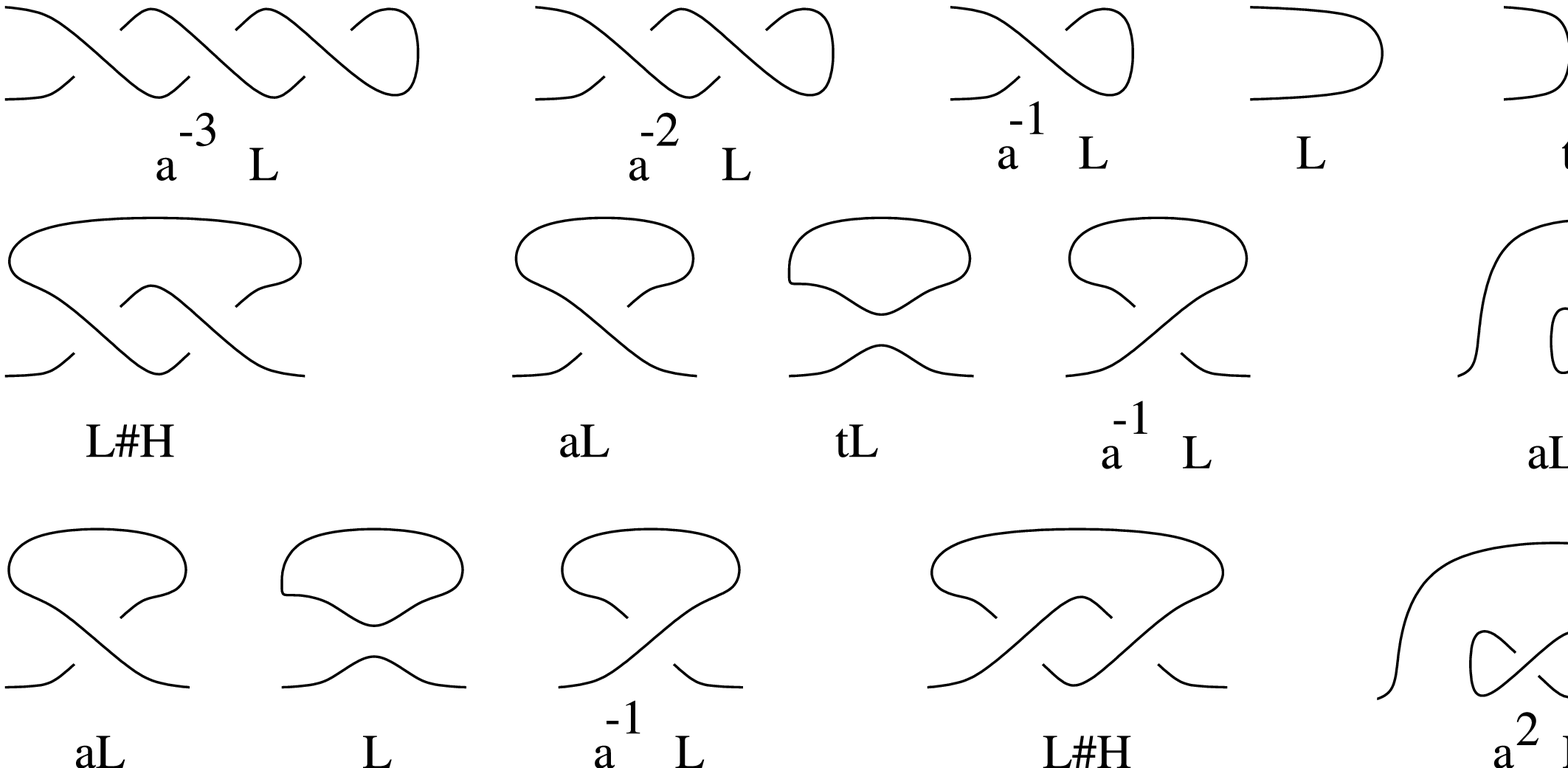}{5.5cm}
$$
\nocolon\caption{}\label{p_fig.2}
\end{figure}

The above conjectures assume that $b_{\infty}=0$ in our
skein relation. Let consider the possibility that
$b_{\infty}$ is invertible in $R$. Using the ``denominator"
of our skein relation (the first line of Figure \ref{p_fig.2})
we get the relation which
allows to compute the effect of adding a trivial component
to a link $L$ (we write $t^n$ for the trivial link $T_n$):
\begin{equation}
\label{p_eq.1}
(a^{-3}b_3 + a^{-2}b_2 + a^{-1}b_1 + b_0 + b_{\infty}t)L=0.
\end{equation}
When considering the ``numerator" of the relation and its mirror image
(Figure \ref{p_fig.2})
we obtain formulas for Hopf link summands, 
and because unoriented Hopf link is amphicheiral 
we can eliminate it from our equations to get the formula (\ref{p_eq.2}):
\begin{align}
& b_3(L\#H) + (a b_2 + b_1t +a^{-1}b_0 + a b_{\infty})L =0. \notag \\
& b_0(L\#H) + (a^{-1}b_1 +b_2t + a b_3 + a^2b_{\infty})L =0. \notag \\
\label{p_eq.2}
& ((b_0b_1 - b_2b_3)t + (a^{-1}b_0^2 -a b_3^2) +
(a b_0b_2 - a^{-1}b_1b_3) + b_{\infty}(a b_0 -a^2b_3))L = 0.
\end{align}

It is possible that (\ref{p_eq.1}) and (\ref{p_eq.2}) 
are the only relations in the module.
Precisely, we ask whether ${\cal S}_{4,\infty}(S^3;R)$
is the quotient ring $R[t]/({\cal I})$ where $t^i$ represents
the trivial link of $i$ components and ${\cal I}$ is the
ideal generated by (\ref{p_eq.1}) and (\ref{p_eq.2}) for $L=t$.
The substitution which realizes the relations is:
$b_0=b_3=a=1$, $b_1=b_2=x$, $b_{\infty}=y$. This may lead to the
polynomial invariant of unoriented links in $S^3$ with values in $\Z[x,y]$
and the skein relation $L_3 +xL_2 + xL_1 + L_0 + yL_{\infty} = 0$.

\begin{prob}[J. Przytycki] 
\label{prob.skeinmod_7col}
For which coefficients of the $(4,\infty)$ skein relation
is the number of Fox $7$-colorings measured by the\index{skein module!$(4,\infty)$ ---}   
$(4,\infty)$ skein module?
\end{prob}


\begin{rem} 
We denote by $Col_p(L)$ the ($\Z/p\Z$)-linear space (for $p$ prime)
of Fox $p$-colorings of a link $L$
(for its definition, see \cite{Przytycki_c})
and $col_p(L)$ denotes the cardinality of the space.
It is known that $Col_p(L)$
can be identified with $H_1 \big( M_2(L); \Z/p\Z \big)$,
where $M_2(L)$ denotes the double cover of $S^3$ branched along $L$.
Since the double covers of tangles
defining $L_0, L_1, \cdots, L_{p-1}, L_\infty$
give all subspaces of $H_1(T^2;\Z/p\Z)$ respectively
(where $T^2$ is the double cover of $(S^2, \mbox{ 4 points})$),
$col_p$ of those links
are equal except for $col_p$ of one link
which is equal to $p$ times the others
\cite{Przytycki_c}. 
This leads to the relation of type $(p,\infty)$.
A relation between the Jones polynomial\index{Jones polynomial!--- and coloring}
 (or the Kauffman bracket)
and $col_3(L)$ has the form: 
$col_3(L) = 3|V_L(e^{\pi \sqrt{-1}/3})|^2$ and
a formula relating the Kauffman polynomial 
and $col_5(L)$ has the form: 
$col_5(L) = 5|F_L(1,e^{2\pi \sqrt{-1}/5} + e^{-2\pi \sqrt{-1}/5})|^2$.
This seems to suggest 
the existence of a similar formula\footnote{\footnotesize
Fran{\c c}ois Jaeger told Przytycki that 
he knew how to get the space of Fox $p$-colorings from 
a short skein relation (of type $(\frac{p+1}{2},\infty)$).
Fran{\c c}ois died prematurely in 1997 and his proof has never been recorded.}
for $col_7(L)$. 
\end{rem}

\subsection{Other problems}


We extend the family $\K$ of oriented knots in a 3-manifold $M$
by singular knots, and resolve a singular crossing by 
$\pict{d-pt1}{0.7cm} 
\!\!=\!\! \pict{d-pt2}{0.7cm} \!\!-\!\! \pict{d-pt3}{0.7cm}$.
These allows us to define the Vassiliev-Goussarov filtration:
$R{\K} = C_0 \supset C_1  \supset C_2  \supset C_3 \cdots$,
where $R$ is a commutative ring with unity and
$C_k$ is generated by knots with $k$ singular points.
Regarding the quotient $W_k(M) = R{\K}/C_{k+1}$
as an invariant of $M$,
we call it the {\it $k$th Vassiliev-Goussarov skein module} of $M$.\index{skein module!Vassiliev-Goussarov ---}  
The completion of the space of knots
with respect to the Vassiliev-Goussarov filtration, $\hat{R{\K}}$,
 is a {\it Hopf algebra} (for $M=S^3$). Functions dual to
Vassiliev-Goussarov skein modules are called 
{\it finite type} or
{\it Vassiliev invariants}\index{Vassiliev invariant!--- of knots in a 3-manifold}
 of knots; see \cite{prz.P-7}.

\begin{prob}
Calculate $W_k(M)$ for each 3-manifold $M$.
\end{prob}

\begin{rem}
When $M=S^3$, and coefficients are from $\Q$ then
the graded space $C_k/C_{k+1}$
can be described by chord diagrams of degree $k$;
see Chapter \ref{sec.Vinv}.
\end{rem}

\begin{prob}
Define a skein module of 3-manifolds,\index{skein module!--- of 3-manifolds} 
and calculate it.
\end{prob}

\begin{rem}
The quantum Hilbert space (or the space of conformal blocks)
of $(S^2, \mbox{ 4 points})$ is known to be finite dimensional.
This is a reason why a quantum invariant
 of links
satisfies a skein relation;
it is a linear relation of tangles bounded by $(S^2, \mbox{ 4 points})$
whose invariants are linearly dependent in the quantum Hilbert space.
The quantum Hilbert space of a closed surface, say, a torus, 
is also known to be finite dimensional.
Hence, a quantum invariant of 3-manifolds satisfies a ``skein relation'';
it should be a linear relation of 3-manifolds bounded by a surface.
A skein module of 3-manifolds might be defined
to be a module spanned by closed oriented 3-manifolds
subject to a suitably chosen ``skein relation'' among 3-manifolds.
It is a problem to define such a skein module which can be calculated.
\end{rem}

\newpage

\section{Quandles}
\label{sec.quandle}

\renewcommand{\thefootnote}{\fnsymbol{footnote}}
\footnotetext[0]{Chapter \ref{sec.quandle} was written by T. Ohtsuki,
following suggestions and comments given by S. Kamada and M. Saito.
Section 5.6 was added by C. Rourke and B.Sanderson.}
\renewcommand{\thefootnote}{\arabic{footnote}}

A {\it quandle}\index{quandle} 
is a set $X$ equipped with a binary operation $\ast$
satisfying the following 3 axioms.
\begin{itemize}
\item[\rm(1)]
$x \ast x = x$ for any $x \in X$.
\item[\rm(2)]
For any $y,z \in X$ 
there exists a unique $x \in X$ such that $z = x \ast y$.
\item[\rm(3)]
$(x \ast y) \ast z = (x \ast z) \ast (y \ast z)$
for any $x,y,z \in X$.
\end{itemize}
The notions of subquandle, homomorphism, isomorphism, automorphism
are appropriately defined.
Each $x$ in a quandle $X$ defines a map $S_x: X \to X$ by $S_x(y) = y \ast x$.
This map is an automorphism of $X$
by the axioms (2) and (3).
The {\it inner automorphism group} 
is a group of automorphisms generated by $S_x$ ($x \in X$).
An orbit under the action of the inner automorphism group on a quandle $X$
is simply called an {\it orbit} of $X$.
This forms a subquandle of $X$.
A quandle is called
{\it connected\/}\footnote{\footnotesize 
We call this property {\it connected} here following \cite{Joyce}.
This is also called {\it weakly homogeneous} in some of the literature.}
if the action of its inner automorphism group is transitive on it
(i.e.\ if $X$ has only one orbit).
A quandle is called {\it simple}
if every surjective homomorphism from the quandle is either
an isomorphism or the constant map to the one-element quandle.
The {\it dual quandle} of $X$
is the set $X$ with the dual binary operation given by
$x \overline{\ast} y = S_y^{-1}(x)$.

The {\it conjugation quandle} of a group
is the group with the binary operation $x \ast y = y^{-1} x y$.
This kind of quandle is a prototype of quandles;
the defining relations of a quandle are relations
satisfied by the conjugation of a group.
Any conjugacy class of a group
is a subquandle of the conjugation quandle of the group.
The {\it dihedral quandle} $R_n$ of order $n$ is 
the subquandle of the conjugation quandle of the dihedral group of order $2n$,
consisting of reflections.
An {\it Alexander quandle} is
a quotient module $\Z[t^{\pm1}]/J$,
where $t$ is an indeterminate and $J$ is an ideal of 
the Laurent polynomial ring $\Z[t^{\pm1}]$,
equipped with the binary operation $x \ast y = t x + (1-t) y$.
The dihedral quandle $R_n$ is isomorphic to $\Z[t^{\pm1}]/(n,t+1)$.

\subsection{Classification of quandles}

It was a classical problem in group theory
to classify the isomorphism classes of groups of order $n$ for each $n$.
The following problem is a corresponding problem for connected quandles.

\begin{prob}
\label{prob.cl-qdl}
Classify\index{quandle!classification of connected ---}
the isomorphism classes of connected quandles of order $n$ 
for each positive integer $n$.
\end{prob}

See Table \ref{tbl.quandle}
for a list of connected quandles of order $n$ for some $n$.

\begin{table}[ht!]
\begin{center}
\begin{tabular}{|c|c|l|l|} \hline 
\dwn{\small $n$} & \dwn{\small $\#$} & 
\multicolumn{2}{c|}{\small Connected quandles of order $n$} \\ 
\cline{3-4}
& &\small  Self-dual &\small  Not self-dual \\
\hline\hline
\small 1 &\small  $1$ &\small  A trivial quandle & \\ \hline
\small 2 &\small  $0$ & & \\ \hline
\small 3 &\small  $1$ &\small  $R_3$ & \\ \hline
\small 4 &\small  $1$ &\small  $\Lambda_2/(t^2+t+1)$ & \\ \hline
\small 5 &\small  $3$ &\small  $R_5$ &\small  
$\Lambda_5/(t-2)$, its dual 
\\ \hline
\small 6 & \small  $2$
& \small 2 subquandles of $\mbox{\small Conj}({\frak S}_4)$ & \\ \hline
\small 7 &\small  $5$ &\small  $R_7$ & 
\small $\Lambda_7/(t-2)$, $\Lambda_7/(t-3)$, their duals
\\ \hline
\dwn{\small 8} & \dwn{\small $\ge 3$} &
\small An abelian extension &
\dwn{\small $\Lambda_2/(t^3+t+1)$, its dual} \\ 
& & of $\Lambda_2/(t^2+t+1)$ & \\ \hline
\dwn{\small 9} & \dwn{\small $8$} &
\small $R_9$, $\Lambda_3/(t^2-t+1)$, & 
\small $\Lambda_9/(t-2)$, $\Lambda_3/(t^2+t-1)$,
\\
& & \small $R_3 \times R_3$, $\Lambda_3/(t^2+1)$ & \small their duals
\\ \hline
\small 10 & \small $\ge 1$ & 
\small A subquandle of $\mbox{\small Conj}({\frak S}_5)$ &
\\ \hline
\small 11 & \small 9 & \small $R_{11}$ &
\small $\Lambda_{11}/(t-a)$ ($a=2,3,\cdots,9$) 
\\ \hline
\dwn{\small 12} & \dwn{\small $\ge 2$} 
&\small $R_3 \times \big(\Lambda_2/(t^2+t+1)\big)$, & \\
& & \small An icosahedral quandle & \\ \hline
\small 13 & \small 11 & \small $R_{13}$ &
\small $\Lambda_{13}/(t-a)$ ($a=2,3,\cdots,11$) \\ \hline
\small 14 & \small $\ge 0$ & & 
\\ \hline
\dwn{\small 15} & \dwn{\small $\ge 4$} & \small $R_3 \times R_5$, 
& \dwn{\small $R_3 \times \big(\Lambda_5/(t-2)\big)$, its dual} \\
& & \small A subquandle of $\mbox{\small Conj}({\frak S}_5)$ & 
\\ \hline
$\vdots$ & & & \\ \hline
\small Prime $p$ &\small  $p-2$ &\small  $R_p$ &\small  
$\Lambda_p/(t-a)$ ($a=2,3,\cdots,p-2$) \\ \hline
\end{tabular}
\end{center}
\caption{\label{tbl.quandle}
A table of some connected quandles.
The second column shows
the numbers of isomorphism classes of connected quandles of order $n$.
We denote $\Z[t^{\pm1}]/(n)$ by $\Lambda_n$.
$\mbox{\footnotesize Conj}({\frak S}_m)$ denotes
the conjugation quandle of the $m$th symmetric group ${\frak S}_m$.
An {\it icosahedral quandle}
is a quandle whose elements are the vertices of an icosahedron
such that $S_x$ of each element $x$ is given by a rotation of 
the icosahedron centered at $x$.}
\end{table}


%

\begin{rem}[\rm (M. Gra\~na)] 
It is shown, in \cite{EGS00} in terms of set theoretical solutions
of the quantum Yang-Baxter equation,
that a connected quandle of prime order $p$ is
isomorphic to the Alexander quandle $\Z[t^{\pm}]/(p,t-a)$ for some $a$.
It is shown in \cite{AnGr02,Nelson} that
two connected Alexander quandles are isomorphic if and only if
they are isomorphic as $\Z[t^{\pm1}]$-modules.
These give the classification of connected quandles of prime order
shown in Table \ref{tbl.quandle}.
\end{rem}

\begin{rem}[\rm (M. Gra\~na)] 
It is shown in \cite{AnGr02} that
a simple quandle of prime power order is an Alexander quandle;
it is a finite field $F$ such that $t$ acts by multiplication by some
primitive element $w$ (i.e.\ $w$ generates $F$ as an algebra).
Further, it is shown in \cite{Grana_p2} that
a connected quandle of prime square order is an Alexander quandle.
This gives the classification of connected quandles of order 9
shown in Table \ref{tbl.quandle}.
\end{rem}

\begin{rem}
S. Yamada made the list of isomorphism classes of 
quandles (and racks) of order $\le 7$ 
by computer search.
The list of connected quandles of order $\le 7$ in Table \ref{tbl.quandle}
follows from it.
S. Nelson \cite{Nelson} classifies the Alexander quandles of order $\le 15$;
connected ones among them are listed in Table \ref{tbl.quandle}.
\end{rem}

\begin{rem}
The following modification of Problem \ref{prob.cl-qdl} 
gives an algorithm to list connected quandles:
classify the isomorphism classes of connected quandles
fixing the conjugacy class of
the union of $S_x$ and an identity map.
For a quandle $X$
we denote by $S_X$ the set of $S_x$ ($x \in X$),
which is regarded as a subset of ${\frak S}_n$
when $X$ is of order $n$.
The map $X \to {\frak S}_n$, taking $x \mapsto S_x$,
is often injective,
though in general the map $X \to S_X$ is a quotient map,
and the order of $S_X$ divides $n$ when $X$ is connected.
Let us investigate this problem in some simple cases.

Let $X$ be a connected quandle of order $n$
whose $S_X$ includes $(12) \in {\frak S}_n$.
Then, for any $i$ there is a sequence $1=a_0, a_1, \cdots, a_k=i$
such that $(a_0 a_1), (a_1 a_2), \cdots \in S_X$
since $X$ is connected.
Further, since $S_X$ is closed with respect to conjugation,
$S_X$ includes $(1 i) \in {\frak S}_n$,
and hence any $(i j) \in {\frak S}_n$.
Therefore, $n=3$, and 
$X$ is isomorphic to the dihedral quandle $R_3$.

Let $X$ be a connected quandle of order $n$
whose $S_X$ includes $(123) \in {\frak S}_n$.
Suppose that $S_X$ further included $(145) \in {\frak S}_n$.
Then, since $S_X$ is closed with respect to conjugation,
$S_X$ would include $(i j k) \in {\frak S}_n$
for any $\{ i,j,k \} \subset \{ 1,2,3,4,5 \}$.
This would contradict,
since the order of $S_X$ is at most $n$.
Hence, $n=4$, and
$X$ is isomorphic to the conjugation subquandle of ${\frak A}_4$
consisting of $(123)$, $(134)$, $(142)$, and $(243)$,
which is isomorphic to $\Z[t^{\pm1}]/(2,t^2+t+1)$.

Let $X$ be a connected quandle of order $n$
whose $S_X$ includes $(1234) \in {\frak S}_n$.
If $S_X$ further included $(1567) \in {\frak S}_n$,
a contradiction would follow
from a similar argument as above.
Hence, it is sufficient to consider the cases that
$S_X$ include $(1234)$ and
either of $(1256)$, $(2156)$, $(1526)$, $(1536)$, or
$(i j k 5)$ for any $\{ i,j,k \} = \{ 1,2,3 \}$.
It follows from some concrete computations that
such a $X$ is isomorphic to
either of the Alexander quandle
$\Z[t^{\pm1}]/(5,t-2)$, its dual quandle, 
or the conjugation subquandle of ${\frak S}_4$ including $(1234)$.
\end{rem}


\subsection{Representations of knot quandles}

Consider the conjugation quandle of 
the fundamental group $\pi_1(S^3-K)$
of the complement of a knot $K$.
The {\it reduced knot quandle} $\hat Q(K)$ is its subquandle
generated by meridians of $K$.
A\index{quandle!knot ---} 
{\it knot quandle\/}\footnote{\footnotesize 
Knot quandle was introduced by Joyce \cite{Joyce} 
and independently by Matveev \cite{Matveev_q};
see \cite{FR_rack} for an exposition.}
$Q(K)$ is a quandle
generated by meridians of $K$
(for its precise definition, see \cite{Joyce})
which is almost equal to $\hat Q(K)$;
to be precise,
there is a surjective (almost, bijective) homomorphism $Q(K) \to \hat Q(K)$.

Homomorphisms to a fixed group/quandle
are often called {\it representations}.
It was said, before quantum invariants
 were discovered, that
to count the numbers of representations of knot groups to a fixed finite group
was a most powerful method to distinguish two given knots.
The following problem is a refinement of it.
A motivation is 
to construct a methodology to count
the number of representations of a knot quandle
to a fixed quandle of finite order.

\begin{prob}
\label{prob.rep_QK}
Describe (the number of) representations\index{quandle!knot ---!representations of ---}
of a knot quandle to a fixed connected quandle of finite order,
say, by using knot invariants known so far,
or by reducing the problem to the case of smaller target quandles.
\end{prob}

\begin{rem}
Since a knot quandle is connected,
the image of a representation to a quandle $X$
is included in an orbit of $X$,
which forms a subquandle of $X$.
Hence, the number of representations to $X$
is equal to
the sum of the numbers of representations to the quandles
which are obtained as orbits of $X$.
Repeating this procedure,
the number of representations to $X$
can be presented by the sum of the numbers of representations
to certain connected quandles.
Hence, it is sufficient to consider this problem
when a target quandle is connected.\footnote{\footnotesize 
This argument is not available for the link case,
since a link quandle is not connected.}
\end{rem}

\begin{rem}
The problem to count
the number of representations of a knot group to a fixed finite group
can be reduced to Problem \ref{prob.rep_QK}.
Because it is equal to the number of representations 
of a knot quandle to the conjugation quandle of the group,
and the problem to count it
can be reduced to Problem \ref{prob.rep_QK} by the above remark.
\end{rem}

\begin{rem}
The number of representations of a knot quandle to an Alexander quandle
can be presented by using the $i$th Alexander polynomials of the knot
\cite{Inoue}.
In particular, 
the number of representations to a dihedral quandle
can be obtained as its corollary.
\end{rem}

\begin{rem}
Let $X$ be a connected finite quandle,
and let $h_X(K)$ denote the number of representations
of the knot quandle of a knot $K$ to $X$.
Then, $h_X$ is multiplicative with respect to connected sum of knots.
It is known (see, for example, \cite{Ohtsuki_book}) that
any $\Q$-valued Vassiliev invariant is equal to a polynomial
in some primitive Vassiliev invariants,
where primitive Vassiliev invariants are additive
with respect to connected sum of knots.
Hence, $h_X$ is not a Vassiliev invariant,
unless it is constant.
(See also \cite{Altsch96} for another proof.)
\end{rem}

\begin{conj}
Let $h_X$ be as above.
Then, $\log h_X$ is not a 
Vassiliev invariant,\index{Vassiliev invariant!detectability by ---}
unless it is constant.
\end{conj}





\subsection{(Co)homology of quandles}

Second cohomology classes of a quandle
are used in order to define quandle cocycle invariants of knots.
They are introduced as follows.
Let $A$ be an abelian group, written additively,
and let $C^n(X;A)$ be the abelian group
consisting of maps $X^n \to A$,
where $X^n$ denotes
the direct product of $n$ copies of $X$.
We put
\begin{align*}
C^1_Q(X;A) &= C^1(X;A), \\*
C^2_Q(X;A) &= \{ f \in C^2(X;A) \ | \ f(x,x) = 0 \mbox{ for any } x \in X
\}, \\*
C^3_Q(X;A) &= \{ g \in C^3(X;A) \ | \ 
g(x,x,y) = 0 \mbox{ and } g(x,y,y) = 0 \mbox{ for any } x,y \in X \}.
\end{align*}
The coboundary operators
$d_i : C^i_Q(X;A) \to C^{i+1}_Q(X;A)$ are given by
\begin{align*}
& d_1 f(x,y) = f(x) - f(x \ast y), \\*
& d_2 g(x,y,z)
= g(x,z) - g(x,y) - g(x \ast y, z)
+ g(x \ast z, y \ast z),
\end{align*}
for $f \in C^1_Q(X;A)$ and $g \in C^2_Q(X;A)$.
We define the second quandle cohomology group by
$H_Q^2(X;A) = ( \mbox{kernel}\, d_2 ) / ( \mbox{image}\, d_1 )$.
It is known that
$H^2_Q(X;A)$ is isomorphic to $\mbox{Hom}\big( H_2^Q(X) ;A \big)$
by the universal coefficient theorem,
noting that $H_1^Q(X)$ is free abelian (see \cite{CJKS_B}).
Here,
$H_2^Q(X)$ denotes the second homology group
of the dual complex of $\{ C_Q^\star(X;\Z), d_\star \}$.
See \cite{CJKS_B} for the definition of the $n$th quandle (co)homology group.
Therefore,
to obtain $H^2_Q(X;A)$ for any $A$,
it is sufficient to compute $H_2^Q(X)$.

\begin{table}[ht!]
\begin{center}
\begin{tabular}{|l|c|l|l|} \hline 
\small Connected quandle $X$ &\small Order &\small  $H_2^Q(X)$ &\small  $H_3^Q(X)$ \\ \hline\hline
\small $R_3$ &\small  $3$ &\small  $0$ &\small  $\Z/3\Z$ \\ \hline
\small $\Z[t^{\pm1}]/(2,t^2+t+1)$ &\small  $4$ &\small  $\Z/2\Z$ &\small  $\Z/2\Z \oplus \Z/4\Z$ \\ \hline
\small $R_5$ & \raisebox{-6pt}[0pt][0pt]{\small $5$}
&\small  $0$ &\small  $\Z/5\Z$ \\ \cline{1-1}\cline{3-4}
\small $\Z[t^{\pm1}]/(5,t-2)$ & &\small  $0$ &\small  $0$ \\ \hline
\small $R_7$ & &\small  $0$ &\small  $\Z/7\Z$ \\ \cline{1-1}\cline{3-4}
\small $\Z[t^{\pm1}]/(7,t-2)$ &\small  $7$ &\small  $0$ &\small  $0$ \\ \cline{1-1}\cline{3-4}
\small $\Z[t^{\pm1}]/(7,t-3)$ & &\small  $0$ &\small  $0$ \\ \hline
\small $\Z[t^{\pm1}]/(2,t^3+t+1)$ &\small  $8$ &\small  $0$ &\small  $\Z/2\Z$ \\ \hline
\small $R_9$ &  &\small  $0$ &\small  $\Z/9\Z$ \\ \cline{1-1}\cline{3-4}
\small $\Z[t^{\pm1}]/(9,t-2)$ &  &\small  $0$ &\small  $\Z/3\Z$ \\ \cline{1-1}\cline{3-4}
\small $\Z[t^{\pm1}]/(3,t^2+1)$ &\small  $9$ &\small  $\Z/3\Z$ &\small  $(\Z/3\Z)^3$ \\ 
\cline{1-1}\cline{3-4}
\small $\Z[t^{\pm1}]/(3,t^2-t+1)$ &  &\small  $\Z/3\Z$ &\small  $\Z/3\Z \oplus \Z/9\Z$
\\ \cline{1-1}\cline{3-4}
\small $\Z[t^{\pm1}]/(3,t^2+t-1)$ &  &\small  $0$ &\small  $0$ \\ \hline\hline
\small $\Z[t^{\pm1}]/(p,t-a)$ & 
\raisebox{-6pt}[0pt][0pt]{\small $p$} & \raisebox{-6pt}[0pt][0pt]{\small $0$} & \\
\small for any prime $p$ and any $a \ne 0,1 \in \Z/p\Z$ & & & \\ \hline
\end{tabular}
\end{center}
\caption{\label{tbl.quandle-homology}
The cohomologies of the quandles, except for the last one,
in the table are due to \cite{LiNe}.
From a table in \cite{LiNe}
we omit one of two dual quandles
and quandles that are not connected
(see remarks on Problem \ref{prob.compute_PaK}).
The 2nd homology of $\Z[t^{\pm1}]/(p,t-a)$ is due to \cite{Mochizuki}.
See \cite{LiNe,Mochizuki}
for computations of cohomology groups of some more quandles.}
\end{table}

\begin{prob}
\label{prob.computeH2QX}
Compute\index{quandle!--- (co)homology!computation of ---}
$H_2^Q(X)$ for each connected quandle $X$.
More generally, 
find a convenient methodology to compute quandle (co)homology groups.
\end{prob}

See Table \ref{tbl.quandle-homology}
for some quandle homology groups given in \cite{LiNe};
see also \cite{Mochizuki}
for computations of quandle cohomology groups of many Alexander quandles.
There are maple programs \cite{JS_maple}
for computing quandle cohomology groups.

\begin{rem}
We consider only connected quandles in this problem,
since computations of quandle cocycle invariants of knots can be reduced to
the cases of connected quandles
(see a remark on Problem \ref{prob.compute_PaK}).
\end{rem}

\begin{prob}[J.S. Carter] 
\label{prob.HiQSnm}
Compute\index{quandle!--- (co)homology!computation of ---} 
$H_i^Q({\frak S}_n^m)$ of ${\frak S}_n^m$
which denotes the quandle of the $n$th symmetric group 
with the binary operation given by $x \ast y = y^{-m} x y^m$.
\end{prob}

\subsection{Quandle cocycle invariant}

The quandle cocycle invariant,
introduced in \cite{CJKLS,CJKLS_Q},
is defined as follows.
For $\alpha \in H^2(X;A)$
we choose a 2-cocycle $\phi$
representing $\alpha$.
Any representation of a knot quandle $Q(K)$ to $X$
is presented by a coloring of a knot diagram of $K$,
where a {\it coloring} of an oriented knot diagram
is a map of the set of over-arcs of it to $X$
satisfying the condition 
depicted in the pictures of (\ref{eq.qw})
at each crossing of the knot diagram.
We define the weight of a crossing of a colored diagram by
\vspace{0.2pc}
\begin{equation}
\label{eq.qw}
W\Big( \begin{picture}(65,20)
\put(5,0){\pict{d-pt2}{1.2cm}}
\put(5,15){\small $x$}
\put(45,-10){\small $x\! \ast \! y$}
\put(45,15){\small $y$}
\end{picture}
\Big) = \phi(x,y) \in A, \qquad
W\Big( \begin{picture}(65,20)
\put(5,0){\pict{d-pt3}{1.2cm}}
\put(5,15){\small $y$}
\put(5,-10){\small $x$}
\put(45,15){\small $x\! \ast \! y$}
\end{picture}
\Big) = \phi(x,y)^{-1} \in A,
\end{equation}
\vspace{0pc}

\noindent
where we write $A$ multiplicatively here.
The {\it quandle cocycle invariant}\index{quandle!--- cocycle invariant} 
of a knot $K$ is defined by
$$
\Phi_{\alpha}(K) = \sum_{\cal C} \prod_{\tau} W(\tau,{\cal C}) \in \Z[A],
$$
where the sum runs over all coloring $\cal C$ of a diagram of $K$,
and the product runs over all crossing $\tau$ of the diagram, 
and $\Z[A]$ denotes the group ring of $A$.
$\Phi_\alpha(K)$ only depends on $K$ and $\alpha$.

\begin{prob}
\label{prob.compute_PaK}
Compute\index{quandle!--- cocycle invariant!computation of ---}  
the quandle cocycle invariant $\Phi_\alpha(K)$ of each knot $K$
for a second cohomology class $\alpha$ of a connected quandle.
\end{prob}

\begin{rem}
When $X = R_4$ (which is not connected),
it is shown as follows
(see also \cite{CJKS_c} for numerical computation)
that $\Phi_\alpha(K) = 4$ for any $K$ and $\alpha$,
though $R_4$ has non-trivial cohomology groups
since $H_2^Q(R_4) = \Z^2 \oplus (\Z/2\Z)^2$.
The quandle $R_4$ has two orbits,
which form subquandles isomorphic to $T_2$,
where $T_n$ denotes the trivial quandle
(i.e.\ $x \ast y = x$ for any $x,y$) of order $n$.
Further, $T_2$ has two orbits,
which form subquandles isomorphic to $T_1$.
Since $Q(K)$ is connected,
any representation of $Q(K)$ to $R_4$ is trivial
(i.e.\ a constant map).
Hence, any coloring is trivial
(i.e.\ colored by a single element of $X$).
Since $\phi(x,x)=0$ for any 2-cocycle $\phi$,
$\Phi_\alpha(K)=4$ by definition.

When $X = \Z[t^{\pm1}]/(9,t-4)$ (which is not connected),
it follows from a similar argument
(see also \cite{CJKS_c} for numerical computation)
that $\Phi_\alpha(K) = 9$
for any $K$ and $\alpha$,
noting that
this $X$ has three orbits,
which form subquandles isomorphic to $T_3$.

In general,
let $X_1, X_2, \cdots$ be the orbits of $X$.
These form subquandles of $X$.
We denote by $i_k : X_k \to X$ the inclusions.
Then, it follows from a similar argument as above that
$\Phi_\alpha(K) = \sum_k \Phi_{i_k^\star \alpha}(K)$.
Repeating this procedure,
the computations of $\Phi_\alpha(K)$ of a knot $K$
can be reduced to those for connected quandles.\footnote{\footnotesize 
This argument is not available for the link case,
since a link quandle is not connected.}
\end{rem}

\begin{rem}
The cohomology group $H^2_Q(\overline{X};A)$
of the dual quandle $\overline{X}$ of a quandle $X$
is isomorphic to $H^2_Q(X;A)$
by an isomorphism taking a 2-cocycle $\overline{\phi}$ to $\phi$
where $\phi(x,y) = \overline{\phi}(x \ast y , y)$.
It follows that $\Phi_{\overline{\alpha}}(K) = \Phi_{\alpha}(\overline{K})$,
where $\overline{K}$ denotes the mirror image of $K$.
Therefore, the computations of quandle cocycle invariants for $\overline{X}$
can be reduced to those for $X$.
\end{rem}

\begin{rem}
When $\alpha=0$,
by definition
$\Phi_\alpha(K)$ is equal to the number of representations $Q(K) \to X$.
In particular, when $X$ is an Alexander quandle,
it can be presented by using the $i$th Alexander polynomials,
as mentioned in a remark of Problem \ref{prob.rep_QK}.
\end{rem}

\begin{rem}[\rm \cite{CJKS_c}] 
When $X = \Z[t^{\pm1}]/(2,t^2+t+1)$,
$H^2_Q(X;\Z/2\Z) = \Z/2\Z$.
For its non-trivial cohomology class $\alpha$,
$$
\Phi_\alpha(K) = \begin{cases}
\mboxsm{$4(1+3u)$} & 
\mboxsm{ for $K = 3_1, 4_1, 7_2, 7_3, 8_1, 8_4, 8_{11}, 8_{13}$, and} \\
& \qquad \mboxsm{$9$ certain knots with 9 crossings}, \\
\mboxsm{$16(1+3u)$} & 
\mboxsm{ for $K = 8_{18}, 9_{40}$,}  \\
\mboxsm{$16$} & 
\mboxsm{ for $K = 8_5, 8_{10}, 8_{15}$, $8_{19}$--$8_{21}$, and} \\
& \qquad \mboxsm{$16$ certain knots with 9 crossings}, \\
\mboxsm{$4$} & 
\mboxsm{ for the other knots $K$ with at most $9$ crossings},
\end{cases}
$$
where $u$ denotes the generator of $\Z/2\Z$.
See \cite{CJKS_c} for details.

When $X = \Z[t^{\pm1}]/(3,t^2+1)$,
$H^2_Q(X;\Z/3\Z) = \Z/3\Z$.
For a non-trivial cohomology class $\alpha$ of it,
$$
\Phi_\alpha(K) = \begin{cases}
\mboxsm{$9(1+4u+4u^2)$} & 
\mboxsm{ if $K = 4_1, 5_2, 8_3, 8_{17}, 8_{18}, 8_{21}$, and} \\
&\qquad \mboxsm{ $9$ certain knots with 9 crossings}, \\
\mboxsm{$27 (11 + 8 u + 8 u^2)$} & \mboxsm{ if $K = 9_{40}$,} \\
\mboxsm{$81$} & 
\mboxsm{ if $K = 6_3, 8_2, 8_{19}, 8_{24}, 9_{12}, 9_{13}, 9_{46}$,} \\
\mboxsm{$9$} & 
\mboxsm{ for the other knots $K$ with at most $9$ crossings},
\end{cases}
$$
where $u$ denotes a generator of $\Z/3\Z$.
See \cite{CJKS_c} for details.
\end{rem}

\begin{rem}
It is known, see \cite{CENS}, that 
each $\alpha \in H^2_Q(X;A)$ gives
an {\it abelian extension} $A \to Y \overset{p}{\to} X$, 
where $Y=A \times X$,
which forms a quandle with the binary operation given by
$(a_1,x_1) \ast (a_2,x_2) = \big( a_1 + \phi(x_1,x_2), x_1 \ast x_2 \big)$
using a 2-cocycle $\phi$ representing $\alpha$,
and $p$ denotes the natural projection.

Let $a_1, a_2, \cdots, a_N$ be a sequence of generators of $Q(K)$
associated with over paths of a diagram of $K$
which are chosen as going around $K$.
Adjacent generators $a_1$ and $a_2$ are related by
$a_1 \ast b = a_2$ (or $a_1 = a_2 \ast b$) for some generator $b$.
Let $\widetilde{\rho(b)}$ be a pre-image of $\rho(b)$
under the projection $p$.
Then, $S_{\widetilde{\rho(b)}}$ (resp. $S^{-1}_{\widetilde{\rho(b)}}$)
induces a map $p^{-1}(a_1) \to p^{-1}(a_2)$,
which does not depend on the choice of a pre-image of $\rho(b)$.
Composing such maps, we have a sequence of maps
$p^{-1}(a_1) \to p^{-1}(a_2) \to \cdots \to p^{-1}(a_N) \to p^{-1}(a_1)$.
The composite map of these maps can be expressed
$a \mapsto a + m(\rho)$ ($a \in A$) for some $m(\rho) \in A$.
Then, the quandle cocycle invariant can be presented by
$\Phi_\alpha(K) = \sum_\rho m(\rho) \in \Z[A]$, 
where the sum runs over all representations $\rho$ of $Q(K)$ to $X$.

In particular, as shown in \cite{CENS}, 
the number of representations $Q(K) \to X$
that can lift to representations $Q(K) \to Y$ is equal to
the coefficient of the unit of $A$ in $\Phi_\alpha(K)$.
For example, when $A = \Z/2\Z$,
it follows that, writing $\Phi_\alpha(K) = a + b t$
(where $t$ is the generator of $\Z/2\Z$),
$a$ is equal to the number of representations $Q(K) \to X$
that can lift to representations $Q(K) \to Y$,
and $b$ is equal to the number of those that do not.

In this way we can compute $\Phi_\alpha(K)$ 
in terms of the abelian extension associated to $\alpha$.
\end{rem}

\begin{prob}
Find relations between
quandle cocycle invariants\index{quandle!--- cocycle invariant!--- and other invariant} 
and knot invariants known so far, 
such as quantum invariants.\index{quantum invariant!--- and quandle cocycle invariant} 
\end{prob}

\begin{rem}
When $\alpha=0$ and $X$ is an Alexander quandle,
$\Phi_\alpha(K)$ can be presented
by using the $i$th Alexander polynomial,
as mentioned in a remark of Problem \ref{prob.compute_PaK}.
\end{rem}

\begin{rem}[\rm (M. Gra\~na)] 
The quandle cocycle invariants can be presented by knot invariants
derived from certain ribbon categories \cite{Grana_qq}.
\end{rem}

A central extension of a group $G$
gives an abelian extension of the conjugation quandle of $G$.
It is known that an abelian extension of a group $G$ 
can be characterized by the cohomology group $H^2(G;A)$ for a $G$-module $A$.
Motivated by this cohomology group
we introduce $H^2_Q(X;A)$ of a quandle $X$
for an ``$X$-module'' $A$ as follows.
We call an abelian group $A$ 
an $X$-{\it module} of a quandle $X$
if there is a map $\rho: X \to \mbox{Aut}(A)$
satisfying that $\rho(x \ast y) = \rho(y)^{-1} \rho(x) \rho(y)$
for any $x,y \in X$.
For simplicity, 
we often write $\rho(x)^{\pm1}a$ as $x^{\pm1}a$ omitting $\rho$.
Let $C^i_Q(X;A)$ be as before.
We give the coboundary operators by
\begin{align*}
& d_1 f( x, y )
= y^{-1} \big( f(x) + x f(y) - f(y) \big) - f(x \ast y), \\*
& d_2 g( x, y, z )
= (y \ast z)^{-1} g(x,z) - z^{-1} g(x,y) 
+ (y \ast z)^{-1} \big( (x \ast z) -1 \big) g(y,z) \\*
& \qquad\qquad\qquad \qquad\qquad
- g(x\ast y, z) + g(x \ast z, y \ast z), 
\end{align*}
for $f \in C^1_Q(X;A)$ and $g \in C^2_Q(X;A)$.
We define the second quandle cohomology group by
$H_Q^2(X;A) = ( \mbox{kernel}\, d_2 ) / ( \mbox{image}\, d_1 )$.

\begin{prob}
Compute\index{quandle!--- (co)homology!computation of ---}   
$H^2_Q(X;A)$ for each $X$-module $A$.
\end{prob}

\begin{rem}
This cohomology group might be isomorphic
to the cohomology group of a quandle space of $X$
(see a remark on Problem \ref{prob.computeH2QX})
with coefficients in the local system
corresponding to the $X$-module $A$.        
\end{rem}

For an $X$-module $A$, 
each $\alpha \in H^2_Q(X;A)$ gives
an extension $A \to Y \to X$, where $Y=A \times X$,
which forms a quandle
with the binary operation given by
$(a_1,x_1) \ast (a_2,x_2) 
= \big( x_2^{-1}(a_1+x_1 a_2 - a_2) + \phi(x_1,x_2), x_1 \ast x_2 \big)$
using a 2-cocycle $\phi$ representing $\alpha$.

\begin{prob}
Let the notation be as above.
Then, extending the definition of the quandle cocycle invariant,
define a knot invariant associated with $\alpha$,
which is, roughly speaking,
an invariant obtained by counting
representations of a knot quandle $Q(K)$ to $X$
with information whether
each representation can lift to a representation $Q(K) \to Y$.
\end{prob}

\subsection{Quantum quandles}

\begin{prob}[M. Polyak] 
\label{prob.qquandle}
Define a quantum quandle.\index{quantum quandle}\index{quandle!quantum ---}
\end{prob}

\begin{rem}
A quantum group is a quantization of a group,
in the sense that it can be regarded as 
a non-commutative perturbation of a (certain) function algebra on a group.
It is a problem to formulate an appropriate quantization of a quandle.
\end{rem}

\subsection{Rack (co)homology}\index{rack}

(C. Rourke, B. Sanderson)\quad A {\it rack} has the same definition as
a quandle, except that axiom (1) is omitted.  Quandles are thus a
special class of racks.  There is a naturally defined classifying
space for a rack (in fact a semi-cubical complex), the {\it rack
space},\index{rack!--- space} and hence both homology and cohomology
groups are defined.  These can be used to define invariants of
(framed) links \cite{FRS_spec, FRS_james, FRS_rack}.  Regarding a
quandle as a rack we have two definitions of cohomology which are
closely related.  Indeed the cochain group for quandle cohomology is a
subgroup of that for rack cohomology with the definition of coboundary
unchanged.\index{rack!--- (co)homology} 

\begin{prob}[C. Rourke, B. Sanderson]\index{quandle!--- space}  
\label{FRSprob1}Is there a natural quandle space whose cohomology groups are the
quandle cohomology groups?
\end{prob}

\begin{rem}The quandle cochain groups do not correspond to setting the
rack cochains to be zero on a subcomplex of the rack space.  Thus the
first guess that the quandle space is obtained by quotienting a
certain subcomplex is false.\end{rem}

There is an anologous problem of classification and computation of
(co)homology groups for racks as for quandles.  One particular
interesting question is the following:

\begin{conj}[R. Fenn, C. Rourke, B. Sanderson] \label{FRSprob2}
$H_3(R_p) \cong \Z\oplus \Z/p\Z$ for $p$ prime.\end{conj}

\begin{rem}The conjecture is equivalent to $H_3^Q(R_p) \cong \Z/p\Z$ for $p$
prime.  This has been verified for $p\le 7$ (see table 5) and the rack
version has been verified (again for $p \le 7$) by Maple calculation
(Rourke and Sanderson).  A direct proof (without using a computer
calculation) has been found for $p=3$.  The conjecture is consistent
(indeed suggested by) the calculation of T. Mochizuki
\cite{Mochizuki} that $H_3^Q(R_p; \Z/p\Z)\cong \Z/p\Z$.\end{rem}

\newpage

\section{Braid group representations}
\label{sec.bg_rep}

\renewcommand{\thefootnote}{\fnsymbol{footnote}}
\footnotetext[0]{Chapter \ref{sec.bg_rep} was written by S.J. Bigelow.}
\renewcommand{\thefootnote}{\arabic{footnote}}

For $n=1,2,\dots$,
the braid group $B_n$\index{braid group}
is the group generated by $\sigma_1,\dots,\sigma_{n-1}$
modulo the relations:
\begin{itemize}
\item
$\sigma_i \sigma_j = \sigma_j \sigma_i$ \quad if $|i-j| > 1$,
\item
$\sigma_i \sigma_j \sigma_i = \sigma_j \sigma_i \sigma_j$ \quad 
if $|i-j| = 1$.
\end{itemize}

\subsection{The Temperley-Lieb algebra}

For $\tau \in \C$,
the {\it Temperley-Lieb algebra}\index{Temperley-Lieb algebra} 
${\mathrm{TL}}_n(\tau)$ is 
the associative $\C$-algebra
generated by $1,e_1,\dots,e_{n-1}$ modulo the relations:
\begin{itemize}
\item
$e_i e_j = e_j e_i$ \quad if $|i-j| > 1$,
\item
$e_i e_j e_i = e_i$ \quad if $|i-j| = 1$,
\item
$e_i e_i = \tau e_i$.
\end{itemize}
We will simply write ${\mathrm{TL}}_n$, where $\tau$ is understood.
There is a map from $B_n$ to ${\mathrm{TL}}_n$ given by
\begin{align*}
\sigma_i &\longmapsto A + A^{-1} e_i, \\
\sigma_i^{-1} &\longmapsto A^{-1} + A e_i,
\end{align*}
where $A \in \C$ is such that $\tau = -A^2-A^{-2}$.

\begin{figure}[ht!]
$$
\pict{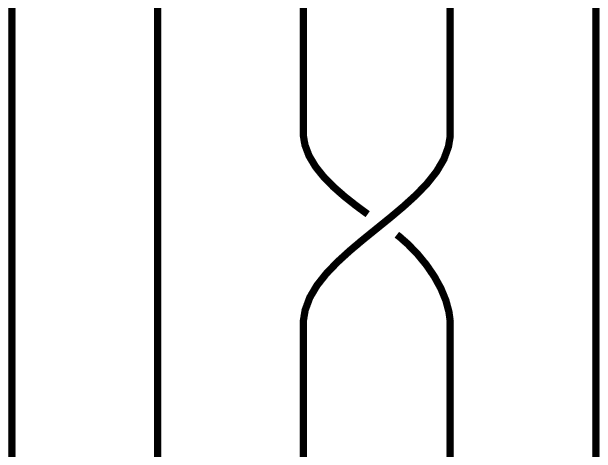}{2cm} \qquad\qquad
\pict{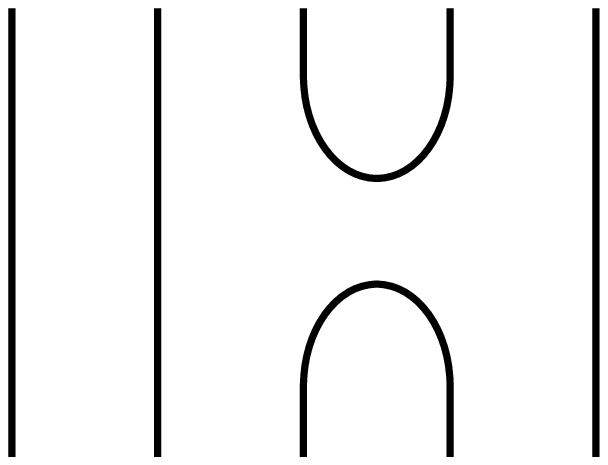}{2cm}
$$
\caption{\label{fig.braid_gen}
The generator $\sigma_3$ of $B_5$ (the left picture) and
the generator $e_3$ of ${\mathrm{TL}}_5$ (the right picture).}
\end{figure}

These definitions can be motivated in terms of
{\em tangle diagrams} in $\R \times I$.
These are similar to knot diagrams,
except that they can include arcs with endpoints on $\R \times \{0,1\}$.
Two tangles are considered the same
if they are related by a sequence of isotopies
and Reidemeister moves of the second and third type.
The generators of $B_n$ and ${\mathrm{TL}}_n$
can be defined to be the tangle diagrams
suggested by Figure \ref{fig.braid_gen}.
The arcs of these diagrams have endpoints
$$\{1,2,\dots,n\} \times \{0,1\}.$$
The product $ab$ of two such diagrams $a$ and $b$
is obtained by placing $a$ on top of $b$
and then shrinking the result vertically to the required height.
The third relation in the Temperley-Lieb algebra
allows one to delete a closed loop
at the expense of multiplying by $\tau$.
Using these definitions,
the map from $B_n$ to ${\mathrm{TL}}_n$ is given by 
resolving all crossings using the Kauffman skein relation.

\begin{prob}[{\cite[Problem 3]{Jones_prob}}] 
\label{prob.TLfaithful}
Is the representation of the braid group\index{braid group!representation of ---!faithful ---}
inside the Temperley-Lieb algebra\index{Temperley-Lieb algebra}  
faithful?
\end{prob}

\begin{rem}
We are mostly interested in the case $\tau$ is a transcendental.
The answer is yes for $n \le 3$,
and unknown for all other values of $n$.

The Jones polynomial\index{Jones polynomial!trivial ---}
 of the closure of a braid $\beta$
is a certain trace function of
the image of $\beta$ inside the Temperley-Lieb algebra.
If $\beta \in B_n \setminus \{1\}$ maps to the identity in ${\mathrm{TL}}_n$,
and $\gamma \in B_n$ is any braid whose closure is the unknot,
then the closure of $\beta \gamma$ would have Jones polynomial one.
It should be easy to arrange for this to be a non-trivial knot. 
Thus a negative answer to Problem \ref{prob.TLfaithful}
would almost certainly lead to a solution to Problem \ref{prob.VK=1}.
\end{rem}

\subsection{The Burau representation}

For $k=0,1,\dots,\lfloor \frac{n}{2} \rfloor$,
let $V^n_{n-2k}$ be the vector space 
spanned by tangle diagrams in $\R \times I$
with no crossings and endpoints
$$\{(1,0),(2,0),\dots,(n-2k,0)\} \cup \{(1,1),(2,1),\dots,(n,1)\}$$
modulo the relations:
\begin{itemize}
\item
a tangle is zero 
if it contains an edge with both endpoints on $\R \times \{0\}$,
\item
a closed loop may be removed at the expense of multiplying by $\tau$.
\end{itemize}
Let ${\mathrm{TL}}_n$ act on $V^n_{n-2k}$
by stacking tangle diagrams in the usual way.
For generic values of $\tau$,
${\mathrm{TL}}_n$ is semisimple
and these are its irreducible representations.

We obtain irreducible representations of $B_n$
by taking its induced action on $V^n_{n-2k}$.
By a result of Long \cite{br.long},
the representation of $B_n$ inside ${\mathrm{TL}}_n$ is faithful
if and only if each of these irreducible representations is faithful.
Note that the action of $B_n$ on the one-dimensional space $V^n_n$
is never faithful for $n>2$.
Also if $n>2$ is even
then the action of $B_n$ on $V^n_0$ is easily shown to be unfaithful.
The action of $B_n$ on $V^n_{n-2}$ is 
the famous {\it Burau representation}.\index{Burau representation}

\begin{prob}
Is\index{braid group!representation of ---!faithful ---} 
the Burau representation\index{Burau representation} 
of $B_4$ faithful?
\end{prob}

\begin{rem}
The Burau representation of $B_n$ is known to be 
faithful for $n \le 3$ and not faithful for $n \ge 5$
\cite{br.bigelow:b5}. 

The representation of $B_4$ in ${\mathrm{TL}}_4$ is faithful
if and only if the Burau representation of $B_4$ is faithful.
\end{rem}

\begin{rem}
The Burau representation of $B_4$ is faithful if and only if 
a certain pair of three-by-three matrices generate a free group.
The matrices given in \cite{br.birman} contain a misprint,
but their description as words in the generators is correct.

The Burau representation of $B_4$ is faithful if and only if
a certain intersection pairing detects intersection of arcs
in the four-times punctured disk \cite{br.bigelow:b5}.

Cooper and Long have explicitly calculated
the kernel of the Burau representation modulo the primes 2, 3 and 5
\cite{br.cooper-long}.
\end{rem}

\begin{prob}[S.J. Bigelow] 
\label{prob.B6V62}
Is\index{braid group!representation of ---!faithful ---} 
the action of $B_6$ on $V^6_2$ faithful?
\end{prob}

\begin{rem}
The Burau representation of $B_6$ is unfaithful \cite{br.long-paton}.
Thus the representation of $B_6$ in ${\mathrm{TL}}_6$ is faithful
if and only if the action of $B_6$ on $V^6_2$ is faithful.

No approach to this problem is known
except for a brute force computer search.
However such a search might find an example more easily
than any of the more subtle approaches
to the Burau representation of $B_4$.
\end{rem}

\begin{rem}
We could also ask whether the action of $B_5$ on $V^5_1$ is faithful.
A computer search of this representation
would be easier because the matrices involved are smaller
(five-by-five instead of nine-by-nine).
On the other hand, this representation is more likely to be faithful,
since if the representation of $B_6$ in ${\mathrm{TL}}_6$ is faithful
then so is the representation of $B_5$ in ${\mathrm{TL}}_5$.
\end{rem}

\subsection{The Hecke and BMW algebras}

We now introduce two algebras
which can be defined in a similar way to the Temperley-Lieb algebra.
The {\it Hecke algebra}\index{Hecke algebra} 
is
the set of formal linear combinations of braids
modulo the relation:
$$
A \pict{bigelow/bmw1}{1.5cm} - A^{-1} \pict{bigelow/bmw2}{1.5cm}
= (A^2-A^{-2}) \pict{bigelow/bmw3}{1.5cm},
$$
where $A \in \C$.
The {\it BMW algebra}\index{BMW algebra} 
is
the set of formal linear combinations of tangles
whose edges have endpoints $\{1,2,\dots,n\} \times \{0,1\}$,
modulo the relations:
\begin{align*}
\pict{bigelow/bmw1}{1.5cm} + \pict{bigelow/bmw2}{1.5cm}
&= m \left( \pict{bigelow/bmw3}{1.5cm} + \pict{bigelow/bmw4}{1.5cm} \right) \\*
\pict{bigelow/RI1}{1.5cm} & = l \pict{bigelow/RI2a}{1.5cm},
\end{align*}
where $m,l \in \C$.
See \cite{br.murakami} and \cite{br.birman-wenzl}
for an analysis of the BMW algebra.
The Temperley-Lieb algebra can be embedded into the Hecke algebra,
which in turn can be embedded into the BMW algebra.

These algebras are semisimple for generic values of their parameters.
The irreducible representations of the BMW algebra
correspond to partitions of $n-2k$
for $k=0,1,\dots,\lfloor \frac{n}{2} \rfloor$.
The irreducible representations of the Hecke algebra
correspond to partitions of $n$.
The irreducible representations of the Temperley-Lieb algebra
correspond to partitions of $n$ into two parts.

Lawrence \cite{br.lawrence} has used a topological construction
to obtain the irreducible representations of the Hecke algebra.
The construction uses the definition of the braid group
as the mapping class group of a punctured disk
to obtain an action on the homology of a related space.

\begin{prob}[S.J. Bigelow] 
\label{prob.BMWalg}
Generalise Lawrence's construction to obtain
the irreducible representations of the BMW algebra.\index{BMW algebra} 
\end{prob}

\begin{rem}
Zinno \cite{br.zinno} has shown how to obtain
the representation of the BMW algebra 
corresponding to the partition of $n-2$ into one part.
\end{rem}

\begin{prob}[S.J. Bigelow] 
\label{prob.repBn_BMW}
Find\index{braid group!representation of ---} 
a larger family of irreducible representations of $B_n$
which includes those coming from the BMW algebra.\index{BMW algebra} 
\end{prob}

\begin{rem}
This might be defined using tangles and some more complicated relations,
or by generalising Lawrence's approach.
\end{rem}

\subsection{Other problems}

\begin{prob}
Classify\index{braid group!representation of ---!irreducible ---} 
all irreducible representations of $B_n$.
\end{prob}

\begin{rem}
This is probably impossibly hard.
However it seems that interesting partial results are possible.
Formanek \cite{br.formanek} has classified 
all irreducible complex representations of $B_n$
having degree at most $n-1$.
\end{rem}

\begin{prob}[S.J. Bigelow] 
\label{prob.repBn_barQ}
Is\index{braid group!representation of ---!faithful ---} 
there a faithful representation of $B_n$
into a group of matrices over $\bar\Q$?
\end{prob}

\begin{rem}
There is a faithful representation of $B_3$ into ${\mathrm{GL}}(2,\Z)$.
The problem is open for all $n \ge 4$.

There is a faithful representation of $B_n$
into a group of matrices over $\Z[q^{\pm 1},t^{\pm 1}]$.
Krammer's proof of this fact \cite{br.krammer}
works when $t$ is assigned any value between $0$ and $1$.
However it is not known whether there is an algebraic value of $q$
for which the representation remains faithful.
\end{rem}

\newpage

\section{Quantum and perturbative invariants of 3-mani\-folds}

\subsection{Witten-Reshetikhin-Turaev invariants and quantum invariants}

Witten \cite{Witten} proposed that,
for a semi-simple compact Lie group $G$ and a positive integer $k$,
a topological invariant of a closed oriented 3-manifold $M$
is given by the path integral\index{Chern-Simons!--- path integral}
\begin{equation}
\label{eq.CSpi}
Z^G_k(M) = \int e^{2\pi\sqrt{-1}k \mboxss{CS$(A)$}} {\cal D}A,
\end{equation}
which is a formal integral
over gauge equivalence classes of connections $A$
on the trivial $G$ bundle on $M$.
Here, the {\it Chern-Simons functional}\index{Chern-Simons!--- functional} 
$\mbox{CS} : \cA \to \R$ is defined by
\begin{equation}
\label{eq.CSfunct}
\mbox{CS}(A) = \frac{1}{8 \pi^2}
\int_M \mbox{ trace} (A \wedge dA + \frac23 A \wedge A \wedge A),
\end{equation}
for a connection $A$,
regarding it as a $\frak g$-valued 1-form on $M$,
where $\frak g$ denotes the Lie algebra of $G$.


Motivated by Witten's proposal,
the {\it quantum $G$ invariant}\index{quantum invariant} 
$\tau^G_r(M)$ has been defined and studied,
first by Reshetikhin and Turaev \cite{RT_inv} and later by other researchers,
where we put $r = k+ h^{{}^\vee}$
with the dual Coxeter number $h^{{}^\vee}$ of $\frak g$.
The quantum invariant
is also called the {\it Witten-Reshetikhin-Turaev 
invariant}.\index{Witten-Reshetikhin-Turaev invariant}
For example, 
when $M$ is obtained from $S^3$ by integral surgery along a framed knot $K$
with a positive framing,
$\tau_r^{SU(2)}(M)$ for $r \ge 3$ and 
$\tau_r^{SO(3)}(M)$ for odd $r \ge 3$ are given by
\begin{align*}
\tau_r^{SU(2)}(M) &= 
\Big( \sum_{n=1}^{r-1} [n] Q^{sl_2;V_n}(U_+) \Big)^{-1}
\sum_{n=1}^{r-1} [n] Q^{sl_2;V_n}(K)
\Big|_{q=\exp(2\pi\sqrt{-1}/r)}, \\
\tau_r^{SO(3)}(M) &= 
\Big( \sum_{\substack{0 < n < r \\ r {\rm \ is \ odd}}}
[n] Q^{sl_2;V_n}(U_+) \Big)^{-1}
\sum_{\substack{0 < n < r \\ r {\rm \ is \ odd}}} [n] Q^{sl_2;V_n}(K)
\Big|_{q=\exp(2\pi\sqrt{-1}/r)},
\end{align*}
where $[n] = (q^{n/2}-q^{-n/2})/(q^{1/2}-q^{-1/2})$, 
and $U_+$ denotes the trivial knot with $+1$ framing,
and $Q^{sl_2;V_n}(K)$ denotes
the quantum invariant
of $K$
associated with the irreducible $n$-dimensional representation
of the quantum group $U_q(sl_2)$;
for details see \cite{KM_inv}
(see also \cite{Ohtsuki_book} for the notation).
It is known \cite{KM_inv} that
$$
\tau_r^{SU(2)}(M) =
\begin{cases} 
\tau_3^{SU(2)}(M) \tau_r^{SO(3)}(M) 
\quad & \mbox{ if } r \equiv 3 \mbox{ mod $4$}, \\
\overline{\tau_3^{SU(2)}(M)} \tau_r^{SO(3)}(M) 
\quad & \mbox{ if } r \equiv 1 \mbox{ mod $4$},
\end{cases}
$$
where $\tau_3^{SU(2)}(M)$ is an invariant
determined by the cohomology ring and the linking pairing of $M$,
which is equal to zero for some $M$ (see (\ref{eq.tauSU2_3})).
For details on quantum $G$ invariants, 
see e.g.\ \cite{Ohtsuki_book} and references therein.

\begin{prob}[see {\cite[Problem 3.108]{Kirby}}] 
\label{prob.quantumS3}
Does there exist a closed 3-manifold $M$, other than $S^3$, such that
$\tau^{SO(3)}_r(M) = \tau^{SO(3)}_r(S^3)$ for all odd $r \ge 3$?
\end{prob}

\begin{rem}[\rm (see {\cite[Remark on Problem 3.108]{Kirby}})] 
Suppose that
$\tau^{SO(3)}_r(M) = \tau^{SO(3)}_r(S^3)$
for a closed 3-manifold $M$ and all odd $r \ge 3$.
If the Betti number of $M$ was positive,
$\tau^{SO(3)}_r(M)$ is divisible by $q-1$.
Hence, $M$ is a rational homology 3-sphere.
We have that
$\tau^{SO(3)}(M) = \tau^{SO(3)}(S^3)$.
Since the leading two coefficients of $\tau^{SO(3)}(M)$
are given by the order of the first homology group 
and Casson invariant of $M$,
$M$ is an integral homology 3-sphere with Casson invariant zero.

Note that
$\tau^{SO(3)}_r \big( L(65,8) \big) = \tau^{SO(3)}_r \big( L(65,18) \big)$
for all odd $r \ge 3$; see \cite{Yamada_lens}.
\end{rem}

\begin{rem}
There is a center in the mapping class group of the closed surface of genus 2,
shown below.
$$
\pict{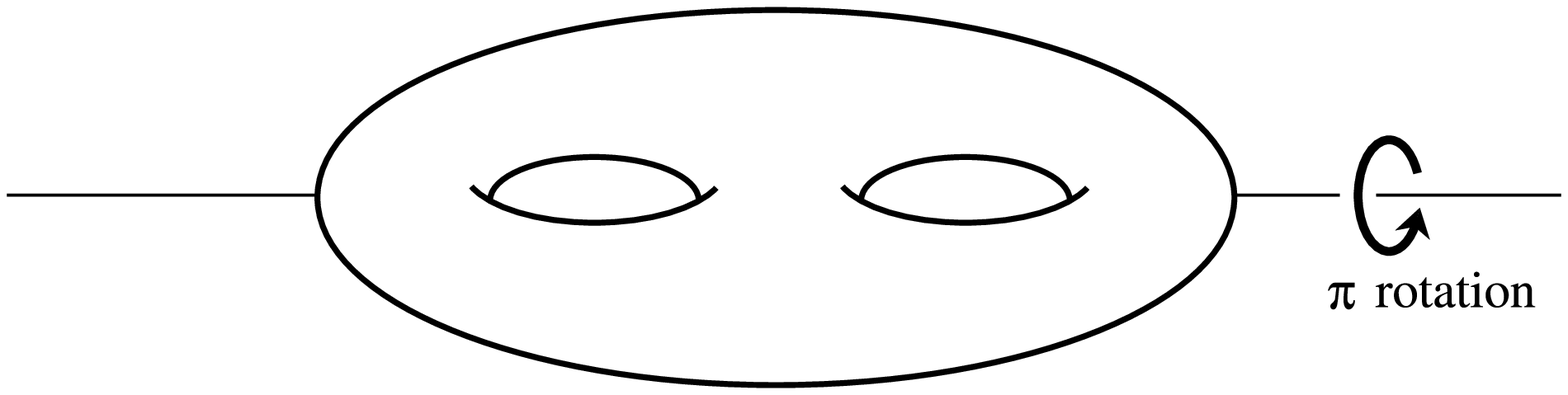}{2cm}
$$
A {\it mutation}\index{mutation} 
of a 3-manifold $M$ 
is defined to be a 3-manifold
obtained from $M$
by cutting along a separating closed surface of genus 2 in $M$
and by gluing again after twisting by the above map.
It is shown in \cite{Kawauchi} that
$\tau_r^{SO(3)}(M)$ does not depend on a change by any mutation of $M$.
\end{rem}

\begin{prob}[S.K. Hansen, T. Takata] 
\label{prob.QHS_tauG_LMO}
Find pairs of non-homeomorphic rational homology 3-spheres 
that can be distinguished
by their quantum $G$ invariants\index{quantum invariant!--- of 3-manifold!strength of ---}  
$\tau_{r}^{G}$ or their quantum $PG$
invariants $\tau_{r}^{PG}$ for some level $r$ and some simply connected
compact simple Lie group $G$ 
but not by their LMO invariants.\index{LMO invariant!strength of ---} 
\end{prob}

\begin{rem}[\rm (S.K. Hansen, T. Takata)] 
For example, the LMO invariants
of the lens spaces\index{lens space} 
$L(25,4)$ and $L(25,9)$
are equal \cite{BaLa}, but
their quantum $SU(2)$ invariants for $r=5$ are not equal.
\end{rem}

\begin{prob}[S.K. Hansen, T. Takata] 
\label{prob.QHS_alltauG}
Do the family of
quantum $G$ invariants\index{quantum invariant!--- of 3-manifold!strength of ---}   
$\tau_{r}^{G}$ or the family of quantum $PG$ invariants
$\tau_{r}^{PG}$, $G$ running through all simply connected compact simple
Lie groups and $r$ running through all allowed levels, 
separate rational homology 3-spheres? 
How well do these families of invariants separate closed
oriented 3-manifolds?
\end{prob}

\begin{rem}[\rm (S.K. Hansen, T. Takata)] 
It is well known that
the LMO invariant is a weak invariant outside the class of
rational homology 3-spheres; 
see the last remark on Problem \ref{prob.calLMO}.
On the contrary there are
3-manifolds with arbitrary high first Betti number and
non-trivial quantum $SU(2)$ invariants as the example of Seifert
manifolds shows. We note that the non-triviality of the invariants
of Seifert manifolds e.g.\ follows from the fact that these
invariants have non-trivial asymptotic expansion in the
limit of large quantum level; 
see \cite{Rozansky_3}, \cite{Hansen}, and Section \ref{sec.aec}.
It is likely to believe, e.g.\ from the asymptotic expansion
conjecture of Andersen, see Conjecture \ref{conj.AEC_Jorgen}, 
that the quantum $G$
invariants are quite strong invariants also outside the class of
rational homology 3-spheres.
It is known, however, that the family of quantum $SU(n)$ invariants,
$n$ running through all integers $>1$, is not a complete invariant,
that is to say that this family of invariants can not separate all
closed oriented 3-manifolds, cf.\ \cite{Lickorish_s}. 
It is still an open
question if this is also the case if we include the quantum invariants
for all the other simply connected compact simple Lie groups.
\end{rem}

\begin{prob}
Find a 3-dimensional topological interpretation of quantum invariants 
of 3-manifolds.\index{quantum invariant!--- of 3-manifold!interpretation of ---} 
\end{prob}

\begin{rem}
Certain special values have some interpretations.
For a closed oriented 3-manifold $M$,
\begin{equation}
\label{eq.tauSU2_3}
\tau^{SU(2)}_3 (M) = 
\begin{cases}
\mboxsm{$0$} \qquad
\mboxsm{if there exists $\alpha \in H^1(M;\Z/2\Z)$ with $\alpha^3 \ne 0$,} \\
\mboxsm{$\sqrt{2}^{{\rm rank}\, H^1(M;\Z/2\Z)}e^{-\beta(M)\pi\sqrt{-1}/4}$} 
\qquad \mboxsm{otherwise,}
\end{cases}
\end{equation}
where $\beta(M)$ denotes the Brown invariant.
Further, for a closed oriented 3-manifold $M$, 
$$
\tau_4^{SU(2)}(M) = \sum_\sigma e^{-\mu(M,\sigma) \cdot 3 \pi \sqrt{-1}/8},
$$
where the sum runs over all spin structures $\sigma$ of $M$
and $\mu(M,\sigma)$ denotes 
the Rokhlin invariant of a spin structure $\sigma$ of $M$.
For details, see \cite{KM_inv}.

It is known \cite{Hitoshi_CW} that,
for any rational homology 3-sphere $M$ and any prime $p > |H_1(M;\Z)|$,
$$
|H_1(M;\Z)| \cdot \tau_p^{SO(3)}(M)
\equiv \left( \frac{|H_1(M;\Z)|}{p} \right)
\Big( 1 + 6 \lambda(M) (\zeta -1) \Big)
$$
mod $(\zeta-1)^2$ in $\Z[\zeta]$, putting $\zeta=e^{2\pi\sqrt{-1}/p}$, 
where $\lambda(M)$ denotes the Casson-Walker invariant of $M$
and $\left( \frac{\cdot}{p} \right)$ denotes the Legendre symbol.
\end{rem}

\begin{rem}
The Chern-Simons path integral
(\ref{eq.CSpi}) by Witten \cite{Witten}
gives a 3-dimensional physical interpretation
of a quantum invariant of 3-manifolds.
Historically speaking, the quantum invariants of 3-manifolds
were introduced, motivated by Witten's Chern-Simons path integral.
\end{rem}

\begin{conj}{\rm\cite{GP98}}\qua
\label{conj.GP98} 
For non-vanishing $\tau^G_r(M)$,
the absolute value\index{quantum invariant!--- of 3-manifold!absolute value of ---} 
$| \tau_r^G(M) |$ only depends on 
the fundamental group $\pi_1(M)$.
\end{conj}


\subsection{The asymptotic expansion conjecture}
\label{sec.aec}

\renewcommand{\thefootnote}{\fnsymbol{footnote}}
\footnotetext[0]{The first version of Section \ref{sec.aec} was written by T. Ohtsuki,
following seminar talks given by J.E. Andersen.
Based on it, J.E. Andersen wrote this section.}
\renewcommand{\thefootnote}{\arabic{footnote}}

The perturbative expansion of the Chern-Simons path integral
(\ref{eq.CSpi}) is given by the semi-classical approximation and
its higher loop perturbations. Roughly speaking, the semi-classical
approximation is obtained from the path integral by ignoring the
contribution from the third order term of the Chern-Simons functional,
and the higher loop perturbation contributions are the
corrections to this semi-classical contribution.

To the best of our knowledge, there is today, no complete
perturbative treatment of the Chern-Simons quantum field theory
available, even from a mathematical physics point of view. In the
following few paragraphs we shall try to outline the main
activities seen so far in this direction.

The the first formula for the semi-classical approximation of the
Chern-Simons path integral
was given by Witten in \cite{Witten},
describing it as a sum of contributions, one for each gauge
equivalence classes of flat connection, involving the Chern-Simons
value, the Reidemeister torsion and a certain
spectral flow for each such gauge equivalence class.
To test this prediction, Freed and Gompf
\cite{FrGo} made for certain Seifert fibered manifolds some
computer studies of the large $k$ behavior of $Z^{SU(2)}_{k}(M)$
and based on these calculations and further discussion of the
semi-classical approximation of the path integral, they proposed
the following formula for the semi-classical approximation
($r=k+2$)
\begin{align*}
Z^{SU(2)}_{k}(M) \underset{r\to\infty}{\sim} 
& e^{- 3 \pi \sqrt{-1} (1+b^1(M))/4 } \\*
& \times \sum_{[A]} e^{2\pi\sqrt{-1}r
\mboxss{CS$(A)$}} r^{(h^1_A - h^0_A)/2} e^{-2\pi\sqrt{-1} (I_A/4 +
h^0_A /8)} \tau_M(A)^{1/2}, 
\end{align*}
where the sum is over the gauge equivalence classes of flat
connections $A$. Let us explain the quantities involved in this
expression and in which cases one can make sense of this
expression as it stands.

For any flat connection $A$, we have the cohomology groups $H^i(M,d_A)$
of the covariant derivative complex
$d_A: \Omega^\star(M;{\frak g}) \to \Omega^{\star +1}(M;{\frak
g})$ given by $d_A f = d f + [A, f]$, and $h^i_A$ is the
dimension of $H^i(M,d_A)$. 
Further associated to this complex we have the Reidemeister torsion 
$\tau_M(A) \in \otimes_i \big(\det H^i(M,d_A)\big)^{(-1)^{i}} \cong 
\Big( \det H^0(M,d_A) \otimes \big(\det H^1(M,d_A)\big)^*\Big)^2$ 
(by Poincar\'e duality).
If one now assumes that
all the gauge equivalence classes of flat connections $A$ are
isolated, in fact Freed and Gompf assumed $H^1(M,d_A) = 0$, so that the 
above sum is finite
and such that the square root of the Reidemeister torsion $\tau_M(A)^{1/2}$ 
is a well-defined
number (once a volume on $H^0(M, d_A)$ has been fixed, but for irreducible 
connections
$H^0(M,d_A)=0$).

The quantity $I_A \in \Z/8\Z$ denotes the spectral flow of the operator 
$$\left(\begin{matrix} \mboxsm{$\star d_{A_t}$} & \mboxsm{$-d_{A_t} \star$} \\
\mboxsm{$d_{A_t} \star$} & \mboxsm{$0$} \end{matrix}\right)$$
on $\Omega^1(M;{\frak g}) \otimes \Omega^3(M;{\frak g})$,
where $A_t$ is a path of connections
running from the trivial connection to $A$. They also looked at
some examples where $H^1(M,d_A) \neq 0$ and checked the overall
growth predicted by the above formula.

Following this Jeffrey \cite{Jeffrey} proposed the following
more general interpretation of the square root of Reidemeister torsion
in the cases where the connections are not isolated: Assume that
the moduli space of flat connections ${\mathcal M}$ on $M$
is smooth and that the tangent space at each equivalence class of
flat connection $A$ equals $H^1(M,d_A)$. Since $H^0(M,d_A)\subset
{\frak g}$ the invariant inner product we have chosen on ${\frak g}$
induces a volume element on $H^0(M,d_A)$. In total this means that
the square root of the Reidemeister torsion induces a measure on
the moduli space when we pair it with the induced volume element on 
$H^0(M,d_A)$
divided by the order of the center of $G$ and one arrives at ($r= k+
h^{{}^\vee}$)
\begin{align*}
Z^G_{k}(M) \underset{r\to\infty}{\sim}
&  e^{- \pi \sqrt{-1} (\mboxss{dim } G)(1+b^1(M))/4 }  \\*
& \hspace{-5pt}\times
\int_{[A] \in {\mathcal M}} e^{2\pi\sqrt{-1}r \mboxss{CS$(A)$}}
r^{(h^1_A - h^0_A)/2} e^{- 2\pi\sqrt{-1} (I_A/4 + (h^0_A + h^1_A)/8)}
\tau_M(A)^{1/2}.
\end{align*}
For some mapping tori of genus $1$ surfaces and lens spaces,\index{lens space} 
Jeffrey verified this form
of the semi-classical approximation.
Garoufalidis \cite{Garoufalidis_rel}
independently proved the semi-classical approximation for lens spaces and
studied in various examples
the growth rate predicted by these approximations.
Rozansky \cite{Rozansky_l}
proposed a further refined version of the above semi-classical approximation,
and offered calculations for a very large class of Seifert fibered
manifolds as evidence. He proposed to divide the volume element on
$H^0(M,d_A)$ by the volume of the stabilizer of $A$, and to use the
resulting quantity paired with the square root of the Reidemeister
torsion as the measure on ${\mathcal M}$ (generalizing the division by the
order of the center above).  This gave a natural
explanation of factors not accounted for in both the work of Freed
and Gompf and the work of Jeffrey. He also proposed corrections
to the formula for the growth rate of the invariant (i.e.\ the power of $r$ 
in the above),
in cases where not all
directions in $H^1(M,d_A)$ are tangent to paths of flat
connections (see \cite{Rozansky_l} and \cite{Rozansky_CMP171}).

Axelrod and Singer \cite{AxSi_I,AxSi_II} (see also \cite{Kontsevich_F})
considered the higher loop contributions in the perturbation expansion and
proposed the following:
\begin{equation}
Z^G_k(M)
\underset{k\to\infty}{\sim}
\int_{[A]\in {\mathcal M}}
\left( \begin{array}{l}
\mbox{semi-classical} \\
\mbox{approximation} \end{array}
\right) \times
\exp \left( \sum_{l=1}^\infty \frac{c^l k^{-l}}{ (2l)! (3l)! } \!\!\!
\sum_{e(\Gamma)=-l}
\!\!\! \frac{Z_\Gamma(M,A)}{ | \mbox{Aut}(\Gamma) | } \right)
\end{equation}
for some scalar $c$,
where the right sum runs over
connected trivalent graphs $\Gamma$ whose Euler number is equal to $-l$,
and $|\mbox{Aut}(\Gamma)|$ denotes
the order of the group of automorphisms of $\Gamma$.
Further, in the case where $A$ is acyclic or when $A\in{\mathcal M}$ is 
contained in a
smooth component, Axelrod and Singer was able to construct
$Z_\Gamma(M,A)$ as a topological invariant of $(M,A)$;
roughly speaking, it is given as follows in the acyclic case.
We identify the set of connection around $A$ with $\Omega^1(M,{\frak g})$.
The second order part of the Chern-Simons functional
gives a bilinear form on $\Omega^1(M,{\frak g})^{\otimes 2}$,
and it determines
a 2-form $L \in \Omega^2(M\times M, {\frak g} \otimes {\frak g})$
and its ``inverse''.
Further, the third order part of the Chern-Simons functional
gives a trilinear form $T$ on $\Omega^1(M,{\frak g})^{\otimes 3}$.
We obtain $Z_\Gamma(M,A)$
by contracting $L^{\otimes (3l)}$ by $T^{\otimes (2l)}$
``along the trivalent graph $\Gamma$''
(roughly regarding $L$ as in $\Omega^1(M,{\frak g})^{\otimes 2}$);
we determine the action of $T^{\otimes (2l)}$ on
$L^{\otimes (3l)} \in \Omega^1(M,{\frak g})^{\otimes (6l)}$
by putting copies of $L$ on $3l$ edges of $\Gamma$
and putting copies of $T$ on $2l$ vertices of $\Gamma$.
For a precise (mathematical) construction (and its topological invariance)
of $Z_\Gamma(M,A)$, see \cite{AxSi_I,AxSi_II}.

{}From the mathematical viewpoint
we regard $Z^G_k(M)$ as
$$
Z^G_k(M)
= \frac{\tau_{k+ h^{{}^\vee}}^G(M)}{\tau_{k+ h^{{}^\vee}}^G(S^1 \times S^2)}
$$
for the quantum $G$ invariant $\tau_r^G(M)$.
Then,
the asymptotic expansion of $Z_k^G(M)$ is predicted by
the semi-classical approximation and its higher loop corrections
stemming from a perturbative expansion of the Chern-Simons path integral,\index{Chern-Simons!--- path integral!perturbative expansion of ---}
explained above in some cases.
This leads us to the following somewhat vague conjecture.

\begin{conj}[The perturbative expansion conjecture] 
\label{conj.AEC}
The\index{quantum invariant!--- of 3-manifold!asymptotic behaviour of ---}
asymptotic expansion of $Z_k^G(M)$ of a closed oriented 3-manifold $M$
is given by
\begin{align*}
Z_k^G(M) \, \underset{k\to\infty}{\sim} \,
& e^{- \pi \sqrt{-1} (\mboxss{\rm dim } G)(1+b^1(M))/4 } \\*
& \times \int_{[A] \in {\mathcal M}} \!\!\!
e^{2\pi\sqrt{-1}r \mboxss{CS$(A)$}}
r^{(h^1_A - h^0_A)/2}
e^{- 2\pi\sqrt{-1} (I_A/4 + (h^0_A + h^1_A)/8)}
\tau_M(A)^{1/2} \\*
& \times \exp \left( \sum_{l=1}^\infty \frac{c^l k^{-l}}{ (2l)! (3l)! } \!\!
\sum_{e(\Gamma)=-l}
\!\! \frac{Z_\Gamma(M,A)}{ | \mbox{\rm Aut}(\Gamma) | } \right),
\end{align*}
putting $r = k + h^{{}^\vee}$,
where the right hand side can be given in the mathematical viewpoint in 
certain cases,
as mentioned above, but which needs further interpretation in general.
\end{conj}

\begin{rem}
The semi-classical approximation stated above
(the upper two lines in the above formula),
has been confirmed
for lens spaces\index{lens space} 
(first partially \cite{FrGo}) and then
by \cite{Jeffrey, Garoufalidis_rel}, 
for certain mapping tori of diffeomorphisms of a torus \cite{Jeffrey}, and for
all finite order mapping
tori of automorphisms
of any closed orientable surface of genus at least 2 \cite{Andersen}.
For a large class of Seifert fibered manifolds \cite{Rozansky_CMP171} and 
\cite{Rozansky_3}
offered
calculations which provided evidence that the phases in the semi-classical 
approximation
is given by the Chern-Simons invariants
and the measure is given by the 
square root of
the Reidemeister torsion as explained above. Also, expressions for the 
higher loop
corrections was offered. Later the necessary
analytic estimates was provided in \cite{Hansen} so as to confirm
this. See also the discussion below. For now, there are no
examples of hyperbolic manifolds, where parts of the above
conjecture has been confirmed.

For other versions of Conjecture \ref{conj.AEC},
see \cite[Problem 3.108]{Kirby}, \cite{Garoufalidis_appl}.

\end{rem}

The formula in Conjecture \ref{conj.AEC}
might not give an exact description of
the asymptotic behavior of $Z_k^G(M)$
even in the semi-classical part, neither is it in all cases well-defined.
Moreover, it might be difficult at present
to calculate the concrete value of the higher loop corrections
in the asymptotic expansion of Conjecture \ref{conj.AEC}
for given $M$, $A$, and $\Gamma$. Nor do we have definitions for
these terms, which has been proven to be well defined topological
invariants in all cases.

The following conjecture offers a kind of reverse viewpoint on Conjecture 
\ref{conj.AEC},
avoiding such ambiguities and difficulties.

\begin{conj}[The asymptotic expansion conjecture, J.E. Andersen 
\cite{Andersen}] 
\label{conj.AEC_Jorgen}
Let $\{ c_0=0, c_1, \cdots, c_m \}$ be
the set of values of the Chern-Simons functional of flat $G$ connections
on a closed oriented 3-manifold $M$.
There exist
$d_j \in \Q$,
$\tilde I_j \in \Q/\Z$,
$v_j \in \R_+$,
and $a^e_j \in \C$
for $j=0,1,\cdots, m$ and $e = 1,2,3,\cdots$ such that ($r = k +
h^{{}^\vee}$)\index{quantum invariant!--- of 3-manifold!asymptotic behaviour of ---}
$$
Z_{k}^{G}(M) \ \underset{r\to\infty}{\sim} \
\sum_{j=0}^m e^{2\pi\sqrt{-1}r c_j}
r^{d_j} e^{\pi\sqrt{-1} \tilde I_j /4} v_j
\big(1 + \sum_{e=1}^\infty a^e_j r^{-e} \big),
$$
that is, for all $E = 0, 1, 2,\ldots$, there exists a constant $c_E$ such that
$$
\Big| Z_{k}^{G}(M) -
\sum_{j=0}^m e^{2\pi\sqrt{-1}r c_j}
r^{d_j} e^{\pi\sqrt{-1} \tilde I_j /4} v_j
\big(1 + \sum_{e=1}^E a^e_j r^{-e}\big) \Big| \le c_E r^{d-E-1}
$$
for all $r=2,3,4,\cdots$.
Here, $d = \mbox{max} \{ d_0, \cdots, d_m \}$.
\end{conj}

\begin{rem}[\rm (J.E. Andersen)] 
If such an expansion in the above conjecture exists,
then $c_j$, $d_j$, $\tilde I_j$, $v_j$, and $a_j^e$
are uniquely determined by $Z_{k+2}^{G}(M)$ for $k=0,1,2,3,4,\cdots$.
\end{rem}

\begin{prob}[J.E. Andersen] 
\label{prob.pexp_O_H}
If such an expansion exists, understand how it is related to the expansion of
Ohtsuki and the expansion of Habiro.
\end{prob}

It will of course be important to establish, that an
expansion of this type exists, however, of far greater importance will be to
give independent topological meaning to the many resulting new
invariants, e.g.\ to prove that the phases are the Chern-Simons values $c_j$.
 From the discussion above on the semi-classical approximation we
derive the following conjecture:

\begin{conj}[Topological interpretations of the $d_j$'s] 
\label{conj.dj}
Let ${\mathcal M}_j$ be the union of components of the moduli
space of flat connections ${\mathcal M}$ which has Chern-Simons
value $c_j$. Then
\[d_j = \frac{1}{2}\max_{A\in {\mathcal M}_j} (h^1_A - h^0_A),\]
where $\max$ here means the maximum value that $(h^1_A - h^0_A)$
assumes on a Zariski open subset of ${\mathcal M}_j$.
\end{conj}

Note that this conjecture might be rather optimistic, and may only
hold in the non-degenerate cases. However, we do not know of any
cases where it fails (see \cite{Garoufalidis_appl}).

\begin{rem}[\rm (J.E. Andersen)] 
The special $\max$ proposed in Conjecture \ref{conj.dj} is certainly
needed, as shown by the example of the mapping torus of the
diffeomorphism $- Id$ of a torus. 
The quantum $SU(2)$ invariant
of
this manifold is easily seen to be $r-1$, since $- Id$ is
represented trivially for all levels, however, there are flat $SU(2)$ 
connections
for which $(h^1_A - h^0_A) > 2$.
\end{rem}

The Conjecture \ref{conj.dj} implies the following growth rate.

\begin{conj}[The growth rate conjecture] 
\label{conj.gr}
Let $d = \max \{d_0, \ldots, d_n\}.$ 
Then\index{quantum invariant!--- of 3-manifold!growth rate of ---}
$|Z_{r}^{G}(M)| = O(r^d).$
\end{conj}

It is well known that the quantum invariants only grows like
$r$ to some power. The power is bounded from above by some simple
function (depending on G) of the Heegaard
genus of the manifold.

\begin{rem}[\rm (J.E. Andersen)]\ \mbox{}

(1)\qua Suppose that $M$ is a closed 3-manifold satisfying that
$\tau^{G}_r(M) = \tau^{G}_r(S^3)$ for all $r$.
If the growth rate conjecture \ref{conj.gr} is true for the group $G$,
then there is no non-central representation of $\pi_1(M)$ to $G$.

(2)\qua  Kronheimer and Mrowka have proposed a program using Seiberg-Witten 
theory and
Floer homology to establish
that any 3-manifold $M$ obtained from $S^3$
by $+1$ surgery along a non-trivial knot $K$ has a non-trivial
(and therefore non-abelian)
representation of $\pi_1(M)$ to $SU(2)$.
Suppose that this is the case and the growth rate conjecture \ref{conj.gr} 
is true.
Then, $J_{K,c} = J_{U,c}$
for all $c = 1,2,\cdots$
if and only if $K$ is the trivial knot $U$,
where $J_{K,c}$ denotes
the colored Jones polynomial\index{Jones polynomial!colored ---} 
of a knot $K$ with a color $c$.
\end{rem}

At this time we do not know of a topological interpretation of the values
of $\tilde I_j$ and $v_j$ which makes sense in all cases. Let us
simply just propose the following
\begin{conj} \label{conj.rt}
There is a construct of the {\em right} measure, say $\tau_M(A)^{1/2}$ for 
$A\in {\mathcal M}_i$, from the
square root of the Reidemeister torsion generalizing the
non-degenerate case explained above and such that
\[e^{\pi\sqrt{-1} \tilde I_j /4} v_j = \int_{A\in {\mathcal M}_i}
e^{\pi\sqrt{-1} (-2 I_A + h^0_A + h^1_A) /4}\tau_M(A)^{1/2}.\]
\end{conj}

Conjectures \ref{conj.AEC_Jorgen} and \ref{conj.dj} together with
Conjecture \ref{conj.rt} were first proved for mapping tori of all
finite order diffeomorphisms of all surfaces of genus at least two
in \cite{Andersen}. Recently, Conjecture
\ref{conj.AEC_Jorgen} was proved for all Seifert fibered spaces in
\cite{Hansen} by supplementing the calculations in \cite{Rozansky_CMP171} 
and \cite{Rozansky_3}
with the need analytic estimates.

\begin{exm}
Let us illustrate the asymptotic behavior
of the quantum $SU(2)$ invariant\index{quantum invariant!--- of 3-manifold!asymptotic behaviour of ---}
of the lens space $L(5,1)$ of type $(5,1)$.\index{lens space}
For simplicity, we let $r$ be an odd prime.
Since $\tau_3^{SU(2)} \big( L(5,1) \big) =1$,
putting $\zeta = \exp(2\pi\sqrt{-1}/r)$, we have that
$$
\tau_r^{SU(2)} \big( L(5,1) \big) =
\tau_r^{SO(3)} \big( L(5,1) \big) =
\left( \frac{5}{r} \right)
\zeta^{- 3 \cdot {5}^\star}
\frac{ \zeta^{{10}^\star} - \zeta^{-{10}^\star} }
      {\zeta^{2^\star}-\zeta^{-2^\star}},
$$
where $k^\star$ denotes the inverse of $k$ in $\Z/r\Z$.
Since $\tau_r^{SU(2)}(S^1 \times S^2) = \sqrt{\frac{r}2} / \sin (\frac\pi{r})$,
we have that
\begin{equation}
\label{eq.ZL1}
Z^{SU(2)}_{r-2} \big( L(5,1) \big)
= \frac{\tau_r^{SU(2)} \big( L(5,1) \big) }
        {\tau_r^{SU(2)} \big( S^1 \times S^2 \big) }
= \sqrt{\frac{2}{r}} \sin \frac{\pi}{r}
\left( \frac{5}{r} \right)
\zeta^{- 3 \cdot {5}^\star}
\frac{ \zeta^{{10}^\star} - \zeta^{-{10}^\star} }
      {\zeta^{2^\star}-\zeta^{-2^\star}}.
\end{equation}
On the other hand,
as in \cite{Jeffrey},
the semi-classical approximation is given as follows.
The lens space $L(5,1)$\index{lens space}
 has three flat connections $A_n$ (n=0,1,2);
each $A_n$ is determined by
the representation of $\pi_1\big( L(5,1) \big) \cong \Z/5\Z$ to $SU(2)$
which takes a generator of $\Z/5\Z$ to
$\left(\begin{matrix} e^{2\pi\sqrt{-1}n/5} & 0 \\
                       0 & e^{-2\pi\sqrt{-1}n/5} \end{matrix}\right)$.
As in \cite{Jeffrey}, we have that
$\mbox{CS}(A_n) = n^2 / 5$,
$h^0_{A_n} = 1$,
$h^1_{A_n} = 0$,
$\tau_M(A_n)^{1/2} = \frac{4\sqrt{2}}{\sqrt{5}} \sin^2 \frac{2\pi n}{5}$,
and $I_n (\mbox{mod $4$}) = 1$ if $n < 5/2$, and $-1$ if $n > 5/2$.
Hence,
\begin{equation}
\label{eq.ZL2}
Z_{r-2}^{SU(2)}\big( \overline{L(5,1)} \big) \ \underset{r\to\infty}{\sim} \
2 \sqrt{\frac{-2}{5 r}}
\sum_{n=0,1,2} e^{2\pi\sqrt{-1} r n^2/5}
\sin^2 \frac{2\pi n}{5},
\end{equation}
noting that
the notatin of lens spaces in \cite{Jeffrey}
is equal to the notation of their mirror images in 
\cite{KM_inv,Garoufalidis_rel}.

The sequence of $\tau_r^{SU(2)}\big( L(5,1) \big)$ for odd primes $r$
splits into four subsequences according to $r \equiv \pm1, \pm3$ mod $10$,
and each subsequence can be approximated by a function of a polynomial order.
Let us describe the subsequence, say,
with $r \equiv -1$ mod 10, as follows.
Since $10^\star = (r+1)/10$,
we calculate (\ref{eq.ZL1}) as
\begin{align*}
Z_{r-2}^{SU(2)}\big( L(5,1) \big)
&= \sqrt{\frac{2}{r}} \sin \frac{\pi}{r}
e^{-6\pi\sqrt{-1}/5r}
\frac{ e^{\pi\sqrt{-1}/5r} - \omega^{-1} e^{-\pi\sqrt{-1}/5r} }
{e^{\pi\sqrt{-1}/r} - e^{-\pi\sqrt{-1}/r} } \\*
& \!\!\! \underset{r\to\infty}{\sim} \
\frac{1 - \omega^{-1}}{\sqrt{-2}} r^{-1/2},
\end{align*}
putting $\omega = \exp(2\pi\sqrt{-1}/5)$.
On the other hand,
the right hand side of (\ref{eq.ZL2}) is calculated as
$$
2 \sqrt{ \frac{-2}{5r} } \Big(
e^{-2\pi\sqrt{-1}/5} \sin^2 \frac{2\pi}{5}
+ e^{2\pi\sqrt{-1}/5} \sin^2 \frac{4\pi}{5} \Big)
= \frac{\omega-1}{\sqrt{-2}} r^{-1/2},
$$
noting that
$\sqrt{5} = 1 + 2 \omega + 2 \omega^{-1}$ (Gaussian sum).
Therefore, it was verified that
the semi-classical approximation is correct
for this subsequence.

This is related to the perturbative invariant
$$
\tau^{SO(3)} \big( L(5,1) \big) =
q^{- 3/5} \frac{q^{1/10} - q^{-1/10}}{q^{1/2}-q^{-1/2}} \in \Q[[q-1]]
$$
as follows.
We regard it as a holomorphic function of $q$ in a suitable domain.
The asymptotic behavior of $\tau_r^{SO(3)}\big( L(5,1) \big)$,
say, for the above mentioned subsequence,
can be presented by using this holomorphic function around $q^{1/5}=\omega$.
\end{exm}

\begin{exm}
It is known, see \cite{LaZa,Le_isp}, that
$$
\tau_r^{SO(3)}\big( \Sigma(2,3,5) \big) = \frac{1}{1-\zeta}
\sum_{n=0}^{r-1} \zeta^n (1-\zeta^{n+1})(1-\zeta^{n+2})\cdots (1-\zeta^{2n+1})
$$
for Poincare homology 3-sphere $\Sigma(2,3,5)$,
where we put $\zeta = \exp(2\pi\sqrt{-1}/r)$.
It is an exercise to compute the asymptotic behaviour of
$Z_{r-2}^{SU(2)}\big( \Sigma(2,3,5) \big)$ as $r\to\infty$
related to Conjecture \ref{conj.AEC_Jorgen},
and to formulate a relation with the perturbative invariant given by
$$
\tau^{SO(3)}\big( \Sigma(2,3,5) \big) = \frac{1}{1-q}
\sum_{n=0}^\infty q^n (1-q^{n+1})(1-q^{n+2}) \cdots (1-q^{2n+1}).
$$
\end{exm}


\subsection{The volume conjecture}

It is known (see Conjecture \ref{conj.gr} and its remark) that
the asymptotic behaviour 
of the quantum $SU(2)$ invariant\index{quantum invariant!--- of 3-manifold!asymptotic behaviour of ---}  
$\tau_N^{SU(2)}(M)$ as $N \to \infty$
is a polynomial growth in $N$.
Nevertheless, this asymptotic behaviour
might be regarded as an exponential growth
in the sense of the following conjecture,
which is a 3-manifold version of
the volume conjecture (Conjecture \ref{conj.vol_conj}).

\begin{conj}[H. Murakami \cite{Hitoshi_opt}] 
\label{conj.vol_olim}
For any closed 3-manifold $M$,\index{volume conjecture!--- for 3-manifolds}
$$
2 \pi \sqrt{-1} \cdot 
\underset{N \to \infty}{\mbox{o-$\lim$}}
\frac{\log \tau_N^{SU(2)}(M)}{N}
= {\rm CS}(M) + \sqrt{-1} \mbox{\rm vol}(M),
$$
where $\mbox{\rm vol}(M)$ and ${\rm CS}(M)$
denote the\index{hyperbolic volume!--- and quantum invariant of 3-manifold}\index{simplicial volume!--- and quantum invariant of 3-manifold}  
hyperbolic volume\footnote{\footnotesize 
When $M$ is not hyperbolic,
we define $\mbox{\rm vol}(M)$ to be $v_3 ||M||$,
where $||M||$ is the simplicial volume and
$v_3$ is the hyperbolic volume of the regular ideal tetrahedron.}
and the Chern-Simons\index{Chern-Simons!--- invariant!--- and quantum invariant of 3-manifold} 
invariant\footnote{\footnotesize 
It is also conjectured (see Problem \ref{prob.topCS})
that there exists an appropriate definition of 
${\rm CS}(M)$ of any closed 3-manifold $M$,
though ${\rm CS}(M)$ is defined only for hyperbolic 3-manifolds $M$ 
at present.}
of $M$ respectively,
and $\mbox{o-$\lim$}$ denotes the 
``optimistic limit''\index{optimistic limit}
introduced in \cite{Hitoshi_opt}.
\end{conj}

\begin{rem}
As mentioned in \cite{Hitoshi_opt}
the ``definition'' of the optimistic limit is not rigorous yet,
because there is some ambiguity in the present definition,
where formal approximation,
such as (\ref{eq.qn_app}) and (\ref{eq.qn_app2}), are used.
It is a problem to find a rigorous formulation of the optimistic limit.
\end{rem}

\begin{rem}
It is shown \cite{Hitoshi_opt}, by using formal approximations,
that Conjecture \ref{conj.vol_olim} is ``true''
for closed 3-manifolds obtained from $S^3$ by surgery 
along the figure-eight knot.
\end{rem}

\begin{rem}
R. Benedetti gave another formulation of the volume conjecture
by using quantum hyperbolic invariants;
see Conjecture \ref{benedetti.conj}.
\end{rem}

\begin{rem}
The statement of Conjecture \ref{conj.vol_olim}
should extend for knot (link) complements $M$,
which should be related to
the volume conjecture for knots (Conjecture \ref{conj.vol_CS}).
\end{rem}

\begin{rem}
By formally applying the (infinite dimensional) saddle point method 
to the Chern-Simons path integral,\index{Chern-Simons!--- path integral!saddle point method on ---}
the value (\ref{eq.CSA=CSM+vol}) appears
at a critical point of the Chern-Simons functional.
This might give a physical explanation of Conjecture \ref{conj.vol_olim}.
Can we justify it in mathematics?
There is an approach,
by using knotted trivalent graphs 
(see Conjecture \ref{conj.KTG_ae}), 
to justify the Chern-Simons path integral mathematically,
which might be helpful to apply
the saddle point method to it rigorously.
\end{rem}

\begin{prob}[H. Murakami] 
\label{prob.olim_Sf}
Calculate $\mbox{o-$\lim$} \frac{\log \tau^{SU(2)}_N(M)}{N}$
for Seifert fibered 3-manifolds $M$.\index{optimistic limit}
\end{prob}

\begin{rem}
When $M$ is a mapping torus of a homeomorphism of a surface,
a quantum invariant\index{quantum invariant!--- of 3-manifold!optimistic limit of ---} 
 of $M$ can be presented by
the trace of the linear map on the quantum Hilbert space
associated to the homeomorphism.
Such a presentation might be useful
to compute the asymptotic behaviour of $\tau_N^{SU(2)}(M)$.
\end{rem}

\begin{rem}
When we choose a simplicial decomposition of $M$,
(the absolute value of) its quantum invariant can be expressed by
using quantum $6j$-symbols.
The computation of the asymptotic behaviour of $\tau_N^{SU(2)}(M)$
might be reduced to the computation of limits of quantum $6j$-symbols.
J. Roberts \cite{vol.R} showed that
a limit of classical $6j$-symbols is given by
the Euclidean volume of a tetrahedron.
Further, J. Murakami and M. Yano \cite{MuYano} recently showed that
a limit of quantum $6j$-symbols\index{6j-symbol!quantum ---!--- and hyperbolic volume} 
is related to
the hyperbolic volume of a tetrahedron
via formal approximation such as (\ref{eq.qn_app}) and (\ref{eq.qn_app2}).
\end{rem}

\begin{prob}[D. Thurston] 
\label{prob.inv_vol}
Find a series of invariants of a 3-manifold
(depending on roots of unity)
that grows as its hyperbolic volume\index{hyperbolic volume!--- and invariants of 3-manifold}   
(or its simplicial volume).\index{simplicial volume!--- and invariants of 3-manifold}  
\end{prob}

\begin{prob}[D. Thurston] 
\label{prob.vol_conj_ncL}
Find a correct generalization of 
the volume conjecture\index{volume conjecture!--- for other Lie groups}
to other non-compact Lie groups.
\end{prob}

\begin{rem}
The volume conjecture is related to
the $SL(2,\C)$ Chern-Simons theory,\index{Chern-Simons!--- theory!--- for $SL(2,\C)$}
which (formally) deduces the hyperbolic volume and the Chern-Simons invariant.
It is a problem to find (or formulate)
such invariants of 3-manifolds
for other non-compact Lie groups.
\end{rem}

The {\it Chern-Simons functional}
${\rm CS}(A) \in \C$
of a $SL(2;\C)$ connection $A$
on a closed 3-manifold $M$ is defined by the formula (\ref{eq.CSfunct}),
where we regard $A$ in the formula as a $sl(2;\C)$-valued 1-form on $M$
in this case.
Since a gauge transformation of $A$ changes ${\rm CS}(A)$ by an integer,
${\rm CS}([A])$ of the gauge equivalence class of $A$
is defined to be in $\C/\Z$.
The {\it Chern-Simons invariant}\index{Chern-Simons!--- invariant} 
${\rm CS}(M) \in \R/\Z$
and the volume $\mbox{\rm vol}(M) \in \R_{>0}$ 
of a closed hyperbolic 3-manifold $M$ is given by\footnote{\footnotesize 
The Chern-Simons invariant
was introduced by
Chern and Simons \cite{ChernSimons}
as an invariant of compact $(4n-1)$-dimensional Riemannian manifolds.
For hyperbolic 3-manifolds,
Meyerhoff \cite{Meyerhoff_h} extended ${\rm CS}(M)$
for $M$ with cusps.
See also \cite{Neumann_H,CGHN00}
for ${\rm CS}(M)$ of hyperbolic 3-manifolds $M$
as a counterpart of $\mbox{\rm vol}(M)$.}
\begin{equation}
\label{eq.CSA=CSM+vol}
{\rm CS}([A_0]) = {\rm CS}(M) + \sqrt{-1} \mbox{\rm vol}(M),
\end{equation}
where $[A_0]$ is the gauge equivalence class 
of a $SL(2;\C)$ flat connection $A_0$
associated to the conjugacy class of a holonomy representation
$\pi_1(M) \to SL(2;\C)$ of the hyperbolic structure on $M$.
Further, 
when $M$ is the complement of a hyperbolic knot (link) in a closed 3-manifold,
${\rm CS}(M)$ can be defined similarly.

\begin{prob}[S. Morita \cite{KojimaNegami}] 
\label{prob.topCS}
Define the Chern-Simons invariant\index{Chern-Simons!--- invariant!topological definition of ---}\index{Chern-Simons!--- invariant!--- of non-hyperbolic 3-manifold}  
${\rm CS}(M)$
as a topological invariant of any closed oriented 3-manifold $M$,
and of any knot (link) complement $M$ in a closed 3-manifold.
\end{prob}

This problem includes two problems: 
to define ${\rm CS}(M)$ 
(topologically or combinatorially) as a topological invariant,
and to define it for non-hyperbolic 3-manifolds.

\begin{rem}
The hyperbolic volume
(which is a counterpart of the Chern-Simons invariant)
has a definition
as a constant multiple of the simplicial volume,
which is combinatorial, 
and can be applied,
not only for hyperbolic 3-manifolds,
but also for any other 3-manifolds.
\end{rem}

\begin{rem}[{\rm (S. Kojima)}] 
The Chern-Simons invariant
${\rm CS}(M)$ of non-hyperbolic 3-manifolds $M$
should be defined satisfying the following two requirements.
One is that
${\rm CS}(-M) = - {\rm CS}(M)$,
where $-M$ denotes $M$ with the opposite orientation.
The other is the requirement explained as follows.
Let $K$ be a hyperbolic knot in a 3-manifold $N$.
Then, it is known that $N_{K;(p,q)}$ has a hyperbolic structure
except for finitely many $(p,q)$,
where $N_{K;(p,q)}$ denotes the 3-manifold obtained from $N$
by Dehn surgery along the slope of type $(p,q)$,
and that such hyperbolic structures can be obtained
in a deformation space of the hyperbolic structures of $N-K$
parameterized by a natural complex parameter,
which can be presented by two real parameters $p$ and $q$.
Moreover, the function 
\begin{equation}
\label{eq.CS+vol}
{\rm CS}(M) + \sqrt{-1} \mbox{\rm vol}(M)
\end{equation}
is a holomorphic function of the complex parameter.
Note that $\mbox{\rm vol}(M)$ can extend for non-hyperbolic 3-manifolds $M$
by redefining it to be a constant multiple of the simplicial volume $||M||$.
${\rm CS}(M)$ should be defined such that, 
for appropriate\footnote{\footnotesize 
These knots should include,
not only all hyperbolic knots, but also other knots.
They might not include the trivial knot.}
knots $K$ in any closed 3-manifold $N$, 
the function (\ref{eq.CS+vol}) on the family $\{ N_{K;(p,q)} \}_{p,q}$
can extend to a holomorphic function
of a complex parameter presented by $p$ and $q$ appropriately.
\end{rem}

\begin{prob}[T. Ohtsuki] 
\label{prob.Cstr_3mfd}
Give a ``complex structure'' to the set of 3-manifolds.
More precisely, 
find an embedding (or, an immersion) of the set of 3-manifolds
to some complex variety such that
its restriction to the set 
$\{ N_{K; (p,q)} \ | \ p^2+q^2 > \!\! > 0 \}$ 
can be extended to a holomorphic map of 
the above mentioned complex parameter
for any (hyperbolic) knot $K$ in any 3-manifold $N$.
\end{prob}

We would expect some structures of the set of 3-manifolds
such as mentioned in Problems \ref{prob.Cstr_3mfd} 
and Problem \ref{prob.product_in_M}.
Such structures would yield
new viewpoints in the study of (the set of, and invariants of) 3-manifolds.

\begin{rem}
As mentioned above,
the set $\{ N_{K; (p,q)} \ | \ p^2+q^2 > \!\! > 0 \}$
can be embedded in $\C$,
on which the function (\ref{eq.CS+vol}) is holomorphic.
In this sense, the infinite family of $N_{K,(p,q)}$
has a ``complex structure'' around the infinity point of $(p,q)$.
The volume conjecture says that
the function (\ref{eq.CS+vol}) would be obtained
as a certain limit of some series of quantum invariants.
This suggests that the above ``complex structure'' would extend to 
the whole set of 3-manifolds.
\end{rem}


\subsection{Quantum hyperbolic invariants of 3-manifolds}
\label{sec.qhypinv}

\renewcommand{\thefootnote}{\fnsymbol{footnote}}
\footnotetext[0]{Section \ref{sec.qhypinv} was written by 
S. Baseilhac and R. Benedetti.}
\renewcommand{\thefootnote}{\arabic{footnote}}

\noindent The main references for this section are
\cite{bene.B,bene.B-B1,bene.B-B2}, a review being \cite{bene.B-B4}. In
\cite{bene.B-B2} the ideas of sections 7-9 in \cite{bene.B-B1} are
developed with some important differences in the way they are
concretized.

\smallskip

\noindent Let $W$ be a compact closed oriented $3$-manifold, $L\subset W$ be a
non-empty link, $\rho$ be a flat principal $B$-bundle on $W$; $B$ is
the upper triangular Borel subgroup of $SL(2,\C)$. In \cite{bene.B-B1}
one constructs a family of 
``quantum hyperbolic invariants'' (QHI)\index{quantum invariant!quantum hyperbolic invariant}
$K_N(W,L,\rho)\in \C$, where $N > 1$ is any odd integer. This 
consists of two main steps:
\begin{itemize}
\item[\rm(1)] For every triple $(W,L,\rho)$, the construction of
so-called $\cal D$-{\it triangulations} ${\cal T}=(T,H,{\cal D})$, where:
$(T,H)$ is a (singular) triangulation of $(W,L)$ such that each edge
has distinct vertices and $H$ contains all the vertices of $T$; the
``decoration'' $\cal D$ is made of a {\it full} simplicial $B$-1-cocycle
representing $\rho$ on $W$, a {\it branching} (for instance one induced by
a total ordering of the vertices of $T$), and an {\it integral
charge}. For these notions, see \cite{bene.B-B4}.
\item[\rm(2)] The proof that a suitable state sum
$H_N({\cal T})$ does
not depend on the choice of the $\cal D$-triangulation ${\cal T}$ up to
multiplication by $N$-th roots of unity, so that
$K_N(W,L,\rho)=K_N({\cal T})=H_N({\cal T})^N$ actually defines an invariant.
\end{itemize}

\noindent The proof of the existence of $\cal D$-triangulations is
difficult essentially due to strong global constraints in $\cal D$. The main
building-blocks of the state sums $H_N({\cal T})$ are the
``quantum-dilogarithm'' $6j$-symbols
 of the $N$-dimensional cyclic
representations of a quantum Borel subalgebra of
$U_{\omega}(sl(2,\C))$, where $\omega = \exp(2\pi i/N)$. Kashaev
proposed in \cite{bene.K1} a conjectural purely topological invariant
$K_N(W,L)$ which should have been expressed by a state sum of this
kind (although in his proposal there were no flat bundles and no
notion of $\cal D$-triangulation); in fact, $K_N(W,L)$ appears as a
special case of $K_N(W,L,\rho)$ when $\rho$ is the trivial flat
$B$-bundle on $W$. The algebraic properties of the $6j$-symbols ensure
the invariance of $K_N({\cal T})$ up to certain elementary moves on
$\cal D$-triangulations.  Then, the proof of the full invariance of
$K_N({\cal T})$ consists in connecting by such elementary moves any
two $\cal D$-triangulations of $(W,L,\rho)$, which is not so easy to
achieve.

\begin{prob}[S. Baseilhac, R. Benedetti] 
\label{bene.0} 
Generalize the construction of the QHI for flat principal $G$-bundles,
for Lie groups $G$ different from $B$.
\end{prob}

\begin{rem}
A basic ingredient of the $B$-QHI is the relationship between the
cyclic representation theory of a quantum Borel subalgebra of $U_{\omega}(sl(2,\C))$, and flat $B$-bundles encoded by simplicial full
$1$-cocycles. This relationship relies on the theory of quantum
coadjoint action of \cite{bene.DCP}, which holds for other Lie groups
such as $G=SL(2,\C)$.
\end{rem}

\begin{prob}[S. Baseilhac, R. Benedetti] 
\label{bene.1} 
Fix $(W,L)$ and vary $\rho$. Study $K_N$ as a function of the bundle,
that is as a function defined on the character variety of $W$ with
respect to $B$: regularity, fibers, and so on.
\end{prob}

\begin{rem}
Denote by $z$ the $B$-$1$-cocycle in $\cal T$ that represents $\rho$. The
state sum $K_N({\cal T})$ is a rational function of the upper diagonal
entries of the whole set of values of $z$. Moreover, the
$6j$-symbols
are rational functions of the moduli of the idealized
triangulation $\widehat{F}({\cal T})$ defined below.
\end{rem}  

\noindent Every $\alpha \in H^1(W;\C)$ leads to two flat $B$-bundles
$\rho_\alpha$ and $\rho'_\alpha$ defined as follows. The first one is
obtained via the natural identification of $(\C,+)$ with the parabolic
subgroup $Par(B)$ of $B$. The second one is obtained by means of the
exponential map of $(\C,+)$ onto the multiplicative $\C^*$, and the
identification of $\C^*$ with the diagonal Cartan subgroup $C(B)$ of
$B$. Similarly, every class in $H^1(W;\Z/p\Z)$ leads to a $B$-bundle
by the natural embedding of $\Z/p\Z$ into the group $S^1\subset
\C^*$.

\begin{prob}[S. Baseilhac, R. Benedetti] 
\label{prob.bene2}
Specialize Problem \ref{bene.1} to bundles coming from the ordinary
cohomology as above. For real additive ones, analyze the behaviour of
the QHI with respect to Thurston's norm. Are they constant on the
faces of the corresponding unit sphere ?
\end{prob}  

\begin{rem}
The ``projective invariance'' property of the QHI (see
\cite{bene.B-B1,bene.B-B4}) implies in particular that they are
constant on the rays of $H^1(W;\R)$.
\end{rem}

\begin{prob}[S. Baseilhac, R. Benedetti] 
\label{prob.bene3}
Understand the `phase factor' (i.e.\ the
ambiguity due to $N$-th roots of unity) of the state sum
$H_N({\cal T})$. Possibly derive from it an invariant for $(W,L,\rho)$
endowed with some extra-structure, thus refining
$K_N(W,L,\rho)$.
\end{prob}

\begin{rem}
The phase factor uniquely depends on the branching and the integral
charge in the decoration $\cal D$. On one hand, it is known
that branchings can be used to encode, for instance, combings,
framings, spin structures and so on. On another hand, combings induce
the extra-structure that allows Turaev's refinement of Reidemeister
torsions.
\end{rem}

\begin{prob}[S. Baseilhac, R. Benedetti] 
\label{prob.bene4} Determine a suitable $(2+1)$ `decorated' cobordism theory 
supporting a (non purely topological) QFT containing the already
defined QHI. Study in particular the behaviour of the QHI with respect
to connected sums.
\end{prob}

\begin{prob}[S. Baseilhac, R. Benedetti] 
\label{prob.bene5} Develop a 4-dimensional theory of QHI based on 
Turaev's shadow theory.
\end{prob}
\begin{rem} A first step should be to determine the right notion of 
$\cal D$-shadow together with a geometric interpretation.  In this direction,
F. Costantino is completing his PhD thesis at Pisa, where he shows in particular that `branched shadows' do encode Spin$^{c}$ structures.
\end{rem}

\begin{prob}[S. Baseilhac, R. Benedetti] 
\label{prob.bene6} Determine the actual relationship between $K_N(S^3,\cdot)$ and the
coloured Jones polynomial\index{Jones polynomial!colored ---} 
$J_N(\cdot)$ (evaluated at
$\omega=\exp(2i\pi/N)$ and normalized by $J_N(\rm{unknot})=1$), as
functions of links.
\end{prob}

\begin{rem}
(1)\qua In \cite{vol.MM} it is shown that $J_N$ may be defined by means of usual
$(1,1)$-tangle presentations (as for the Alexander polynomial), using an enhanced Yang-Baxter operator whose $R$-matrix is derived
from the quantum-dilogarithm $6j$-symbols.
This suggests that there
could be a relationship between $K_N(S^3,\cdot)$ (necessarily associated to the
trivial flat $B$-bundle on $S^3$) and $J_N(\cdot)^N$. The most immediate guess
would be that $K_N(S^3,L)=J_N(L)^N$ for each $L$. In fact, one can
give an $R$-matrix formulation of $K_N(S^3,\cdot)$ involving $R$-matrices
depending on parameters. These parameters are specified in terms of
the decorations of special $\cal D$-triangulations adapted to planar link
diagrams \cite{bene.B-B3}. So $K_N(S^3,\cdot)$ can be computed by
using suitably decorated link diagrams, and the decoration must
satisfy non trivial global constraints. In this setup, $(1,1)$-tangle
presentations do not play any role. On another side, the constant
$R$-matrix used for $J_N$ corresponds to one {\it fixed} particular
choice in the parameters. This is not enough to confirm the above
guess.

\noindent (2)\qua A motivation of Problem \ref{prob.bene6} is also to make working for $J_N$ a theory of scissors congruence classes, as described below for the QHI.  
\end{rem}

\noindent The so-called {\it Volume Conjectures} concern the
asymptotic behaviour of the invariants constructed on the base of the
quantum dilogarithm $6j$-symbols,
that is of $K_N(W,L,\rho)$ or
$J_N(L)$ (for $L \subset S^3$), when $N \to \infty$. They are
originally motivated by the asymptotic behaviour of the quantum
dilogarithm $6j$-symbols, whose dominant term involves dilogarithm
functions that may be used to compute the volume of oriented ideal
hyperbolic tetrahedra. In the case of $J_N(L)$ there are also some
numerical computations (sometimes using formal manipulations) - see
for instance the first section of the present volume for details. In
the case of QHI, we develop in \cite{bene.B-B1,bene.B-B2} (see also
\cite{bene.B-B4}) a theory of scissors congruence classes
for triples $(W,L,\rho)$ which gives a natural framework for a
formulation of a volume conjecture.

\smallskip

\noindent This goes roughly as follows. One constructs a `Bloch-like'
group ${\cal P}({\cal D})$ based on $\cal D$-decorated tetrahedra, which
maps via an explicit {\it idealization} map $\widehat{F}$ onto an
enriched version ${\cal P}({\cal I})$ of the classical Bloch group,
built on hyperbolic ideal tetrahedra. Any $\cal D$-triangulation
${\cal T}$ of $(W,L,\rho)$ leads to elements
${\frak c}_{\cal D}(W,L,\rho)\in {\cal P} ({\cal D})$ and
${\frak c}_{\cal I}(W,L,\rho) =
\widehat{F}({\frak c}_{\cal D}(W,L,\rho)) \in {\cal P} ({\cal
I})$. They are respectively called the $\cal D$- and $\cal I$-scissors
congruence classes of $(W,L,\rho)$. The QHI essentially depend on the
$\cal D$-class, and for any given $\cal D$-triangulation $\cal T$ the
$6j$-symbols occurring in $H_N({\cal T})$ depend on the moduli of the
hyperbolic tetrahedra of the idealization $\widehat{F}({\cal T})$ of
$\cal T$. By using the classical Rogers dilogarithm one can also define a
{\it dilogarithmic invariant}\index{dilogarithm!dilogarithmic invariant} 
$R(W,L,\rho)$ which only depends on
the $\cal I$-class.

\begin{conj}[S. Baseilhac, R. Benedetti] (Real Volume Conjecture for QHI)\qua
\label{benedetti.conj}
\noindent For any triple $(W,L,\rho)$ one 
has:\index{volume conjecture!--- for quantum hyperbolic invariant}
$$ \lim_{N\to \infty} \ (2 \pi/N^2) \log
(\vert K_N(W,L,\rho) \vert)= {\rm Im} \ R \bigl({\frak c}_{\cal I}(W,L,\rho) \bigr)
\ .$$
\end{conj}

\begin{rem}
From the explicit formula of $H_N({\cal T})$ one easily shows that the
left-hand side of Conjecture \ref{benedetti.conj}, if it exists, only
depends on the moduli of the hyperbolic tetrahedra of
$\widehat{F}({\cal T})$. A natural problem is to find a geometric
interpretation of the dilogarithmic invariant. Indeed, for scissors
congruence classes built with ideal triangulations of genuine
(non-compact finite volume) hyperbolic $3$-manifolds $M$, a similar
dilogarithmic invariant gives $i ({\rm Vol}(M) +i{\rm CS}(M))$, where Vol
is the Volume and CS is the Chern-Simons invariant (see \cite{bene.N}).
\end{rem}

\noindent In \cite{bene.B-B4} one proposes a complex version of
Conjecture \ref{benedetti.conj}, for the whole $K_N$ (not only its
modulus).

\subsection{Perturbative invariants}

The {\it perturbative $SO(3)$ invariant}\index{perturbative invariant} 
(or the {\it Ohtsuki series})\index{Ohtsuki series}
$\tau^{SO(3)}(M) = \sum_{d=0}^\infty \lambda_d$ $(q-1)^d \in \Q[[q-1]]$ 
of a rational homology 3-sphere $M$
is the invariant characterized by the property that
$\sum_{d=0}^k \lambda_d (e^{2\pi\sqrt{-1}/r}-1)^d$ for any $k$
is congruent to
$\left( \frac{|H_1(M;\Z)|}r \right) \tau_r^{SO(3)}(M)$
modulo $r$
for infinitely many primes $r$;
for a detailed definition see \cite{Ohtsuki_rat,Ohtsuki_book}.
(It is known, see \cite{Rozansky_pWRT,Habiro_expand}, that
$\tau^{SO(3)}(M) \in \Z[[q-1]]$ for any integral homology 3-sphere $M$.)
The perturbative $PG$ invariant $\tau^{PG}(M)$
of a rational homology 3-sphere $M$,
say, for $G=SU(N)$, 
is defined in $\Q[[q-1]]$ similarly,
related to the quantum invariant $\tau_r^{PG}(M)$;
see \cite{Le_PSUn,Le_isp}.

\begin{prob}
\label{prob.cal_tauSO3}
For each rational homology 3-sphere $M$,
calculate\index{perturbative invariant!calculation of ---} 
$\tau^{SO(3)}(M)$ and $\tau^{PSU(N)}(M)$ for all degrees.
\end{prob}

\begin{rem}
The value of $\tau_r^{SO(3)}\big( L(a,b) \big)$
of the lens space $L(a,b)$\index{lens space}
is concretely calculated in \cite{Jeffrey,Garoufalidis_rel}.
It follows from those values that
$$
\tau^{SO(3)}\big( L(a,b) \big)
= q^{- 3 s(b,a)} \frac{ q^{1/2a} - q^{-1/2a} }{q^{1/2}-q^{-1/2}},
$$
where we regard it as in $\Q[[q-1]]$
and $s(b,a)$ denotes the Dedekind sum.

Concrete presentations of $\tau^{SO(3)}(M)$
for Seifert fibered 3-manifolds $M$
are given in \cite{LaRo99}.

Lawrence \cite{Lawrence_h}
has given holomorphic expression for
the perturbative $SO(3)$ invariants of rational homology 3-spheres 
obtained by integral surgery along $(2,n)$ torus knot.

Habiro's expansion (\ref{eq.Habiro_expansion})
gives a presentation of $\tau^{SO(3)}(M)$.
See examples of Problem \ref{prob.ch_He},
for presentations of
$\tau^{SO(3)}\big( \Sigma(2,3,5) \big)$
and $\tau^{SO(3)}\big( \Sigma(2,3,7) \big)$,
which are due to \cite{Le_isp}.
See also \cite{LaZa}
for a computation of $\tau^{SO(3)}\big(\Sigma(2,3,5)\big)$.
\end{rem}


\begin{rem}
From the value of $\tau_r^{PSU(N)}\big( L(a,b) \big)$
of the lens space\index{lens space} $L(a,b)$ calculated in \cite{Takata_lens},
we obtain
$$
\tau^{PSU(N)}\big( L(a,b) \big)
= q^{ - N(N^2-1) s(b,a) /2}
\frac{[1/a]^{N-1}[2/a]^{N-2}\cdots[(N-1)/a]}{[1]^{N-1} [2]^{N-2} \cdots [N-1]},
$$
where we regard it as in $\Q[[q-1]]$
putting $[\alpha] = (q^{\alpha/2}-q^{-\alpha/2})/(q^{1/2}-q^{-1/2})$.

Takata \cite{Takata_S} computed the quantum $PSU(N)$ invariant
of Seifert fibered manifolds $M$.
Concrete presentations of $\tau^{PSU(N)}(M)$
might follow from the computation.
\end{rem}

\begin{rem}
$\tau^{PSU(N)}(M)$ is recovered from the LMO invariant by
$$
\tau^{PSU(N)}(M) =
| H_1(M;\Z) |^{- n(n-1)/2} \hat{W}_{sl_n} \big(\hLMO(M) \big).
$$
In particular, noting $PSU(2) = SO(3)$, 
$$
\tau^{SO(3)}(M)
= | H_1(M;\Z) |^{-1}  \hat{W}_{sl_2} \big(\hLMO(M) \big).
$$
For details see \cite{Ohtsuki_book}.
In this sense Problem \ref{prob.cal_tauSO3} is related to
Problem \ref{prob.calLMO}.
\end{rem}

\begin{prob}[J. Roberts]  
\label{prob.roberts14}
Explain\index{perturbative invariant!modular form in  ---}  
the appearance of modular forms in the Witten invariants. 
\end{prob}

\begin{rem}[{\rm (J. Roberts)}]  
Lawrence and Zagier discovered in \cite{LaZa} that the
perturbative series for the Poincar\'e homology sphere was close to a
modular form. Is this a random coincidence, or is there a more
systematic explanation? Does such a relation ever hold for a
{\it hyperbolic} $3$-manifold?
\end{rem}

\begin{prob}
\label{prob.ch_tauSO3}
Characterize\index{perturbative invariant!image of ---}  
those elements of $\Z[[q-1]]$ of the form $\tau^{SO(3)}(M)$ 
of integral homology 3-spheres $M$.
\end{prob}

\begin{rem}
The degree $\le d$ part of $\tau^{SO(3)}(M)$
can have any value in the degree $\le d$ part of $\Z[[q-1]]$.
Hence, it is meaningful to consider this problem
for the form $\tau^{SO(3)}(M)$ for all degrees.
\end{rem}

\begin{rem}
Problem \ref{prob.ch_tauSO3} is related to Problem \ref{prob.ch_He},
which is on the characterization of
Habiro's expansion (\ref{eq.Habiro_expansion}).
See examples of Problem \ref{prob.ch_He},
for calculations of Habiro's expansions of
$\tau^{SO(3)}\big( \Sigma(2,3,5) \big)$
and $\tau^{SO(3)}\big( \Sigma(2,3,7) \big)$.
\end{rem}


Let $q$ be an indeterminate, and let $\zeta$ be an $r$-th root of unity.
Set
$$
R_1 = \varprojlim_n \Z[q,q^{-1}]/ \big( (q-1)(q^2-1)\cdots(q^n-1) \big).
$$
For an integral homology 3-sphere $M$,
relations between $\tau_r^{SU(2)}(M)$
(which equals $\tau_r^{SO(3)}(M)$ for odd $r$, in this case)
and $\tau^{SO(3)}(M)$
can be described in the following commutative diagram.
\begin{equation}
\label{eq.CDforR1}
\begin{CD}
\quad \qquad\qquad\qquad I^{sl_2}(M)\ \in\ @. R_1  @>{\mboxss{injection}}>> \Z[[q-1]] 
@. \ \subset \ \Q[[q-1]] \ \ni\ \tau^{SO(3)}(M) \\
@. @V{\mboxss{put $q=\zeta$}}VV 
@VV{\mboxss{put $q=\zeta$}}V @. \\
\tau_r^{SU(2)}(M) = \tau_r^{SO(3)}(M) \ \in\ @. \Z[\zeta]
@>{\mboxss{injection}}>> \Z_r[\zeta]
\end{CD}
\end{equation}
Here,
the two horizontal maps are defined to be natural injections,
and the two vertical maps are defined by substituting $q=\zeta$.

It was conjectured by Lawrence \cite{Lawrence_a},
and proved by Rozansky \cite{Rozansky_pWRT}, 
that
$\tau^{SO(3)}(M) \in \Z[[q-1]]$
for any integral homology 3-sphere $M$, 
and that the images of $\tau^{SO(3)}(M)$ and $\tau_r^{SO(3)}(M)$
coincide in $\Z_r[\zeta]$ in the above diagram
for any odd prime power $r$.
See \cite{Rozansky_p_adic} for their numerical examples.

Habiro \cite{Habiro_expand} showed\footnote{\footnotesize 
Hence, $\tau^{SO(3)}(M)$ is as powerful as
the set of $\tau_r^{SU(2)}(M)$ for any integer $r \ge 3$,
and as powerful as the set of $\tau_r^{SO(3)}(M)$ for any odd $r \ge 3$,
for any integral homology 3-sphere $M$.
Further, the LMO invariant dominates
$\tau_r^{SU(2)}(M)$ for any integer $r \ge 3$.}
that there exists an $R_1$-valued invariant $I^{sl_2}(M)$
of an integral homology 3-sphere $M$
whose images in $\Q[[q-1]]$ and $\Z[\zeta]$ in the above diagram
are equal to $\tau^{SO(3)}(M)$ and $\tau_r^{SU(2)}(M)$ respectively
for any positive integer $r$.
(Here we set $\tau_r^{SU(2)}(M)=1$ for $r=1,2$.)
This gives another proof of the above mentioned conjecture of Lawrence
for integral homology 3-spheres.
This also implies that $\tau^{SO(3)}(M)$ can be presented by
\begin{equation}
\label{eq.Habiro_expansion}
\tau^{SO(3)}(M) = \sum_{n=0}^\infty \lambda'_n
(q-1)(q^2-1)\cdots(q^n-1)
\end{equation}
with some $\lambda'_n \in \Z[q,q^{-1}]$ (in the above sense)
such that
$$
\tau_r^{SU(2)}(M) = \sum_{0\le n < r} \lambda'_n
(\zeta-1)(\zeta^2-1)\cdots(\zeta^n-1).
$$
Note that the presentation (\ref{eq.Habiro_expansion}) is not unique.

\medskip

\noindent
{{\namae}(K. Habiro)}\quad
Let $\frak g$ be a finite dimensional simple complex Lie algebra.  
Let $d\in\{1,2,3\}$ be such that $d=1$ in the $ADE$ cases, 
$d=2$ in the $BCF$ cases and $d=3$ in the $G_2$ case.  
If $M$ is a closed 3-manifold and 
if $\zeta $ is a root of unity of order $r$ divisible by $d$, 
then the quantum $\frak g$ invariant\index{quantum invariant!--- of 3-manifold}  
$\tau_\zeta^{\frak g}(M)\in\Q[\zeta ]$ of $M$ at $\zeta $ is defined.

\begin{conj}[K. Habiro, T. Le] 
\label{conj.HL_IgM}
For each $\frak g$ as above, 
there is a (unique) invariant $I^{\frak g}(M)\in R_1$ 
of an integral homology 3-sphere $M$ such that 
for each root of unity $\zeta$ of order $r$ divisible by $d$ 
we have\index{perturbative invariant!--- and quantum invariant} 
$$
I^{\frak g}(M) \big|_{q=\zeta} = \tau^{\frak g}_\zeta (M).
$$
\end{conj}

\begin{rem}
When $(r, \det(a_{ij}))=1$, 
where $(a_{ij})$ is the Cartan matrix of the Lie algebra $\frak g$, 
the projective $\frak g$-invariant $\tau_\zeta^{P{\frak g}}(M)$
can be defined \cite{Le_isp}. 
Then Habiro and Le also conjecture that 
$I^{\frak g}(M)|_{q=\zeta}=\tau_\zeta^{P{\frak g}}(M)$, 
if $(r,\det(a_{ij}))=1$.  
Note that for an integral homology 3-sphere, 
$\tau_\zeta^{P{\frak g}}(M)=\tau_\zeta^{\frak g}(M)$ 
when both are defined 
(i.e.\ when $r$ is divisible by $d$ and $(r,\det(a_{ij}))=1$). 
If this is the case, then we have
$$
i \big(I^{\frak g}(M) \big) = \tau^{\frak g}(M) 
$$
where $\tau^{\frak g}(M)\in\Q[[q-1]]$ is 
the perturbative ${\frak g}$ invariant of $M$ \cite{Le_isp}, 
and $i : R_1\to \Z[[q-1]]$ is the upper injection in (\ref{eq.CDforR1}). 
\end{rem}

\begin{rem}
The above conjecture implies that the quantum ${\frak g}$ invariant
$\tau_\zeta^{{\frak g}}(M)$ of an integral homology sphere $M$ takes values
in the ring of cyclotomic integers $\Z[\zeta]$, and also that 
the perturbative invariant $\tau^{\frak g}(M)$ takes values in $\Z[[q-1]]$.
\end{rem}

\begin{update}
Habiro and Le \cite{HabiroLe_Ig} proved Conjecture \ref{conj.HL_IgM}.
\end{update}

\begin{conj}[K. Habiro] 
\label{conj.IslM}
Suppose that Conjecture \ref{conj.HL_IgM} would hold.  
For a new indeterminate $t$, set
$$
  R_1' = \varprojlim_{n}R_1[t]/((t-q)(t-q^2)\cdots(t-q^n))
$$
Then there exists an invariant $I^{sl}(M)\in R_1'$ of an integral
homology $3$-sphere $M$ such that $I^{sl}(M)|_{t=q^n}=I^{sl_n}(M)$
for any $n\ge1$, where we set $I^{sl_1}(M)=1$.
\end{conj}

\begin{prob}
\label{prob.ch_He}
Characterize\index{perturbative invariant!image of ---}  
those elements of Habiro's expansion (\ref{eq.Habiro_expansion}) 
of \newline
$\tau^{SO(3)}(M)$ of integral homology 3-spheres $M$.
\end{prob}

\begin{exm}
For the Poincare homology 3-sphere $\Sigma(2,3,5)$
(obtained by surgery on a left-hand trefoil with framing $-1$)
and the Brieskorn sphere $\Sigma(2,3,7)$
(obtained by surgery on a right-hand trefoil with framing $-1$),
it is computed in \cite{Le_isp} that
\begin{align*}
&\tau^{SO(3)}\big(\Sigma(2,3,5)\big) = \frac{1}{1-q}
\sum_{n=0}^\infty q^n (1-q^{n+1})(1-q^{n+2}) \cdots (1-q^{2n+1}), \\
&\tau^{SO(3)}\big( \Sigma(2,3,7) \big) = \frac{1}{1-q}
\sum_{n=0}^\infty q^{-n(n+2)} (1-q^{n+1})(1-q^{n+2}) \cdots (1-q^{2n+1}).
\end{align*}
See also \cite{LaZa}
for a computation of $\tau^{SO(3)}\big(\Sigma(2,3,5)\big)$.
\end{exm}

\begin{rem}
Such an infinite sum as (\ref{eq.Habiro_expansion})
would be interesting from the number theoretical viewpoint.
For example, 
$$
1 + \sum_{n=1}^\infty q^n(q-1)(q^2-1) \cdots (q^n-1)
= \sum_{\substack{k \in \Z \\ k \ne 0}} (-1)^{k+1} 
  q^{\frac32 k^2 - \frac12 k -1}.
$$
A similar infinite sum appears in (\ref{eq.Zag}); see also \cite{Sikora}.
\end{rem}

\newpage

\section{Topological quantum field theory}
\label{sec.TQFT}

\renewcommand{\thefootnote}{\fnsymbol{footnote}}
\footnotetext[0]{The first version of 
the introductory part of Chapter \ref{sec.TQFT} 
and Sections \ref{sec.TQFT.1}--\ref{sec.TQFT.4}
was written by T. Ohtsuki,
following seminar talks given by G. Masbaum.
Based on it, G. Masbaum wrote this introductory part and these sections.
Section \ref{sec.Kerler} was written by T. Kerler.}
\renewcommand{\thefootnote}{\arabic{footnote}}

The notion of topological quantum field theory (TQFT) was introduced
in \cite{Atiyah_t,AHLS},
motivated by the operator formalism
of a partition function in a quantum field 
       theory
which does not depend on the metric of the space.
In the mathematical viewpoint,
any quantum invariant\index{quantum invariant!--- of 3-manifold} 
 of 3-manifolds
can be formulated by a TQFT,
which 
enables
us to compute the invariant by the cut-and-paste method.

A {\it TQFT}\index{TQFT} 
is a functor which 
takes an oriented closed surface $\Sigma$
to a finite dimensional complex vector space $V(\Sigma)$,
and takes an oriented compact 3-manifold $M$ with boundary $\Sigma$
to a vector $Z(M) \in V(\Sigma)$, 
satisfying the following 5 axioms.
\begin{itemize}
\item[\rm(1)]
$V(-\Sigma) = V(\Sigma)^\star$,
where $-\Sigma$ denotes $\Sigma$ with the opposite orientation
and $V(\Sigma)^\star$ denotes the dual vector space of $V(\Sigma)$.
\item[\rm(2)]
$V(\Sigma_1 \sqcup \Sigma_2) = V(\Sigma_1) \otimes V(\Sigma_2)$,
where $\Sigma_1 \sqcup \Sigma_2$ denotes
the disjoint union of two surfaces $\Sigma_1$ and $\Sigma_2$.
\item[\rm(3)]
$V(\emptyset) = \C$,
where $\emptyset$ denotes the empty surface.
\item[\rm(4)]
For 3-cobordisms $M_1$ and $M_2$
with $\partial M_1 = (-\Sigma_1) \sqcup \Sigma_2$
and $\partial M_2 = (-\Sigma_2) \sqcup \Sigma_3$
we have that
$Z(M_1 \underset{\Sigma_2}\cup M_2) = 
       Z(M_2) \circ Z(M_1)$
as linear maps\footnote{\footnotesize 
For a 3-cobordism $M$ with $\partial M = (-\Sigma_1) \sqcup \Sigma_2$
the vector $Z(M)$ belongs to
$V(-\Sigma_1 \sqcup \Sigma_2)
= V(-\Sigma_1) \otimes V(\Sigma_2)
= V(\Sigma_1)^\star \otimes V(\Sigma_2)$
by the axioms (1) and (2).
Hence, $Z(M)$ can be regarded as a linear map
$V(\Sigma_1) \to V(\Sigma_2)$.}
$V(\Sigma_1) \to V(\Sigma_3)$.
\item[\rm(5)]
$Z(\Sigma \times I)$ is equal to the identity map of $V(\Sigma)$.
\end{itemize}


To be precise,
in many (but not in all) examples
we need 
``extended
3-manifolds'' 
instead of 3-manifolds to formulate a TQFT,
where 
an {\it extended 
3-manifold} is a 3-manifold $M$ equipped with some kind of
framing, e.g.\ a $p_1$-structure $\alpha$ on $M$
(see \cite{BHMV95})\footnote{\footnotesize 
There is another formulation of a ``framing'' of a 3-manifold
using signature cocycle; see \cite{Turaev_book}.}.
Namely, we extend the above definition of TQFT
to a functor from the category of 
extended 
3-cobordisms
in an appropriate way (see \cite{BHMV95}).
Then, each quantum invariant\index{quantum invariant!--- of 3-manifold} 
 can be formulated
as a TQFT of the category of 
extended
3-cobordisms.
In the remaining part of this section
we call such a TQFT simply a TQFT.

\subsection{Classification and characterization of TQFT's}
\label{sec.TQFT.1}

To understand TQFT's is an important problem
in order to investigate the 3-cobordism category,
similarly as the representation theory is important
in order to investigate groups and algebras.

\begin{prob}
\label{prob.classify_TQFT}
Find (and classify) all TQFT's.\index{TQFT!classification of ---}
\end{prob}

\begin{rem}
The operator formalism of the Chern-Simons path integral
suggests the existence of many TQFT's.
It is known, see \cite{Turaev_book,BaKi01},
that a modular category
is derived from a quantum group
 at a root of unity
and a TQFT is derived from a modular category.
The underlying 3-manifold invariant is called 
the {\it Reshetikhin-Turaev invariant}.
Some other TQFT's might be obtained from 
quantum groupoids \cite{NTV00}.\index{quantum groupoid!--- and TQFT}
A TQFT for the LMO invariant is discussed in \cite{MuOh97}.

Another major construction of TQFT's
is derived from sets of $6j$-symbols;
for the construction see \cite{TV,BarrWest96}.   
When a set of $6j$-symbols arises from a subfactor,
the underlying 3-manifold invariant is called 
the {\it Turaev-Viro-Ocneanu invariant}
(see Section \ref{sec.TVOinv}). 
Further, when a set of $6j$-symbols
comes from a quantum group,
such a TQFT is isomorphic to a tensor product of 
two TQFT's derived from the quantum group
\cite{Turaev_book}.
See Problems in Chapter \ref{sec.state-sum-inv}
for concrete problems for TQFT's derived from $6j$-symbols.

There are TQFT's derived from finite groups, 
whose invariants are called the {\it Dijkgraaf-Witten invariants} 
\cite{DiWi90}.
Such TQFT's can alternatively be formulated
by using certain sets of $6j$-symbols.

It is known \cite{Atiyah_g} that
the vector space $V(\Sigma)$ of a TQFT $(V,Z)$
derived from a quantum group
is isomorphic to the space of conformal blocks
of a conformal field theory (CFT) of the Wess-Zumino-Witten model.
Some other (possibly, ``new'') TQFT's might be obtained from
the orbifold construction of CFT. 
It is a problem to understand TQFT's derived from
the Rozansky-Witten invariant (see \cite{rw.RSW}); 
their isomorphism types might be described by known TQFT's, 
or they might be ``new'' TQFT's.
\end{rem}

%
%

The following problem is a part of Problem \ref{prob.classify_TQFT}
in the sense that
some TQFT's are derived from modular categories,
as mentioned in a remark after Problem \ref{prob.classify_TQFT}.

\begin{prob}
\label{prob.classify_modu_cat}
Find (and classify) all 
modular categories.\index{tensor category!modular ---!classification of ---} 
\end{prob}

For a TQFT $(V,Z)$, put
$P_{(V,Z)}(t) = \sum_{g=0}^\infty \big( \mbox{dim} V(\Sigma_g) \big) t^g$,
where $\Sigma_g$ denotes a closed surface of genus $g$.
The following problem is a refinement of Problem \ref{prob.classify_TQFT}.

\begin{prob}
\mbox{}\newline
\noindent
(1)\qua
Characterize the power series of the form $P_{(V,Z)}(t)$.\newline
\noindent
(2)\qua
For each power series $P(t)$
(satisfying the characterization of (1)),\index{TQFT!classification of ---}
classify all TQFT's $(V,Z)$ such that $P_{(V,Z)}(t) = P(t)$.
\end{prob}

\begin{rem}
A concrete form of such a power series for a TQFT derived from a quantum group
is given by Verlinde formula \cite{Verlinde}.\index{Verlinde formula}
For example, such a power series of the TQFT derived from
$U_q(sl_2)$ at level $k$ is presented by

$$
\sum_{g=0}^\infty t^g \Big( \frac{k+2}2 \Big)^{g-1} 
\sum_{j=1}^{k+1} \Big( \sin \frac{\pi j}{k+2} \Big)^{2-2g}.
$$
\end{rem}

%
%

\subsection{Spin TQFT's}

There are some refinements of TQFT's.

A {\it spin TQFT} is a TQFT 
on
the category of spin 3-cobordisms,
whose invariants depend on spin structures;
such a TQFT can be formulated
by extending the definition of a usual TQFT (see \cite{BlMa96}).
It is shown \cite{BlMa96} that
a spin TQFT can be obtained from the modular category
of $U_q(sl_2)$ at level $k \equiv 2$ (mod 4).

\begin{prob}
Find other spin TQFT's.\index{TQFT!spin ---}
\end{prob}

\begin{rem} Some examples of spin TQFT's can be constructed from 
the refined quantum invariants
 of \cite[Theorem 6.2]{BeBl01}.
\end{rem}

\begin{rem}

A spin TQFT is expected to be a refinement of a usual TQFT
in the sense that
a spin TQFT $(V^s,Z^s)$ should be related to a usual TQFT $(V,Z)$
such that
$V(\Sigma)$ 
for connected $\Sigma$
can be described by
the direct sum of $V^s(\Sigma,\sigma_\Sigma)$ 
over the spin structures $\sigma_\Sigma$ on $\Sigma$
(see \cite{BlMa96})
and $Z(M)$ of 
a closed manifold $M$
can be described by the sum of $Z^s(M,\sigma_M)$
over the spin structures $\sigma_M$ on $M$.
\end{rem}

A spin${}^c$ TQFT should be a TQFT
on the category of spin${}^c$ 3-cobordisms,
whose invariants depend on spin${}^c$ structures.

\begin{prob}
Formulate and find spin${}^c$ TQFT's.\index{TQFT!spin${}^c$ ---}
\end{prob}

\begin{rem}
The Seiberg-Witten invariant 
(for its exposition see, e.g.\ \cite{Marcolli})
and the torsion invariant $\tau$ (see \cite{Turaev_SW})
are defined for closed 3-manifolds
with spin${}^c$ structures.
Are there TQFT's which are related to these invariants?
\end{rem}

\subsection{Homotopy QFT's}

V. Turaev \cite{Turaev_hft2,Turaev_hft3} introduced and developed
{\it HQFT} (homotopy QFT)\index{homotopy QFT} 
with a target space $X$ in dimension $d+1$.

\begin{prob}[V. Turaev] 
\label{prob.HQFT_spin}
\begin{itemize}
\item[\rm(1)]
Extend HQFT's to spin and spin${}^c$ settings.\index{homotopy QFT!spin ---}\index{homotopy QFT!spin${}^c$ ---}
\item[\rm(2)]
Find algebra structures behind spin and spin${}^c$ HQFT's in dimension $1+1$.
\end{itemize}
\end{prob}

\begin{prob}[V. Turaev] 
\label{prob.HQFT_spin2}
Study (spin and spin${}^c$) HQFT's 
with the target space $K(H,2)$ in dimensions $1+1, 2+1$, and $3+1$  
for $H = \Z^N$.
\end{prob}

\begin{rem}
It is shown by V. Turaev that
HQFT's with the target space $K(\pi,1)$ in dimension $1+1$
can be described by crossed $\pi$-algebras,
and that any modular $G$-category gives rise to a HQFT
with the target space $K(G,1)$ in dimension $2+1$ \cite{Turaev_hft3}.
HQFT's with the target space $K(H,2)$ in dimension $1+1$ 
were studied and classified by M. Brightwell and P. Turner \cite{BrTu00}.
\end{rem}


\subsection{Geometric construction of TQFT's}
\label{sec.TQFT.4}

Assume that the surface $\Sigma$ is equipped with the structure
of a smooth algebraic curve over $\C$.
We denote by
$H^0({\cal M}_\Sigma, {\cal L}^{\otimes k})$
the space of 
sections of 
${\cal L}^{\otimes k}$ on ${\cal M}_\Sigma$,
where ${\cal M}_\Sigma$ is the moduli space of
semi-stable
rank $N$ bundles
with trivial determinant over $\Sigma$,
and $\cal L$ is the determinant line bundle on ${\cal M}_\Sigma$.
It is known 
that
$H^0({\cal M}_\Sigma, {\cal L}^{\otimes k})$
is 
isomorphic to 
$V(\Sigma)$ of a TQFT $(V,Z)$ derived from 
the quantum group $U_q(sl_N)$ at a $(k+N)$-th root of unity.
In this sense, $H^0({\cal M}_\Sigma, {\cal L}^{\otimes k})$
gives a geometric construction of such a $V(\Sigma)$.

\begin{prob}
Find a geometric construction of 
a TQFT\index{TQFT!geometric construction of ---}
using $H^0({\cal M}_\Sigma, {\cal L}^{\otimes k})$.
Namely, find a geometric way to associate a vector in 
$H^0({\cal M}_\Sigma, {\cal L}^{\otimes k})$
to a 3-manifold $M$ with $\partial M = \Sigma$.
\end{prob}

\begin{rem}
In physics such a vector
is obtained by applying
an infinite dimensional formal analogue of
the geometric invariant theory and the symplectic quotient
to the Chern-Simons path integral; see \cite{Atiyah_g}.
It is a problem to justify this argument
in some mathematical sense.
\end{rem}

Here is a concrete problem which may be of interest in studying 
the relationship 
between $V(\Sigma)$ and $H^0({\cal M}_\Sigma, {\cal L}^{\otimes k})$. 
The group $J^{(N)}(\Sigma)$ of $N$-torsion points  on the Jacobian $J(\Sigma)$
acts on ${\cal M}_\Sigma$ by tensoring. This gives an action of 
a central extension ${\mathcal E}(\Sigma)$ of $J^{(N)}(\Sigma)$ 
on $H^0({\cal M}_\Sigma, {\cal L}^{\otimes k})$.

\begin{prob}[G. Masbaum] 
\label{prob.EGonHMSLk}
Study this action of  the finite group ${\mathcal E}(\Sigma)$ on 
$H^0({\cal M}_\Sigma, {\cal L}^{\otimes k})$, and describe the induced 
decompositions of this vector space according to the characters of  ${\mathcal E}(\Sigma)$. Also  relate these decompositions 
to decompositions of $V(\Sigma)$ for the TQFT $(V,Z)$ 
derived from\index{TQFT!decomposition of ---} 
the quantum group $U_q(sl_2)$ at a $(k+N)$-th root of unity.
\end{prob}

\begin{rem} This was done for $N=2$ in \cite{AnMa99}. 
\end{rem}
\begin{rem} The group $J^{(N)}(\Sigma)$ is isomorphic to $H^1(\Sigma;\Z/N)$ 
and the extension ${\mathcal E}(\Sigma)$ is described using the Weil 
pairing, which corresponds to the intersection form on $H_1(\Sigma;\Z/N)$. 
For $N=2$, an action of ${\mathcal E}(\Sigma)$ on the vector space 
$V(\Sigma)$ is described 
in \cite[Section 7]{BHMV95}, and it was shown in \cite{AnMa99} that 
$V(\Sigma)$ and $H^0({\cal M}_\Sigma, {\cal L}^{\otimes k})$ are 
isomorphic as representations of ${\mathcal E}(\Sigma)$; here the 
torsion points on the Jacobian $J(\Sigma)$  correspond to simple closed
curves on the surface $\Sigma$. For example, if $k\equiv 2$ mod $4$, one 
obtains decompositions indexed by spin structures (theta-characteristics) 
on $\Sigma$. For $N\geq 3$, the action of ${\mathcal E}(\Sigma)$ and the 
spin decompositions of $V(\Sigma)$ were constructed in \cite{Bl01}.
\end{rem}

Let ${\frak M}_g$ denote the mapping class group 
of a closed surface $\Sigma_g$ of genus $g$,
and let $\tilde{\frak M}_g$ denote its central extension 
(see \cite{Atiyah_f,MaRo95})
arising in the category of 
extended 
3-cobordisms.

\begin{prob}
For a given TQFT $(V,Z)$,\index{TQFT!image of mapping class in ---}
determine whether the image of $\tilde {\frak M}_g$ 
in $\mbox{End}\big( V(\Sigma_g) \big)$ is finite.
\end{prob}

\begin{rem}
Using physical arguments, Bantay \cite{Ban01} (see also 
references therein)  showed that
for every CFT the image of $\tilde {\frak M}_1$ in
$\mbox{End}\big( V(S^1\times S^1) \big)$ is finite.
This had been rigorously proved by Gilmer \cite{Gilmer99} for 
the $SU(2)$ case.

In higher genus, it is known \cite{Funar99,Masbaum99}  that
the image of ${\frak M}_g$ ($g \ge 2$) is infinite in general.
\end{rem}

\begin{prob}[G. Masbaum] 
\label{prob.NT_TQFT}
Is there a relation between the Nielsen-Thurs\-ton classification
of mapping classes of $\Sigma_g$\index{TQFT!image of mapping class in ---}
and their images on $V(\Sigma_g)$ for TQFT's $(V,Z)$?
\end{prob}

\begin{rem}
The Nielsen-Thurston classification says that
any mapping class of a surface
is either finite order, reducible, or pseudo-Anosov 
(see, e.g.\ \cite{CaBl88}).
It is known that
a Dehn twist is taken to a matrix of finite order
by any TQFT derived from a modular category
of a quantum group.
On the other hand, it is shown
in  \cite{Masbaum99} that a certain product of two non-commuting
Dehn twists is taken to a matrix of infinite order in the $SU(2)$ TQFT 
at level $k$ unless $k = 1,2,4,8$.
\end{rem}



\subsection{Half-projective and homological TQFT's}
\label{sec.Kerler}


In \cite{Gilmer01} it is shown that, 
for a restricted set of cobordisms, the Reshetikhin-Turaev TQFT
at a prime $p$-th root of unity $\zeta_p$ can be defined, at least
abstractly,  as a functor 
${\cal V}_p:\ThKeCob\to \Z[\zeta_p]\mbox{-mod}$, meaning the category of
{\em free} $\Z[\zeta_p]$-modules. 
Note that there is a well defined ring epimorphism 
$\Z[\zeta_p]\ThKeonto{10}\, {\FF}_p[\mbox{\sf y}]/\mbox{\sf y}^{p-1}$, 
which sends $\zeta_p\mapsto 1+\mbox{\sf y}$ and maps integer
coefficients canonically onto the finite field $\FF_p=\Z/p\Z$. 
Thus an endomorphism, which 
for a choice of basis of the free $\Z[\zeta_p]$-modules is given by a
matrix with entries in $\Z[\zeta_p]$, will be represented by the same matrix 
with reduced coefficients now in  ${\FF}_p[\mbox{\sf y}]/\mbox{\sf y}^{p-1}$.
Collecting the coefficients for each degree we can thus reexpress such a 
matrix as a sum of matrices over $\FF_p$ multiplied with powers of $\mbox{\sf y}$, or,
more succinctly, use 
${\rm Mat}({\FF}_p[\mbox{\sf y}]/\mbox{\sf y}^{p-1})=
{\rm Mat}({\FF}_p)[\mbox{\sf y}]/\mbox{\sf y}^{p-1}$. 
This means that in the  ring-reduction the TQFT assigns to 
cobordisms a polynomial $\overline{\cal V}_p(M)=
\sum_{j=0}^{p-2}\mbox{\sf y}^j\cdot {\cal V}_p^{[j]}(M)$, where each 
${\cal V}_p^{[j]}(M)$ is a matrix over $\FF_p$ and 
is well defined for given bases.

Recall also the notion of a half-projective
TQFT with respect to an element $\mbox{\sf x}\in\mbox{\sf R}$ in the base ring, 
introduced in \cite{Kerler98}. It is 
defined, by perturbing functoriality into 
${\cal V}(N\circ M)=\mbox{\sf x}^{\mu(M,N)}{\cal V}(N){\cal V}(M)$, where 
$\mu(M,N)={\rm rank}\bigl(H_1(N\circ M)\stackrel{\delta}{\to}H_0(N\cap M)
\bigr)$.

\begin{prob}[T. Kerler] 
\label{prob.kerler1}
{\rm [Cyclotomic integer TQFT's]}\index{TQFT!cyclotomic integer ---}
\begin{enumerate}
\item Find explicit/computable bases for the ${\cal V}_p(\Sigma_g)$ as free 
modules over $\Z[\zeta_p]$. 
\item Show that ${\cal V}_p$ can be extended to all cobordisms as a 
half-projective TQFT with $\mbox{\sf x}=(\zeta_p-1)^{\frac {p-3}2}\in \mbox{\sf R}
=\Z[\zeta_p]$. 
\item Determine the structure of the ${\cal V}_p^{[j]}(M)$ and in how far
they have lifts from $\FF_p$ to $\Z$, analogous to the Ohtsuki invariants
for closed 3-manifolds.
\item Find a universal TQFT that combines all ${\cal V}_p$, at least
perturbatively, into one. 
\end{enumerate}
\end{prob}

In the case of $p=5$ the program for   items (1)--(3) has been mostly
carried out  in \cite{Kerler02b}, for primes $p\geq 7$ not much is 
known though. Some explicit bases have been found for genus $g=1$ by Gilmer,
but the situation for higher genera $g\geq 3$ is unknown.
An immediate application
of item (2) is that the quantum order, as introduced in \cite{CoMe01},
is also an upper bound 
for the cut-number of a 3-manifold. A closely related statement for 
(2) would also yield a very different proof for the fact that the 
Ohtsuki invariants are of finite type. In item (3) the ``lift''
must depend on $p$ since the dimensions of the vector spaces do, 
and must also involve further quotients 
that arise since the irreducible TQFT's over $\Z$ 
do not match the required dimensions either, but they
become reducible when reduced to $\FF_p$. Item (4) is rather vague at this 
point, indicating for some sort of infinite filtered space with finite 
graded components. 
\medskip

Any TQFT ${\cal V}:\ThKeCob\to \mbox{\sf R}\mbox{-mod}$ implies a sequence of representation
${\cal V}_{[g]}:\Gamma_g\to 
{\rm GL}_{\mbox{\scriptsize\sf R}}({\cal V}(\Sigma_g))$ 
of the mapping class groups. We say
that a TQFT is {\em homological} if each of these representations factors
through the quotient $\Gamma_g\ThKeonto{15} {\rm Sp}(2g,\Z)$ (given by the action
on $H_1(\Sigma_g)$), and we say it is 
{\em strictly homological} if each of the ${\rm Sp}(2g,\Z)$-representations
is algebraic, i.e.\ either faithful or zero. A particular example of
strictly homological TQFT's over $\mbox{\sf R}=\Z$ are the Lefschetz 
components ${\cal V}^{(j)}$ of the Frohman-Nicas TQFT, 
see \cite{FrNi92,Ke00}.
{From} these we can generate a larger family $\ThKeQc^0$ of such TQFT's 
by taking all direct sums of  ${\cal V}^{(j)}$'s. For example all the TQFT's
constructed in \cite{Don99} lie in $\ThKeQc^0$. An even larger family 
$\ThKeQc^*$ is found by taking also tensor products and their irreducible 
summands.

\begin{prob}[T. Kerler]  
\label{prob.kerler2}
{\rm [Homological TQFT's]}\index{TQFT!homological ---}
\begin{enumerate}
\item Find the irreducible components and ring structure 
 (w.r.t $\oplus$ and
$\otimes$) of $\ThKeQc^*$. 
\item Determine whether all strictly homological TQFT's lie in $\ThKeQc^*$.
\item Identify the homological TQFT's that arise from the gauge theory of 
higher rank groups (such as $PSU(n)$ in \cite{FrNi94}) with elements in  
$\ThKeQc^*$. 
\item Identify the irreducible factors of the constant orders 
${\cal V}^{[0]}_p$ of the cyclotomic
integer expansion of the Reshetikhin-Turaev theory
with elements in  $\ThKeQc^*$.
\end{enumerate}
\end{prob}

The first item is in some sense about finding the representation ring
of ${\rm Sp}(2,\Z)\times {\rm Sp}(4,\Z)\times \ldots\times {\rm Sp}(2g,\Z)
\times \ldots$ equipped with further generators and relations given by the 
standard handle attachments.
The constraints given by the latter may be just good enough
to ensure that the answer to item (2) is positive. The application
of (3) is a better understanding and possibly a closed form for the
polynomials from \cite{FrNi94} that express the $PSU(n)$-invariants in terms
of the coefficients of the Alexander polynomial. Evidence seems to
suggest that the TQFT's from (4) stem from $\frac {p-3}2$-fold
symmetric products of elements in $\ThKeQc^0$. A plausible corollary would be 
that for a closed manifold with  $b_1(M)\geq 1$ we have 
\begin{equation}\label{eq-LesExp}
{\cal V}_p(M)\;=\;(\zeta_p-1)^{\frac {p-3}2}
P_{\frac {p-3}2}(\lambda_{CWL}(M))\,\;\;+\,\;\;{\cal O}((\zeta_p-1)^{\frac {p-1}2})\;,
\end{equation}
where $\lambda_{CWL}$ is the Casson-Walker-Lescop invariant, and $P_j$ is a
polynomial of degree $j$ with integer coefficients. (Note our normalization ${\cal V}_p(S^3)=1$). 
As remarked in  \cite{Kerler02} the identity in (\ref{eq-LesExp})
is true for $p=5$ and general $M$ with $b_1(M)\geq 1$. Moreover,
work in progress shows that (\ref{eq-LesExp}) holds also for general $p$ if $M$ is a
torus-bundle over a circle. 
\medskip

The homological TQFT's are the starting point for a more general,
perturbative view point on TQFT's that should parallel and extend that of the
finite type theory of homology-3-spheres. At least for fixed $p$ one can understand, for
example, the Reshetikhin-Turaev theory
as deformation of the 
$\ThKeQc^*$-theories. 
The notion that is somewhat parallel to that of finite type 
for closed 3-manifolds is what we shall call {\em finite length}.
More precisely, the representations 
${\cal V}_{[g]}:\Gamma_g\to 
{\rm GL}_{\mbox{\scriptsize\sf R}}({\cal V}(\Sigma_g))$ 
of the mapping class groups extend 
linearly to homomorphisms ${\cal V}_{[g]}:\Z[\Gamma_g]\to 
{\rm End}_{\mbox{\scriptsize\sf R}}({\cal V}(\Sigma_g))$. 
Denote by 
$I{\cal I}_g\subset\Z[\Gamma_g]$ the augmentation ideal of the Torelli
group. The {\em length} of ${\cal V}$ is the maximal $L\in\N$ such that
${\cal V}_{[g]}((I{\cal I}_g)^{L+1})=0$. Clearly, the $L=0$-theories are
just the homological ones. The $L=1$-theories can   be 
thought of as   elements of some ${\rm Ext}({\cal V},{\cal W})$ with
${\cal V},\,{\cal W}\in\ThKeQc^*$. Restricted to representations of the
$\Gamma_g$'s they factor (in $char\neq 2$)
through the Johnson-Morita-homomorphism
$\Gamma_g\to \ext 3H_1(\Sigma_g)\rtimes {\rm Sp}(2g,\Z)$, for which
such extension are explicitly constructible \cite{Kerler01b}.

\begin{prob}[T. Kerler]  
\label{prob.kerler3}
{\rm [Length = 1  TQFT's]}\index{TQFT!length 1 ---}
\begin{enumerate}
\item Describe and construct 
algebraic $L=1$-extensions of $\Gamma_g$-representations
to TQFT's, preferably as  ``simple''  generalizations 
of the Frohman-Nicas-$U(1)$-theory.
\item Produce a classification of $L=1$-TQFT's in the sense of an extension
theory of $\ThKeQc^*$. 
\item Identify the 
$\Gamma_g$-representations on relative $SU(2)$-moduli space 
from \cite{CLM00} with these TQFT's, and find similar, higher rank theories. 
\item Identify the ${\cal V}^{[0]}_p$ as $L=1$-theories, if possible. 
\end{enumerate}
\end{prob}

The conceivable generalizations of the TQFT construction of Frohman and Nicas
described in (1) include using different, possibly non-compact
 gauge groups  instead of $U(1)$ and using more refined versions of
intersection homologies for stratified moduli spaces. 
Given the theory for ${\ThKeQc^*}$ the solution to item (2) 
will lead to well defined problems in ${\frak s}{\frak p}$-invariant theory.  
Constructions of $L=1$-theories follow the schemes from (1) and (3). 
The identification in (4) is carried out for $p=5$ in \cite{Kerler01b}. 

The notion of {\em finite length} can be refined into the notion
of {\em $q/l$-solvable} introduced in \cite{Kerler02}, indicating a TQFT over 
$\mbox{\sf R}={\Bbb M}[\mbox{\sf y}]/{\mbox{\sf y}^{l+1}}$  such that
the constant order TQFT over the ground ring $\Bbb M$ is of length $q$. 
This, clearly,
defines a special case of a TQFT of length $\leq (q\cdot l+q+l)\,$. 
Murakami's result \cite{Hitoshi_CW} can be restated as saying that the Reshetikhin-Turaev theory
gives rise to a 1/1-solvable TQFT ${\cal V}_p^{[\leq 1]}$
with ground ring $\FF_p$ 
(i.e.\ a TQFT of length 3 over $\FF_p[\mbox{\sf y}]/\mbox{\sf y}^2$) such that
\begin{equation}\label{eq-MurId}
{\cal V}_p^{[\leq 1]}(M)=\,1\,+\,\mbox{\sf y}\frac 16 \lambda_{CWL}(M)
\end{equation}
for any closed homology sphere $M$.  Following Ohtsuki's work
Murakami's identity (with some extra renormalizations by the order of $H_1(M)$)
 extends also to rational homology spheres. Let us call a theory  with this property
a TQFT of {\em Casson type}.

Recall, that the similar
relation  (\ref{eq-LesExp}) 
for $\lambda_{CWL}$ for manifolds with $b_1(M)\geq 1$ 
is already contained in the information of a homological ($L=0$) TQFT, and is
indeed a special evaluation of the Turaev-Milnor Torsion, see \cite{Kerler02}.  Given
the richer structure of a  1/1-solvable TQFT we will expect new invariants $\Xi$ that are
refinements of  $\lambda_{CWL}$ and the torsion invariants. 

To be more precise, note that for a
pair $(M,\varphi)$, where $\varphi:\pi_1(M)\to\!\!\!\!\to\Z$ 
defines a cyclic cover, any  TQFT ${\cal V}$ yields an invariant   
${\cal V}(M,\varphi)=trace({\cal V}(C_{\Sigma}))$ where 
$C_{\Sigma}=\overline{M-\Sigma}:\Sigma\to\Sigma$ and $\Sigma\subset M$
is any surface dual to $\varphi$. In this way the Frohman Nicas
theories ${\cal V}^{(j)}$ yields the coefficients of the Alexander Polynomial,
and, as shown in  \cite{Kerler02}, thus also  $\lambda_{CWL}$. 

A more refined
invariant, which, roughly speaking, generalizes the Alexander module, 
is the Turaev-Viro module\index{Turaev-Viro!--- module} 
${\cal M}_{TV}(M,\varphi)$. It is
described by Gilmer in  \cite{Gilmer97}. 
${\cal M}_{TV}(M,\varphi)$ is given, up to conjugacy, 
by $\raise  .3ex\hbox{${\cal V}(\Sigma)$}\big/
\raise  -.3ex\hbox{ker$({\cal V}(C_{\Sigma})^N)$}$ (with $N$ large enough) 
together with the action of ${\cal V}(C_{\Sigma})$ on it. The traces of 
${\cal V}(C_{\Sigma})$ or its powers are the most obvious well defined 
numerical invariants of  ${\cal M}_{TV}(M,\varphi)$. The dimension of 
the module is yet another such invariant.

For a 1/1-solvable theory ${\cal V}$ the invariant ${\cal V}(M,\varphi)$
takes values in
${\mathbb M}[\mbox{\sf y}]/\mbox{\sf y}^2$ and can hence be written
as ${\cal V}(M,\varphi)=\lambda_{\varphi}^{\cal V}(M)\,+\,\mbox{\sf y}\cdot\Xi_{\varphi}^{\cal V}(M)$,
where $\lambda^{\cal V}$ and $\Xi^{\cal V}$ are now ${\mathbb M}$-valued invariants. 
If $\mbox{\sf y}$ coincides with the half projective parameter $\lambda^{\cal V}$ does not
depend on $\varphi$, and we expect it to be some function of $\lambda_{CWL}$. 
Moreover, if ${\cal V}$ descends from a 1/2-solvable TQFT with the 
same property also  $\Xi^{\cal V}$ would be independent of $\varphi$. 

For the modular TQFT over $\FF_5[\mbox{\sf y}]/\mbox{\sf y}^2$ obtained from the Reshetikhin
Turaev theory this invariant has already been defined in  \cite{Kerler02}, and we
may expect it to lift, similarly, to an invariant   $\Xi_{\Z}$ over $\Z$. For
$p>5$ we expect, as in the case of $\lambda_{CWL}$, the next order terms in
the expansions (\ref{eq-LesExp}) of the 
Reshetikhin Turaev theories to be 
polynomial expressions in $\lambda_{CWL}$ and $\Xi_{\Z}$.

\begin{prob}[T. Kerler]  
\label{prob.kerler4}
{\rm [$q/l$-solvable and Casson  TQFT's]}\index{TQFT!Casson ---}\index{TQFT!solvable ---}
\begin{enumerate}
\item  
Lift the 1/1-solvable TQFT's of Casson type over $\FF_p$\index{Casson invariant!--- and TQFT}
to a universal 1/1-solvable TQFT's of Casson type over $\Z$.
\item 
Describe the {\em resulting } invariant
$\Xi_{\Z}$ for 3-manifolds with $b_1(M)\geq 1$.
\item   
Develop a  perturbation theory for general $q/l$-solvable TQFT's. 
\item 
Relate those with the various, standard resolutions of $\Gamma_g$. 
\item 
Relate them also to the traditional finite type theory for closed 3-manifolds.
\item 
Describe the Reshetikhin-Turaev theories
in this pattern.  
\end{enumerate}
\end{prob}

Preparations for item (1) can be found in \cite{Kerler02} in which formulae
for the Casson invariant over $\Z$ are derived that have the same form
as general TQFT formulae. Item (2) is immediate from the preceding
discussion.  
The remaining items are logical continuations. 

\medskip

The category of 3-dim cobordisms $\ThKeCob^{\bullet}$ between compact, oriented
surfaces with one 
boundary component has a natural structure   of 
a braided tensor category.
Another,
category ${\cal A}lg$ can be defined  entirely algebraically in terms of 
generators and relations with respect to a tensor product and a composition
product. On the level of objects it has exactly one generator, say $A$, so 
that all other objects  are of the form $A^{\otimes g}$ with $1=A^{\otimes 0}$.
 The morphisms are given by all words that can be generated by taking 
composition and tensor products of elementary morphisms 
$m:A\otimes A\to A$, $\Delta:A\to A\otimes A$, $e:1\to A$, 
$\varepsilon: A\to 1$, \dots, that appear in the definition of a braided,
ribbon Hopf algebra with integrals and a non-degenerate pairing. 
For example, in \cite{Kerler01} 
a surjective functor ${\cal A}lg\ThKeonto{10}\ThKeCob^{\bullet}$ 
is constructed, which,
in the genus one restriction in fact an isomorphism. 

\begin{prob}[T. Kerler]  
\label{prob.kerler5}
{\rm [3-dim cobordisms from Hopf algebras]}
\begin{enumerate}
\item 
Find further relations on ${\cal A}lg$, besides the ones arising from
the axiomatics of Hopf algebras, 
that would make ${\cal A}lg\to\ThKeCob^{\bullet}$ an isomorphism. 
\item 
Find  relations on ${\cal A}lg$ such that the maps 
$${\rm Aut}_{{\cal A}lg}(A^{\otimes g})\to \Gamma_g
\cong {\rm Aut}_{\ThKeCob^{\bullet}}(\Sigma_{g,1}) $$ are isomorphisms. 
\item 
Relate this to obstructions, such as Steinberg and Whitehead 
groups, via stratified function spaces. 
\item 
What are the analogous algebraic structures in higher dimensions. 
\end{enumerate}
\end{prob}

The first problem is easily stated, 
but presumably very difficult as it implies a faithful translation of 
3-dimensional topology into an algebraic gadget. 
In this respect it is vaguely parallel to 
the geometrization and Poincar\'e conjectures. 
The easier problem stated in item (2) can,
in theory, be attacked head-on, given the known presentations of
the mapping class groups. The third point hints to the fact that the
generators  in ${\cal A}lg$ correspond to Morse-theoretically elementary 
cobordisms, and the relations can be interpreted, similarly,
 in terms of handle slides and cancellation. This is, thus, reminiscent of
the definitions of, e.g.\ Steinberg groups of 3-manifolds. The problem in 
item (4) is, again, easily stated but even in 4 dimensions lingers in 
almost complete  total darkness. It is not hard to understand that higher 
category theory has  to be invoked and not just one ``object'' $A$ suffices 
as a ``generator''. Any partial answers may open the possibility  of
 constructing  functorial 4-manifold invariants by ``linear representation'' of such 
structures.

\medskip

In \cite{KeLy01}
ETQFT's $\ThKeV$ are defined as double functors  from the double category
of relative, 2-framed 1+1+1-dim cobordisms $\ThKeCob^*$
to the double category of 
linear, abelian categories over a perfect field.  (The ``E'' 
stands for ``extended to surfaces with boundaries''). Applied to a single 
circle, thought of as a 0-object in $\ThKeCob^*$, it yields an abelian
category ${\cal C}_{\ThKev}\,=\,\ThKeV(S^1)$, which we call the associated {\em circle
category}. The main result of \cite{KeLy01} is a construction of a $\ThKeV_{\cal C}$, for
each given {\em modular} 
tensor category
$\cal C$ (meaning a bounded, ribbon, braided
tensor category with some additional properties)
such that ${\cal C}_{\ThKev_{\cal C}}={\cal C}$. The construction is made
for all semisimple $\cal C$, and is extended, in the case of non-semisimple $\cal C$,
to both to the situation of connected surfaces with boundary  as well as disconnected,
closed surfaces using the previously mentioned notion of half-projective TQFT's.

\begin{prob}[T. Kerler]  
\label{prob.kerler6}
{\rm [Extended and half-projective  TQFT's]}\index{TQFT!extended ---}\index{TQFT!half-projective ---}
\begin{enumerate}
\item 
Describe in how far an ETQFT $\ThKeV$ with circle category $\cal C$
can differ from $\ThKeV_{\cal C}$, thus introducing a equivalence notion that
would establish a bijective correspondence between the class of
ETQFT's and the class of modular tensor categories.\index{tensor category!modular ---!--- and TQFT}  
\item 
Find an extended notion of half-projectivity that includes also
surfaces that are both disconnected  and have boundary. 
\item 
Find constructions and axioms of ETQFT's 
that apply to more relaxed notions of boundedness or modularity. 
\end{enumerate}
\end{prob}

The functor ${\cal A}lg\to\ThKeCob^{\bullet}$ already imposes that a
circle category ${\cal C}_{\ThKev}$
must fulfill about all axioms of a modular tensor category, and 
contain a Hopf algebra object with properties. Given some rigidity 
assumption it actually must be the same chosen in the construction
of $\ThKeV_{\cal C}$. What may still differ is the choice of algebra 
structures of the same object in the same category, which is thus
the main source of possible ambiguities. Already in \cite{KeLy01} it is clear 
that there are several choices. The correct axiomatics for item (2)
should follow from a careful analysis of the  double composition laws
for surgery tangles from \cite{KeLy01} and generalization of \cite{Kerler98}. 
Item (3) is relevant to include 
more general notions of TQFT's as they would be of interest in the
theory of finite type invariants. 

\medskip

The Reshetikhin-Turaev theory
typically starts with non-semisimple 
modular category
${\cal C}$, typically the representation category of a 
non-semisimple quantum groups $U_q({\frak g})$,
and then considers a 
canonical semisimple sub-quotient $\overline{\cal C}$, see \cite{Kerler92}. 
Thus 
${\cal V}_{\overline{\cal C}}$ yields a semisimple TQFT. It is known that
this is different from the non-semisimple TQFT ${\cal V}_{\cal C}$,
which in the case of a quantum group is obtained via the Hennings 
algorithm.

 TQFT's can also be 
generated from a rigid, monoidal category
${\cal B}$ without any braiding. 
One way is to take the Drinfel'd double $D({\cal B})$, which is then
a modular category for some choice of ribbon element,  and use 
${\cal V}_{D({\cal B})}$. For semisimple ${\cal B}$ one can also extract 
the 6j-symbol data and follow the Turaev-Viro construction to obtain
a TQFT ${\cal W}_{\cal B}$. 

\begin{prob}[T. Kerler] 
\label{prob.kerler7}
{\rm [Non-semisimple vs.\ semisimple TQFT's, the double conjecture]}\index{TQFT!semisimple ---}
\begin{enumerate}
\item  
Clarify the difference in the content of ${\cal V}_{\overline{\cal C}}$
and ${\cal V}_{\cal C}$! \ 
Are there homological TQFT's $\cal H$  such that ${\cal V}_{\cal C}$ is 
in some essential way equivalent to 
${\cal V}_{\overline{\cal C}}\otimes{\cal H}$?  
\item 
Find a  construction of  ${\cal W}_{\cal B}$  that generalizes 
the Turaev-Viro TQFT's\index{Turaev-Viro!--- TQFT}
to non-semisimple ${\cal B}$'s, similar to the
way \cite{KeLy01} generalized the Reshetikhin-Turaev construction.  
In the case of quantum groups
and closed 3-mani\-folds this should reproduce
a version of the Kuperberg invariant. 
\item  
What is the relation between ${\cal W}_{\cal B}$ and 
${\cal V}_{D({\cal B})}$? Are they in some sense isomorphic TQFT's? 
\end{enumerate}
\end{prob}

For the case of $U_q({\frak s}{\frak l}_2)$ there is evidence from
the genus=1 case that such an $\cal H$ is indeed given by the 
Frohman-Nicas-$U(1)$-theory. Item (2) is rather natural as a problem. 
As is apparent in \cite{Kuperberg96} one may expect technical challenges
requiring ``minimal'' cell decompositions of cobordisms, as opposed
to general triangulations, as well as ``combings'' instead of framings. 

The last  conjecture appears also as Question~5 in \cite{Kerler97}
which was motivated by works of and discussions with D. Kazhdan and
S. Gelfand in 1994. Since it is a rather nearby conjecture from a formal
point of view it may have been posed already earlier. For categories
arising from subfactors
and closed manifolds results answering this
conjecture have been obtained in \cite{KSW02}. 
As outlined in \cite{Kerler97} further, more general results in this direction 
should yield a deeper understanding of both TQFT constructions involved 
as well as entail a topological picture for the Drinfel'd double construction.

\newpage

\section{The state-sum invariants of 3-manifolds derived from $6j$-symbols}
\label{sec.state-sum-inv}

\renewcommand{\thefootnote}{\fnsymbol{footnote}}
\footnotetext[0]{The introductory part of 
each section of Chapter \ref{sec.state-sum-inv} was written by T. Ohtsuki,
following suggestions given by Y. Kawahigashi and J. Roberts.}
\renewcommand{\thefootnote}{\arabic{footnote}}

Turaev and Viro \cite{TV} introduced a formulation of 
a state-sum invariant of 3-manifolds
as a state-sum on triangulations of 3-manifolds
derived from certain $6j$-symbols.
After that, Ocneanu gave a general formulation of this state-sum
for general $6j$-symbols
and constructed 3-manifold invariants from subfactors
based on this formulation.
This general formulation was also given by 
Barrett and Westbury \cite{BarrWest96}.

\subsection{Monoidal categories, $6j$-symbols, and subfactors}

Consider a collection, $\{ V_i \}_{i \in I}$,
of (irreducible) modules over $\C$
(of a quantum group or a subfactor)
which is closed under tensor product,
i.e.\ for any $i, j \in I$,
$V_i \otimes V_j \cong \oplus_{k \in I} \cH^k_{i,j} \otimes V_k$
for some $N^k_{i,j}$ dimensional vector space $\cH^k_{i,j}$, 
which expresses the multiplicity of $V_k$ in $V_i \otimes V_j$.
Such a collection (with a certain property)
is called a {\it monoidal category},\index{tensor category!monoidal ---}\index{monoidal category|see {tensor category}} 
where each $V_i$ is called a {\it simple object} of the category
(for details see \cite{BaKi01}).
A monoidal category is provided by
a certain set of representations of a quantum group
(see, e.g.\ \cite{Kassel}),
and also by a certain set of $N$-$N$ bimodules 
arising from a subfactor $N \subset M$
(as explained below).
The algebra spanned by $I$ with the multiplication given by
$a \cdot b = \sum_{c \in I} N^c_{a,b} c$ for $a,b \in I$
is called the {\it fusion rule algebra}.\index{fusion rule algebra}

Let $\{ V_i \}_{i \in I}$ be a monoidal category
(with a finite set $I$)
provided by a quantum group (at a root of unity)
or a subfactor (of finite depth).
Fix the above mentioned isomorphism
$V_i \otimes V_j \cong \oplus_{k \in I} \cH^k_{i,j} \otimes V_k$
for each $i, j$.
Then, we have two bases of the vector space
$\mbox{Hom}(V_l, V_i \otimes V_j \otimes V_k)$
for each $i,j,k,l$ as follows.
Consider the maps
$$
V_l \longrightarrow V_n \otimes V_k \longrightarrow
(V_i \otimes V_j) \otimes V_k
$$
determined by basis vectors
$A \in \cH^l_{n,k}$ and $B \in \cH^n_{i,j}$.
The composition of these maps gives a vector of
$\mbox{Hom}(V_l, V_i \otimes V_j \otimes V_k)$.
Thus, we obtain a basis of this vector space
consisting of vectors labeled by triples $(n,A,B)$.
Moreover, we obtain another basis consisting of
vectors labeled by triples $(m,C,D)$, 
where $C \in \cH_{i,m}^l$, $D \in \cH_{j,k}^m$,
by considering the following maps,
$$
V_l \longrightarrow V_i \otimes V_m \longrightarrow
V_i \otimes (V_j \otimes V_k).
$$
The collection of the entries of the matrix which relates these two bases
is a typical example of a set of $6j$-symbols,
where a {\it set of $6j$-symbols}\index{6j-symbol} 
is defined to be
a solution of certain polynomial equations:
the tetrahedral symmetry, the unitarity, and the pentagon relation.
Each $6j$-symbol is labeled by $i,j,k,l,m,n \in I$, and
$A \in \cH^l_{n,k}$, $B \in \cH^n_{i,j}$, 
$C \in \cH_{i,m}^l$, and $D \in \cH_{j,k}^m$.
This $6j$-symbol will be associated to a tetrahedron
\pict{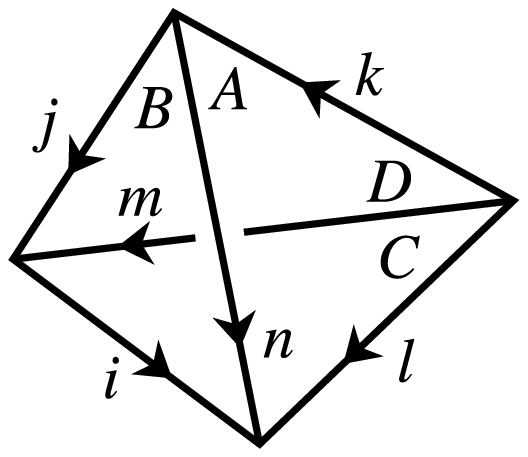}{2cm}
labeled by them.

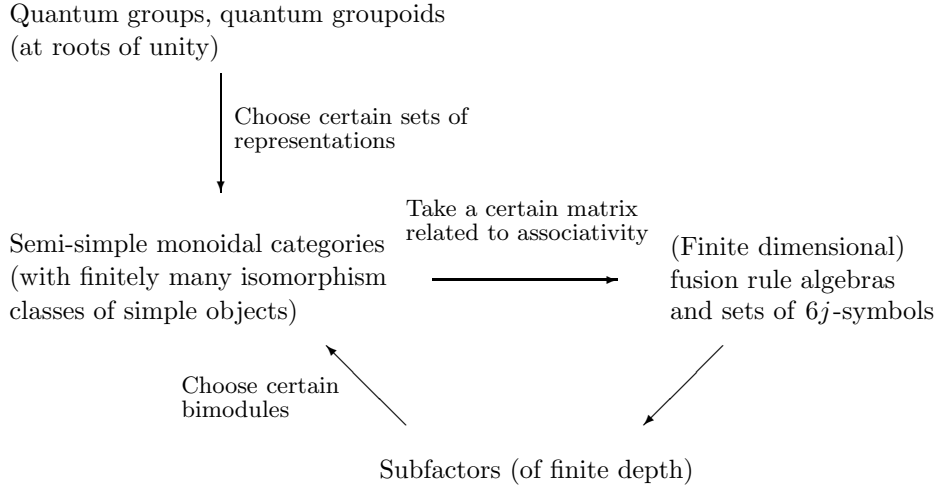
\begin{figure}[ht!]
\begin{center}
\begin{picture}(350,190)
   \put(-10,160){\shortstack[l]{\small Quantum groups, quantum groupoids \\
                           \small (at roots of unity)}}
   \put(70,153){\vector(0,-1){45}}
   \put(75,125){\shortstack[l]{\footnotesize Choose certain sets of \\ 
                               \footnotesize representations}}
   \put(-10,60){\shortstack[l]{\small Semi-simple monoidal categories \\
        \small (with finitely many isomorphism \\ \small classes of simple objects)}}
   \put(150,75){\vector(1,0){70}}
   \put(140,90){\shortstack[l]{\footnotesize Take a certain matrix \\ 
                               \footnotesize related to associativity}}
   \put(240,60){\shortstack[l]{\small (Finite dimensional) \\
                            \small fusion rule algebras \\ 
                            \small and sets of $6j$-symbols}}
   \put(140,20){\vector(-1,1){30}}
   \put(55,23){\shortstack[l]{\footnotesize Choose certain \\ 
                              \footnotesize bimodules}}
   \put(260,50){\vector(-1,-1){30}}
   \put(265,23){\shortstack[l]{}}
   \put(130,0){\shortstack[l]{\small Subfactors (of finite depth)}}
\end{picture}
\end{center}
\caption{\label{fig.mono_6j_subfac}
$6j$-symbols and related objects}
\end{figure}

A {\it subfactor}\index{subfactor} 
is a pair of infinite dimensional algebras $N$ and $M$
with an inclusion relation $N \subset M$
satisfying some property.
A major class of subfactors is a class of 
{\it WZW model subfactors} of level $k = 1,2, \cdots$,
which are related to quantum groups.
Another well-known class is a class of subfactors of the Jones index $<4$;
they are classified to be of types $A_n, D_{2n}, E_6$, or $E_8$.
A left $X$ right $Y$ module $Z$ is called
a $X$-$Y$ {\it bimodule}, and is written ${}_X Z_Y$.
For a subfactor $N \subset M$,
consider irreducible $N$-$N$ bimodules
appearing as direct summands of $N$-$N$ bimodules
in the following sequence,
$$
{}_N N_N, \qquad
{}_N M_M \underset{M}{\otimes} {}_M M_N, \qquad
{}_N M_M \underset{M}{\otimes} {}_M M_N \underset{N}{\otimes}
{}_N M_M \underset{M}{\otimes} {}_M M_N, \qquad \cdots.
$$
The collection of (isomorphism classes of) such irreducible modules
provides a monoidal category
$\{ V_i \}_{i \in I}$.
It is known that
$I$ is a finite set
when the subfactor is of finite depth
(this always holds when its index $< 4$).
For a fusion rule algebra with a set of $6j$-symbols
there exists a subfactor
(if quantum dimensions are positive)
such that the diagram in Figure \ref{fig.mono_6j_subfac} commutes.
For details of this paragraph see \cite{GHV_opalg,sato.EvKa_OUP}.

\index{quantum group}\index{quantum groupoid}\index{6j-symbol}\index{subfactor}\index{tensor category!monoidal ---}

Thus, the following classification problems are almost equivalent.
Each of them is fundamental, but probably impossibly hard.
(See also Problem \ref{prob.classify_modu_cat}.)

\begin{prob}
\mbox{}
\begin{itemize}
\item[\rm(1)]
Find (and classify) all semi-simple monoidal categories\index{tensor category!monoidal ---!classification of ---}  
(with finitely many isomorphism classes of simple objects).
\item[\rm(2)]
Find (and classify) (finite dimensional) 
fusion rule algebras\index{fusion rule algebra!classification of ---} 
and sets of $6j$-symbols.\index{6j-symbol!classification of ---}
\item[\rm(3)]
Find (and classify) all 
subfactors (of finite depth).\index{subfactor!classification of ---}
\end{itemize}
\end{prob}

\begin{rem}
Major sets of $6j$-symbols,
what we call {\it quantum $6j$-symbols},\index{6j-symbol!quantum ---}
are the sets of $6j$-symbols 
derived from quantum groups,
resp. WZW model subfactors.
Another class of $6j$-symbols is derived from finite groups;\index{finite group} 
for a 3-cocycle $\alpha$ of a finite group $G$,
a set of $6j$-symbols is given by
$$
W\Big(\pict{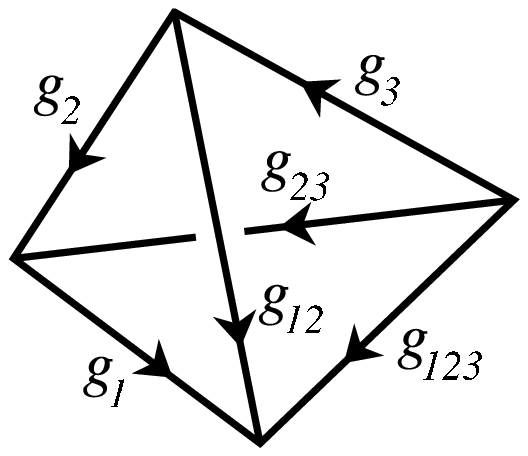}{2cm}\Big) = \begin{cases}
\mboxsm{$\alpha(g_1,g_2,g_3)$}  & 
\mboxsm{if $g_{12} {=} g_1 g_2, \, g_{23} {=} g_2 g_3$, and 
$g_{123} {=} g_1 g_2 g_3$,} \\
\mboxsm{$1$}  & \mboxsm{otherwise,} \end{cases}
$$
where the tetrahedra is given a trivial face coloring.
There are still other infinitely many sets of $6j$-symbols
arising from subfactors; see Table \ref{tbl.subfactor}.
These $6j$-symbols might have a universal presentation
given by a tetrahedron in the theory of knotted trivalent graphs
(see Section \ref{sec.ktg}).
\end{rem}

\begin{table}[ht!]
\begin{center}
\begin{tabular}{|l|l|l||l|l|} \hline 
\multicolumn{3}{|l||}{\dwn{\small subfactor}}
& \small monoidal & \dwn{\small $S$-matrix} \\ 
\multicolumn{3}{|l||}{} & \small category & \\ \hline\hline
\multicolumn{3}{|l||}{\small WZW model subfactors of level $k=1,2,\cdots$} & 
\dwn{\small braided}
& \dwn{\small non-degenerate} \\
\multicolumn{3}{|l||}{\small $SU(N)_k$, $SO(N)_k$, $Sp(N)_k$, $\cdots$} 
& & \\ \hline
\multicolumn{2}{|l|}{\dwn{\small subfactors of}}
& \small type $A_n$ ($= SU(2)_{n}$) & \small braided 
& \small non-degenerate \\ \cline{3-5}
\multicolumn{2}{|l|}{\dwn{\small index $<4$}} & \small type $D_{2n}$ & 
\small braided\footnotemark\addtocounter{footnote}{-1} 
& \small non-degenerate\footnotemark \\ \cline{3-5}
\multicolumn{2}{|l|}{} & 
\small type $E_6, E_8$ & \small not braided & \small none \\ \hline
\multicolumn{2}{|l|}{\small subfactors of}
& \small (generalized) Haagerup, & \dwn{\small not braided} 
& \dwn{\small none} \\
\multicolumn{2}{|l|}{\small index $> 4$: }
& \small Asaeda-Haagerup, $\cdots$ & & \\ \cline{3-5}
\multicolumn{2}{|l|}{\small exotic subfactors,} 
& \small quantum doubles of & \dwn{\small braided} 
& \dwn{\small non-degenerate} \\
\multicolumn{2}{|l|}{\small $\cdots$} 
& \small Haagerup subfactor, $\cdots$ & & \\ \hline
\small subfactors
& \multicolumn{2}{|l||}{\small 3-cocycles of finite groups}
& \small not braided & \small none 
\\ \cline{2-5}
\small from ---
& \multicolumn{2}{|l||}{\small representations of finite groups}
& \small braided\footnotemark & \small degenerate \\ \hline
\end{tabular}\end{center}
\caption{\label{tbl.subfactor}\index{subfactor}\index{tensor category!monoidal ---}
Subfactors, their monoidal categories, and $S$-matrices}
\end{table}\addtocounter{footnote}{-1}\footnotetext{
To be precise, the even part of the subfactor of type $D_{2n}$ is braided, 
and its $S$-matrix is non-degenerate.}\addtocounter{footnote}{1}\footnotetext{
This is trivially braided.}

\subsection{Turaev-Viro invariants and 
the state-sum invariants derived from monoidal categories}

A state-sum invariant
of 3-manifolds
is defined by using such a set of $6j$-symbols
with a monoidal category
$\{ V_i \}_{i \in I}$, as follows.
Choose a simplicial decomposition of a closed 3-manifold $M$,
and fix a total order of its vertices,
which induces orientations of edges.
Further, choose an {\it edge coloring} $\lambda$,
which is a map of the set of edges to $I$,
and choose a {\it face coloring} $\varphi$,
which is a collection of such assignments that
a basis vector of $\cH^k_{i,j}$
is assigned to a triangle \pict{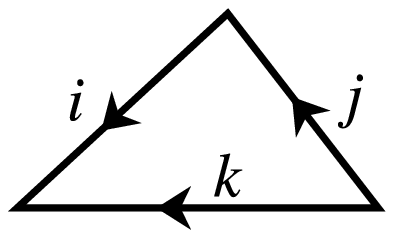}{1cm}
with an edge coloring.
To a tetrahedron $\sigma$ 
with an edge coloring $\lambda$ and a face coloring $\varphi$, 
we associate the above mentioned $6j$-symbol,
which we denote by $W(\sigma;\lambda,\varphi)$.
Then, a {\it state-sum invariant}\index{state-sum!--- invariant}   
of $M$ is defined by
\begin{equation}
\label{eq.defTVO}
Z(M) = w^{-v} \sum_{\lambda} \Big( \prod_E \mu_{\lambda(E)} \Big)
\sum_\varphi \prod_\sigma W(\sigma; \lambda,\varphi),
\end{equation}
where the sums of $\lambda$ and $\varphi$ run over
all edge colorings and all face colorings,
and the products of $E$ and $\sigma$
run over all edges and all tetrahedra
of the simplicial decomposition of $M$, and
$\mu_i$ is a constant, which corresponds to a ``quantum dimension'', and
$w = \sum_{i \in I} \mu_i^2$, and
$v$ is the number of vertices of the simplicial decomposition.
It is known (see \cite{BarrWest96}, \cite[Chapter 12]{sato.EvKa_OUP}) that
the invariant (\ref{eq.defTVO}) is a topological invariant of $M$.
The definition of the invariant (\ref{eq.defTVO})
can naturally be extended to an invariant of 3-manifolds with boundaries,
and a TQFT can be formulated based on it.

In particular, for the set of $6j$-symbols
arising from representations of 
the quantum group $U_q(Sl_2)$ at a root of unity,
the invariant (\ref{eq.defTVO}) is called
{\it the Turaev-Viro invariant} \cite{TV}.\index{Turaev-Viro!--- invariant}
In its definition
it is not necessary to introduce face colorings
(because $N_{i,j}^k$ is always equal to $0$ or $1$ for any $i,j,k$ 
in this case)
and orientations of edges
(because each representation of $U_q(sl_2)$ is self-dual).


The monoidal category of a set of quantum $6j$-symbols
is a modular category,
and 
we can construct the Reshetikhin-Turaev invariant
from it 
(see Section \ref{sec.ssinv_ribbon}).
The square of the absolute value of the invariant
is equal to the value of the state-sum invariant
derived from these $6j$-symbols.

The state-sum invariant
derived from the set of $6j$-symbols
given by a 3-cocycle $\alpha$ of a finite group\index{finite group}  
$G$
is called the {\it Dijkgraaf-Witten invariant} \cite{DiWi90}.
In particular, when $\alpha = 1$, it is equal to 
the number of conjugacy classes of representations $\pi_1(M) \to G$.
It is further equal to the state-sum invariant
derived from the set of $6j$-symbols
obtained from the representations of the finite group $G$.

When a set of $6j$-symbols arises from a subfactor,
the state-sum invariant
derived from these $6j$-symbols is called
the {\it Turaev-Viro-Ocneanu invariant}.\index{Turaev-Viro!Turaev-Viro-Ocneanu invariant}
There are infinitely many subfactors other than the above cases
as shown in Table \ref{tbl.subfactor}.
The Turaev-Viro-Ocneanu invariants derived from 
such subfactors might be new invariants of 3-manifolds.

\begin{prob}[Y. Kawahigashi] 
\label{kawah_prob.4}
Suppose we have a three-dimensional TQFT.\index{TQFT!--- from 6j-symbol}
Can we determine whether it arises from a 
fusion rule algebra\index{fusion rule algebra!TQFT arising from ---}
and $6j$-symbols?\index{6j-symbol!TQFT arising from ---}
If yes, can we describe all fusion rule algebras with $6j$-symbols 
producing the TQFT?
\end{prob}

\begin{rem}[\rm (Y. Kawahigashi)] 
By a result of Ocneanu, we have at most only finitely many
such fusion rule algebras with $6j$-symbols, up to equivalence
of $6j$-symbols.
\end{rem}

\begin{prob}[Y. Kawahigashi] 
\label{kawah_prob.5}
Suppose we have two fusion rule algebras with 
$6j$-symbols\index{6j-symbol!TQFT arising from ---}\index{TQFT!--- from 6j-symbol}\index{fusion rule algebra!TQFT arising from ---}
and that two TQFT's arising from them are isomorphic.  
What relation do we have for the two sets of $6j$-symbols?  
\end{prob}

\begin{rem}[\rm (Y. Kawahigashi)] 
Are they equivalent in the sense of \cite{k_S}?
\end{rem}

\begin{prob}[Y. Kawahigashi] 
\label{kawah_prob.6}
Suppose we have a TQFT arising from a fusion rule algebra with 
$6j$-symbols.\index{6j-symbol!TQFT arising from ---}\index{TQFT!--- from 6j-symbol}\index{fusion rule algebra!TQFT arising from ---}  
Using a fusion rule subalgebra and $6j$-symbols restricted on it, 
we can construct another TQFT.  
What relation do we have for these TQFT's?  
\end{prob}

\begin{rem}[\rm (Y. Kawahigashi)] 
How about the case where the fusion rule subalgebra 
arises from $\alpha$-induction?  
The $\alpha$-induction produces a fusion rule algebra with
$6j$-symbols from a  semisimple ribbon category
with finitely many isomorphism classes of simple objects
and a specific choice of an object satisfying certain
axioms.  See \cite{k_BEK}, \cite{k_KO} and their references.
If the original ribbon category is modular, we have some
answer in \cite{k_BEK}, so it is particularly interesting when
the $S$-matrix is not invertible.
\end{rem}

\subsection{The state-sum invariants derived from ribbon categories}
\label{sec.ssinv_ribbon}

A {\it ribbon category}\index{tensor category!ribbon ---}\index{ribbon category|see {tensor category}} 
is a monoidal category $\{ V_i \}_{i \in I}$
equipped with a {\it braiding} $V \otimes W \to W \otimes V$
and a {\it twist} $V \to V$ for any objects $V$ and $W$
which are maps satisfying certain properties.
We obtain an invariant of framed links from a ribbon category
by associating a braiding to a crossing of a link diagram
and a twist to a full-twist of a framing of a link.
A monoidal category is called {\it semi-simple}
if any object is isomorphic to a direct sum of simple ones.
The {\it S-matrix} $S = ( S_{ij} )_{i,j \in I}$
of a semi-simple ribbon category $\{ V_i \}_{i \in I}$
is defined by putting
$S_{ij}$ to be the invariant of the Hopf link
whose components are associated with $V_i$ and $V_j$.
A {\it modular category}\index{tensor category!modular ---}\index{modular category|see {tensor category}}
is a semi-simple ribbon category
with finitely many isomorphism classes of simple objects
whose S-matrix is invertible.
We obtain the Reshetikhin-Turaev invariant
of 3-manifolds and its TQFT from a modular category
by using surgery presentations of the 3-manifolds.
See \cite{BaKi01} for details of this paragraph.

\begin{figure}[ht!]
\begin{center}
\begin{picture}(380,140)
   \put(0,120){\small Monoidal category}
   \put(90,120){\small $\Longrightarrow$ \ State-sum invariant of 3-manifolds and its TQFT}
   \put(30,110){\vector(0,-1){35}}
   \put(35,85){\shortstack[l]{\footnotesize $+$ braiding \\ 
                               \footnotesize $+$ twist}}
   \put(0,60){\small Ribbon category}
   \put(90,60){\small $\Longrightarrow$ \ Invariant of framed links}
   \put(30,50){\vector(0,-1){35}}
   \put(35,20){\shortstack[l]{\footnotesize $+$ semi-simple \\
                              \footnotesize $+$ finiteness of $I$ \\
                              \footnotesize $+$ invertibility of $S$}}
   \put(0,0){\small Modular category}
   \put(90,0){\small $\Longrightarrow$ \ Reshetikhin-Turaev invariant of 3-manifolds and its TQFT}
\end{picture}
\end{center}
\caption{\label{fig.mono_rib_modu}
Monoidal, ribbon, modular categories and their consequences}
\end{figure}\index{Reshetikhin-Turaev!--- invariant}\index{state-sum!--- invariant}\index{tensor category!ribbon ---}\index{tensor category!monoidal ---}\index{tensor category!modular ---}

The quantum $6j$-symbols are typical $6j$-symbols
which induce modular categories.
The square of the absolute value of the Reshetikhin-Turaev invariant
derived from a modular category
is equal to the value of 
the state-sum invariant derived from the category.
It is suggested by Ocneanu that
the monoidal category of the quantum double
of each of such subfactors
would be braided, and that 
the Reshetikhin-Turaev invariant derived from this quantum double
would be equal to the Turaev-Viro-Ocneanu invariant
derived from the original subfactor.

\begin{prob}[Y. Kawahigashi] 
\label{kawah_prob.3}
Suppose we have a semisimple 
ribbon category\index{tensor category!ribbon ---!--- with degenerate S-matrix} 
$C$
with finitely many isomorphism classes of simple objects.
If the $S$-matrix is invertible, we can construct the
Reshetikhin-Turaev invariant
and the state-sum invariant\index{state-sum!--- invariant!--- from degenerate S-matrix} 
from $C$ and the latter is the square of the
absolute value of the former.  If the $S$-matrix is
not invertible, do we still have a similar description 
of the state-sum invariant?
\end{prob}

\begin{rem}
See also Problem \ref{prob.sato3}
for a similar problem for the Turaev-Viro-Ocneanu invariants.
\end{rem}

\begin{prob}[Y. Kawahigashi] 
\label{kawah_prob.1}
Suppose we have a semisimple 
ribbon category $C_1$\index{tensor category!ribbon ---!--- with degenerate S-matrix}
with finitely many isomorphism classes of simple objects, 
but the $S$-matrix is not invertible.  
Then we can construct a new modular category\index{tensor category!modular ---!--- as extension of ribbon category} 
$C_2$ containing $C_1$ 
as a full subcategory by the ``quantum double'' construction\index{quantum double} 
\cite{k_O1,k_O2,k_I}, but
there may be another extension of $C_1$ to a modular category.  
Theorem 2.13 in \cite{k_O1} claims
that we have a ``minimal'' extension in an ``essentially unique'' way.  
Do we indeed have existence and certain
uniqueness of such an extension?  If so, what is the relation
between the two TQFT's\index{TQFT!--- from tensor category} 
arising from $C_1$ and its minimal extension?
\end{prob}

\begin{prob}[Y. Kawahigashi] 
\label{kawah_prob.2}
Suppose we have a semisimple ribbon category $C_1$\index{tensor category!ribbon ---!--- with degenerate S-matrix} 
with
a degenerate $S$-matrix as in Problem \ref{kawah_prob.1}.
By the method in \cite{k_M}, we can also make a modular
tensor category\index{tensor category!modular ---!--- as extension of ribbon category} 
$C_2$ from $C_1$.  What is the relation
between the two TQFT's\index{TQFT!--- from tensor category}
arising from $C_1$ and $C_2$?
\end{prob}

\begin{prob}[Y. Kawahigashi] 
\label{kawah_prob.7}
There are some fusion rule algebras with 
$6j$-symbols\index{6j-symbol!non-quantum ---}\index{fusion rule algebra!non-quantum ---}
that do not seem to arise from quantum groups
 in \cite{k_AH}
and more conjectured candidates of such examples in \cite{k_H}.
What are the corresponding TQFT's?\index{TQFT!non-quantum ---}
Especially if the series conjectured in \cite{k_H} does exist,
it would give a  parametrized family of TQFT's.  Does a
differentiation by a parameter (after a certain reparametrization)
give a more interesting invariant,
possibly of Vassiliev type?
\end{prob}

\subsection{Turaev-Viro-Ocneanu invariants} 
\label{sec.TVOinv}

The state-sum invariant of 3-manifolds derived from $6j$-symbols
is called the {\it Turaev-Viro-Ocneanu invariant}\index{Turaev-Viro!Turaev-Viro-Ocneanu invariant} 
when the set of $6j$-symbols arises from a subfactor.
There are infinitely many subfactors
other than those derived from quantum groups or finite groups.
The Turaev-Viro-Ocneanu invariants derived from 
such subfactors might be new invariants of 3-manifolds.

\medskip
\noindent
{{\namae}(N. Sato)}
\quad
\noindent
The Haagerup subfactor\index{subfactor!Haagerup ---} 
of Jones index $\frac{5+\sqrt{13}}{2}$ has the 
smallest index among finite depth subfactors with Jones index 
bigger than 4 and it is expected to have some \lq\lq exotic" properties from 
the subfactor theoretical viewpoint. However, it does not seem so sensitive to 
classify 3-manifolds. 
The Turaev-Viro-Ocneanu invariant\index{Turaev-Viro!Turaev-Viro-Ocneanu invariant} 
constructed from 
the Haagerup subfactor cannot distinguish lens spaces\index{lens space} $L(5,1)$ and $L(5,2)$, 
as well as $L(7,1)$ and $L(7,2)$. On the other hand, generalized 
$E_6$-subfactors with the group symmetries $\Bbb{Z}/3\Bbb{Z}$ and 
$\Bbb{Z}/5\Bbb{Z}$ can distinguish $L(3,1)$ and $L(3,2)$, $L(5,1)$ and 
$L(5,2)$, respectively. 

\begin{prob}[N. Sato] 
\label{prob.sato1}
Find a subfactor\index{subfactor!--- distinguishing lens spaces} 
which can distinguish lens spaces\index{lens space} $L(7,1)$ and $L(7,2)$. 
Moreover, find a subfactor to classify 3-manifolds as well as possible. 
\end{prob}

In the lattice field theory, 
Ponzano and Regge \cite{sato.PoRe68} constructed 
a state sum model 
for $SU(2)$ and investigated an asymptotic behavior of the model. 

Some infinite depth subfactors are manageable in the sense of growth rate 
(amenability). 
Such subfactors are called {\it strongly amenable}.\index{subfactor!strongly amenable ---} 
The strong 
amen\-ab\-ility condition might be enough to control the asymptotic behavior 
of the state sum model constructed from a strongly amenable subfactor.  

\begin{prob}[N. Sato] 
\label{prob.sato2}
Construct a well-defined 
state sum type invariant\index{state-sum!--- invariant!--- from strongly amenable subfactor}
from a strongly amenable 
subfactor. 
\end{prob}

Note that, unlike the Ponzano-Regge model,  we do not have an asymptotic 
description of the quantum $6j$-symbols
in general.
(Recall that $6j$-symbols 
of $SU(2)$ have an asymptotic description.)

Let us consider the Turaev-Viro-Ocneanu invariant for a closed 3-manifold 
constructed from a subfactor.
Then, this invariant can be considered as a 
Reshetikhin-Turaev type invariant
constructed from a subfactor by passing 
the initial subfactor through the Longo-Rehren construction. If we start with 
a subfactor which has a non-degenerate braiding in particular, 
then this Turaev-Viro-Ocneanu invariant splits into a Reshetikhin-Turaev 
invariant and its complex conjugate. The following question will open a way 
to establish a theory of the minimal non-degenerate extension of a degenerate 
braiding.   

\begin{prob}[N. Sato] 
\label{prob.sato3}
Let us consider the Turaev-Viro-Ocneanu invariant\index{Turaev-Viro!Turaev-Viro-Ocneanu invariant} 
from a subfactor\index{subfactor!--- with degenerate braiding} 
with a degenerate braiding. Then, find a description of this 
invariant as a Reshetikhin-Turaev invariant.\index{Reshetikhin-Turaev!--- invariant!--- and Turaev-Viro-Ocneanu invariant} 
\end{prob}

\begin{rem}
See also Problem \ref{kawah_prob.3}
for a similar problem for the state-sum invariants
derived from ribbon categories.
\end{rem}

\newpage

\section{Casson invariant and finite type invariants of 3-manifolds}

\subsection{Casson and Rokhlin invariants}

It is known as Rokhlin theorem that
the signature of a spin smooth closed 4-manifold is divisible by 16,
which deduce the following definition of the Rokhlin invariant.
For a closed 3-manifold $M$ and a spin structure $\sigma$ on $M$,
the {\it Rokhlin invariant} $\mu(M,\sigma) \in \Z/16\Z$
is defined to be
the signature of any smooth compact spin 4-manifold
with spin boundary $(M,\sigma)$.
In particular, for a $\Z/2\Z$ homology 3-sphere $M$,
the {\it Rokhlin invariant}\index{Rokhlin invariant} 
$\mu(M) \in \Z/16\Z$ is defined to be
the signature of any smooth compact spin 4-manifold with boundary $M$,
noting that there exists a unique spin structure on such a $M$.
The Casson invariant\index{Casson invariant} 
is a $\Z$-valued lift of
the Rokhlin invariant of integral homology 3-spheres.
Further, it is known \cite{Walker_Cinv} that
$$
\mu(M) \equiv 4 | H_1(M;\Z) |^2 \lambda_{\rm CW}(M) 
\equiv 8 | H_1(M;\Z) | \lambda_{\rm CWL}(M) 
\qquad \mbox{ (mod 16)}
$$
for any $\Z/2\Z$ homology 3-sphere, where
$\lambda_{\rm CW}$ denotes the Casson-Walker\index{Casson invariant!Casson-Walker invariant} 
invariant\footnote{\footnotesize 
The normalization here is that
$\lambda_{\rm CW}(M) =  2 \lambda_{\rm C}(M)$ 
for an integral homology 3-sphere $M$.}
\cite{Walker_Cinv} and
$\lambda_{\rm CWL}$ denotes the Casson-Walker-Lescop\index{Casson invariant!Casson-Walker-Lescop invariant} 
invariant\footnote{\footnotesize 
The normalization here is that
$\lambda_{\rm CWL}(M) = \big( | H_1(M;\Z) | /2 \big) 
\lambda_{\rm CW}(M)$ 
for a rational homology 3-sphere $M$.}
\cite{Lescop96}.
For an exposition of the Casson and Rokhlin invariants,
see \cite{KM_inv,Lescop96,Saveliev}.

\begin{prob}
Can the Casson invariant\index{Casson invariant!--- and signature of 4-manifold}
 of an integral homology 3-sphere $M$
be characterized by the signature\index{signature!--- of 4-manifold}
of a certain 4-manifold bounded by $M$?
\end{prob}

\begin{rem}
It is shown in \cite{FMM} that
the Casson invariant of the Seifert fibered homology 3-sphere
$\Sigma(\alpha_1,\cdots,\alpha_n)$
is equal to $1/8$ times the signature of its Milnor fiber.
\end{rem}

The Casson-Walker-Lescop invariant 
of closed 3-manifolds with positive Betti number
can be computed from the torsion invariant $\tau$ of V. Turaev.
He \cite{Turaev_SW} gave a surgery formula for $\tau$,
which implies a surgery formula for the Casson-Walker-Lescop invariant.

\begin{prob}[V. Turaev] 
\label{prob.sf_CWL}
Relate this surgery formula for 
the Casson-Wal\-ker-Lescop invariant\index{Casson invariant!Casson-Walker-Lescop invariant!--- and Reidemeister-Turaev torsion}\index{Reidemeister-Turaev torsion!--- and Casson-Walker-Lescop invariant}
with that of Lescop \cite{Lescop96}.
\end{prob}

\noindent
{{\namae}(C. Lescop)}\quad
In 1984, Casson defined his invariant of integral homology 3-spheres 
as an integer that ``counts'' the $SU(2)$-representations 
of their fundamental group in an appropriate way (see \cite{Akbulut-McCarthy,Guillou-Marin}).
Cappell, Lee and Miller \cite{Cappell-Lee-Miller}
showed that the Casson way of counting $SU(2)$-representations 
of the $\pi_1$ works for any compact Lie group and 
provides other invariants of integral homology spheres.

\begin{quest}[C. Lescop] 
\label{quest.SUn_Cinv}
Are the Cappell-Lee-Miller 
Casson-type $SU(n)$-invariants\index{Cappell-Lee-Miller invariant}
of\index{finite type invariant!--- of homology spheres!--- and Cappell-Lee-Miller invariant} 
finite type?
If so, what are their degrees and their\index{weight system!--- of Cappell-Lee-Miller invariant} 
 weight systems?
\end{quest}

\begin{prob}[M. Polyak] 
\label{prob.spinCWL}
Define an invariant $\lambda$ of a pair $(M,\sigma)$
of a closed 3-manifold $M$ and a spin structure $\sigma$ on $M$
such that
$$
\lambda_{\rm CWL}(M) = \sum_\sigma \lambda(M,\sigma)
$$
for any closed 3-manifold $M$,\index{Casson invariant!Casson-Walker-Lescop invariant!refinement of ---}
where the sum runs over all spin structures $\sigma$ on $M$.
\end{prob}

Note that the set of spin structures on $M$
is a torsor over $H^1(M;\Z/2\Z)$
in the sense that
differences of spin structures can be detected by
cohomology classes in $H^1(M;\Z/2\Z)$,
while the set of spin${}^c$ structures on $M$
is a torsor over $H^1(M;\Z)$
in a similar sense.

\begin{rem}
It is shown \cite{OzSz00} that
there exists an invariant $\hat\theta$ of a rational homology 3-sphere $M$
associated with a spin${}^c$ structure $\alpha$ on $M$ such that
$$
\frac12 | H_1(M,\Z) | \lambda_{\rm CW}(M) = \sum_\alpha \hat\theta(M,\alpha)
$$
for any rational homology 3-sphere $M$,
where the sum runs over all spin${}^c$ structures $\alpha$ on $M$.
It is conjectured \cite{OzSz00} that $\hat\theta$ is equal to
Seiberg-Witten invariant for all rational homology 3-spheres.
\end{rem}

\begin{rem}[\rm (M. Polyak)] 
The Casson invariant is a lift of The Rokhlin invariant.
We expect that
$\lambda(M,\sigma)$ of Problem \ref{prob.spinCWL}
should be a lift of $\mu(M,\sigma)$.
How is 
$\sum_\sigma \mu(M,\sigma) \in \Z/16\Z$
related to $\lambda_{\rm CWL}(M)$?

It is known that
this sum vanishes in $\Z / 16 \Z$ when $b_1(M) > 3$,
while it is known \cite{Lescop96} that
$\lambda_{\rm CWL}(M) = 0$ when $b_1(M) > 3$.
\end{rem}

\begin{rem}[\rm (C. Lescop)] 
Let $M$ be the 3-manifold obtained by surgery along a framed link $L$,
and let $W$ be the 4-manifold associated to the surgery presentation.
Then,
$ 24 \lambda_{\rm CWL}(M) - 3 |H_1(M;\Z)| \mbox{sign} W$
can be presented by a formula of Alexander polynomial coefficients and 
linking numbers of $L$
\cite[Formula 6.3.1]{Lescop96}, which might be helpful.

Note that
the list of $\mu(M,\sigma)$ for a given $M$ is richer than
their sum $\sum_\sigma \mu(M,\sigma)$.
For example, $\mu(\R P^3, \sigma) = 1$, $-1$ and
$\mu (\R P^3 \# \mbox{(Poincare sphere)},\sigma) = 7$, $9$,
while their sums are equal in $\Z / 16 \Z$.
\end{rem}

\begin{rem}
The invariant of Problem \ref{prob.spinCWL} should be related to 
the Goussarov-Habiro theory for spin 3-manifolds \cite{Massuyeau}.
Recall that the Rokhlin and Casson invariants can be characterized
as invariants under $Y_2$-equivalence and $Y_3$-equivalence among {\ZHSs}
respectively.
It was shown \cite{Massuyeau} that
Rokhlin invariant of spin closed 3-manifolds
is the invariant under spin $Y_2$-equivalence
among spin closed 3-manifolds.
What is the invariant under spin $Y_3$-equivalence?
\end{rem}

\begin{rem}
The Casson-Walker invariant can be characterized as
the first coefficient of the perturbative expansion of
the quantum $SO(3)$ invariant
 $\tau^{SO(3)}(M)$ \cite{Hitoshi_CW}.
We have a spin refinement
$\tau_r^{SU(2)}(M,\sigma)$
of the quantum $SU(2)$ invariant $\tau_r^{SU(2)}(M)$ for $r \equiv 0$ mod 4
such that
$$
\tau_r^{SU(2)}(M) = \sum_\sigma \tau_r^{SU(2)}(M,\sigma), 
$$
where the sum runs over all spin structures $\sigma$ on $M$ \cite{KM_inv}.
We expect that
$\lambda(M,\sigma)$ of Problem \ref{prob.spinCWL}
should be related to the first coefficient of
the perturbative expansion of $\tau_r^{SU(2)}(M,\sigma)$.

For $r \equiv 2$ mod 4,
we have another refinement $\tau_r^{SU(2)}(M, \xi)$
for $\xi \in H_1(M;\Z/2\Z)$
such that
$$
\tau_r^{SU(2)}(M) = \sum_\xi \tau_r^{SU(2)}(M,\xi),
$$
where the sum runs over all cohomology classes in $H^1(M;\Z/2\Z)$.
The first coefficient of the perturbative expansion
of $\tau_r^{SU(2)}(M,\xi)$ was discussed in \cite{Hitoshi_tfc,Hitoshi_ncc}.
It might be a problem
to find a refinement $\lambda(M,\xi)$ of $\lambda_{\rm CW}(M)$
for some cohomology class $\xi$.
\end{rem}


\begin{rem}
Problem \ref{prob.spinCWL} is related to Problem \ref{prob.spinLMO},
which is a problem to find a spin refinement of the LMO invariant,
noting that the first coefficient of the LMO invariant
is given by the Casson-Walker-Lescop invariant.
\end{rem}

\begin{quest}[M. Polyak] 
\label{quest.Rohlin_spin}
Is there a ``Rokhlin invariant''\index{Rokhlin invariant!--- of spin${}^c$ 3-manifold}
of a pair $(M, \alpha)$
of a closed 3-manifold $M$ and a spin${}^c$ structure $\alpha$ on $M$?
(See Question \ref{quest.lift_g}.)
\end{quest}

\begin{prob}[M. Polyak] 
\label{prob.Gdf_ftinv}
By presenting 3-manifolds by surgery along\break framed links in $S^3$,
we can regard an invariant of 3-manifolds as an invariant of framed links.
Establish a Gauss diagram formula
for the link invariant derived from
each finite type invariant of 3-manifolds.\index{finite type invariant!--- of homology spheres!--- by Gauss diagram formula} 
\end{prob}

\begin{rem}[\rm (M. Polyak)] 
The first step is to find a Gauss diagram formula for the
Casson invariant.\index{Casson invariant!Gauss diagram formula for ---}
The Casson-Walker invariant
as an invariant of 2-component links
is studied in \cite{KiLi97}.

If we would obtain a Gauss diagram formula 
for the Casson-Walker-Lescop invariant,
then a spin refinement of it (of Problem \ref{prob.spinCWL})
would be obtained
by decorating the Gauss diagram formula
by characteristic sublinks,
noting that
the spin structures on 
the 3-manifold obtained by surgery along a framed link $L$
can be presented by characteristic sublinks of $L$
(see \cite{KM_inv}).
\end{rem}

\subsection{Finite type invariants}

A link in an integral homology 3-sphere is called
{\it algebraically-split}
if the linking number of any pair of its components vanishes,
and is called {\it boundary}
if all its components bound disjoint surfaces.
A framed link is called {\it unit-framed}
if the framings of its components are $\pm1$.
Let $\M$ be the set of 
(homeomorphism classes of) oriented integral homology 3-spheres,
and let $R$ be a commutative ring with $1$.
For an algebraically-split unit-framed link $L$
in an integral homology 3-sphere $M$, we put
$$
[M,L] = \sum_{L' \subset L} (-1)^{\# L'} M_{L'} \in R \M, 
$$
where the sum runs over all sublinks $L'$ of $L$,
and $\# L'$ denotes the number of components of $L'$,
and $M_{L'}$ denotes the 3-manifold
obtained from $M$ by surgery along $L'$.
Let $\cF_d^{\rm as}(R \M)$ \cite{Ohtsuki_f}
(resp. $\cF_d^{\rm b}(R \M)$ \cite{Garoufalidis_ftI})
denote the submodule of $R \M$ spanned by $[M,L]$ such that
$M$ is an integral homology 3-sphere
and $L$ is a unit-framed algebraically-split link $L$
with $d$ components in $M$
(resp. a unit-framed boundary link $L$ in $M$).
Let $\cF_d^{\rm Y}(R \M)$ \cite{GGP} denote
the submodule of $R \M$ spanned by 
$[M,G]$
such that 
$M$ is an integral homology 3-sphere and 
$G$ is a collection of $d$ disjoint Y-graphs (see Figure \ref{fig.surgeryY})
in $M$,
where $[M,G]$ is defined similarly as $[M,L]$ (see \cite{GGP}).\footnote{\footnotesize 
$\cF_d^{\rm Y}(R \M)$ can alternatively be defined
by using blinks \cite{GaLe_blink}; see \cite{GGP}.}
A homomorphism $v: R \M \to R$ is called 
a {\it finite type invariant}\index{finite type invariant!--- of homology spheres} 
of $\cF_\star^{\rm as}$-degree $d$
(resp. $\cF_\star^{\rm b}$-degree $d$, or $\cF_\star^{\rm Y}$-degree $d$)
if $v$ vanishes on $\cF_{d+1}^{\rm as}(R \M)$
(resp. $\cF_{d+1}^{\rm b}(R \M)$, or $\cF_{d+1}^{\rm Y}(R \M)$).
It is known \cite{GaOh98} that
$$
\cF_{3d}^{\rm as}(\Q \M) = \cF_{3d-1}^{\rm as}(\Q \M)
= \cF_{3d-2}^{\rm as}(\Q \M)
$$
and that there is an isomorphism
$$
\cA(\emptyset;\Q)^{(d)} \longrightarrow
\cF_{3d}^{\rm as}(\Q \M) / \cF_{3d+3}^{\rm as}(\Q \M)
$$
between vector spaces \cite{GaOh98,Le_uf}.
It is known \cite{GGP} that
\begin{align*}
& \cF_d^{\rm b}(\Z \M) \supset \cF_{2d}^{\rm Y}(\Z \M), \qquad
\cF_{3d}^{\rm as}(\Z \M) \supset \cF_{2d}^{\rm Y}(\Z \M), \\
& \cF_{3d}^{\rm as}(R \M) = \cF_d^{\rm b}(R \M) = \cF_{2d}^{\rm Y}(R \M), \\*
& \cF_{2d-1}^{\rm Y}(R \M) = \cF_{2d}^{\rm Y}(R \M)
\end{align*}
if $1/2 \in R$.

\subsubsection{Torsion and finite type invariants}\index{finite type invariant!--- of homology spheres!torsion and ---} 

\begin{conj}
$\cF^{\rm as}_d(\Z\M) / \cF^{\rm as}_{d+1}(\Z\M)$ 
(resp. $\cF^{\rm b}_d(\Z\M) / \cF^{\rm b}_{d+1}(\Z\M)$)
is torsion free\index{torsion!--- free}   
for each $d$.
\end{conj}

\begin{rem}[\rm (K. Habiro)] 
The group $\cF^{\rm Y}_d(\Z\M) / \cF^{\rm Y}_{d+1}(\Z\M)$
has 2-torsion for each $d>0$.
\end{rem}

\begin{conj}
\label{conj.Aphi_tf}
$\cA(\emptyset;\Z)$ is\index{Jacobi diagram!space of ---!torsion of ---}   
torsion free.\index{torsion!--- free}  
\end{conj}

\subsubsection{Do finite type invariants distinguish homology 3-spheres?}

\begin{conj}
\label{conj.Fdist_hS3}
Finite type invariants\index{finite type invariant!--- of homology spheres!strength of ---}  
distinguish integral homology\break 3-spheres.
(See Conjecture \ref{conj.LMOdist_hS3}.)
\end{conj}

\subsubsection{Dimensions of spaces of finite type invariants}

A finite type invariant
$v$ is called {\it primitive}
if $v( M_1 \# M_2 ) = v(M_1) + v(M_2)$
for any integral homology 3-spheres $M_1$ and $M_2$.
We denote by $\cA(\emptyset;R)\conn$ the submodule of $\cA(\emptyset;R)$
spanned by Jacobi diagrams with connected trivalent graphs.
As a graded vector space
$\cA(\emptyset;\Q)$ is isomorphic to
the symmetric tensor algebra of $\cA(\emptyset;\Q)\conn$.

\begin{prob}
\label{prob.dimF}
Determine the dimension
of the space of primitive 
finite type invariants\index{finite type invariant!--- of homology spheres!dimension of ---}
of integral homology 3-spheres of each degree $d$.
Equivalently,
determine the dimension of the space\index{Jacobi diagram!space of ---!dimension of ---}  
$\cA(\emptyset;\Q)\conn^{(d)}$
for each $d$.
\end{prob}
\eject

\begin{table}[ht!]
\begin{center}
\begin{tabular}{|c||c|c|c|c|c|c|c|c|c|c|c|c|c|c|c|}
\hline
\small $d$ &\small 0 &\small 1 &\small 2 &\small 3 &\small 4 &\small 5 &\small 6 &\small 7 &\small 8 &\small 9 &\small 10 \\
\hline\hline
\small prime diag. &\small 
0 &\small 1 &\small 0 &\small 0 &\small 1 &\small 0 &\small 1 &\small 1 
&\small 1 &\small 1&\small 2 \\
\hline
\small dim $\cA(\emptyset)\conn^{(d)}$ &\small 
0 &\small 1 &\small 1 &\small 1 &\small 2 &\small 2 &\small 3 &\small 4 
&\small 5 &\small 6 &\small 8 \\
\hline
\small dim $\cA(\emptyset)^{(d)}$ &\small 
1 &\small 1 &\small 2 &\small 3 &\small 6 &\small 9 &\small 16 &\small 25 
&\small 42 &\small 65 &\small 105 \\
\hline
\end{tabular}
\end{center}
\vspace{0.2pc}
\begin{center}
\begin{tabular}{|c||c|c|c|c|}
\hline
\small $d$ &\small 11 &\small 12 &\small 13 &\small 14 \\
\hline\hline
\small prime diag. &\small 1&&& \\
\hline
\small dim $\cA(\emptyset)\conn^{(d)}$
&\small 9 &\small \!$\ge \! 11$\! &\small \!$\ge \! 13$\! &\small \!$\ge \! 15$\! \\
\hline
\small dim $\cA(\emptyset)^{(d)}$ 
&\small 161 &\small \!$\ge \! 254$\! &\small \!$\ge \! 386$\! &\small \!$\ge \! 595$\! \\
\hline
\end{tabular}
\end{center}
\caption{\label{tbl.dimF}
Some dimensions for Problem \ref{prob.dimF}}
\end{table}

\begin{rem}
$\cA(\emptyset;\Q)\conn^{(d)}$ is isomorphic to $\cB\conn^{(d+1,2)}$
mentioned in a remark of Problem \ref{prob.dimV},
by the isomorphism
taking a trivalent graph to a uni-trivalent graph
obtained from the trivalent graph
by cutting a middle point of an edge.
Hence, the dimension of $\cA(\emptyset;\Q)\conn^{(d)}$
is equal to the dimension $\beta_{d+1,2}$ of $\cB\conn^{(d+1,2)}$.
Therefore, we obtain the row of $\cA(\emptyset;\Q)\conn^{(d)}$ 
in Table \ref{tbl.dimF} from a column of Table \ref{tbl.bdu}.
\end{rem}

\begin{rem}
$\cA(\emptyset)\conn$ is an algebra
with the product given by connected sum of Jacobi diagrams.
Let us look for prime diagrams
with respect to the connected sum;
they generate the algebra $\cA(\emptyset)\conn$.
By the AS and IHX relations,
we can remove a triangle, and
we can break a polygon with odd edges.
Hence, prime diagrams are given by
\begin{align*}
& p_1=\pict{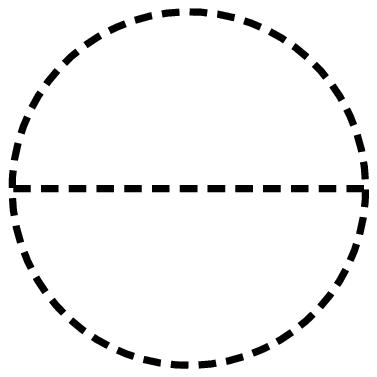}{1.5cm}, \qquad
p_4=\pict{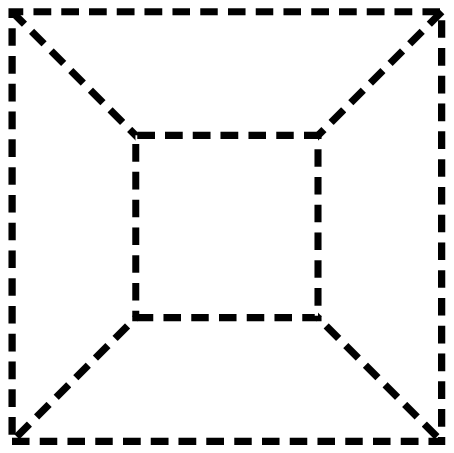}{1.7cm}, \qquad
p_6=\pict{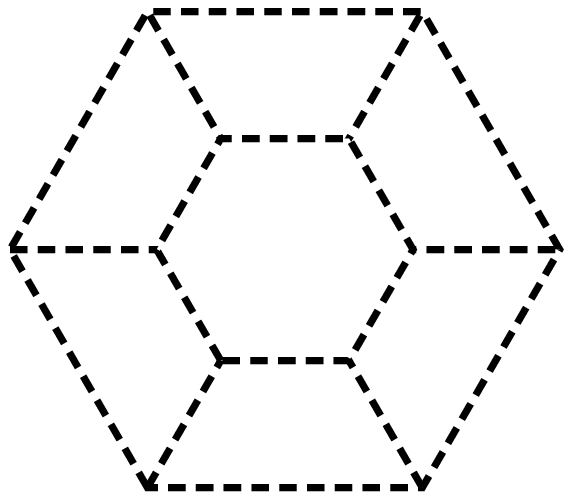}{1.9cm}, \\*
& p_7=\pict{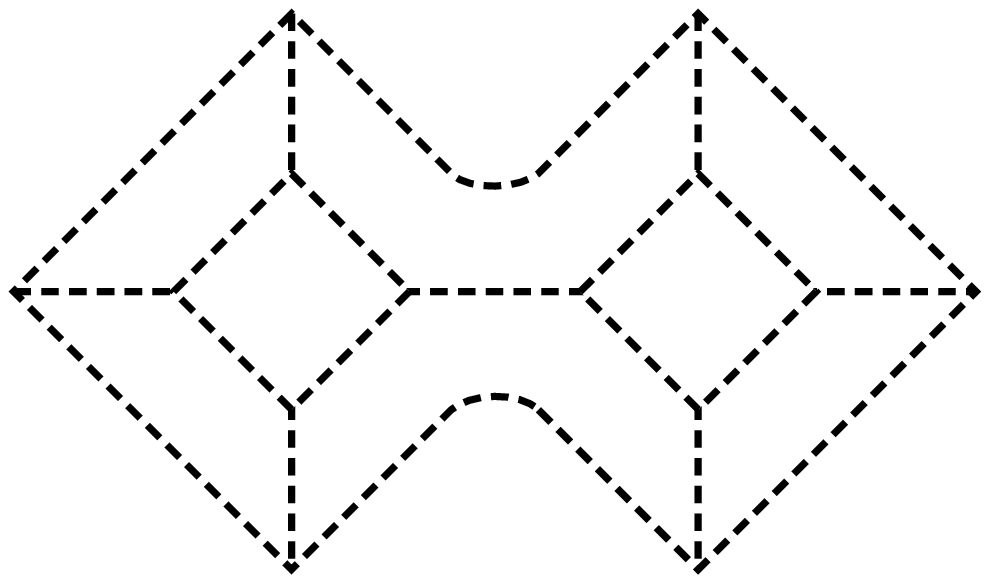}{2cm}, \qquad
p_8=\pict{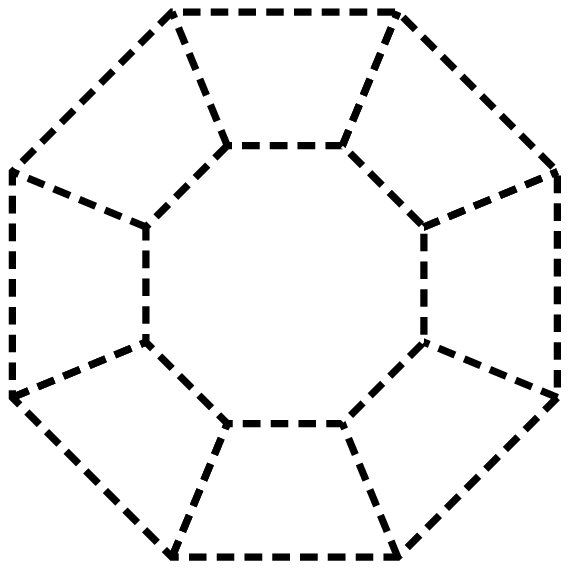}{2.3cm}, \qquad \cdots.
\end{align*}
They have the relation $p_1 p_7 = p_4^2$,
since $p_1 p_7 = p_1 (x_3 p_4) = (x_3 p_1) p_4 = p_4^2$,
where $x_3$ is the element of Vogel's algebra $\Lambda$ given 
in (\ref{eq.Lambda_gen}) below,
which acts on $\cA(\emptyset)\conn$.
It is a problem to find a complete list of
generators and relations of the algebra $\cA(\emptyset)\conn$.
\end{rem}

\begin{rem}
$\cA(\emptyset;\Q)\conn$ is a $\Lambda$-algebra,
where $\Lambda$ is Vogel's algebra given below,
whose generators and relations have been known
in degree $\le 10$;
see a remark on Problem \ref{prob.Lambda_gen}.
It is a problem to find
generators and relations of $\cA(\emptyset;\Q)\conn$ as a $\Lambda$-algebra.
\end{rem}

\begin{update}
The prime diagrams of degree $\le 11$ are given in \cite{ChDuKa}.
\end{update}

\subsubsection{Vogel's algebra}

Vogel's algebra \cite{Vogel_alg}
is defined as follows.
For fixed 3 points,
we denote by $\cA(\mbox{3 points})\conn$
the module over $\Q$
spanned by vertex-oriented connected uni-trivalent graphs
whose univalent vertices are
the fixed 3 points
subject to the AS and IHX relations.
The symmetric group ${\frak S}_3$ acts on $\cA(\mbox{3 points})\conn$
by permutation of 3 points.
The module $\Lambda$ is defined to be 
the submodule of $\cA(\mbox{3 points})\conn$
consisting of all elements $u$ satisfying that
$\sigma(u) = \mbox{sgn}(\sigma) \cdot u$
for any $\sigma \in {\frak S}_3$.
It is well defined to insert $u \in \Lambda$ 
in a vertex-oriented trivalent vertex as
$$
\pict{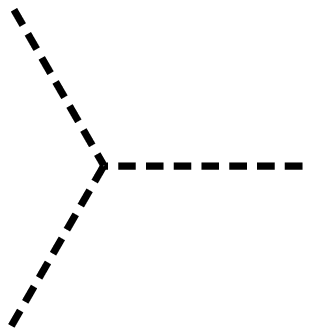}{1.65cm}
\longmapsto
\pict{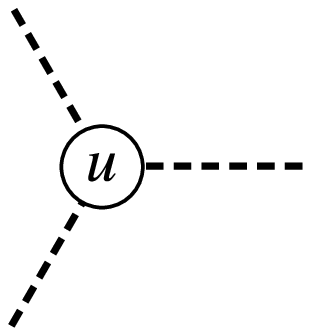}{1.65cm}.
$$
Moreover, this insertion is independent, modulo the AS and IHX relations,
of a choice of a trivalent vertex as follows.
By the AS and IHX relations,
\begin{align*}
\pict{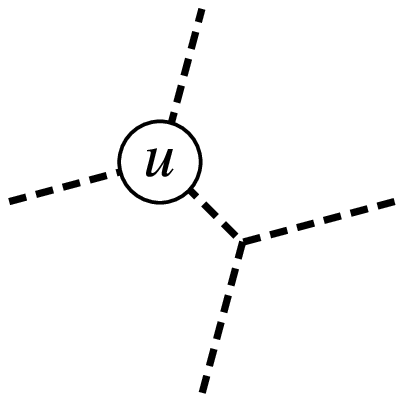}{1.85cm} - \pict{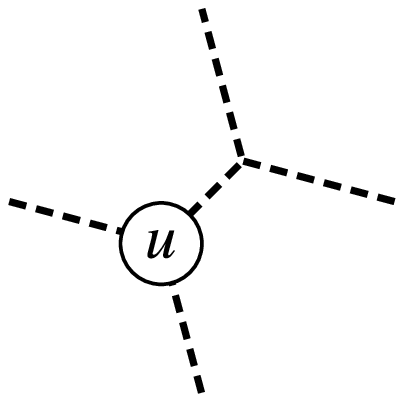}{1.85cm}
& = \pict{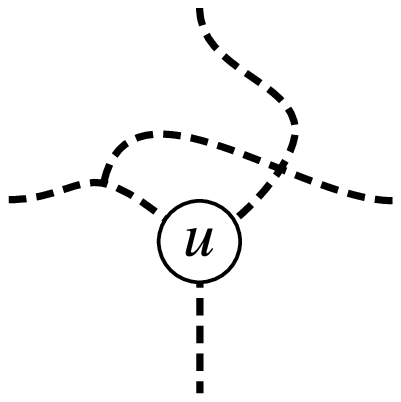}{1.85cm} = - \pict{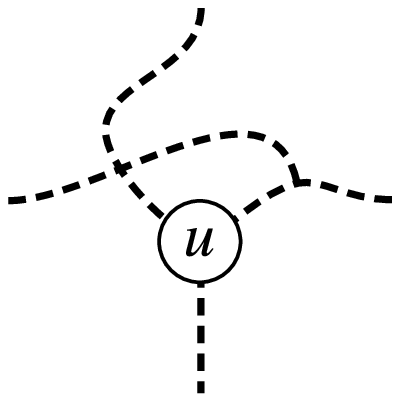}{1.85cm} \\*
& = \pict{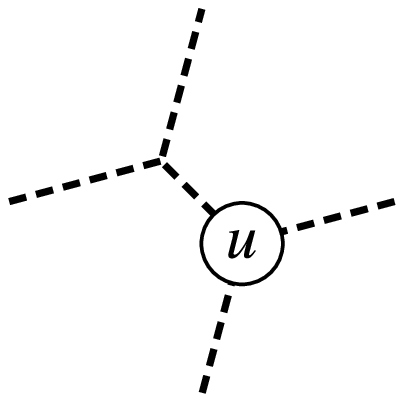}{1.85cm} - \pict{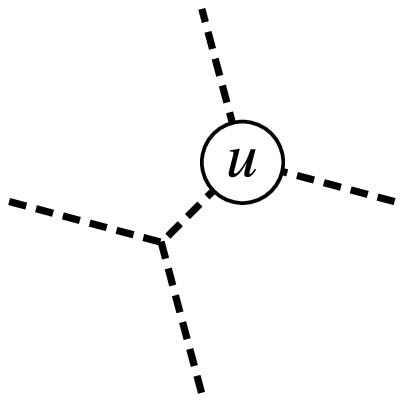}{1.85cm},
\end{align*}
where the middle equality is derived from the anti-symmetry of $u$.
By the $\pi/4$ and $-\pi/4$ rotations of the above formula,
we have that
\begin{align*}
\pict{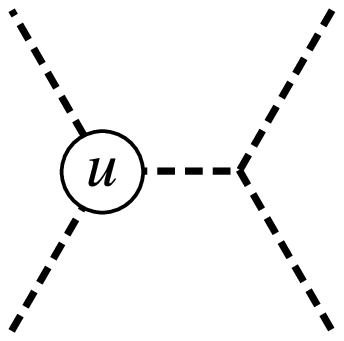}{1.65cm} - \pict{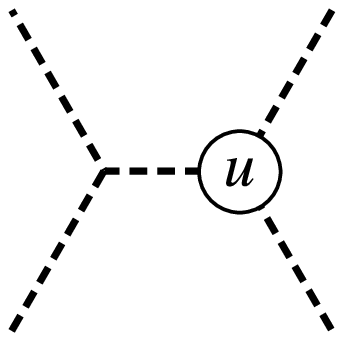}{1.65cm}
& = \pict{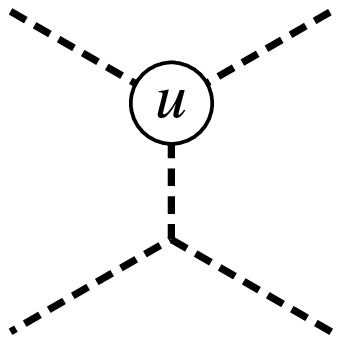}{1.65cm} - \pict{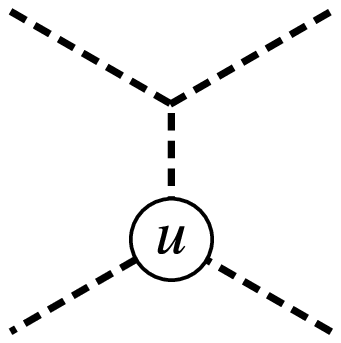}{1.65cm} \\*
& = \pict{v/vog2}{1.65cm} - \pict{v/vog1}{1.65cm}.
\end{align*}

Hence, the left hand side of the above formula is equal to $0$.
This implies that
the insertion of $u$ is independent of
a choice of a trivalent vertex.
The module $\Lambda$ is an algebra,
called {\it Vogel's algebra},\index{Vogel's algebra}
whose product of $x,y \in \Lambda$ is defined to be
the element of $\Lambda$ obtained by inserting $x$ 
in a trivalent vertex of $y$.
It is a commutative algebra.
Some generators of $\Lambda$ in low degrees are given by
\begin{equation}
\label{eq.Lambda_gen}
1 = \pict{v/vog-g1}{1.65cm}, \qquad
t = \pict{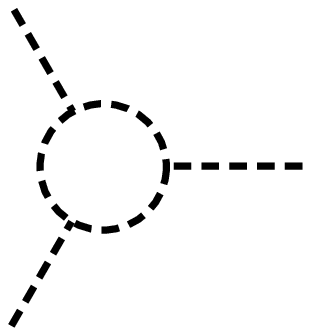}{1.65cm}, \qquad
x_3 = \pict{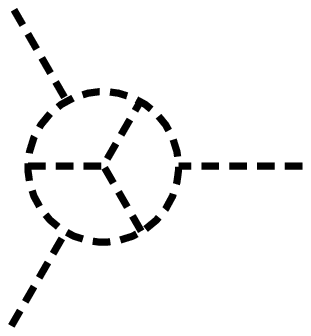}{1.65cm},
\end{equation}
and further,
$$
x_n = \pict{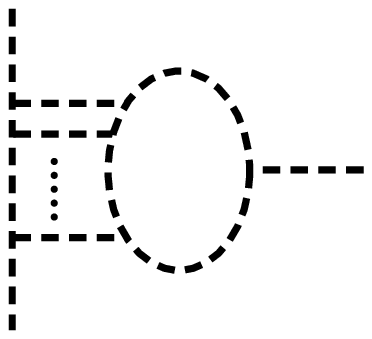}{1.65cm},
$$
having $n$ horizontal lines between the vertical line and the circle.
It is known that the even $x_n$'s can be presented by odd $x_n$'s.

\begin{prob}
\label{prob.Lambda_gen}
Describe Vogel's algebra $\Lambda$,\index{Vogel's algebra!presentation of ---}
say, by giving complete sets of generators and relations of $\Lambda$.
\end{prob}

\begin{rem}
Vogel \cite{Vogel_grenoble} conjectured that
the homomorphism $\varphi: R_0 \to \Lambda$ given in \cite{Vogel_grenoble}
was bijective,
where $R_0$ is the subalgebra of a polynomial algebra in 3 variables,
generated by elements given in \cite{Vogel_grenoble}.
As mentioned in \cite{Vogel_grenoble}, 
$\varphi$ has been known to be
bijective in degree $\le 10$,
and injective in degree $\le 15$.

Recently (in June, 2001),
Vogel found a polynomial in $R_0$ whose image in $\Lambda$ vanishes;
this implies that $\varphi$ is not injective.
Surjectivity of $\varphi$
(which implies that $\Lambda$ is generated by $t, x_3, x_5, x_7, \cdots$)
is still an open problem.

Vogel \cite{Vogel_0div} further found
a divisor of zero in $\Lambda$.
It is given as follows.
Putting
$$
U = \pict{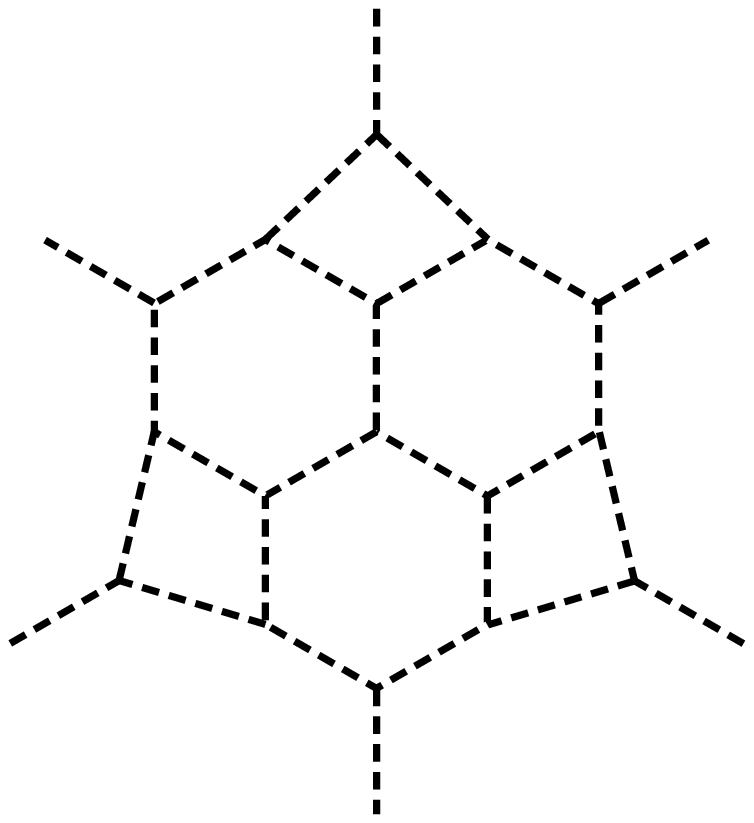}{3.5cm},
$$
we define $W \in \cA(\emptyset)$ and $\lambda \in \Lambda$ by
$$
W = \sum_\sigma \mbox{sign}(\sigma) \pict{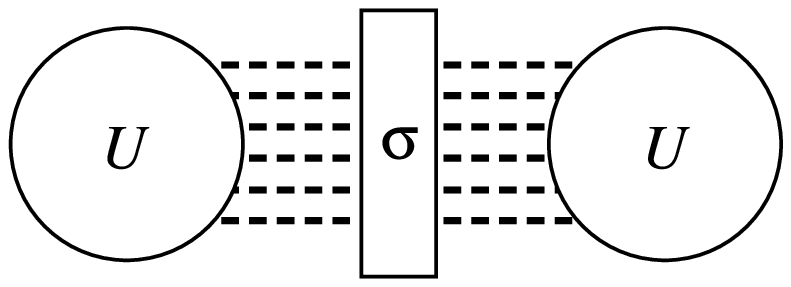}{1.5cm}
= \pict{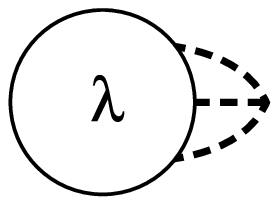}{1cm},
$$
where the sum of $\sigma$ runs over all permutations $\sigma \in {\frak S}_6$,
and $\lambda$ is obtained from $W$ by removing a neighborhood 
of a trivalent vertex.
Vogel showed that
$t \lambda = 0 \in \Lambda$ and $\lambda \ne 0 \in \Lambda$.
\end{rem}

\subsubsection{Other problems}

\begin{prob}
Find a constructive combinatorial presentation of each 
finite type invariant\index{finite type invariant!--- of homology spheres!constructive presentation of ---}   
of integral homology 3-spheres, and, in particular, 
of the Casson invariant,\index{Casson invariant!--- by counting configurations}
by localizing 
configuration space integrals.\index{configuration space!--- integral!localization of ---} 
\end{prob}

\begin{rem}
The perturbative expansion of the path integral
of the Chern-Simons field theory
suggests that
each Vassiliev invariant of knots can be obtained
as a mapping degree of a certain map on a configuration space,
whose localization deduces
a Gauss diagram formula of this Vassiliev invariant;
see comments before Problem \ref{prob.khs}.
In the 3-manifold case
G. Kuperberg and D. Thurston \cite{Kuperberg-Thurston} 
gave a presentation of each finite type invariant
by using configuration space integrals,
whose localization might deduce a combinatorial formula,
similarly as a Gauss diagram formula.
It would be a difficult point of such localization
to deal with ``hidden strata'' (anomaly faces).
\end{rem}

\begin{prob}[J. Roberts]  
\label{prob.roberts15}
What is the space of $3$-manifolds? 
\end{prob}

\begin{rem}[{\rm (J. Roberts)}]  
Vassiliev invariants are usually characterised in purely
combinatorial terms, but it is worth remembering that Vassiliev was
led to this definition by considering the natural stratification of
the space of smooth maps $S^1 \rightarrow \R^3$. The combinatorial
theory of 
finite type invariants\index{finite type invariant!--- of homology spheres!--- by space of 3-manifolds}   
of homology spheres is now equally
well-developed but there remains no natural justification for
considering the relations introduced by Ohtsuki, other than that these
turn out to interact very well with the perturbative expansion of the
Witten invariants. One would like to find a stratified space of
integer homology spheres, in which crossing a codimension $1$ stratum
corresponds to doing $\pm 1$ surgery on a knot. Now the space of
smooth maps $f: S^{n+3} \rightarrow S^n$ is a natural choice for a ``space of
framed $3$-manifolds'', via the Pontrjagin-Thom construction (take the
preimage of a fixed point in $S^n$). But this space gives the wrong
filtration,\eject and it's not clear how to alter it to implement (for
example) constraints on the homology of the preimages. See Shirokova
\cite{Shi00}.
\end{rem}

\subsection{Goussarov-Habiro theory}

\subsubsection{Goussarov-Habiro theory for 3-manifolds}

Related to finite type invariants of 3-manifolds,
equivalence relations among 3-manifolds
have been studied by Goussarov \cite{Gou_kg,Gou_fn} and Habiro \cite{Habiro_GT},
which is called the Goussarov-Habiro theory for 3-manifolds.\index{Goussarov-Habiro theory!--- for 3-manifolds}
These equivalence relations are helpful for us
to study structures of the set of 3-manifolds.

The\index{equivalence relation!$Y_d$-equivalence}\index{equivalence relation!$d$-equivalence}     
{\it $Y_d$-equivalence}\footnote{\footnotesize 
The $Y_d$-equivalence is also called
the {\it $(d-1)$-equivalence} (due to Goussarov) in some literatures.}
among oriented 3-manifolds is the equivalence relation
generated by either of the following relations,
\begin{itemize}
\item[\rm(1)]
surgery on a tree clasper with $d$ trivalent vertices \cite{Habiro_GT},
\item[\rm(2)]
Goussarov's $d$-variation 
(which generates Goussarov's notion of $(d-1)$-equivalence)
\cite{Gou_kg,Gou_fn},
\item[\rm(3)]
surgery by an element in the $d$th lower central series subgroup of
the Torelli group of a compact connected surface.
\end{itemize}
It is known \cite{Habiro_GT} that
these relations generate the same equivalence relation among {\ZHSs}.
Two closed 3-manifolds $M$ and $M'$ are
$Y_1$-equivalent if and only if
there is an isomorphism $H_1(M;\Z) \to H_1(M';\Z)$
which induces an isomorphism between their linking pairings \cite{Matveev87}.

It is known \cite{Habiro_GT} that
$\{ \mbox{integral homology 3-spheres (\ZHSs)} \} / \simsubY{2} \cong \Z/2\Z$
and that $\{ \ZHSs \} / \simsubY{3} \cong \Z$,
which deduce the Rokhlin and Casson invariants
respectively.
Further, it is known \cite{Habiro_GT} that
$\{ M \simsubY{2d-1} \! S^3 \} / \simsubY{2d} = 0$ for $d>1$ and that
there exists a natural surjective homomorphism
\begin{equation}
\label{eq.AtoMsim}
\cA(\emptyset;\Z)\conn^{(d)} \longrightarrow 
\{ M \ \simsubY{2d} S^3 \} / \simsubY{2d+1} 
\end{equation}
such that the tensor product of this map and $\Q$ is an isomorphism.
In particular, 
$\{ M \ \simsubY{2d} S^3 \} / \simsubY{2d+1}$ forms
an abelian group
with respect to the connected sum of \ZHSs,
and hence, so does $\{ \ZHSs \} / \simsubY{2d+1}$.

\begin{conj}
The map (\ref{eq.AtoMsim}) is an isomorphism.
\end{conj}
\eject

This conjecture might be reduced to Conjecture \ref{conj.Aphi_tf} and
the following conjecture.

\begin{conj}
\label{conj.ZHS_Yd}
$\{ M \ \simsubY{2d} S^3 \} / \simsubY{2d+1}$ is 
torsion free\index{torsion!--- free}   
for each $d$.
\end{conj}

\begin{rem}
Conjecture \ref{conj.Aphi_tf} implies this conjecture,
since the surjective homomorphism (\ref{eq.AtoMsim})
gives a $\Q$-isomorphism.
\end{rem}

\begin{rem}[\rm (K. Habiro)] 
It is also a problem to describe the graded set
$\{ M \ \simsubY{d} M_0 \} / \simsubY{d+1}$
for an arbitrarily given 3-manifold $M_0$.
For $d=0$, the quotient set $\{ \mbox{3-manifolds} \} / \simsubY{1}$
can be identified with the set of isomorphism classes
of $H_1(M;\Z)$ and their linking pairings (as mentioned above).
For $d>0$, there is a surjective map to this graded set
from a certain module of Jacobi diagrams
(subject to the AS and IHX relations).
\end{rem}

\begin{prob}[T. Ohtsuki] 
\label{prob.product_in_M}
Define a product $M_1 \circ M_2$ of 
integral homology 3-spheres $M_1$ and $M_2$
which is related, by (\ref{eq.AtoMsim}),
to the product of Jacobi diagrams given by their connected sum.
\end{prob}

\begin{rem}
$\cA(\emptyset)\conn$ is an algebra
with the product given by connected sum of Jacobi diagrams.
The connected sum of Jacobi diagrams on $\emptyset$ is well defined
by the AS and IHX relations.
The sum of $\cA(\emptyset)\conn$ corresponds, by (\ref{eq.AtoMsim}),
to the connected sum of integral homology 3-spheres.
The problem is to define a product among integral homology 3-spheres
corresponding to the product of $\cA(\emptyset)\conn$ by (\ref{eq.AtoMsim}).
\end{rem}

It is known \cite{Gou_kg,Habiro_GT} that
two integral homology 3-spheres
$M$ and $M'$ are $Y_{d}$-equivalent if and only if
$v(M) = v(M')$ for any $A$-valued finite type invariant\footnote{
For an abelian group $A$,
a homomorphism $v:\Z\M \to A$ is called a {\it finite type invariant}
of $\cF_\star^{\rm Y}$-degree $d$
if $v$ vanishes on $\cF_{d+1}^{\rm Y}(\Z\M)$.}
$v$ of $\cF_\star^{\rm Y}$-degree $< d$
for any abelian group $A$.
In fact, a natural quotient map
$\{ \ZHSs \} \to \{ \ZHSs \} / \simsubY{d}$
is a finite type invariant of $\cF_\star^{\rm Y}$-degree $< d$, 
which classifies $Y_{d}$-equivalence classes of integral homology 3-spheres.

For an oriented compact surface $F$,
a {\it homology cylinder} over $F$
is a homology $F \times I$ 
whose boundary is parameterized by $\partial (F \times I)$.

\begin{conj}[M. Polyak, {\rm see \cite[``Theorem 4'']{Gou_fn}}] 
\label{conj.hc_Yd}
Let $F$ be an oriented compact surface.
Two homology cylinders $C$ and $C'$ over $F$ are 
$Y_{d}$-equivalent\index{equivalence relation!$Y_d$-equivalence}  
if and only if
$v(C) = v(C')$ for any $A$-valued\index{finite type invariant!--- of homology cylinders}
finite type invariant $v$ 
of $\cF_\star^{\rm Y}$-degree $< d$
for any abelian group $A$.
\end{conj}

\begin{rem}[\rm (M. Polyak)] 
The corresponding assertion for closed 3-manifolds does not hold;
note that $\{ \mbox{closed 3-manifolds} \} / \simsubY{d}$ 
does not (naturally) form a group.
Recall that $\{ \ZHSs \} / \simsubY{d}$ forms an abelian group,
which guarantees the corresponding assertion for \ZHSs,
as mentioned above.
The set $\{ \mbox{homology cylinders}$ $\mbox{on $F$} \} / \simsubY{d}$ 
forms a group with respect to the composition of homology cylinders, 
though it is not abelian.
\end{rem}

%
%

\subsubsection{Goussarov-Habiro theory for spin and spin${}^c$ 3-manifolds}
\label{sec.GHspin}

\renewcommand{\thefootnote}{\fnsymbol{footnote}}
\footnotetext[0]{The first version of Section \ref{sec.GHspin}
was written by T. Ohtsuki, following a report of F. Deloup.
Based on it, F. Deloup wrote this section.}
\renewcommand{\thefootnote}{\arabic{footnote}}

As\index{Goussarov-Habiro theory!--- for spin 3-manifolds}\index{Goussarov-Habiro theory!--- for spin${}^c$ 3-manifolds} 
shown in \cite{Massuyeau},
we have a natural spin (resp. spin${}^c$ structure)
on the 3-manifold
obtained from a spin (resp. spin${}^c$) 3-manifold
by surgery along a Y graph
(or a tree clasper).
We define the {\it $Y_d^s$-equivalence} 
({\it spin $Y_d$-equivalence})\index{equivalence relation!$Y_d$-equivalence!spin ---}
(resp. {\it $Y_d^c$-equivalence} 
({\it spin${}^c$ $Y_d$-equivalence}))\index{equivalence relation!$Y_d$-equivalence!spin${}^c$ ---}
to be the equivalence relation
among spin (resp. spin${}^c$) 3-manifolds
given by the $Y_d$-equivalence.
It is known \cite{Massuyeau} that the quotient set
$\{ \mbox{spin closed 3-manifolds} \}/ \!\!\underset{{}^{Y_{1}^s}}{\sim}$
can be identified with the isomorphism classes of
pairs of $H_1(M;\Z)$ and
certain quadratic forms $\phi_{M,\sigma} : \mbox{Tor} H_1(M;\Z) \to \Q/\Z$,
or equivalently,
the isomorphism classes of triples of
$H_1(M;\Z)$ and linking pairings
$\lambda_M: \big( \mbox{Tor} H_1(M;\Z) \big)^{\otimes 2} \to \Q/\Z$
and the mod 8 reduction of the Rokhlin invariant $\mu(M,\sigma)$.
Further, it is known \cite{DeMa_q}
the quotient set
$\{ \mbox{spin${}^c$ closed 3-manifolds} \}/ 
\!\!\underset{{}^{Y_{1}^c}}{\sim}$
can be identified with
the set of the isomorphism classes of pairs of
$H_1(M;\Z)$ and certain quadratic forms $q_\sigma$.
This set would be well described 
by the classification of the following problem.

\begin{prob}[F. Deloup] 
\label{prob.q_sigma}
Classify the monoid (for orthogonal sum)
of isomorphism classes of quadratic forms $q_\sigma$.
\end{prob}

\begin{rem}
The quotient set
$\{ \mbox{closed 3-manifolds} \}/\simsubY{1}$
can be identified with the set of 
the isomorphism classes of pairs of $H_1(M;\Z)$ and linking pairings.
This set can be well described by 
the classification of linking pairings given in \cite{KaKo80}.
\end{rem}

\begin{prob}[G. Massuyeau] 
\label{prob.spin_Yd}
Describe the quotient set \newline
$\{ \mbox{spin closed 3-manifolds} \}/ \!\!\underset{{}^{Y_{d}^s}}{\sim}$,
in particular, 
for $d=2,3$.\index{Goussarov-Habiro theory!--- for spin 3-manifolds} 
\end{prob}



\begin{prob}[F. Deloup, G. Massuyeau] 
\label{prob.spinc_Yd}
Describe the quotient set \newline
$\{\mbox{spin${}^c$ closed 3-manifolds}\}/ \!\!\underset{{}^{Y_{d}^c}}{\sim}$,
in particular, 
for $d=2,3$.\index{Goussarov-Habiro theory!--- for spin${}^c$ 3-manifolds} 
\end{prob}

\begin{rem}
There is a unique spin (resp spin${}^c$) structure on a \ZHS.
Hence, $\{ \mbox{spin \ZHSs} \} / \!\!\underset{{}^{Y_{d}^s}}{\sim}$
(resp. $\{ \mbox{spin${}^c$ \ZHSs} \} / \!\!\underset{{}^{Y_{d}^s}}{\sim}$)
is equal to $\{ \ZHSs \} / \simsubY{d}$.
This quotient set can be described by
Jacobi diagrams (see Conjecture \ref{conj.ZHS_Yd}).
\end{rem}

\begin{rem}
The above two problems are related to
spin and spin${}^c$ refinements of the Casson-Walker-Lescop invariant;
see Problem \ref{prob.spinCWL}.
\end{rem}

Deloup and Massuyeau \cite{DeMa_q} obtained
a complete system of invariants for quadratic functions 
on finite abelian groups
which involves the Gauss-Brown invariant
$\gamma(q) = \sum_{x \in G} e^{2\pi \sqrt{-1} q(x)}$
of a quadratic form $q$.
In the case $q_\sigma$ comes from a usual spin structure, 
$q_\sigma$ is homogeneous\footnote{\footnotesize 
A quadratic function $q$ is a a map such that
$q(x+y)-q(x)-q(y)$ is bilinear in $x$ and $y$.
It is called {\it homogeneous} 
if $q(n x) = n^{2}q(x)$ for any $n \in \Z$ and $x \in G$. 
In fact, there is a canonical map
$\sigma \mapsto q_\sigma$ from spin${}^c$ structures to quadratic functions 
and $q_\sigma$ is homogeneous 
if and only if $\sigma$ actually comes from a spin structure.
Note that not all spin${}^c$ structures come from spin structures.}
and the argument of $\gamma(q_\sigma)$ is 
just the mod 8 reduction of the Rokhlin invariant. 
(Here we take the classical Rokhlin invariant of a spin structure on $M$ 
to be the signature mod 16 of an oriented smooth simply-connected $4$-manifold
bounded by $M$.)
Thus, in general, ${\rm arg}\, \gamma(q_\sigma) \in \Q/\Z$
may be viewed as mod 8 generalization of Rokhlin invariant
for spin${}^c$ structures.
In the context of spin Goussarov-Habiro theory,
Massuyeau proved that
the Rokhlin invariant is a finite type invariant of degree 1.
This suggests the following question.

\begin{quest}[F. Deloup] 
\label{quest.lift_g}
Is there a lift of ${\rm arg}\, \gamma(q_\sigma)$
to a mod 16 invariant?
This would give a finite type 
invariant\index{finite type invariant!--- of spin${}^c$ 3-manifolds}
of degree 1 in the spin${}^c$ Goussarov-Habiro 
theory.\index{Goussarov-Habiro theory!--- for spin${}^c$ 3-manifolds} 
\end{quest}

\newpage

\section{The LMO invariant}

The LMO invariant\index{LMO invariant}  $\LMO(M) \in \cA(\emptyset)$ 
of closed oriented 3-manifolds was introduced in \cite{LMO}.
The LMO invariant of rational homology 3-spheres
was reformulated by Aarhus integral \cite{BGRT_A}.
The LMO invariant
is a universal perturbative invariant of
rational homology 3-spheres
(see \cite{Ohtsuki_rec,BGRT_A,Ohtsuki_book}),
and a universal finite type invariant
of integral homology 3-spheres \cite{Le_uf}.

\subsection{Calculation of the LMO invariant}

\begin{prob}
\label{prob.calLMO}
For each rational homology 3-sphere $M$,
calculate $\LMO(M)$ for all degrees.\index{LMO invariant!calculation of ---} 
\end{prob}

{\bf Remark}\qua
Bar-Natan and Lawrence \cite{BaLa} showed
a rational surgery formula for the LMO invariant.
By using it, they obtained
\begin{equation}
\label{eq.LMO_lens}
\hLMO\big( L(p,q) \big) = 
\langle \Omega_x, \Omega_x^{-1} \Omega_{x/p} \rangle_x
\exp \frac{-s(q,p)}{48} \theta
\end{equation}
for the lens space\index{lens space} $L(p,q)$ of type $(p,q)$,
where $s(q,p)$ denotes the Dedekind sum.
For the notation $\langle \Omega_x, \Omega_x^{-1} \Omega_{x/p} \rangle_x$
see \cite{BaLa}.

{\bf Remark}\qua
The degree 1 part of $\LMO(M)$ is given by the
Casson-Walker invariant of $M$ \cite{LMO}.
Further,
the degree $\le d$ part of $\LMO(M)$ of
integral homology 3-spheres
are given by finite type invariants of degree $\le d$.
Hence, it is algorithmically possible
to calculate the degree $\le d$ part of
$\LMO(M)$ of an integral homology 3-sphere for each $d$.
It is meaningful to calculate $\LMO(M)$ for all degrees.

{\bf Remark}\qua
It is meaningful to calculate $\LMO(M)$
when $M$ is a rational homology 3-sphere.
Otherwise, it is known that
$\LMO(M)$ can be given by some ``classical'' invariants.
When $b_1(M)=1$, 
the value of $\LMO(M)$ can be presented by using
the Alexander polynomial of $M$ \cite{GH,Lieberum}.
When $b_1(M)=2$, 
the value of $\LMO(M)$ can be presented by using
the Casson-Walker-Lescop invariant of $M$ \cite{HB}.
When $b_1(M)=3$,
the value of $\LMO(M)$ can be presented by using
the cohomology ring of $M$ \cite{Habegger_c}.
When $b_1(M)>3$,
we always have that $\LMO(M)=1$ \cite{Habegger_c}.

\subsection{Does the LMO invariant distinguish integral homology 3-spheres?}

\begin{conj}
\label{conj.LMOdist_hS3}
The LMO invariant distinguishes 
integral homology 3-spheres.\index{LMO invariant!strength of ---} 
(See Conjecture \ref{conj.Fdist_hS3}.)
\end{conj}
\eject

\begin{rem}
Bar-Natan and Lawrence \cite{BaLa} showed 
(as a corollary of their calculation (\ref{eq.LMO_lens}))
that
the LMO invariant does not separate lens spaces.\index{lens space}
They also showed in \cite{BaLa} that
the LMO invariant
separates integral homology Seifert fibered spaces.
\end{rem}

\begin{prob}
Does there exist an integral/rational homology 3-sphere $M$
such that $\LMO(M) = \LMO(S^3)$?\index{LMO invariant!strength of ---} 
\end{prob}

\subsection{Characterization of the image of the LMO invariant}

\begin{prob}
Characterize those elements of $\hat\cA(\emptyset)\conn$ which are of the form 
$\log \LMO (M)$
for integral/rational homology 3-spheres.\index{LMO invariant!image of ---} 
\end{prob}

\begin{rem}
Since $\tau^{SO(3)}(M)$ can be obtained from $\LMO(M)$
by applying the weight system $W_{sl_2}$,
some characterization of this problem
might be obtained from 
the characterization of the form $\tau^{SO(3)}(M)$
(Problem \ref{prob.ch_tauSO3}),
say, from the integrality of the coefficients of $\tau^{SO(3)}(M)$
for integral/rational homology 3-spheres $M$.
Some other characterization of this problem might by obtained from
the loop expansion of the Kontsevich invariant.
\end{rem}

\subsection{Variations of the LMO invariant}

\begin{prob}
Construct the LMO invariant with coefficients in a finite field.\index{LMO invariant!--- in a finite field}\index{finite field} 
\end{prob}

\begin{rem}
If the Kontsevich invariant with coefficients in a finite field
would be constructed (see Problem \ref{prob.Kinv_ff}),
then it would be helpful for this problem.
\end{rem}

\begin{prob}
Construct the LMO invariant\index{LMO invariant!--- in arrow diagrams}\index{arrow diagram!LMO invariant in ---}
(or the theory of finite type 
invariants)\index{finite type invariant!--- of homology spheres!--- in arrow diagrams}
in arrow diagrams.
\end{prob}

\subsection{Refinements of the LMO invariant}

\noindent {{\namae}(T. Le)}\quad 
As mentioned in a remark in Problem
\ref{prob.calLMO}, the LMO invariant is a weak invariant when
$b_1(M) > 0$; in particular, $\LMO(M) = 1$ when $b_1(M) > 3$. The
following two problems might give refinements of $\LMO(M)$ which
would be stronger than $\LMO(M)$, in particular, when $b_1(M) >
0$.

\begin{prob}[T. Le, V. Turaev] 
\label{prob.spinLMO}
Define the LMO invariant\index{LMO invariant!refinement of ---}  
$\LMO(M,\sigma)$ of the pair of a closed
3-manifold $M$ and a spin structure $\sigma$ of $M$ such that
$\LMO(M) = \sum_\sigma \LMO(M,\sigma)$, where the sum runs over
all spin structures on $M$. There is also a similar problem for
spin${}^c$ structures.
\end{prob}
\eject

\begin{rem}
The quantum $SU(2)$ invariant
 of $(M,\sigma)$ satisfies that
$\tau_r^{SU(2)}(M) = \sum_\sigma \tau_r^{SU(2)}(M,\sigma)$ for
$r$ divisible by $4$ (see \cite{KM_inv}). The $\LMO(M,\sigma)$
should be defined such that $\tau_r^{SU(2)}(M,\sigma)$ can be
recovered from $\LMO(M,\sigma)$ in an appropriate sense, and such
that the coefficients of $\LMO(M,\sigma)$ are ``finite type
invariants'' of $(M,\sigma)$ under an appropriate definition of
finite type invariants of $(M,\sigma)$.
\end{rem}

The set of spin structure is a torsor over $H^1(M;\Z/2\Z)$ in the
sense that the difference of two spin structures is an element in
$H^1(M;\Z/2Z)$, and every element of $H^1(M;\Z/2Z)$ is the
difference of some spin structure and a fixed one. Similarly, the
set of all spin${}^c$ structure is a torsor over $H^1(M,\Z)$. In this
sense the previous problem might be related to the following problem.

\begin{prob}[T. Le, V. Turaev] 
\label{prob.xiLMO}
For every element $\xi\in H^1(M,\Z)$ construct an extension of
$\LMO(M,\xi)$ of the LMO invariant\index{LMO invariant!refinement of ---}  
such that when $\xi=0$ one
recovers the usual LMO invariant.
\end{prob}

The idea is that the usual LMO invariant corresponds only to the
trivial cohomology class, and for manifolds with high Betti
number, it is equal to 0. K. Habiro has an extension of the LMO
invariant that might be a solution to this problem.

\begin{rem}
For a finite abelian group $A$ and $\xi \in H^1(M,A)$,
let $\tau(M,\xi)$ be the invariant of $(M,\xi)$,
defined from a modular $A$-category,\index{tensor category!modular ---!modular A-category}
and let $\tau(M)$ be the invariant of $M$
derived from a modular category forgetting $A$-grading.
Then, $\tau(M) = \sum_\xi \tau(M,\xi)$.
(For details, see \cite{LeTuraev}.)
The $\LMO(M,\xi)$ should be defined such that
a suitable $\tau(M,\xi)$ can be recovered from $\LMO(M,\xi)$
in an appropriate sense,
and such that the coefficients of $\LMO(M,\xi)$ are
``finite type invariants'' of $(M,\xi)$
under an appropriate definition of finite type invariants of $(M,\xi)$.
\end{rem}

\subsection{Other problems}

\begin{quest}
\label{quest.KT_LMO}
\begin{itemize}
\item[\rm(1)]
Find a surgery formula for the 
Kuperberg-Thurston invariant\index{Kuperberg-Thurston invariant}
\cite{Kuperberg-Thurston}  
in terms of the Chern-Simons series
of Question \ref{quest.Kinv=CS}
\item[\rm(2)]
Compare the Kuperberg-Thurston invariant to the LMO invariant.\index{LMO invariant!--- and Kuperberg-Thurston invariant} 
\end{itemize}
\end{quest}

\begin{prob}[D. Thurston] 
\label{prob.tor_KTinv}
Do configuration spaces of \cite{Kuperberg-Thurston}  
have torsion\index{torsion}\index{configuration space!--- of Kuperberg-Thurston invariant}    
in $\Z$-homology?
Does such torsion deduce a torsion invariant of homology 3-spheres?
\end{prob}


\newpage

\section{Other problems}

\subsection{(Pseudo) Legendrian knot invariants}
\label{sec.leknotinv}

\renewcommand{\thefootnote}{\fnsymbol{footnote}}
\footnotetext[0]{Section \ref{sec.leknotinv} was written by R. Benedetti.}
\renewcommand{\thefootnote}{\arabic{footnote}}

\noindent
Let $W$ be a compact closed oriented $3$-manifold. $(K,v)$ is said a
{\it pseudo Legendrian pair} in $W$ if $K\subset W$ is a knot, $v$ is
a non singular vector field on $W$ and $K$ is transverse to $v$. $K$
is simply said a (pL)-knot. $(K_t,v_t)$, $t\in [0,1]$, is a {\it
pseudo Legendrian isotopy} if $K_t$ is an ambient isotopy of knots,
$v_t$ is a homotopy of fields and $(K_t,v_t)$ is a (pL)-pair for every
$t\in [0,1]$.  Every (pL)-knot is naturally a {\it framed} knot, and
every (pL)-isotopy is in particular a framed knots isotopy.  If $\xi$
is a transversely oriented contact structure on and $K$ is
$\xi$-Legendrian in the classical sense, then $K$ is a (pL)-knot
w.r.t. any field $v$ which is positively transverse to $\xi$. Every
Legendrian isotopy between $\xi$-Legendrian knots induces a
(pL)-isotopy.  So we have $3$ categories of knots, related by natural
forgetting maps:
$$ 
\{ \mbox{Legendrian knots} \}\stackrel{f_1}{\rightarrow} 
\{ \mbox{(pL)-knots} \} \stackrel{f_2}{\rightarrow} 
\{ \mbox{framed knots} \}.
$$
\noindent Note that, for each one of these categories, $\cC$ say, also
the $\cC$-{\it homotopy immersion class} of any $\cC$-knot is
naturally defined, this contains the $\cC$-isotopy class and is
preserved by the forgetting maps.

In \cite{benedetti.B-P1} one has introduced the 
{\it Reidemeister-Turaev torsions}\index{Reidemeister-Turaev torsion!--- of Legendrian knot} 
of (pL)-knots; one has realized that torsions include a
correct lifting to the (pL)-category of the classical Alexander
invariant; moreover, in many cases (for instance when $W$ is a
$\Z$-homology sphere), they can distinguish (pL)-knots which are
isotopic as framed knots.

\begin{quest}[R. Benedetti] 
\label{benedetti.1} 
Are torsions actually sensitive only to the (pL)-{\rm homotopy
immersion classes} of (pL)-knots?
\end{quest}

If one fix a $\cC$- homotopy immersion class of knots, say $\alpha$ , then
one can define the set of 
{\it finite type invariants} $\cF (\alpha)$
of the $\cC$-isotopy classes contained in $\alpha$. If $\alpha_0$ is a
class of Legendrian knots, one can take $\alpha_1 = f_1(\alpha_0)$ and
$\alpha_2 = f_2(\alpha_1)$; a finite type invariant for $\alpha_i$
lifts to a finite type invariant for $\alpha_{i-1}$. So one has
natural maps
 $$ \cF (\alpha_2) \stackrel{f_2^*}{\rightarrow}  \cF(\alpha_1)
\stackrel{f_1^*}{\rightarrow}  \cF (\alpha_0) .$$
\noindent It is known \cite{benedetti.F-T} that, under certain hypotheses on $W$
(for instance when $W$ is a $\Z$-homology sphere), $f_1^*\circ f_2^*$ is a
bijection. On the oder hand, one can find in \cite{benedetti.T} examples where
$f_1^*\circ f_2^*$ is not surjective and Legendrian finite type invariants
can eventually distinguish some Legendrian knots which are isotopic as
framed knots. In fact one can realize that for these examples $f_2^*$
is already not surjective and that (pL)-finite type invariants can
eventually distinguish some (pL)-knots which are isotopic as framed
knots. The following conjecture is not in contradiction with all these
known results on the subject.

\begin{conj}[R. Benedetti] 
\label{benedetti.2} 
For every $W$, for every (pL)-class $\alpha_1$ 
as above, $f_1^*$ is an isomorphism. This means, in
particular, that 
finite type invariants\index{finite type invariant!--- of Legendrian knots}  
of Legendrian knots should be
definitely not sensitive to geometric (rigid) properties of the
contact structures like ``tightness''.
\end{conj}

\noindent See also \cite{benedetti.B-P2} for a more detailed discussion and
related questions.

\subsection{Knots and finite groups}
\label{sec.kfingr}

\renewcommand{\thefootnote}{\fnsymbol{footnote}}
\footnotetext[0]{Section \ref{sec.kfingr} was written by H.R. Morton.}
\renewcommand{\thefootnote}{\arabic{footnote}}

Knot groups are known to be residually finite,
that is, any non-trivial element can be
detected by a homomorphism to some finite group.

Now by Dehn's lemma and the loop theorem a knot is trivial if and only if 
its longitude represents the trivial element of the knot group. 
Consequently for each non-trivial knot
there is a homomorphism to {\em some} finite group\index{finite group}  
which carries the longitude to a non-trivial element.

\begin{prob}[H.R. Morton] 
\label{prob.morton}
From a knot diagram find an explicit such
homomorphism to some permutation group or
establish that the knot is trivial.
\end{prob}

\noindent{\bf Refinements.}

(1)\qua Give an upper bound in terms of the diagram 
for the order of the permutation groups which need to be considered.

(2)\qua See what happens if the meridians (which are all conjugate) 
are restricted to map to permutations of some specified cycle type, 
for example, single transpositions.

\begin{rem}
Every finite group is the subgroup of a permutation group, 
so no restrictions are implied here.

The language of quandles could be adopted for {2} 
when referring to the chosen meridian conjugacy class.

It is possible to represent some knot groups
onto a finite non-cyclic group with the longitude mapping trivially. 
This always happens when $n$-colouring a knot, 
as the knot group is mapped onto the dihedral group $D_n$,
and the longitude goes into its commutator subgroup.  
The problem here focusses on the stronger question 
of representing the longitude non-trivially.
\end{rem}

\subsection{The numbers of 3-, 5-colorings and some local moves}

\begin{figure}[ht!]
\begin{align*}
\mboxsm{The 3-move}: & \quad
\pict{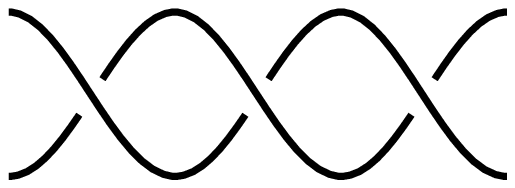}{1cm} \longleftrightarrow \pict{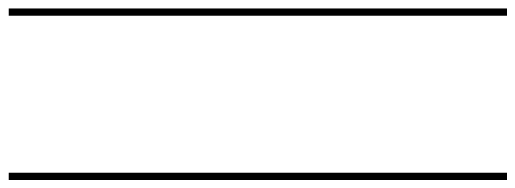}{1cm} \\
\mboxsm{The (2,2)-move}: & \quad\quad
\pict{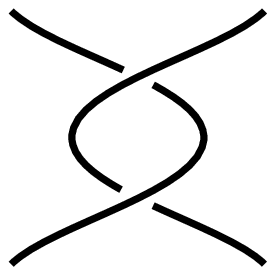}{1.5cm} \longleftrightarrow \pict{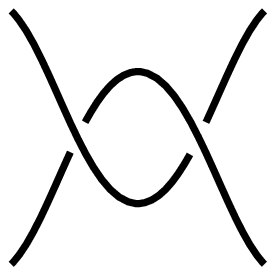}{1.5cm} 
\end{align*}
\caption{\label{fig.3-22-move}
The 3-move and the (2,2)-move}\index{local move!(2,2)-move}\index{local move!3-move}
\end{figure}

A {\it $p$-coloring} of a link $L$
is a homomorphism of the link quandle of $L$
to the dihedral quandle $R_p$ of order $p$
(or, alternatively,
a homomorphism of $\pi_1(S^3-L)$
to the dihedral group of order $2p$
which takes each meridian to a reflection).\footnote{\footnotesize 
The original definition of a 3-coloring by Fox
(see \cite[Chapter VI, Exercise 6]{CrFo63})
is (an equivalent notion of) a {\it non-trivial} homomorphism
of the link quandle of $L$ to the dihedral quandle $R_3$.
Przytycki \cite{Przytycki_c} studied the number of 3-colorings.
His definition allows trivial homomorphisms.}
Let $Col_p(L)$ denote the number of $p$-colorings of $L$
(see the remark of Problem \ref{prob.skeinmod_7col}).
The following conjecture implies that
the 3-move (see Figure \ref{fig.3-22-move}) would topologically characterize
the partition of the set of links given by $Col_3(L)$;
note that $Col_3(L)$ is unchanged under the 3-move.

\begin{conj}[3-move conjecture, Y. Nakanishi \cite{Nakanishi_guo}]
\label{conj.3-move}
Any link can be related to a trivial link
by a sequence of 3-moves.\index{local move!3-move}
\end{conj}

\begin{rem}
$Col_3(L)$ is equal to $3^{n+1}$,
where $n$ is the rank of $H_1(M_2(L);\Z/3\Z)$
and $M_2(L)$ denotes the double cover of $S^3$ branched along $L$.
Further, $Col_3$ of the trivial link with $n$ components is equal to $3^n$.
Hence, if a link $L$ is related to a trivial link by 3-moves,
then such a trivial link has $\log_3 Col_3(L)$ components.
\end{rem}

\begin{rem}[{\rm \cite[Remark on Conjecture 1.59 (1)]{Kirby}}] 
Since $B_n / \langle \sigma_i^3 \rangle$ is finite for $n \le 5$,
the proof of this conjecture for closures of braids of at most 5 strands
is reduced to verifying finitely many cases.
According to Y. Nakanishi, the smallest known obstruction of this conjecture
is the 2-parallel of a set of Borromean rings.
\end{rem}

\begin{rem}[{\rm \cite{Stoimenow_g2}}]
This conjecture is true for weak genus two knots.
\end{rem}

\begin{update}
Dabkowski and Przytycki \cite{DaPr_3m} showed that
some links cannot be reduced to trivial links by 3-moves,
which are counterexamples to this conjecture.
\end{update}

It is shown in \cite{HaUc93} that
$Col_5(L)$ is invariant under the (2,2)-move (see Figure \ref{fig.3-22-move}).
The following conjecture implies that
the (2,2)-move would topologically characterize
the partition of the set of links given by $Col_5(L)$.

\begin{conj}[{Y. Nakanishi, T. Harikae \cite[Conjecture 1.59 (6)]{Kirby}}] 
\label{conj.22-move}
Any link can be related to a trivial link
by a sequence of (2,2)-moves.\index{local move!(2,2)-move}
\end{conj}

\begin{rem}
This conjecture holds for algebraic links;
see \cite[Conjecture 1.59 (6)]{Kirby}, \cite{Przytycki_c},
and references therein.
\end{rem}

\subsection{Knotted trivalent graphs}
\label{sec.ktg}

\renewcommand{\thefootnote}{\fnsymbol{footnote}}
\footnotetext[0]{Section \ref{sec.ktg} was written by T. Ohtsuki,
following seminar talks given by D. Thurston.}
\renewcommand{\thefootnote}{\arabic{footnote}}

D. Bar-Natan and D. Thurston
\cite{BaTh_KTG1,BaTh_KTG2,Thurston_lect}
developed the theory of knotted trivalent graphs and their algebra,
related to shadow surfaces of V. Turaev \cite{Turaev_book}
and Lie groups/algebras.

\begin{figure}[ht!]
$$
\begin{array}{lrcl}
\mboxsm{Connected sum}: & 
\pict{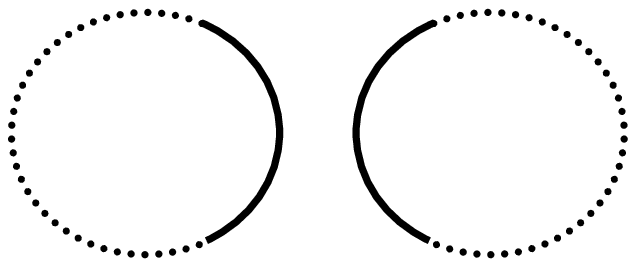}{1.3cm} &\longrightarrow& \pict{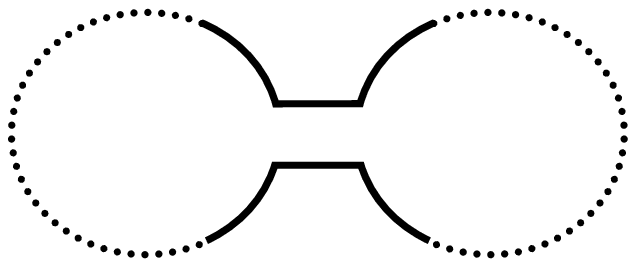}{1.3cm} \\
\mboxsm{Unzip}: & 
\pict{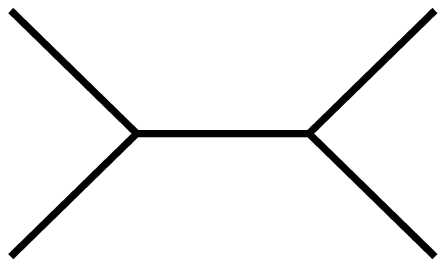}{1.3cm} &\longrightarrow& \pict{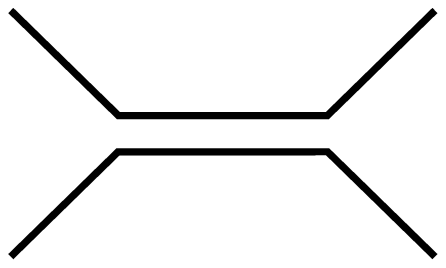}{1.3cm} \\
\mboxsm{Bubbling}: & 
\pict{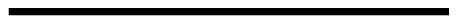}{0.06cm} &\longrightarrow& \pict{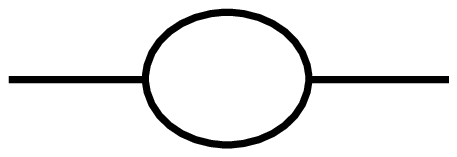}{0.95cm} \\
\mboxsm{Unknot}: & 
\qquad \emptyset \qquad &\longrightarrow& \pict{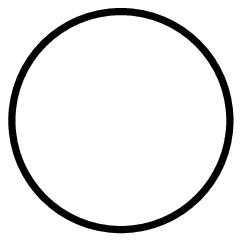}{1.3cm}
\end{array}
$$
\caption{\label{fig.KTGop}
Four operations of KTG's \cite{BaTh_KTG2}.
The left hand side of the connected sum denotes a disjoint union
of two separate graphs.}
\end{figure}

A {\it knotted trivalent graph}\index{knotted trivalent graph} 
({\it KTG})
is a (framed) embedding of
a (ribbon) trivalent graph $\Gamma$ into $S^3$,
where framing is an integer of a half integer
(hence, the ribbon of a trivalent graph is not necessarily orientable).
There are four operations of KTG's:
connected sum, unzip, bubbling and unknot; see Figure \ref{fig.KTGop}.
Any KTG (in particular, any link) can be obtained from copies of 
tetrahedron and M\"obius strip with $\pm1/2$ framing
by applying KTG operations.
Further, two sequences of KTG operations give the same KTG,
if and only if they are related by 
certain (finitely many) relations including
the pentagon and hexagon relations (see \cite{BaTh_KTG2}). 
Thus, the theory of KTG's is finitely presented in this sense.

The Kontsevich invariant of framed links have an extension for KTG's 
(see \cite{MuOh97})\index{Kontsevich invariant!--- of knotted trivalent graphs}
and the extended Kontsevich invariant is well-behaved under the KTG operations
such that they give another construction of the Kontsevich invariant
starting from the invariants of tetrahedron and M\"obius strip.

\begin{prob}
Find a new proof of the existence of a universal Vassiliev invariant of knots,
presenting them by KTG's and their operations.\index{Vassiliev invariant!universal ---!--- by knotted trivalent graphs}\index{knotted trivalent graph!universal Vassiliev invariant by ---}
\end{prob}

\begin{conj}[D. Bar-Natan, D. Thurston] 
\label{conj.KTG_ae}
For each compact Lie group $G$, level $k$, and every KTG $K: \Gamma \to \R^3$,
there exists a collection of measures $\mu_{{}_K}$ 
on the space of gauge equivalence classes of $G$-connections on $\Gamma$
satisfying the following conditions.
\begin{itemize}
\item
It is well-behaved under KTG operations.
\item
It is ``localized'' near connections that extend to $S^3-K$.
\item
A half-twist framing change acts by $e^{\sqrt{-1} H \hbar /2}$,
where $H$ is the Schr\"odinger operator on $G$.
\item
It recovers quantum invariants by\index{quantum invariant!--- from knotted trivalent graph}\index{knotted trivalent graph!quantum invariant from ---}
$$
I_R(K) = \int h_R(A) d \mu_{{}_K}(A),
$$
where $h_R(A)$ denotes the holonomy of $A$ in $R$.
Here, $R$ is a set of representations of $G$ associated to edges of $\Gamma$
and appropriate intertwiners associated to vertices of $\Gamma$.
\end{itemize}
\end{conj}

\begin{rem}[\rm \cite{BarNatan_ic,BaTh_KTG1}] 
The physical presentation of the quantum invariant of a knot $K$
associated with a representation $R$ of $G$
is given by the Chern-Simons path integral,\index{Chern-Simons!--- path integral}
$$
Z_k(S^3,K) = \int h_R(A) e^{2\pi\sqrt{-1} k \, {\rm CS}(A)} {\cal D} A,
$$
where ${\rm CS}(A)$ denotes the Chern-Simons functional of $A$
and the integral is a formal integral
over the infinite dimensional space of all $G$ connections on $S^3$.
It is a motivation of Conjecture \ref{conj.KTG_ae}
that a collection of $\mu_{{}_K}$ should play a role
of $e^{2\pi\sqrt{-1} k \, {\rm CS}(A)} {\cal D} A$.
It is expected \cite{BarNatan_ic,BaTh_KTG1} that
the collection of measures $\mu_{{}_K}$ of Conjecture \ref{conj.KTG_ae}
would prove the asymptotic expansion conjecture (Conjecture \ref{conj.AEC}).
\end{rem}

\begin{prob}
Construct an invariant of KTG's 
from configuration space integrals 
in a natural way.\index{configuration space!--- integral!--- for knotted trivalent graphs}\index{knotted trivalent graph!configuration space integral for ---}
\end{prob}

Turaev \cite{Turaev_book} introduced a presentation of 3-manifolds
as $S^1$-bundles over ``sha\-dow surfaces'', as follows
(for details see \cite{Turaev_book,BaTh_KTG2,Thurston_lect}).
A {\it fake surface} is a singular surface
such that a neighborhood of each point is homeomorphic to 
an open subset of the cone over a tetrahedron.
A $S^1$-bundle over a fake surface can appropriately be defined
and its isomorphism class is determined by the Chern number,
which is an integer or half-integer associated to each face;
we call the Chern number the {\it gleam}.
A {\it shadow surface} is a fake surface with gleams associated to the faces.
Every (closed) 3-manifold can be presented by a $S^1$-bundle 
over a (closed) shadow surface.
The pentagon and hexagon relations
(see \cite[Figure 1.1 of Chapter VIII]{Turaev_book})
are moves among shadow surfaces
which present a homeomorphic 3-manifold,
though they are not enough to characterize 
a homeomorphism class of 3-manifolds.

\begin{exe}
Find a complete set of moves among shadow surfaces\index{shadow surface}
which present a homeomorphic 3-manifold.
\end{exe}

\begin{figure}[ht!]
\begin{center}
\begin{picture}(320,100)
 \put(0,70){\small Framed links}
 \put(63,73){\vector(1,0){30}}
 \put(60,80){\footnotesize exterior}
 \put(100,70){\small Framed link exteriors}
  \put(203,73){\vector(1,0){30}}
 \put(200,80){\footnotesize surgery}
 \put(240,70){\small Closed 3-manifolds} 
 \put(-25,40){\footnotesize presentation}
 \put(30,25){\vector(0,1){40}}
 \put(160,35){\shortstack[l]{\footnotesize by making \\
                             \footnotesize $S^1$-bundle}}
 \put(155,25){\vector(0,1){40}}
 \put(285,35){\shortstack[l]{\footnotesize by making \\
                             \footnotesize $S^1$-bundle}}
 \put(280,25){\vector(0,1){40}}
 \put(-20,0){\shortstack[l]{\small Certain sequences \\
                          \small of KTG's}}
 \put(63,13){\vector(1,0){30}}
 \put(65,-5){\shortstack[l]{\footnotesize time \\ 
                           \footnotesize evolution}}
 \put(120,0){\shortstack[l]{\small Collapsible \\ 
                            \small shadow surfaces}}
 \put(203,13){\vector(1,0){30}}
 \put(200,0){\footnotesize cap off $\partial$}
 \put(240,10){\small Closed shadow surfaces}
\end{picture}
\end{center}
\caption{\label{fig.l3Ks}
Links, 3-manifolds, KTG's, and shadow surfaces}
\end{figure}
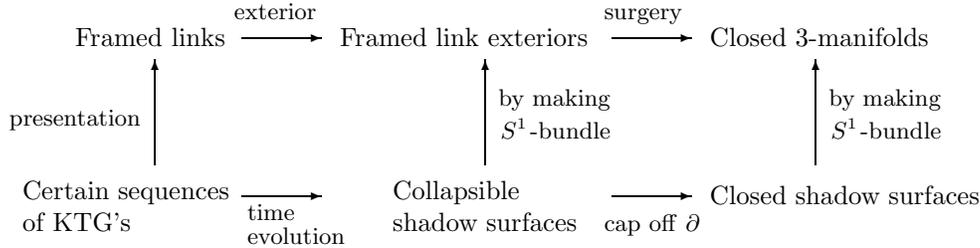

We obtain a shadow surface as a time evolution of a sequence of KTG's
given by KTG operations.
Thus, we have relations among links, 3-manifolds, KTG's and shadow surfaces
as in the commutative diagram in Figure \ref{fig.l3Ks};
for detailed statements see \cite{BaTh_KTG2,Thurston_lect}.

Motivated by a complexity of 3-manifolds discussed in 
\cite{Matveev_cplx,MaPe01,MaPe_new},
D. Thurston introduced the shadow number of 3-manifolds.
The {\it shadow number} is defined to be
the minimal number of vertices of a shadow surface.
All graph manifolds have shadow number 0 and all surgeries
on the Borromean rings have shadow number 1.
The volume conjecture\index{volume conjecture!--- and shadow number} 
might be related to the following conjecture.

\begin{conj}[D. Thurston] 
\label{conj.shadow_num}
The shadow number of a 3-manifold is quasi-linear in its Gromov norm.
That is, there exist constants $c_1$ and $c_2$ such that
$$
c_1 || M || \le \mbox{(shadow number of $M$)} \le c_2 || M ||
$$
for any 3-manifold $M$,
where $|| M ||$ denotes the Gromov norm\index{simplicial volume!--- and shadow number}   
of $M$.
\end{conj}

\begin{rem}[\rm (D. Thurston)] 
It is easy to bound the Gromov norm in terms of the shadow number
(i.e.\ to prove the left inequality for some $c_1$).
\end{rem}

\begin{rem}[\rm (D. Thurston)] 
It is shown by W. Thurston that
the hyperbolic volume of a hyperbolic 3-manifold is quasi-linear 
in the minimal number of ideal tetrahedra in a ``spun triangulation'' 
(i.e.\ the minimal number of ideal tetrahedra 
in some link complement in the 3-manifold).
It is shown by J. Brock \cite{Brock_WP} that
the volume of a mapping torus is quasi-linear in the pants translation
distance (for fixed genus).
\end{rem}

Lackenby \cite{Lackenby_v} showed
that alternating knot diagrams give good information about the
hyperbolic volume.  Knot diagrams are a special case of shadow
diagrams, but shadow diagrams can be much more efficient.  This
suggests the following problem:

\begin{prob}[D. Thurston] 
\label{prob.shadow_vol}
Find a condition on shadow diagrams which is satisfied
by shadow diagrams from alternating knots; and gives a lower bound
on the hyperbolic volume.\index{hyperbolic volume!--- and shadow number}
\end{prob}

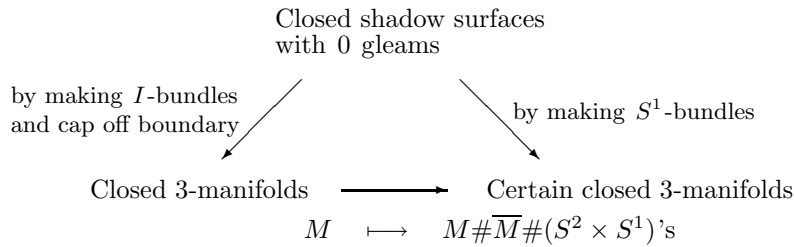
\begin{figure}[ht!]
\begin{center}
\begin{picture}(250,100)
 \put(70,70){\shortstack[l]{\small Closed shadow surfaces \\
                             \small with $0$ gleams}}
 \put(80,60){\vector(-1,-1){30}}
 \put(-30,40){\shortstack[l]{\footnotesize by making $I$-bundles \\
                            \footnotesize and cap off boundary}}
 \put(140,60){\vector(1,-1){30}}
 \put(160,45){\footnotesize by making $S^1$-bundles}
 \put(0,15){\small Closed 3-manifolds}
 \put(150,15){\small Certain closed 3-manifolds}
 \put(95,18){\vector(1,0){40}}
 \put(80,0){\small $M \quad\longmapsto\quad M \# \overline{M} \# (S^2 \times S^1)$'s}
\end{picture}
\end{center}
\caption{\label{fig.ss_3mfds}
Two ways to obtain 3-manifolds from shadow surfaces}
\end{figure}

The Reshetikhin-Turaev invariant and
the Turaev-Viro-Ocneanu invariant can be described
in terms of the KTG algebra,
via $I$-bundles and $S^1$-bundles over shadow surfaces respectively.
The relation between the two invariants is derived from
the relation between the two construction of 3-manifolds 
shown in Figure \ref{fig.ss_3mfds}.

\begin{prob}
Construct a universal Reshetikhin-Turaev invariant\index{Reshetikhin-Turaev!universal --- invariant} 
and
a universal Turaev-Viro-Ocneanu invariant\index{Turaev-Viro!Turaev-Viro-Ocneanu invariant!universal ---} 
of closed 3-manifolds,
in terms of the KTG algebra.
\end{prob}

\begin{rem}
The LMO invariant
and the even degree part of it
might be a universal Reshetikhin-Turaev invariant
and a universal Turaev-Viro-Ocneanu invariant
of rational homology 3-spheres, respectively.
\end{rem}

\subsection{Quantum groups}

\begin{prob}[J. Roberts]  
\label{prob.roberts_wqg}
What are quantum groups?\index{quantum group!interpretation of ---} 
\end{prob}

\begin{rem}[{\rm (J. Roberts)}]  
A naive answer is to simply define them by means of generators
and relations, but this is appallingly unsatisfying. Better is
Drinfel'd's original construction \cite{Drinfeld_ICM86}, 
which begins with the
geometric construction of {\em quasi-quantum groups} using the
monodromy of the KZ equation. He then uses completely algebraic
results about uniqueness of deformations to obtain from each one a
quantum group, whose category of representations is equivalent to that
of the quasi-quantum group, though the first has a trivial associator
and a complicated $R$-matrix, the second vice versa. (In particular,
the braid group representation associated to a quantum group is {\em
local} in the sense that the $R$-matrix implementing the action of a
braid generator on a tensor product of representations of the quantum
group involves only the tensor factors associated to the two strings
concerned. This is certainly not true for the KZ equation. Is there
any way to understand this using geometry?)

These constructions are very subtle and complicated. What really is a
quantum group, in fact? I believe that algebraists have some
reasonably geometric descriptions of pieces of them in terms of
perverse sheaves, etc., but I do not pretend to understand
these. Atiyah made the very interesting suggestion that quantum groups
might be in some sense the ``quaternionifications'' of compact Lie
groups. Literal quaternionification does not make sense, but
substitutes might be available, in the sense that hyperk\"ahler
geometry provides a working substitute for the non-existent
quaternionic version of complex manifold theory.  Some evidence for
this point of view is presented in Atiyah and Bielawski \cite{AtBi01}. 
\end{rem}

\subsection{Other problems}

\begin{prob}[N.~Askitas] 
\label{prob.gs}
Can a knot of 4-genus $g_s$ always be sliced (made into a slice
knot) by $g_s$ crossing switches?
\end{prob}

\begin{rem}[\rm (A. Stoimenow)] 
Clearly (at least) $g_s$ crossing switches are needed, but sometimes
more are needed to {\em unknot} the knot.
\end{rem}

\begin{update}
Livingston \cite{Livingston_sn} showed that
the knot $7_4$ provides a counterexample to this problem;
$g_s(7_4)=1$ but no crossing change results in a slice knot.
\end{update}

\begin{prob}[M. Boileau {\cite[Problem 1.69 (C)]{Kirby}}] 
\label{prob.mutant_un}
Are there mutants\index{mutation!--- and unknotting number} 
of distinct unknotting numbers?
\end{prob}

\begin{rem}[\rm (A. Stoimenow)] 
There are mutants of distinct genera (Gabai \cite{Gabai_f}) 
and slice genera (Livingston \cite{Livingston_nc}).
\end{rem}

Let $G$ be the graph such that
its vertices are isotopy classes of unoriented knots,
and two vertices are adjacent if
the corresponding knots differ by a single crossing change.

\begin{conj}[{X.-S. Lin \cite{Lin_website}}] 
\label{conj.autoG}
Any automorphism of $G$ is either the identity or the mirror map,
that is,\index{graph of knots}
any automorphism of $G$ is induced by a diffeomorphism of the ambient space.
\end{conj}

\begin{prob}[{X.-S. Lin \cite{Lin_website}}] 
\label{prob.long_rope}
What is the homotopy type of the space\index{space of knots} 
$L(K)$ of long ropes 
(as shown in the picture below) with the fixed knot type $K$?
$$
\pict{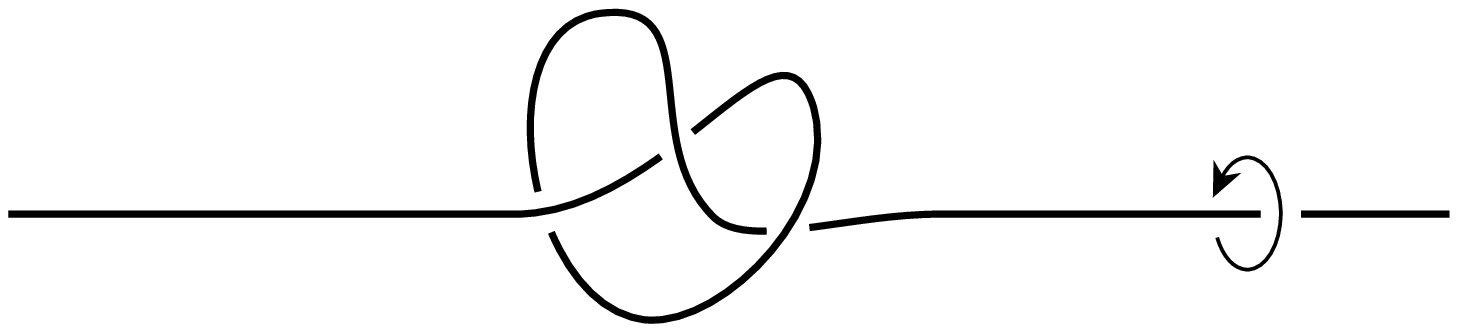}{1.5cm}
$$
\end{prob}

\begin{rem}[{\rm \cite{Lin_website}}] 
A conjecture would be that, if $K$ is a prime knot, 
$L(K)$ is homotopy equivalent to the circle 
if and only if $K$ is non-trivial, with the fundamental group
generated by the obvious loop in $L(K)$ shown in the above picture. 
This question is motivated by the paper \cite{Mostovoy}.
If the conjecture holds, the homotopy type of the space of short
ropes studied by Mostovoy would be clear. 
A paper of Hatcher \cite{Hatcher}
seems to be related with this problem.
\end{rem}

\begin{prob}[J. Roberts]  
\label{prob.roberts11}
Extend Kuperberg's work on webs. 
\end{prob}

\begin{rem}[{\rm (J. Roberts)}]  
Kuperberg posed in \cite{Kup96} the question of giving a
presentation, as a tensor category,\index{tensor category!presentation of ---} 
of the representation category of
a compact Lie group or quantum group.
The generators should be
(roughly) the fundamental modules and their bilinear and trilinear
invariants; more complicated morphisms in the category can be built
out of these according to a graphical calculus (essentially Penrose's
tensor calculus) of ``webs''. The first main problem is to describe a
set of elementary linear relations (skein relations) among such
pictures which generates all the relations among morphisms in the
category. The second is to describe a canonical basis of any invariant
space in terms of canonical pictures in the disc. Kuperberg solved
both these problems for groups of ranks one (in which case the
pictures are just Temperley-Lieb diagrams) and two and, with Khovanov
in \cite{KhKu99}, made tantalising but imprecise conjectures about how
in the higher-rank case the pictures might be related to the geometry
of the weight lattice. These ideas are closely related to the work of
Vaughan Jones \cite{Jon00} on planar algebra, which is a similar kind
of calculus describing the category of bimodules over a subfactor.
(Aside: Is it possible to find a bimodule category whose
intertwining rules are described by quasiperiodic Penrose tiles?)
\end{rem}


\begin{prob}[J. Roberts]  
\label{prob.roberts12}
Extend the theory of measured laminations\index{measured lamination}
to higher rank groups. 
\end{prob}

\begin{rem}[{\rm (J. Roberts)}]  
Let $\Sigma$ be a closed oriented surface of genus $g$, and let
$C(\Sigma)$ be its set of multicurves (isotopy classes of collections
of disjoint simple closed curves). Let $T(\Sigma)$ be its
Teichm\"uller space; that is, the space of hyperbolic structures,
considered up to diffeomorphisms isotopic to the
identity. Topologically, $T(\Sigma)$ is an open ball of dimension
$6g-6$.

Each of $C(\Sigma), T(\Sigma)$ has a natural embedding in the space of
functions $C(\Sigma) \rightarrow \R_{\geq 0}$: one sends a multicurve to its
associated minimal geometric intersection number function, and a
metric to its associated geodesic length function. It is a remarkable
fact that the $\R_+$-projective boundaries of these sets
coincide. They define the space of {\em measured laminations}, which
compactifies $T(\Sigma)$ into a closed ball and is of great importance
in Thurston's theory of surface automorphisms. For further details see
for example Penner and Harer \cite{PeHa92}.

Now $T(\Sigma)$ may also be described algebraically as a certain
component of the space of flat $SL(2,\R)$ connections on $\Sigma$
(that is, homomorphisms $\pi_1(\Sigma) \rightarrow SL(2, \R)$), and in this
context the geodesic length function is replaced by a
trace-of-holonomy function. Is there a generalisation of the above
picture to a higher rank group such as $SL(n,\R)$?

Hitchin \cite{Hit92} proves that in fact the space of flat $SL(n,\R)$
connections has a special ``Teichm\"uller component'', which is
topologically an open ball, so we have a candidate for
$T(\Sigma)$.(Aside: he asks whether there is an interpretation of the
points of the Teichm\"uller component in terms of some kind of
geometric structures on $\Sigma$. Choi and Goldman showed that for
$n=3$ they parametrise convex real projective structures, but no
general answer is known.)

A candidate for $C(\Sigma)$ might be the set of Kuperberg-style
(closed) webs drawn on the surface, for there is then a natural
holonomy-type map $T(\Sigma) \times C(\Sigma) \rightarrow \R$ which is a
substitute for the geodesic length function. (In the $SL(2)$ case,
this $C(\Sigma)$ is just the set of multicurves, as it should be.)
What might replace the geometric intersection number, and lead to some
notion of ``measured lamination'' for higher-rank groups, is unclear.
\end{rem}


\begin{prob}[J. Roberts]  
\label{prob.roberts13}
What is the generating function for $q$-spin net\index{q-spin net} 
evaluations? 
\end{prob}

\begin{rem}[{\rm (J. Roberts)}]  
A $q$-spin net is a trivalent planar graph whose edges are
labelled by irreducible representations of $SU(2)$. By placing
idempotents from the Temperley-Lieb algebra on its edges and joining
up their external strings in a planar fashion at the vertices, one
forms an evaluation in ${\Z}[q^{\pm 1}]$. The goal is to find a
power series in variables associated to the edges which serves as a
generating function for the evaluations corresponding to all possible
labellings of a given graph. Such a formula is known for any graph at
the classical value $q=1$, and Westbury \cite{Wes98} found a
generating function for the tetrahedral graph (the quantum
$6j$-symbol).\index{6j-symbol!quantum ---!--- and q-spin net} 
A general formula is, however, unknown, and Westbury
also shows that the naive guess (simply replacing factorials in the
$q=1$ formula by quantum factorials) is wrong.
\end{rem}

\begin{prob}[Y. Shinohara \cite{sto.Shi}] 
\label{prob.det_sign}
If $n=4k+1$ with $k>0$, is there a knot with 
determinant\index{determinant (of knot)} 
$n$ and
signature\index{signature!--- of knot}  
$4$?
\end{prob}

\begin{rem}[\rm (A. Stoimenow)] 
The form $4k+1$ follows from Murasugi \cite{Murasugi_cn}, 
and the condition $k\ne 0$
from a signature theorem for even unimodular quadratic forms over
$\Z$. If a counterexample for $n>1$ exists, then all prime divisors of
$n$ are of the form $24k+1$ and not smaller than $2857$. If
$\sigma_{4+8l,8l+5}$ is the elementary symmetric polynomial of degree
$4+8l$ in $8l+5$ variables, then all values of $\sigma_{4+8l,8l+5}$
on {\em positive odd} arguments are no
counterexamples, so the problem could ``reduce''
to showing that some of the $\sigma_{4+8l,8l+5}$ realizes almost all $n$
on positive odd arguments. This appears number theoretically hard,
however.
\end{rem}

The set of concordance\index{concordance}
classes of 2-strand string links
forms a group ${\cal C}_2$.
Stanford showed that ${\cal C}_2$ is not nilpotent, in
particular not abelian.

\begin{prob}[T. Stanford] 
\label{prob.C2_solv}
Is ${\cal C}_2$ solvable?
Does ${\cal C}_2$ contain a free group?
\end{prob}

\begin{prob}[A. Stoimenow] 
\label{prob.sign_chi}
Do positive links\index{positive link} 
of given signature\index{signature!--- of link} 
$\sigma$ have bounded (below) maximal
Euler characteristic $\chi$?
\end{prob}

\begin{rem}[\rm (A. Stoimenow)] 
So far for general positive links only $\sigma>0$ is known 
\cite{sto.Rud,sto.CG}, and for
positive knots $\sigma\ge 4$ if $2g=1-\chi\ge 4$
(it follows from \cite{Taniyama_pk}).
For positive
braid links the answer is positive, and also for special alternating
links by Murasugi \cite{Murasugi_cn}.
\end{rem}

\begin{prob}[A. Stoimenow] 
\label{prob.pri_ach_sli}
If a prime knot $K$ can be transformed into its mirror image by one
crossing change, is $K$ achiral\index{achiral} 
or (algebraically?) slice?
\end{prob}

\begin{rem}[\rm (A. Stoimenow)] 
Smoothing out this crossing gives a link of 
zero Tristram-Levine-signatures \cite{Tristram,Levine_kcg} 
and zero Alexander polynomial.
Many such links are slice, and then
$K$ would be slice also. But unlikely.
\end{rem}

\begin{prob}[A. Stoimenow] 
\label{prob.sto.sq}
Let $n$ be an odd natural number, different from 1, 9, and 49, such that
$n$ is the sum of two squares. 
Is there a prime alternating achiral\index{achiral} 
knot of determinant\index{determinant (of knot)} 
$n$?
\end{prob}

\begin{rem}[\rm (A. Stoimenow)] 
If there is an achiral knot  of determinant $n$, then $n$ is the odd
sum of two squares \cite{sto.Har}. 
The converse is also true, and the 
achiral knot of determinant $n$ can be chosen to be alternating
{\em or} prime, but {\em not} always both. For $n=1$, 9, and 49, 
there is no prime alternating achiral knot of determinant $n$.
If there is another such $n$, then $n>2000$ and $n$ is not
a square. See \cite{Stoimenow_ss}.
\end{rem}

\begin{conj}[V. Turaev] 
\label{conj.turaev.Alex}
A pair (a finitely generated abelian group $H$ of rank 1, 
an element  $\Delta(t) \in \Z [H/\mbox{Tors}\, H]=\Z[t^{\pm 1}]$)  
(where $t$ is a generator of $H/\mbox{Tors}\, H$) 
can be realized as the pair 
$(H_1(M)$, the Alexander polynomial $ \Delta_M$ of $M$)\index{Alexander polynomial!realization as ---}  
for a closed connected oriented 3-manifold $M$ if and only if   
$\Delta(t) = t^k \Delta(t^{-1})$ with even $k\in \Z$ and 
$\Delta(1) =\pm \vert \mbox{Tors}\, H\vert$.
\end{conj}

\begin{rem}[\rm (V. Turaev)]  
Both conditions are known to be necessary. 
They are presumably sufficient. 
This is known for $H=\Z$ and for $H=\Z \times (\Z /n\Z)$ with $n\geq 2$. 
When $M$ is obtained from $S^3$ by $0$-surgery along a knot $K$, 
$H_1(M)=\Z$ and $\Delta_M(t) = \Delta_K(t)$.
It is known that a Laurent polynomial $f(t) \in \Z[t^{\pm1}]$
is realized as the Alexander polynomial of a knot
if and only if $f(t)= t^k f(t^{-1})$ with even $k$ and $f(1)=1$. 
Using surgery on a 2-component link in $S^3$ with linking number 0 
and framing numbers $0, n$, respectively, 
one can prove (cf. \cite{Levine-67}) 
the conjecture for $H=\Z \times (\Z /n\Z)$.
\end{rem}

\newpage


\vspace{2pc}
{\small Received:\quad 28 December 2001\hfill
Revised:\quad 4 December 2002 -- 8 April 2004}

\newpage

\twocolumn[
\section*{Index of contributors}
]
\addcontentsline{toc}{section}{Index of contributors}

\small

Andersen, JE, \ref{conj.AEC_Jorgen}, \ref{prob.pexp_O_H}, \ref{conj.dj}
\vspace{0.2pc}\newline
Bar-Natan, D, \ref{prob.HN_Vinv}, \ref{conj.arrow_Ktriv}, \ref{conj.KTG_ae}
\vspace{0.2pc}\newline
Baseilhac, S,
\ref{bene.0}, \ref{bene.1}, \ref{prob.bene2}, \ref{prob.bene3},
\ref{prob.bene4}, \ref{prob.bene5}, \ref{prob.bene6}, \ref{benedetti.conj}
\vspace{0.2pc}\newline
Benedetti, R,
\ref{bene.0}, \ref{bene.1}, \ref{prob.bene2}, \ref{prob.bene3},
\ref{prob.bene4}, \ref{prob.bene5}, \ref{prob.bene6}, \ref{benedetti.conj}, 
\ref{benedetti.1}, \ref{benedetti.2}
\vspace{0.2pc}\newline
Bigelow, S,
\ref{prob.B6V62}, \ref{prob.BMWalg}, \ref{prob.repBn_BMW}, 
\ref{prob.repBn_barQ}
\vspace{0.2pc}\newline
Carter, JS, \ref{prob.HiQSnm}
\vspace{0.2pc}\newline
Deloup, F, \ref{prob.q_sigma}, \ref{prob.spinc_Yd}, \ref{quest.lift_g}
\vspace{0.2pc}\newline
Dunfield, N, \ref{prob.vol_VK-1}
\vspace{0.2pc}\newline
Fenn, R, \ref{FRSprob2}, 
\vspace{0.2pc}\newline
Ferrand, E, \ref{prob.sign_PL_FL}
\vspace{0.2pc}\newline
Garoufalidis, S, \ref{conj.hair_map_inj}
\vspace{0.2pc}\newline
Habiro, H,
\ref{conj.Cd_string_link}, \ref{prob.MK_HL}, \ref{conj.HL_IgM}, 
\ref{conj.IslM}
\vspace{0.2pc}\newline
Hansen, SK, \ref{prob.QHS_tauG_LMO}, \ref{prob.QHS_alltauG}
\vspace{0.2pc}\newline
Harikae, T, \ref{conj.22-move}
\vspace{0.2pc}\newline
Haviv, A,
\ref{conj.arrow_Ktriv}
\vspace{0.2pc}\newline
Jeong, M-J, \ref{quest.An_Vn}
\vspace{0.2pc}\newline
Kawahigashi, Y,
\ref{kawah_prob.4}, \ref{kawah_prob.5}, \ref{kawah_prob.6}, 
\ref{kawah_prob.3}, \ref{kawah_prob.1}, \ref{kawah_prob.2}, 
\ref{kawah_prob.7}
\vspace{0.2pc}\newline
Kerler, T,
\ref{prob.kerler1}, \ref{prob.kerler2}, \ref{prob.kerler3},
\ref{prob.kerler4}, \ref{prob.kerler5}, \ref{prob.kerler6}, \ref{prob.kerler7}
\vspace{0.2pc}\newline
Kidwell, M, \ref{prob.PWK=FK}
\vspace{0.2pc}\newline
Kohno, T,
\ref{kohno_prob.1}, \ref{prob.qPalg}, \ref{prob.qCMG}, 
\ref{prob.ker_tau}, \ref{prob.hol_cb},\vspace{0.2pc}
\ref{prob.PnS2AnS}
\newline
Kricker, A,
\ref{conj.p_loop_move}, \ref{prob.pres_Theta}, \ref{conj.hair_map_inj}
\vspace{0.2pc}\newline
Le, TTQ, \ref{conj.HL_IgM}, \ref{prob.spinLMO}, \ref{prob.xiLMO}
\vspace{0.2pc}\newline
Lescop, C, \ref{quest.Kinv=CS}, \ref{quest.SUn_Cinv}
\vspace{0.2pc}\newline
Lin, X.-S,
\ref{prob.desc_0_VK}, \ref{conj.torsion_ws}, \ref{prob.sign_Vinv}, 
\ref{conj.autoG}, \ref{prob.long_rope}
\vspace{0.2pc}\newline
Masbaum, G, \ref{prob.EGonHMSLk}, \ref{prob.NT_TQFT}
\vspace{0.2pc}\newline
Massuyeau, G, \ref{prob.spin_Yd}, \ref{prob.spinc_Yd}
\vspace{0.2pc}\newline
Morita, S, \ref{prob.topCS}
\vspace{0.2pc}\newline
Morton, HR, \ref{prob.morton}
\vspace{0.2pc}\newline
Murakami, H,
\ref{conj.vol_conj}, \ref{conj.vol_CS}, \ref{prob.CS_torus_knot},
 \ref{conj.vol_olim}, \ref{prob.olim_Sf}
\vspace{0.2pc}\newline
Murakami, J, \ref{conj.vol_conj}, \ref{conj.vol_CS}
\vspace{0.2pc}\newline
Nakanishi, Y, \ref{conj.3-move}, \ref{conj.22-move}
\vspace{0.2pc}\newline
Ohtsuki, T, \ref{prob.Cstr_3mfd}, \ref{prob.product_in_M}
\vspace{0.2pc}\newline
Ohyama, Y, \ref{prob.D_link_htpc}
\vspace{0.2pc}\newline
Okamoto, M, \ref{conj.vol_CS}
\vspace{0.2pc}\newline
Okuda, N, \ref{prob.Okuda_fish}
\vspace{0.2pc}\newline
Park, C-Y, \ref{quest.An_Vn}
\vspace{0.2pc}\newline
Polyak, M,
\ref{prob.top_pres_muinv}, \ref{conj.inj_A2vA}, 
\ref{conj.Vinv_ckvk}, \ref{prob.GH_vknot}, \ref{prob.Kinv_vk}, \ref{prob.khs},
 \ref{prob.qquandle}, \ref{prob.spinCWL}, \ref{quest.Rohlin_spin}, 
\ref{prob.Gdf_ftinv}, \ref{conj.hc_Yd}
\vspace{0.2pc}\newline
Przytycki, J,
\ref{prob.prz.it_tS2}, \ref{prob.prz.cS2},
\ref{p_4.1}, \ref{3.2}, \ref{prob.skeinmod_7col}
\vspace{0.2pc}\newline
Roberts, J,
\ref{prob.roberts1}, \ref{prob.roberts2}, \ref{prob.roberts3},
\ref{prob.roberts4}, \ref{prob.roberts5}, \ref{prob.roberts_tpK},
\ref{drinfeldprob}, \ref{prob.roberts6}, 
\ref{prob.roberts14}, 
\ref{prob.roberts15}, \ref{prob.roberts_wqg},
\ref{prob.roberts11}, \ref{prob.roberts12},\vspace{0.2pc}
\ref{prob.roberts13}
\newline
Rourke, C, \ref{FRSprob1}, \ref{FRSprob2}, 
\vspace{0.2pc}\newline
Sanderson, B, \ref{FRSprob1}, \ref{FRSprob2}, 
\vspace{0.2pc}\newline
Sato, N, \ref{prob.sato1}, \ref{prob.sato2}, \ref{prob.sato3}
\vspace{0.2pc}\newline
Stanford, T, \ref{quest.2tor_cA}, \ref{quest.hG_Vinv}, 
\ref{quest.Milnorinv_Vinv}, \ref{prob.C2_solv}
\vspace{0.2pc}\newline
Stoimenow, A, \ref{prob.VK_span}, \ref{prob.QKunit}, \ref{prob.PWK=FK}, 
\ref{prob.sign_PL_FL}, \ref{prob.Conway_cr}, \ref{stoi.prm}, 
\ref{prob.sign_chi}, \ref{prob.pri_ach_sli}, \ref{prob.sto.sq}
\vspace{0.2pc}\newline
Takata, T, \ref{conj.vol_CS}, \ref{prob.QHS_tauG_LMO}, \ref{prob.QHS_alltauG}
\vspace{0.2pc}\newline
Thurston, D, \ref{prob.csi_vknot}, \ref{prob.inv_vol}, 
\ref{prob.vol_conj_ncL}, \ref{prob.tor_KTinv}, \ref{conj.KTG_ae}, 
\ref{conj.shadow_num}, \ref{prob.shadow_vol}
\vspace{0.2pc}\newline
Turaev, V, \ref{prob.HQFT_spin}, \ref{prob.HQFT_spin2}, \ref{prob.sf_CWL}, 
\ref{prob.spinLMO}, \ref{prob.xiLMO},\vspace{0.2pc}
\ref{conj.turaev.Alex}
\newline 
Willerton, S, \ref{conj.bound_v3}
\vspace{0.2pc}\newline
Yokota, Y, \ref{conj.vol_CS}

\vspace{0.5pc}
{\bf\normalsize Quotation}

Askitas, N, \ref{prob.gs}
\vspace{0.2pc}\newline
Boileau, M, \ref{prob.mutant_un}
\vspace{0.2pc}\newline
Bott, R, \ref{prob.bott}
\vspace{0.2pc}\newline
Goussarov, M, \ref{conj.Vinv_ckvk}
\vspace{0.2pc}\newline
Guadagnini, E, \ref{conj.GP98} 
\vspace{0.2pc}\newline
Jones, V, \ref{prob.VK=1}, \ref{prob.char_VK}, \ref{prob.TLfaithful}
\vspace{0.2pc}\newline
Kashaev, R, \ref{conj.vol_conj}
\vspace{0.2pc}\newline
Pilo, L, \ref{conj.GP98} 
\vspace{0.2pc}\newline
Rozansky, L, \ref{conj.hair_map_inj}
\vspace{0.2pc}\newline
Shinohara, Y, \ref{prob.det_sign}
\vspace{0.2pc}\newline
Viro, O, \ref{conj.Vinv_ckvk}
\vspace{0.2pc}\newline

\onecolumn

\begin{theindex}
\addcontentsline{toc}{section}{Index}

  \item 2-loop polynomial, \hyperpage{439}
    \subitem --- of knots of genus 1, \hyperpage{439}
    \subitem topological construction of ---, \hyperpage{439}
  \item 6j-symbol, \hyperpage{507}
    \subitem classification of ---, \hyperpage{507}
    \subitem non-quantum ---, \hyperpage{512}
    \subitem quantum ---, \hyperpage{507}
      \subsubitem --- and hyperbolic volume, \hyperpage{482}
      \subsubitem --- and q-spin net, \hyperpage{540}
    \subitem TQFT arising from ---, \hyperpage{510}

  \indexspace

  \item achiral, \hyperpage{541}
  \item Alexander polynomial, \hyperpage{391}
    \subitem categorification of ---, \hyperpage{390}
    \subitem realization as ---, \hyperpage{542}
  \item anomaly, \hyperpage{431}
  \item arrow diagram, \hyperpage{409}
    \subitem Kontsevich invariant in ---, \hyperpage{427}
    \subitem LMO invariant in ---, \hyperpage{528}
    \subitem space of ---, \hyperpage{411}
      \subsubitem dimension of ---, \hyperpage{409}, \hyperpage{411}
  \item associator, \hyperpage{432}
    \subitem --- in rational homotopy theory, \hyperpage{435}
    \subitem presentation of ---, \hyperpage{433}

  \indexspace

  \item BMW algebra, \hyperpage{469}
  \item braid group, \hyperpage{466}
    \subitem --- of $\Sigma$, \hyperpage{444}
    \subitem representation of ---, \hyperpage{469}
      \subsubitem faithful ---, \hyperpage{467, 468}, \hyperpage{470}
      \subsubitem irreducible ---, \hyperpage{469}
  \item Burau representation, \hyperpage{467, 468}

  \indexspace

  \item Cappell-Lee-Miller invariant, \hyperpage{515}
  \item Casson invariant, \hyperpage{385}, \hyperpage{514}
    \subitem --- and TQFT, \hyperpage{502}
    \subitem --- and signature of 4-manifold, \hyperpage{514}
    \subitem --- by counting configurations, \hyperpage{522}
    \subitem Casson-Walker invariant, \hyperpage{514}
    \subitem Casson-Walker-Lescop invariant, \hyperpage{514}
      \subsubitem --- and Reidemeister-Turaev torsion, \hyperpage{514}
      \subsubitem refinement of ---, \hyperpage{515}
    \subitem Gauss diagram formula for ---, \hyperpage{517}
  \item categorification, \hyperpage{384}, \hyperpage{388}, 
		\hyperpage{390}
  \item Chern-Simons
    \subitem --- functional, \hyperpage{471}
    \subitem --- invariant, \hyperpage{483}
      \subsubitem --- and colored Jones polynomial, \hyperpage{396}
      \subsubitem --- and quantum invariant of 3-manifold, 
		\hyperpage{481}
      \subsubitem --- of non-hyperbolic 3-manifold, \hyperpage{483}
      \subsubitem topological definition of ---, \hyperpage{483}
    \subitem --- path integral, \hyperpage{378}, \hyperpage{471}, 
		\hyperpage{535}
      \subsubitem perturbative expansion of ---, \hyperpage{477}
      \subsubitem saddle point method on ---, \hyperpage{482}
    \subitem --- theory
      \subsubitem --- for $SL(2,\C)$, \hyperpage{483}
  \item concordance, \hyperpage{541}
  \item configuration space, \hyperpage{429, 430}
    \subitem --- integral, \hyperpage{378}, \hyperpage{429, 430}
      \subsubitem --- for knotted trivalent graphs, \hyperpage{535}
      \subsubitem --- for virtual knots, \hyperpage{427}
      \subsubitem hidden strata of ---, \hyperpage{428}
      \subsubitem Kontsevich invariant by ---, \hyperpage{429}, 
		\hyperpage{431}
      \subsubitem localization of ---, \hyperpage{522}
    \subitem --- of Kuperberg-Thurston invariant, \hyperpage{529}
    \subitem homology of ---, \hyperpage{436}
  \item conformal block, \hyperpage{443, 444}
  \item Conway polynomial
    \subitem --- and crossing number, \hyperpage{392}
  \item crossing number, \hyperpage{392}, \hyperpage{403}, 
		\hyperpage{405}

  \indexspace

  \item deformation quantization, \hyperpage{443}
  \item determinant (of knot), \hyperpage{540, 541}
  \item dilogarithm
    \subitem --- function, \hyperpage{394}
    \subitem dilogarithmic invariant, \hyperpage{488}
    \subitem quantum ---, \hyperpage{393}

  \indexspace

  \item equivalence relation
    \subitem $C_d$-equivalence, \hyperpage{418, 419}
    \subitem $Y_d$-equivalence, \hyperpage{523}, \hyperpage{525}
      \subsubitem spin ---, \hyperpage{525}
      \subsubitem spin${}^c$ ---, \hyperpage{525}
    \subitem $d$-equivalence, \hyperpage{418}, \hyperpage{523}
    \subitem homology $d$-loop equivalence, \hyperpage{420}
    \subitem homotopy $d$-loop equivalence, \hyperpage{421}

  \indexspace

  \item finite field, \hyperpage{384}, \hyperpage{426}, \hyperpage{528}
  \item finite group, \hyperpage{402}, \hyperpage{508}, \hyperpage{510}, 
		\hyperpage{531}
  \item finite type invariant
    \subitem --- by local move, \hyperpage{412, 413}
    \subitem --- of Legendrian knots, \hyperpage{531}
    \subitem --- of homology cylinders, \hyperpage{525}
    \subitem --- of homology spheres, \hyperpage{518}
      \subsubitem --- and Cappell-Lee-Miller invariant, \hyperpage{515}
      \subsubitem --- by Gauss diagram formula, \hyperpage{517}
      \subsubitem --- by space of 3-manifolds, \hyperpage{522}
      \subsubitem --- in arrow diagrams, \hyperpage{528}
      \subsubitem constructive presentation of ---, \hyperpage{522}
      \subsubitem dimension of ---, \hyperpage{518}
      \subsubitem strength of ---, \hyperpage{518}
      \subsubitem torsion and ---, \hyperpage{518}
    \subitem --- of knots, \see {Vassiliev invariant}{398}
    \subitem --- of links, \hyperpage{422}
    \subitem --- of links in $\Sigma\times I$, \hyperpage{443}
    \subitem --- of spin${}^c$ 3-manifolds, \hyperpage{526}
    \subitem --- of string links, \hyperpage{400}, \hyperpage{419}
    \subitem --- of virtual knots, \hyperpage{409}, \hyperpage{412}, 
		\hyperpage{427}
      \subsubitem universal ---, \hyperpage{427}
    \subitem loop ---, \hyperpage{416}, \hyperpage{421}
  \item fusion rule algebra, \hyperpage{506}
    \subitem classification of ---, \hyperpage{507}
    \subitem non-quantum ---, \hyperpage{512}
    \subitem TQFT arising from ---, \hyperpage{510}

  \indexspace

  \item Goussarov-Habiro theory
    \subitem --- for 3-manifolds, \hyperpage{523}
    \subitem --- for knots, \hyperpage{418}
    \subitem --- for loop finite type invariant, \hyperpage{421}
    \subitem --- for spin 3-manifolds, \hyperpage{525, 526}
    \subitem --- for spin${}^c$ 3-manifolds, \hyperpage{525, 526}
    \subitem --- for virtual knots, \hyperpage{420}
  \item graph cohomology, \hyperpage{435}
  \item graph of knots, \hyperpage{538}
  \item Gromov-Witten theory, \hyperpage{385}

  \indexspace

  \item Hecke algebra, \hyperpage{384}, \hyperpage{468}
  \item HOMFLY polynomial, \see {skein polynomial}{390}
  \item homotopy QFT, \hyperpage{496}
    \subitem spin ---, \hyperpage{496}
    \subitem spin${}^c$ ---, \hyperpage{496}
  \item hyperbolic volume
    \subitem --- and $\log V_K(-1)$, \hyperpage{387}
    \subitem --- and colored Jones polynomial, \hyperpage{393}, 
		\hyperpage{396}
    \subitem --- and invariants of 3-manifold, \hyperpage{483}
    \subitem --- and quantum invariant of 3-manifold, \hyperpage{481}
    \subitem --- and shadow number, \hyperpage{536}

  \indexspace

  \item Jacobi diagram, \hyperpage{398}
    \subitem oriented ---, \see {arrow diagram}{420}
    \subitem space of ---
      \subsubitem dimension of ---, \hyperpage{405}, \hyperpage{518}
      \subsubitem torsion of ---, \hyperpage{400}, \hyperpage{518}
  \item Jones polynomial, \hyperpage{378}, \hyperpage{380}
    \subitem --- and Gromov-Witten theory, \hyperpage{385}
    \subitem --- and coloring, \hyperpage{453}
    \subitem --- of given span, \hyperpage{391}
    \subitem categorification of ---, \hyperpage{390}
    \subitem colored ---, \hyperpage{393}, \hyperpage{479}, 
		\hyperpage{487}
    \subitem image of ---, \hyperpage{382}
    \subitem interpretation of ---, \hyperpage{382--385}
    \subitem special value of ---, \hyperpage{384}
    \subitem trivial ---, \hyperpage{381}, \hyperpage{467}
    \subitem zeros of ---, \hyperpage{386}

  \indexspace

  \item Kauffman bracket, \hyperpage{380}, \hyperpage{447}
    \subitem --- skein module, \see {skein module}{446}
  \item Kauffman polynomial, \hyperpage{390}
    \subitem --- and Whitehead double, \hyperpage{391}
    \subitem --- and genus, \hyperpage{392}
    \subitem --- and signature, \hyperpage{392}
    \subitem --- and unknotting number, \hyperpage{392}
    \subitem --- of given span, \hyperpage{391}
  \item knotted trivalent graph, \hyperpage{533}
    \subitem configuration space integral for ---, \hyperpage{535}
    \subitem quantum invariant from ---, \hyperpage{534}
    \subitem universal Vassiliev invariant by ---, \hyperpage{534}
  \item Kontsevich integral, \see {Kontsevich invariant}{378}
  \item Kontsevich invariant, \hyperpage{378}, \hyperpage{424}
    \subitem --- and 2-loop polynomial, \hyperpage{439}
    \subitem --- and Jones polynomial, \hyperpage{381--383}
    \subitem --- by configuration space integral, \hyperpage{429}, 
		\hyperpage{431}
    \subitem --- in a finite field, \hyperpage{426}
    \subitem --- in arrow diagrams, \hyperpage{427}
    \subitem --- of knotted trivalent graphs, \hyperpage{533}
    \subitem --- of the trivial knot, \hyperpage{424}, \hyperpage{427}
    \subitem --- of virtual knots, \hyperpage{427}
    \subitem calculation of ---, \hyperpage{424}, \hyperpage{433}
    \subitem image of ---, \hyperpage{425}
    \subitem interpretation of ---, \hyperpage{426}, \hyperpage{437}
    \subitem loop expansion of ---, \hyperpage{437}
    \subitem strength of ---, \hyperpage{425}
  \item Kuperberg-Thurston invariant, \hyperpage{529}

  \indexspace

  \item lens space, \hyperpage{446}, \hyperpage{472}, \hyperpage{475}, 
		\hyperpage{477}, \hyperpage{480}, \hyperpage{489}, 
		\hyperpage{513}, \hyperpage{527, 528}
  \item linking pairing, \hyperpage{418}
  \item LMO invariant, \hyperpage{378}, \hyperpage{527}
    \subitem --- and Kuperberg-Thurston invariant, \hyperpage{529}
    \subitem --- in a finite field, \hyperpage{528}
    \subitem --- in arrow diagrams, \hyperpage{528}
    \subitem calculation of ---, \hyperpage{527}
    \subitem image of ---, \hyperpage{528}
    \subitem refinement of ---, \hyperpage{528, 529}
    \subitem strength of ---, \hyperpage{472}, \hyperpage{527, 528}
  \item local move, \hyperpage{412, 413}
    \subitem $C_d$-move, \hyperpage{419}
    \subitem $\Delta$ move, \hyperpage{413}, \hyperpage{415}
    \subitem (2,2)-move, \hyperpage{533}
    \subitem 3-move, \hyperpage{532, 533}
    \subitem doubled delta move, \hyperpage{413}, \hyperpage{416}
    \subitem finite type invariant by ---, \hyperpage{412, 413}
      \subsubitem loop finite type invariant, \hyperpage{416}
    \subitem loop move, \hyperpage{416}
    \subitem mod $p$ loop move, \hyperpage{418}
  \item loop expansion, \see {Kontsevich invariant}{437}

  \indexspace

  \item measured lamination, \hyperpage{539}
  \item Milnor invariant, \hyperpage{408, 409}
  \item modular category, \see {tensor category}{511}
  \item monoidal category, \see {tensor category}{506}
  \item mutation, \hyperpage{381}, \hyperpage{472}
    \subitem --- and unknotting number, \hyperpage{538}

  \indexspace

  \item Nullstellensatz, \hyperpage{422}

  \indexspace

  \item Ohtsuki series, \hyperpage{488}
  \item optimistic limit, \hyperpage{481, 482}

  \indexspace

  \item path integral, \see {Chern-Simons}{378}
  \item perturbative invariant, \hyperpage{488}
    \subitem --- and quantum invariant, \hyperpage{491}
    \subitem calculation of ---, \hyperpage{489}
    \subitem image of ---, \hyperpage{490}, \hyperpage{492}
    \subitem modular form in  ---, \hyperpage{490}
  \item positive link, \hyperpage{541}

  \indexspace

  \item Q polynomial, \hyperpage{390}
    \subitem --- of given span, \hyperpage{391}
    \subitem values of ---, \hyperpage{391}
  \item q-spin net, \hyperpage{540}
  \item quandle, \hyperpage{455}
    \subitem --- (co)homology
      \subsubitem computation of ---, \hyperpage{460}, \hyperpage{464}
    \subitem --- cocycle invariant, \hyperpage{461}
      \subsubitem --- and other invariant, \hyperpage{463}
      \subsubitem computation of ---, \hyperpage{461}
    \subitem --- space, \hyperpage{465}
    \subitem classification of connected ---, \hyperpage{456}
    \subitem knot ---, \hyperpage{458}
      \subsubitem representations of ---, \hyperpage{458}
    \subitem quantum ---, \hyperpage{464}
  \item quantum double, \hyperpage{512}
  \item quantum group, \hyperpage{378}, \hyperpage{507}
    \subitem interpretation of ---, \hyperpage{537}
  \item quantum groupoid, \hyperpage{507}
    \subitem --- and TQFT, \hyperpage{494}
  \item quantum invariant, \hyperpage{378}, \hyperpage{471}
    \subitem --- and Gromov-Witten theory, \hyperpage{385}
    \subitem --- and quandle cocycle invariant, \hyperpage{463}
    \subitem --- from knotted trivalent graph, \hyperpage{534}
    \subitem --- of 3-manifold, \hyperpage{491}, \hyperpage{493}
      \subsubitem absolute value of ---, \hyperpage{474}
      \subsubitem asymptotic behaviour of ---, \hyperpage{396}, 
		\hyperpage{477, 478}, \hyperpage{480, 481}
      \subsubitem growth rate of ---, \hyperpage{479}
      \subsubitem interpretation of ---, \hyperpage{473}
      \subsubitem optimistic limit of ---, \hyperpage{482}
      \subsubitem strength of ---, \hyperpage{472, 473}
    \subitem quantum hyperbolic invariant, \hyperpage{485}
  \item quantum quandle, \hyperpage{464}

  \indexspace

  \item rack, \hyperpage{464}
    \subitem --- (co)homology, \hyperpage{465}
    \subitem --- space, \hyperpage{465}
  \item rational Z invariant, \hyperpage{442}
  \item Reidemeister-Turaev torsion
    \subitem --- and Casson-Walker-Lescop invariant, \hyperpage{514}
    \subitem --- of Legendrian knot, \hyperpage{530}
  \item Reshetikhin-Turaev
    \subitem --- invariant, \hyperpage{511}
      \subsubitem --- and Turaev-Viro-Ocneanu invariant, 
		\hyperpage{513}
    \subitem universal --- invariant, \hyperpage{537}
  \item reversed knot, \hyperpage{402}
  \item ribbon category, \see {tensor category}{511}
  \item Rokhlin invariant, \hyperpage{514}
    \subitem --- of spin${}^c$ 3-manifold, \hyperpage{517}

  \indexspace

  \item shadow surface, \hyperpage{535}
  \item signature
    \subitem --- of 4-manifold, \hyperpage{514}
    \subitem --- of knot, \hyperpage{403}, \hyperpage{540}
    \subitem --- of link, \hyperpage{392}, \hyperpage{541}
  \item simplicial volume
    \subitem --- and colored Jones polynomial, \hyperpage{393}
    \subitem --- and invariants of 3-manifold, \hyperpage{483}
    \subitem --- and quantum invariant of 3-manifold, \hyperpage{481}
    \subitem --- and shadow number, \hyperpage{536}
  \item skein homology, \hyperpage{447}
  \item skein module, \hyperpage{445}
    \subitem $(4,\infty)$ ---, \hyperpage{450--452}
    \subitem $sl_3$ ---, \hyperpage{447}
    \subitem --- of 3-manifolds, \hyperpage{454}
    \subitem Homflypt ---, \hyperpage{447, 448}
    \subitem homotopy ---, \hyperpage{449, 450}
    \subitem Kauffman ---, \hyperpage{448, 449}
    \subitem Kauffman bracket ---, \hyperpage{446, 447}
    \subitem Vassiliev-Goussarov ---, \hyperpage{454}
  \item skein polynomial, \hyperpage{390}
    \subitem --- and signature, \hyperpage{392}
    \subitem --- of Whitehead double, \hyperpage{391}
    \subitem --- of given span, \hyperpage{391}
  \item space of knots, \hyperpage{538}
  \item state-sum
    \subitem --- invariant, \hyperpage{509}, \hyperpage{511}
      \subsubitem --- from degenerate S-matrix, \hyperpage{512}
      \subsubitem --- from strongly amenable subfactor, \hyperpage{513}
  \item subfactor, \hyperpage{507}, \hyperpage{509}
    \subitem --- distinguishing lens spaces, \hyperpage{513}
    \subitem --- with degenerate braiding, \hyperpage{513}
    \subitem classification of ---, \hyperpage{507}
    \subitem Haagerup ---, \hyperpage{513}
    \subitem strongly amenable ---, \hyperpage{513}

  \indexspace

  \item Temperley-Lieb algebra, \hyperpage{466, 467}
  \item tensor category
    \subitem modular ---, \hyperpage{511}
      \subsubitem --- and TQFT, \hyperpage{504}
      \subsubitem --- as extension of ribbon category, \hyperpage{512}
      \subsubitem classification of ---, \hyperpage{494}
      \subsubitem modular A-category, \hyperpage{529}
    \subitem monoidal ---, \hyperpage{506, 507}, \hyperpage{509}, 
		\hyperpage{511}
      \subsubitem classification of ---, \hyperpage{507}
    \subitem presentation of ---, \hyperpage{539}
    \subitem ribbon ---, \hyperpage{511}
      \subsubitem --- with degenerate S-matrix, \hyperpage{511, 512}
  \item torsion, \hyperpage{401}, \hyperpage{446}, \hyperpage{451}, 
		\hyperpage{529}
    \subitem --- free, \hyperpage{400}, \hyperpage{419}, 
		\hyperpage{446}, \hyperpage{518}, \hyperpage{524}
  \item TQFT, \hyperpage{378}, \hyperpage{493}
    \subitem --- and skein module, \hyperpage{445}
    \subitem --- from 6j-symbol, \hyperpage{510}
    \subitem --- from tensor category, \hyperpage{512}
    \subitem Casson ---, \hyperpage{502}
    \subitem classification of ---, \hyperpage{494, 495}
    \subitem cyclotomic integer ---, \hyperpage{498}
    \subitem decomposition of ---, \hyperpage{497}
    \subitem extended ---, \hyperpage{504}
    \subitem geometric construction of ---, \hyperpage{496}
    \subitem half-projective ---, \hyperpage{504}
    \subitem homological ---, \hyperpage{499}
    \subitem image of mapping class in ---, \hyperpage{497, 498}
    \subitem length 1 ---, \hyperpage{500}
    \subitem non-quantum ---, \hyperpage{512}
    \subitem semisimple ---, \hyperpage{504}
    \subitem solvable ---, \hyperpage{502}
    \subitem spin ---, \hyperpage{495}
    \subitem spin${}^c$ ---, \hyperpage{496}
  \item Turaev-Viro
    \subitem --- TQFT, \hyperpage{505}
    \subitem --- invariant, \hyperpage{509}
    \subitem --- module, \hyperpage{501}
    \subitem Turaev-Viro-Ocneanu invariant, \hyperpage{510}, 
		\hyperpage{512, 513}
      \subsubitem universal ---, \hyperpage{537}

  \indexspace

  \item Vassiliev invariant, \hyperpage{378}, \hyperpage{398}
    \subitem --- and crossing number, \hyperpage{403}
    \subitem --- and polynomial invariants, \hyperpage{423}
    \subitem --- by configuration space, \hyperpage{428}
    \subitem --- of knots in a 3-manifold, \hyperpage{454}
    \subitem detectability by ---, \hyperpage{402, 403}, 
		\hyperpage{459}
    \subitem dimension of the space of ---, \hyperpage{405}
    \subitem extension for virtual knots, \hyperpage{412}
    \subitem strength of ---, \hyperpage{401}
    \subitem universal ---, \hyperpage{424}
      \subsubitem --- by knotted trivalent graphs, \hyperpage{534}
      \subsubitem --- in a finite field, \hyperpage{426}
    \subitem weight system of ---, \hyperpage{400}
  \item Verlinde formula, \hyperpage{495}
  \item virtual knot, \hyperpage{409}
    \subitem finite type invariant of ---, \hyperpage{412}, 
		\hyperpage{427}
      \subsubitem universal ---, \hyperpage{427}
    \subitem Goussarov-Habiro theory for ---, \hyperpage{420}
  \item Vogel's algebra, \hyperpage{521}
    \subitem presentation of ---, \hyperpage{521}
  \item volume conjecture
    \subitem --- and shadow number, \hyperpage{535}
    \subitem --- for 3-manifolds, \hyperpage{481}
    \subitem --- for knots, \hyperpage{393}
    \subitem --- for links, \hyperpage{396}
    \subitem --- for other Lie groups, \hyperpage{483}
    \subitem --- for quantum hyperbolic invariant, \hyperpage{488}

  \indexspace

  \item weight system, \hyperpage{398}, \hyperpage{409}
    \subitem --- and configuration space, \hyperpage{436}
    \subitem --- from Lie bialgebra, \hyperpage{411}
    \subitem --- of Cappell-Lee-Miller invariant, \hyperpage{515}
    \subitem realization by Vassiliev invariant, \hyperpage{400}
  \item Witten-Reshetikhin-Turaev invariant, \hyperpage{471}

\end{theindex}

\end{document}